\newif\ifarxiv
\arxivtrue

\ifarxiv

\documentclass{article}

\usepackage[colorlinks=True,allcolors=blue]{hyperref}
\usepackage{amsthm}
\usepackage[margin=1in]{geometry}
\usepackage{wrapfig}

\usepackage{graphicx}
\usepackage{amsmath}
\usepackage{amssymb}
\usepackage{mathrsfs}
\usepackage{bm}
\usepackage{cancel}
\usepackage[bbgreekl]{mathbbol}
\usepackage{bbm}
\usepackage{bbold}

\makeatletter
\def\amsbb{\use@mathgroup \M@U \symAMSb}
\makeatother
\renewcommand{\mathbb}[1]{\amsbb{#1}}

\numberwithin{equation}{section}

\def\namedlabel#1#2{\begingroup
    #2%
    \def\@currentlabel{#2}%
    \phantomsection\label{#1}\endgroup
}
\makeatletter
\newcommand{\leqnomode}{\tagsleft@true\let\veqno\@@leqno}
\makeatother

\makeatletter
\newcommand{\customlabel}[2]{%
\protected@write\@auxout{}{\string \newlabel{#1}{{#2}{\thepage}{#2}{#1}{}}}%
\hypertarget{#1}{#2}
}
\makeatother

\newcommand\OCP[1]{\hyperref[OCP]{$\textbf{OCP}_{#1}$}\xspace}
\newcommand\SOCP[1]{\hyperref[SOCP]{$\textbf{SOCP}_{#1}$}\xspace}
\newcommand\BVP[1]{\hyperref[BVP]{$\textbf{BVP}_{#1}$}\xspace}
\newcommand\PMP{\hyperref[PMP]{$\textbf{PMP}$}\xspace}

\usepackage{tabularx}
\newcolumntype{C}{>{\centering\arraybackslash}p{19mm}}
\newcolumntype{G}{>{\centering\arraybackslash}p{4mm}}
\newcolumntype{S}{>{\centering\arraybackslash\scriptsize}p{4mm}}

\usepackage[font=small,labelfont=bf]{caption}
\usepackage{color}

\usepackage{cleveref} %

\usepackage{xspace}

\usepackage{enumitem} %

\usepackage{cite} %

\usepackage{amssymb}%
\usepackage{algorithm}
\makeatletter
\renewcommand*{\ALG@name}{Alg.}
\makeatother
\usepackage[noend]{algpseudocode}

\usepackage{xpatch} %

\newcounter{myexamplecounter}[myexamplecounter]
\newcommand\blue[1]{\textcolor{blue}{#1}} 

\newcommand\rev[1]{#1} 
\newcommand\mydots{\hbox to 1em{.\hss.\hss.}}

\newcommand{\ocp}{\hyperref[OCP]{\textbf{OCP}}\xspace}

\newcommand{\pmp}{\hyperref[PMP]{\textbf{PMP}}\xspace}

\newcommand{\olocp}{\hyperref[OLOCP]{\textbf{OL-OCP}}\xspace}
\newcommand{\fbocp}{\hyperref[FBOCP]{\textbf{FB-OCP}}\xspace}

\newcommand{\dd}{\textnormal{\textrm{d}}}

\newcommand{\Int}{\textrm{Int}}

\newcommand{\M}{\mathcal{M}}

\newcommand{\Prob}{\mathbb{P}}

\newcommand{\nablax}{\nabla_{\hspace{-1pt}x}}
\newcommand\nablaof[1]{\nabla_{\hspace{-1pt}#1}}

\newcommand{\B}{\mathcal{B}}
\newcommand{\C}{\mathcal{C}}
\newcommand{\sC}{\mathscr{C}}
\newcommand{\sD}{\mathscr{D}}

\newcommand{\bX}{\mathbb{X}}

\newcommand{\bZ}{\mathbb{Z}}
\newcommand{\mbX}{\mathbf{X}}

\newcommand{\mbZ}{\mathbf{Z}}
\newcommand{\bB}{\mathbb{B}}
\newcommand{\mbB}{\mathbf{B}}
\renewcommand{\L}{\mathcal{L}}

\newcommand{\cT}{\mathcal{T}}
\newcommand{\mbT}{\boldsymbol{\mathcal{T}}}

\newcommand{\bJ}{\mathbb{J}}
\newcommand{\mbJ}{\boldsymbol{J}}
\usepackage{stackengine,xcolor}
\input pdf-trans
\newbox\qbox
\def\usecolor#1{\csname\string\color@#1\endcsname\space}
\newcommand\bordercolor[1]{\colsplit{1}{#1}}
\newcommand\fillcolor[1]{\colsplit{0}{#1}}
\newcommand\colsplit[2]{\colorlet{tmpcolor}{#2}\edef\tmp{\usecolor{tmpcolor}}%
  \def\tmpB{}\expandafter\colsplithelp\tmp\relax%
  \ifnum0=#1\relax\edef\fillcol{\tmpB}\else\edef\bordercol{\tmpC}\fi}
\def\colsplithelp#1#2 #3\relax{%
  \edef\tmpB{\tmpB#1#2 }%
  \ifnum `#1>`9\relax\def\tmpC{#3}\else\colsplithelp#3\relax\fi
}
\newcommand\outline[1]{\leavevmode%
  \def\maltext{#1}%
  \setbox\qbox=\hbox{\maltext}%
  \boxgs{Q q 2 Tr \thickness\space w \fillcol\space \bordercol\space}{}%
  \copy\qbox%
}
\newcommand\mathcalbb[2][1]{\outline{$\mathcal{#2}$}}
\bordercolor{black}
\fillcolor{white}
\def\thickness{.1}%
\newcommand{\bT}{\mathcalbb{T}}

\newcommand{\E}{\mathbb{E}}
\newcommand{\R}{\mathbb{R}}
\newcommand{\N}{\mathbb{N}}

\newcommand{\F}{\mathcal{F}}

\NewDocumentCommand{\opnorm}{sO{}m}{%
  \IfBooleanTF{#1}{%
    \left|\opnormkern\left|\opnormkern\left|
    #3
    \right|\opnormkern\right|\opnormkern\right|
  }{
    \mathopen{#2|\opnormkern #2|\opnormkern #2|}
    #3
    \mathclose{#2|\opnormkern #2|\opnormkern #2|}
  }%
}
\newcommand{\opnormkern}{\mkern-1.5mu\relax}%

\usepackage{multirow}

\newcommand\SmallMatrix[1]{{%
  \tiny\arraycolsep=0.5\arraycolsep\ensuremath{\begin{bmatrix}#1\end{bmatrix}}}}

\usepackage{xfrac}
\usepackage{booktabs}
\newtheorem{proposition}{Proposition}[section]

\newtheorem{definition}[proposition]{Definition}
\newtheorem{assumption}[proposition]{Assumption}
\newtheorem{lemma}[proposition]{Lemma}
\newtheorem{theorem}[proposition]{Theorem}
\newtheorem{corollary}[proposition]{Corollary}
\newtheorem{remark}[proposition]{Remark}
\title{\LARGE \bf Rough Stochastic Pontryagin Maximum Principle\\and an Indirect Shooting Method} 
\author{Thomas Lew\footnote{Toyota Research Institute. Email: thomas.lew@tri.global}}

\else

\documentclass[review,hidelinks,onefignum,onetabnum]{siamart250211}
\input{shared}
\usepackage{xfrac}
\usepackage{booktabs}
\headers{Rough Stochastic Pontryagin Maximum Principle}{T. Lew}
\title{Rough Stochastic Pontryagin Maximum Principle\\and an Indirect Shooting Method\thanks{Submitted to the editors on March 15, 2025. Revision submitted on October 3, 2025.}
}
\author{Thomas Lew\thanks{Toyota Research Institute, Los Altos, CA 
  (\email{thomas.lew@tri.global}).}}

\fi

\ifarxiv
\else
\fi

\usepackage{xcolor}
\newif\ifshowrevisions

\showrevisionsfalse

\newcommand{\revproof}[1]{%
  \ifshowrevisions
    {\color{blue}#1}
  \else
    #1%
  \fi
} 
\newcommand{\revproofnote}[1]{%
  \ifshowrevisions
    {\textit{\color{blue!70!cyan}\textit{Remark}:  #1}}%
  \fi
}

\begin{document}

\maketitle

\begin{abstract}
We derive first-order Pontryagin optimality conditions for stochastic optimal control with deterministic control inputs for systems modeled by rough differential equations (RDE) driven by Gaussian rough paths. 
This   Pontryagin Maximum Principle (PMP) 
applies to systems following stochastic differential equations (SDE) 
driven by Brownian motion, yet %
it does not rely on forward-backward SDEs and involves the same Hamiltonian as the deterministic PMP. 
The proof consists of first deriving various  integrable error bounds for solutions to nonlinear and linear RDEs by leveraging recent results on Gaussian rough paths. 
The PMP then   
follows using standard techniques based on needle-like variations. %
As an application, we propose the first indirect shooting method for nonlinear stochastic optimal control. Numerical experiments  on a  stabilization problem show that it converges $10\times$ faster than a direct method.
\end{abstract}

\revproofnote{Revisions are listed below.
\begin{itemize}[itemsep=1pt]
\item The maximality condition is $u_t
\in\mathop{\arg\max}_{v\in U}\,
\E\left[
H(t,x_t,v,p_t,\mathfrak{p}_0)
\right]$, not $u_t
=\mathop{\arg\max}_{v\in U}\,
(\dots)$.
\item Section \ref{sec:linearization}: Needle-like variations are defined using a countable dense subset of controls $u\in U$. Error bounds now state an $o(\eta)$ bound instead of one in $\eta^2$.
\item Appendix: Clarified the proof of Corollary \ref{cor:needle_like_error:etas}. Fixed a sign error in the Milstein scheme for the adjoint equation in Appendix \ref{apdx:proofs:example}.
\end{itemize}
\vspace{-5mm}
}

\ifarxiv

\setcounter{tocdepth}{2}
\tableofcontents
 \newpage

\fi

\section{Introduction and main results}\label{sec:introduction} 
Stochastic optimal control has found numerous applications  %
such as in finance \cite{Guasoni2006}, aerospace \cite{Leparoux2024}, robotics \cite{Blackmore2011}, automotive \cite{LewMPC2024}, and biology \cite{Berret2020}. %
Stochastic optimal control problems 
typically %
involve %
a  dynamical system described by a %
stochastic differential equation (SDE)
\begin{align}
\label{eq:SDE} 
\dd x_t = b(t,x_t,u_t)\dd t + \sigma(t,x_t)\circ\dd B_t, \quad t\in[0,T],
\end{align}
 in Stratonovich or It\^o form, 
where $x_t$ is the state of the system at time $t$, 
$u_t$ is the control input, 
$b$ is the drift, $\sigma$ is the diffusion, $B$ is a Brownian motion, $T$ is the final time, 
and consist of optimizing an objective $\E[\int_0^T f(t,x_t,u_t)\dd t+g(x_T)]$  over a set of control input trajectories subject to state and control constraints.

By now, a rich literature on stochastic optimal control is available, with optimality conditions characterized by the dynamic programming principle as Hamilton-Jacobi-Bellman (HJB) partial differential equations (PDEs) \cite{Lions1983,Peng1992,Yong1999}, and by the Pontryagin Maximum Principle (PMP) as forward-backward stochastic differential equations (FBSDEs)  \cite{Peng1990,Yong1999,Frankowska2018,Bonalli2023}. 
For  %
problems with linear dynamics and linear-quadratic costs, both approaches lead to tractable  solutions characterized by stochastic Riccati equations \cite{Bismut1976,Peng1992,Tang2003}. 
However, for general nonlinear problems, solving HJB-PDEs or FBSDEs remains computationally challenging for high-dimensional state spaces, despite recent progress  \cite{Kushner2001,Gobet2016,BonalliLewESAIM2022,E2017}.  
In practice, %
an effective approach consists of optimizing over a class of solutions $u_t^\theta$ %
parameterized by finitely-many parameters $\theta\in\R^k$ \cite{Gobet2005,Massaroli2022} (see \cite{LiSDEs2020,Kidger2022} for machine learning applications). However, restricting %
solutions to a finite-dimensional space may obscure the structure of solutions and lead to suboptimality. 
For example, in deterministic optimal control, the PMP can provide closed-form expressions for optimal controls, such as bang-bang controls \cite{Leparoux2022}, that drastically reduce the search space %
and guide the design of indirect methods \cite{Trelat2012,Bonalli2018} for efficient and accurate numerical resolution. 
Thus, in this work, our main motivation  is to derive %
a stochastic PMP that is as close as possible to the deterministic PMP. In particular, we seek optimality conditions that are interpreted pathwise and do not rely on FBSDEs, %
to guide  the  future development of efficient algorithms inspired by deterministic optimal control techniques.

Rough path theory \cite{Lyons2002,Lyons2007,Friz2010,Friz2020,Allan2021} provides a deterministic framework of  \textit{pathwise} integration against irregular signals such as sample paths of Brownian motion and has been recognized as a robust tool for stochastic calculus. 
Pathwise stochastic optimal control \cite{HADavis1992,Rogers2007,Bhauryal2024} has been studied using rough path theory in \cite{Diehl2016,Allan2020}, but this formulation results in anticipative controls. %
In contrast, we focus on optimizing objectives averaged over random realizations of the driving rough path under deterministic open-loop controls\rev{, which include feedback controls with a fixed parameterization (see Remark \ref{remark:feedback})}. 
 However, using rough path theory in this classical setting presents two main challenges. 
First, we cannot directly use the techniques for proving the PMP from \cite{Diehl2016} as  state constraints are not accounted for.  
Second, the pathwise error bounds for solutions to \rev{rough differential equations (RDE)}  in \cite{Diehl2016,Allan2020}  %
are too crude to be integrable in general, %
so they  cannot be used in our setting that requires error bounds with integrable constants.  
To address the first challenge, we adapt 
standard arguments for proving the deterministic PMP via needle-like variations and Brouwer's fixed point theorem \cite{LeeMarkus1967,Agrachev2004,Bonnard2005,BonalliLewESAIM2022}. %
To address the second challenge, we derive finer error bounds %
by leveraging greedy partitions and favorable integrability properties of Gaussian rough paths \cite{Cass2013,Friz2013}.

\textbf{Problem setting}. We %
 consider 
stochastic optimal control problems (\ocp) with deterministic controls%
\begin{equation}\label{OCP}%
\tag{\ocp}
\begin{cases}
\min\limits_{u\in L^\infty([0,T],U)} \qquad 
	&\E\left[
\int_0^T f(t,x_t,u_t)\dd t+g(x_T)\right]
\\
\ \,  \textrm{\rev{subject to}} \quad\ \  
&
x_t = x_0+\int_0^t b(s,x_s,u_s)\dd s +\int_0^t  \sigma(s,x_s)\dd\mbB_s,
\quad t\in[0,T],
\\[1mm]
&\E\left[h(x_T)\right]=0,
\end{cases}
\end{equation}
where the  differential 
equation $\dd x_t = b(t,x_t,u_t)\dd t +\sigma(t,x_t)\dd\mbB_t$  is a random rough differential equation (RDE) %
defined pathwise  \cite{Lyons2002,Lyons2007,Friz2010,Friz2020,Allan2021} as in Theorem \ref{thm:rdes:integrable}, and the sample paths of the stochastic process $\mbB$ are Gaussian rough paths  \cite{Cass2013,Friz2013,Friz2020} as in Theorem \ref{thm:gaussian_rough_paths}, see Assumption \ref{assum:pmp} for definitions. 
In particular, if $\mbB=(B,\bB)$ is the Stratonovich lift of a Brownian motion $B$, then the RDE in \ocp %
is equivalent to the Stratonovich SDE in \eqref{eq:SDE}  \cite[Theorem 9.1]{Friz2020} and  \ocp is a classical optimal control problem.

\textbf{Main contributions}. 
First, we prove the well-posedness and various regularity and integrability results for solutions to (random) RDEs  under regularity assumptions on $(b,\sigma,\mbB)$, see for example Proposition \ref{prop:rdes:error_bound:entire_interval} and Theorem \ref{thm:rdes:integrable}. Such RDEs are slightly outside the scope of classical results in rough path theory due to the presence of the  control input $u$ which is potentially irregular, and previous error bounds in \cite{Diehl2016,Allan2020} that are derived pathwise for a deterministic rough path $\mbB$ are not integrable in general, even if $\mbB$ is the Stratonovich lift of a Brownian motion. %

Second, we derive the first-order necessary optimality conditions for \ocp stated below. 
Notations are described in Section \ref{sec:preliminaries}. The assumptions are stated in Assumption \ref{assum:pmp} in Section \ref{sec:pmp_proof}.
\begin{theorem}[Rough Stochastic Pontryagin Maximum Principle (\PMP)]\label{thm:pmp}\label{PMP} 
Define $T,\ell,U,\Omega,\Prob,b,\sigma,f$, $g,h,x_0$, the  enhanced Gaussian process $\mbB$ %
as in  Assumption \ref{assum:pmp}, %
and 
 the Hamiltonian 
\begin{align}\label{eq:hamiltonian}
H:[0,T]\times\R^n\times U\times\R^n\times\R\to\R,\  (t,x,u,p,\mathfrak{p}_0)\mapsto p^\top b(t,x,u)+\mathfrak{p}_0f(t,x,u).
\end{align}
Let $(x,u)\in L^\ell(\Omega,C([0,T],\R^n))\times L^\infty([0,T],U)$ be an optimal solution to \ocp, where  $x$ solves the random RDE in \ocp  in the pathwise sense of Theorem \ref{thm:rdes:integrable}. 

Then, there \rev{exist} a stochastic process $p\in L^\ell(\Omega,C([0,T],\R^n))$, called adjoint vector, 
and 
non-trivial Lagrange multipliers $(\mathfrak{p}_0,\dots,\mathfrak{p}_r)\in\{-1,0\}\times \R^r$ %
such that:
\begin{enumerate}[label=(\roman*)]
\item \textbf{Adjoint equation}: for some initial conditions $p_0\in L^\ell(\Omega,\R^n)$, the adjoint vector solves the random RDE %
\footnote{Each component of the  term $\rev{\nablaof{x}\sigma}(t,x)^\top p\in\R^{n\times d}$ in  the adjoint equation \eqref{eq:spmp:pmp_equations} is %
$\big[\rev{\nablaof{x}\sigma}(t,x)^\top p\big]^{ij}
=
\sum_{k=1}^n 
\frac{\partial\sigma^{kj}}{\partial x^i}(t,x) [p]^k$.}
\begin{align}\label{eq:spmp:pmp_equations}
p_t
=
p_0
-\int_0^t
\rev{\nablaof{x}H}
(s,x_s,u_s,p_s,\mathfrak{p}_0)\dd s
-
\int_0^t
\rev{\nablaof{x}\sigma}(s,x_s)^\top p_s\dd\mbB_s,
\end{align}
\rev{for almost every $t\in[0,T]$,} 
in the pathwise sense of Theorem \ref{thm:rdes:integrable}. 
\item \textbf{Transversality condition}: almost surely, the final value of the adjoint vector satisfies
\begin{gather}
\label{eq:spmp:transversality_condition:pT} 
p_T
=
\mathfrak{p}_0
\rev{\nabla g}(x_T)
+
\sum_{i=1}^r
\mathfrak{p}_i
\rev{\nabla h_i}(x_T).
\end{gather}
\item \textbf{Maximality condition}: for almost every $t\in[0,T]$, the optimal control  satisfies
\begin{align}\label{eq:spmp:maximizality_condition}
u_t
\revproof{\,\in\,}
\mathop{\arg\max}_{v\in U}\,
\E\left[
H(t,x_t,v,p_t,\mathfrak{p}_0)
\right].
\end{align} 
\end{enumerate}
\end{theorem}
The similarity to the deterministic setting is striking. The Hamiltonian $H$ remains unchanged, the adjoint equation is interpreted pathwise and has the same drift term $-\rev{\nablaof{x}H}$, 
the transversality condition is identical (pathwise almost surely), and 
the only difference in the maximality condition 
is an expected value. 
The adjoint equation \eqref{eq:spmp:pmp_equations} is \textit{not} an FBSDE, which is the key to unlock a practical indirect shooting method.

\begin{figure}[t]
\centering
\rotatebox{90}{\quad\quad\hspace{0pt}
\fbocp
\quad\hspace{-7pt}
\olocp}
\hspace{1mm}
\includegraphics[width=0.95\textwidth,trim={5mm 6mm 5mm 5mm},clip]{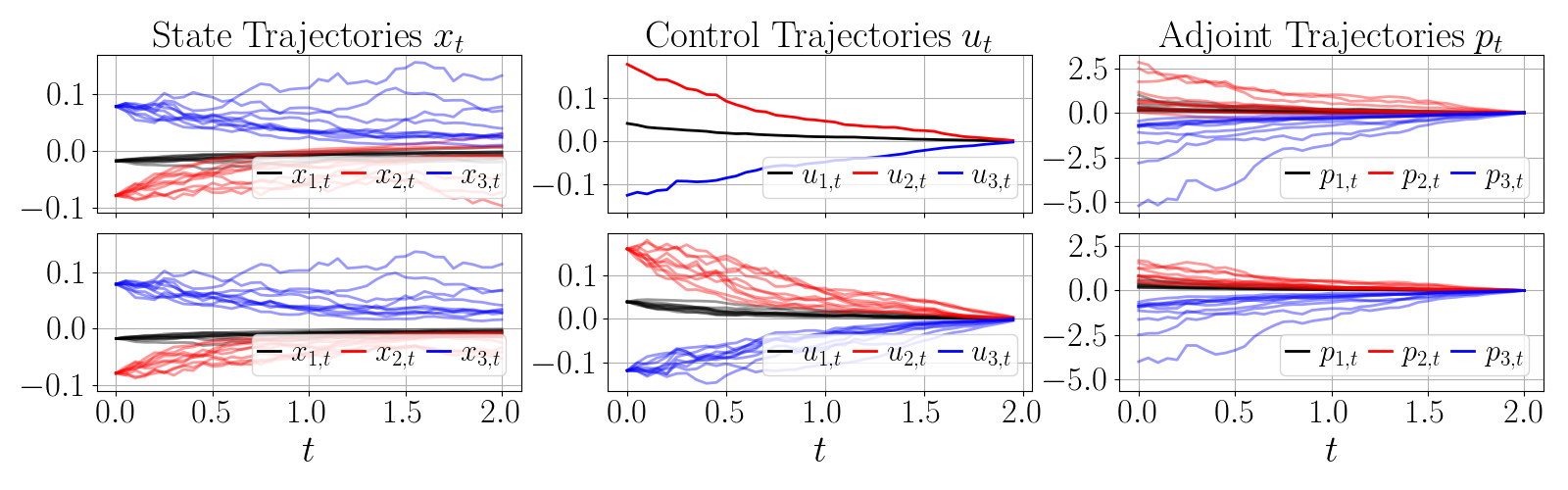}
\vspace{-1mm}
\caption{Solutions to the open-loop (\olocp) and feedback (\fbocp) optimal control problems in Section \ref{sec:example} computed using the \texttt{Indirect} method. For \fbocp, we plot the closed-loop control trajectories $u_t^i=K_tx_t^i$.}
\vspace{-2mm}
\label{fig:openloop_feedback}
\end{figure}

\newpage
\textbf{Indirect shooting method}. These optimality conditions %
inform the design of %
an indirect  method for nonlinear stochastic optimal control. That is, %
if we approximate all expectations in \ocp and \pmp using sample average (Monte Carlo) estimates for a sample size $M$, the search for (approximate) solutions to \ocp amounts to finding multipliers $(\mathfrak{p}_j)_{j=1}^r\in\R^r$ and initial values of the adjoint vector  $(p_0^i)_{i=1}^M\in\R^{Mn}$ satisfying
\begin{align*}
\begin{bmatrix}
\frac{1}{M}\sum_{i=1}^Mh(x_T^i)
\\[1mm]
\mathfrak{p}_0
\rev{\nabla g}(x_T^1)
\,{+}
\sum_{j=1}^r
\mathfrak{p}_j\hspace{-2pt}
\rev{\nabla h_j}(x_T^1)
\\
\vdots
\\
\mathfrak{p}_0
\rev{\nabla g}(x_T^M)
\,{+}
\sum_{j=1}^r\hspace{-2pt}
\mathfrak{p}_j
\rev{\nabla h_j}(x_T^M)
\end{bmatrix}
\hspace{-2pt}{=}\hspace{-2pt}
\begin{bmatrix}
0
\\[2pt]
p_T^1
\\[2pt]
\vdots
\\[2pt]
p_T^M
\end{bmatrix}
\text{where}
\begin{cases}
x_T^i = x_0^i+
\int_0^T\hspace{-1pt}
b(t,x_t^i,u_t^M)\dd t +
\int_0^T\hspace{-1pt}
\sigma(t,x_t^i)\dd\mbB_t^i
\\[2mm]
p_T^i
=
p_0^i
{-}\hspace{-1pt}
\int_0^T\hspace{-2pt}
\rev{\nablaof{x}H}
(t,x_t^i,u_t^M,p_t^i,\mathfrak{p}_0)\dd t
{-}\hspace{-1pt}
\int_0^T\hspace{-2pt}
\rev{\nablaof{x}\sigma}(t,x_t^i)^\top\hspace{-1pt} p_t^i\dd\mbB_t^i
\\[1mm]
u_t^M=
\mathop{\arg\max}\limits_{v\in U}
\frac{1}{M}\sum\limits_{i=1}^M
H(t,x_t^i,v,p_t^i,\mathfrak{p}_0)
\text{ for a.e. }
t\in[0,T].
\end{cases}
\hspace{-4mm}
\end{align*} 
Assuming that the maximality condition gives a closed-form expression of the control $u_t^M$ as a function of $(x_t^i,p_t^i)_{i=1}^M$ (e.g., as is often the case for control-affine systems), the final values $(x_T^i,p_T^i)_{i=1}^M$ can be computed as a function of $(p_0^i)_{i=1}^M$ by integrating the corresponding RDEs pathwise. Then, solutions $\big((\mathfrak{p}_j)_{j=1}^r,(p_0^i)_{i=1}^M\big)\in\R^{r+Mn}$ to this system of $(r+Mn)$ equations can be efficiently found via  a root-finding Newton method.

This approach is %
known as an indirect shooting method in deterministic optimal control \cite{Trelat2012}, and is a natural extension to the stochastic setting using \pmp and a sample average approximation  \cite{Phelps2016,Shapiro2021,Lew2024,Melnikov2024}. 
To our knowledge, this method has not appeared in the literature yet, %
as previous optimality conditions  rely on FBSDEs that introduce greater complexity. 
Indeed, previous indirect methods use deep learning to solve the FBSDEs from the classical stochastic PMP \cite{Pereira2019,Carmona2021,E2017}, whereas this indirect method does not require training a
neural network and uses a Newton method instead.  
 While we do not derive guarantees for this method (such as asymptotic optimality) %
 and thus treat it as a heuristic informed by \pmp, we evaluate it on an example in Section \ref{sec:example} and show that it is substantially faster ($10\times$ speedup) %
 than  a direct method \cite{Lew2024} solving the sample average approximation of \ocp via sequential quadratic programming. 

Figure \ref{fig:openloop_feedback} shows the solutions returned by this indirect method for two problems without final cost and state constraints ($g=h=0$), see Section \ref{sec:example}. The sample paths of the adjoint vector start from different initial conditions $p^i_0$ and are all zero at the final time ($p^i_T=0$) to satisfy the transversality condition \eqref{eq:spmp:transversality_condition:pT}.

\textbf{Connections to classical stochastic optimal control}. 
If $B$ is a Brownian motion with filtration $(\F_t)_{t\in[0,T]}$, the initial conditions $x_0$ are $\F_0$-measurable, and  $\mbB$ is the Stratonovich lift of $B$, then \ocp is equivalent to the classical stochastic optimal control problem
\begin{align*}
\min\limits_{u\in L^\infty([0,T],U)} \,
\E\bigg[
\int_0^T \hspace{-1mm}
f(t,x_t,u_t)\dd t+g(x_T)
\bigg]
\  \textrm{s.t.} \   \,
x_t = x_0+\int_0^t b(s,x_s,u_s)\dd s +\int_0^t  \sigma(s,x_s){\circ}\dd B_s,
\ 
\E\left[h(x_T)\right]=0
\end{align*}
with  a Stratonovich SDE. %
Thus, Theorem \ref{thm:pmp} applies to this standard setting %
with the following observations.  
\begin{itemize}
\item \textit{Pathwise adjoint equation}:  The  adjoint equation \eqref{eq:spmp:pmp_equations} is still understood pathwise in the rough path sense, as the adjoint vector $p$ is defined pathwise with initial conditions $p_0$ that depend on the entire path of $\mbB$ (i.e., $p_0$ is not $\F_0$-measurable, so \eqref{eq:spmp:pmp_equations} cannot be replaced by a Stratonovich SDE), 
see Section \ref{sec:pmp_proof}. %
Rough path theory is a natural framework to make sense of this equation. In contrast,  the classical stochastic PMP \cite{Peng1990,Yong1999,Bonalli2023} relies on It\^o calculus and FBSDEs.  %
\item 
\textit{Gaussian rough paths}: 
Our results hold for a large class of Gaussian rough paths (see Assumption \ref{assum:Gaussian_lift} and the many examples in \cite{Friz2016}). %
For example, fractional Brownian motion (fBM) satisfies our assumptions \cite{Bayer2016} and has  applications in finance \cite{Guasoni2006}, yet fBM (except for Brownian motion) is not a semimartingale \cite{Bayer2016} and thus cannot be handled via It\^o (or Stratonovich) integration.  
\item \textit{Regularity of the diffusion $\sigma$}: %
Our results rely on stronger regularity assumptions on $\sigma$ than the classical PMP. %
Informally, 
rough path theory cannot distinguish between It\^o and Stratonovich integration, and an It\^o SDE can be written as a Stratonovich SDE with a correction term involving derivatives of $\sigma$, so we expect smoothness assumptions on $\sigma$ to be stronger than if using It\^o calculus \cite[Section 4.1]{Perkowski2016}. On the other hand, these assumptions unlock stronger pathwise regularity results with respect to the driving process $\mbB$ (for example, Proposition \ref{prop:rdes:error_bound:entire_interval}) that may be of independent interest. %

\item \textit{Independence of the diffusion $\sigma$ on the control $u$}: Our results rely on the independence of $\sigma$ on  $u$, as allowing a dependence on the control $u$ may lead to a degenerate formulation  with irregular controls as described in \cite{Diehl2016,Allan2020}, and rough path theory relies on coefficients that are smooth-enough in time. 
\end{itemize}

\begin{remark}[Open-loop and feedback control]\label{remark:feedback}
\rev{
We optimize over deterministic open-loop controls $u$ as in \cite{Phelps2016,Lew2024,Berret2020,BonalliLewESAIM2022,LewMPC2024,Blackmore2011,Stannat2021}, which ensures that solutions are non-anticipative. %
\pmp also informs the search of  closed-loop controls with fixed parameterization, see Section \ref{sec:example} for a feedback optimization example. 
Solving such stochastic optimal control problems remains computationally challenging today. %
In practice, a feedback controller is often pre-specified and  open-loop controls are recomputed in real-time %
via model predictive control \cite{Houska2018,LewMPC2024}, which yields feedback through online replanning and is often  effective. 
However, considering open-loop controls or parameterized feedback controls is more restrictive than optimizing over %
stochastic non-anticipative controls, which remains challenging in nonlinear settings. We leave this generalization to future work, e.g., via pathwise approaches \cite{Diehl2016,Allan2020} or suitable non-anticipative parameterizations \cite{Leparoux2024,Gobet2005,Massaroli2022,Bank2024}. }
\end{remark}

\textbf{Sketch of proof and paper outline}. We adapt the proof of the deterministic PMP based on needle-like variations   \cite{LeeMarkus1967,Agrachev2004,Bonnard2005,BonalliLewESAIM2022} to the stochastic setting by leveraging recent results in rough path theory \cite{Cass2013,Friz2013} to enable a pathwise analysis. %
The main steps of the proof of \pmp and the paper outline are below.
\begin{itemize}
\item In Section \ref{sec:preliminaries}, we review concepts in rough path theory  and Gaussian rough paths. %
Given a rough path $\mbX$, we will partition the interval $[0,T]$ with a greedy partition \cite{Cass2013,Friz2013} of $N_{\alpha}(\mbX)$ increments whose size is a function of the $p$-variation of the rough path $\mbX$. Importantly, for Gaussian rough paths $\mbX=\mbB(\omega)$, the number of increments $N_{\alpha,[0,T]}(\mbB)$ enjoys favorable integrability properties. %
\item In Section \ref{sec:rdes}, we first prove existence and unicity of solutions to nonlinear and linear RDEs, used to describe the evolution of the dynamical system in \ocp and the adjoint equation \eqref{eq:spmp:pmp_equations}. 
Such results are slightly outside the scope of classical results in rough path theory due to the presence of the control $u$ (results in \cite{Diehl2016,Allan2020} do not directly apply since $(b,\sigma)$ are time-varying and the map $t\mapsto b(t,\cdot,u_t)$ is not smooth enough  to append time $t$ to the state $x$ and use previous  results). 
Then, we derive error bounds for RDEs using  greedy partitions and the quantity $N_{\alpha}(\mbB)$, %
and 
 obtain the integrability of solutions to random RDEs driven by Gaussian rough paths and of their Jacobian (Theorem \ref{thm:rdes:integrable}), and integrable error bounds for  solutions to RDEs with different control inputs  and initial conditions (Corollary \ref{cor:rdes:integrable:error_bound}).

\item In Section \ref{sec:pmp}, 
we state the assumptions for \pmp (Assumption \ref{assum:pmp}) and prove the result. The proof uses a standard technique based on needle-like variations and a separation hyperplane argument using Brouwer's fixed point theorem \cite{LeeMarkus1967,Agrachev2004,Bonnard2005}. The main differences with classical proofs of the PMP of It\^o type \cite{BonalliLewESAIM2022} are the pathwise use of It\^o's Lemma for rough paths (Lemma \ref{lem:rough_path:ito_formula}) and defining the adjoint vector pathwise using a rough  differential equation instead of  an FBSDE. 
\end{itemize}
We implement the indirect shooting method  on an example in Section \ref{sec:example} and conclude in Section \ref{sec:conclusion}.  
\ifarxiv
Additional
proofs are provided in the appendix. 
\fi

\section{Preliminaries}\label{sec:preliminaries}
We use the following standard notations  \cite{Bugini2024}. 
Given $a,b\in\R\cup\{\infty\}$, $a\wedge b:=\min(a,b)$. 
Given $a,b\in\R^n$, we write $a^\top b=\sum_{i=1}^na_ib_i$ for the inner product.  
The Kronecker product is denoted by $\otimes$. 
Let   $E,\tilde{E},F$  be  Banach spaces (usually $E=\R^n$).  
The space of linear and continuous functions from $E$ to $F$ is denoted by $\L(E,F)$ and is 
endowed with the norm $\|A\|=\sup_{x\in E,\|x\|\leq 1}\|A(x)\|$ such that $\|Ax\|\leq\|A\|\|x\|$ for any $A\in\L(E,F)$ and $x\in E$. We use the usual identifications $\L(\R^n,\R^m)\cong\R^n\otimes\R^m\cong\R^{n\times m}$ where $\R^{n\times m}$ is the matrix space, 
and 
$\L(E\otimes\tilde{E},F)\cong\L(E,\L(\tilde{E},F))\cong\L^2(E\times\tilde{E},F)$ is the space of bilinear and continuous maps from $E\times\tilde{E}$ to $F$. 
Given a function $f:E\to F$, we write $\|f\|_\infty:=\sup_{x\in E}\|f(x)\|$, and say that $f$ is bounded if $\|f\|_\infty<\infty$. 
Given $g:E\to\R$, we write $f(x)=o(g(x))$ if $\|f(x)\|/\|g(x)\|\to 0$ as $x\to 0$. 
A function $f:E\to F$ is said to be continuously differentiable  (in Fr\'echet sense) if there exists
a continuous map $\nabla f: E\to \L(E, F)$ such that $f(y) - f(x) - \nabla f(x)(y - x) = o(\|y - x\|)$. Partial derivatives \rev{$\nablaof{x}f$} \rev{(also denoted by $\frac{\partial f}{\partial x}$ if $f$ is scalar valued)} and higher order derivatives $\nabla^k f$ are defined as usual.  
Given $n\in\N$, %
$C^n=C^n(E,F)$  denotes the space of continuous functions $f:E\to F$ that are $n$-times continuously differentiable, 
and $C_b^n=C_b^n(E,F)$ denotes 
the space of bounded functions $f\in C^n(E,F)$   with bounded derivatives. 
The space $C_b^n$ is 
endowed with the norm
$\|f\|_{C_b^n}:=\|f\|_\infty+\|\nabla f\|_\infty+\dots+\|\nabla^n f\|_\infty<\infty$. 

Let $T>0$. For any interval $I=[s,t]\subseteq[0,T]$, we write $|I|:=|t-s|$ \rev{and $I^2:=I\times I$}. Given a path $X:[0,T]\to\R^n$, its increments are denoted by $X_{s,t}:=X_t-X_s$ for any $s,t\in[0,T]$. 
Let $C([0,T],\R^n)$ be the set of continuous maps $x:[0,T]\to\R^n$ and $(\Omega,\F,\Prob)$ be a probability space. %
For $\ell\in\N$, we denote by 
$L^\ell(\Omega,\R^n)$ the set of random variables $x:\Omega\to\R^n$ such that $\E[\|x\|^\ell]<\infty$, and 
by $L^\ell(\Omega,C([0,T],\R^n))$ the set of stochastic processes with continuous sample paths $x:\Omega\to C([0,T],\R^n)$ such that $\E[\|x\|_\infty^\ell]<\infty$. 
Given $U\subseteq\R^n$, we denote by  $L^\infty([0,T],U)$ the set of measurable maps $u:[0,T]\to U$ with $\|u\|_{L^\infty([0,T],U)}:=\inf\{C : \|u_t\|\leq C\text{ for almost every }t\in[0,T]\}<\infty$. 
\rev{Given two maps $\sigma:[0,T]\times E\to F$ and $X:[0,T]\to E$, 
$\sigma(\cdot,X_\cdot)$ denotes the map $t\mapsto\sigma(t,X_t)$.}
Throughout derivations, we denote by $C_a$  a constant that depends only on $a$ and can change line by line.

\subsection{Rough paths, controlled rough paths, and rough integration}\label{sec:preliminaries:rough_paths}

We recall concepts in rough path theory \cite{Friz2010,Friz2020,Allan2021}. 
The main tool we use for  quantifying the regularity of a path $X:[0,T]\to\R^n$ is the notion of $p$-variation, which  bounds the sum  of increments $\|X_{s,t}\|^p$ of $X$ over arbitrary partitions of $[0,T]$. For example, Lipschitz continuous paths $X$ have finite $1$-variation, and sample paths $B_{s,t}(\omega)$ of Brownian motion have finite $p$-variation for $p>2$ almost surely. 
Rough path theory can also be studied using  $\frac{1}{p}$-H\"older continuity properties (note that a path that is $\frac{1}{p}$-H\"older continuous has finite $p$-variation), but the resulting analysis gives bounds that are generally not integrable (see Section \ref{sec:preliminaries:greedy_gaussian_paths} for further discussion), which motivates using a rough path analysis via $p$-variation properties.

\begin{definition}[$p$-variation]
Let $p\geq 1$ and $T>0$. The $p$-variation  of a path $X:[0,T]\to\R^n$  is defined as
$$
\|X\|_p
:=
\|X\|_{p,[0,T]},
\ \ \text{where }\ \,
\|X\|_{p,[s,t]}
:=
\bigg(
\sup_{\pi\in\mathcal{P}([s,t])}\sum_{[u,v]\in\pi}\|X_{u,v}\|^p
\bigg)^{\frac{1}{p}}
\  
\text{for any $[s,t]\subseteq[0,T]$},
$$
where  $\mathcal{P}([s,t])$ denotes the set of all partitions of $[s,t]$, 
and the supremum is over all partitions $\pi$ of $[s,t]$.  

The set $\C^p=\C^p([0,T],\R^n)$ denotes the space of $\R^n$-valued continuous paths of finite $p$-variation, that is, continuous paths  $X:[0,T]\to\R^n$ such that $\|X\|_p<\infty$. 
\end{definition}

  \begin{lemma}[Stitching a partition of an interval {\cite[Lemma 2.3]{Allan2020}}]\label{lem:pvar:intervals}
 Let $p\geq 1$, $T>0$,  $0=t_0<t_1<\dots<t_n=T$ be a partition of $[0,T]$, and  $X:[0,T]\to\R^n$ be a path. Then, $
\|X\|_{p,[0,T]}\leq n\big(
\sum_{i=1}^n\|X\|_{p,[t_{i-1},t_i]}^p
\big)^{1/p}$. 
\end{lemma}

\begin{lemma}[Inequalities for  $p$-variations]
\label{lem:pvariation:inequalities}
Let $T>0$, $p\geq 1$, and $X\in\C^p([0,T],\R^d)$. Then,
\begin{equation}
\label{eq:path_finite_var:infty_ineq}
\|X\|_\infty\leq
\|X_0\|+\|X\|_p.
\end{equation}
Let $X:[0,T]\to\R^d$, $Y^1,\dots,Y^n\in\C^p([0,T],\R^d)$, 
and $c\geq 0$. 
Then, there \rev{is} a constant $C_{\rev{n,}p}\geq 1$ such that
\begin{equation}
\label{eq:path_finite_var:sum_pvars}
\|X_{s,t}\|\leq \sum_{i=1}^n\|Y^i_{s,t}\|
+
c|t-s|
\ \forall s,t\in[0,T]
\implies 
\|X\|_p\leq C_{\rev{n,}p}
\left(
\sum_{i=1}^n\|Y^i\|_p+c\,T
\right).
\end{equation}
Let $p\geq 2$, $X:[0,T]\to\R^d$, $Y^i,\rev{\widetilde{Y}}^i\in\C^p$ for $i=1,\dots,n$, $Z^j\in\C^\frac{p}{2}$ for $j=1,\dots,m$. %
Then, $\|X\|_p\leq\|X\|_\frac{p}{2}$, and there exists a constant $C_{\rev{m,n,}p}\geq 1$ such that
{\small
\begin{equation}\label{eq:path_finite_var:sum_p/2vars}
\|X_{s,t}\|\,{\leq}\sum_{i=1}^n\|Y^i_{s,t}\|\|\rev{\widetilde{Y}}^i_{s,t}\|
+\sum_{j=1}^m\|Z^j_{s,t}\|
+
c|t-s|
\ \forall s,t\in[0,T]
\hspace{-2pt}
\implies 
\hspace{-2pt}
\|X\|_\frac{p}{2}\leq C_{\rev{m,n,}p}
\hspace{-1pt}
\left(\sum_{i=1}^n\|Y^i\|_p\|\rev{\widetilde{Y}}^i\|_p\,{+}\sum_{j=1}^m\|Z^j\|_\frac{p}{2}\,{+}\,c\,T
\hspace{-1pt}
\right)\hspace{-2pt}.
\end{equation}
}%
Let $p\geq 1$, $\sigma\in C^1_b([0,T]\times\R^n,\R^m)$, and $X\in\C^p([0,T],\R^n)$.  Then, there exists  $C_p\geq 1$ such that
\begin{equation}\label{eq:sigma(.,X):pvar}
\|\sigma(\cdot,X)\|_p\leq C_p\|\sigma\|_{C^1_b}(\|X\|_p+T).
\end{equation}
Moreover, if $\sigma\in C^2_b\rev{([0,T]\times\R^n,\R^m)}$ and $X,\rev{\widetilde{X}}\in\C^p\rev{([0,T],\R^n)}$, then there exists  $C_p\geq 1$ such that
\begin{equation}\label{eq:Delta_sigma(.,X):pvar}
\|\sigma(\cdot,X_\cdot)-\sigma(\cdot,\rev{\widetilde{X}}_\cdot)\|_p\leq 
C_p\|\sigma\|_{C^2_b}
(
1 + \|X\|_p+\|\rev{\widetilde{X}}\|_p
+
T
)
(\|X_0-\rev{\widetilde{X}}_0\|+\|X-\rev{\widetilde{X}}\|_p).
\end{equation}
 \end{lemma}
\ifarxiv
The results in Lemma \ref{lem:pvariation:inequalities} are standard and are proved in the appendix.
\fi

 \begin{definition}[Rough path]
 Let $p\in[2,3)$ and $T>0$. A \rev{rough} path is a pair $\mbX=(X,\bX)$ that consists of a path $X:[0,T]\to\R^d$ and its enhancement $\bX:\rev{[0,T]^2}\to\R^{d\times d}$ that satisfy Chen's relation	 
\begin{equation}\label{eq:chen's_relation}
\bX^{ij}_{s,t}=\bX^{ij}_{s,r}+\bX^{ij}_{r,t}+X^i_{s,r} X^j_{r,t}
\end{equation}
for all $1\leq i,j\leq d$ and $\rev{s,r,t\in[0,T]}$, %
and that has finite inhomogeneous $p$-variation rough path norm:
$$
\|\mbX\|_p:=\|X\|_p+\|\bX\|_{\frac{p}{2}}<\infty,
\quad\text{where}\quad
\|\bX\|_{\frac{p}{2},[s,t]}:=\bigg(\sup_{\pi\in\mathcal{P}([s,t])}\sum_{[u,v]\in\pi}
\|\bX_{u,v}\|^{\frac{p}{2}}\bigg)^{\frac{2}{p}}
\  
\text{for any $[s,t]\subseteq[0,T]$},
$$
and $\|\bX\|_{\frac{p}{2}}:=\|\bX\|_{\frac{p}{2},[0,T]}$. We also write $\|\mbX\|_{p,[s,t]}:=\|X\|_{p,[s,t]}+\|\bX\|_{\frac{p}{2},[s,t]}$ for any $[s,t]\subseteq[0,T]$. 

A geometric \rev{rough} path is a  \rev{rough} path $\mbX=(X,\bX)$ that additionally satisfies
\begin{equation}\label{eq:rough_path:integration_by_parts}
\bX^{ij}_{s,t}+\bX^{ji}_{s,t}=X^i_{s,t}X^j_{s,t}
\end{equation}
 for all $1\leq i,j\leq d$ and $\rev{s,t\in[0,T]}$. %

The sets $\sC^p=\sC^p([0,T],\R^d)$ and $\sC^p_g=\sC^p_g([0,T],\R^d)$ denote the sets of \rev{rough} paths and of geometric \rev{rough} paths, respectively. 

Given two rough paths $\mbX=(X,\bX)\in\sC^p$ and $\widetilde{\mbX}=(\rev{\widetilde{X}},\rev{\widetilde{\bX}})\in\sC^p$, we define $\Delta X:=X-\rev{\widetilde{X}}$, $\Delta\bX:=\bX-\rev{\widetilde{\bX}}$, and %
\vspace{-2mm}
\begin{align*}
\|\Delta\mbX\|_p:=\|\Delta X\|_p+\|\Delta\bX\|_{\frac{p}{2}}.
\end{align*}
 \end{definition}

Clearly, $\sC^p_g\subset\sC^p$. %
As an example, the sample paths $\mbX=\mbB(\omega)$ of the Stratonovich lift $\mbB=(B,\bB)$ of a Brownian motion $B$, where $\bB$ is defined by the Stratonovich integrals $\bB_{s,t}^{ij}:=\int_s^tB_{s,u}^i\circ\dd B_u^j$, are geometric rough paths ($\mbB(\omega)\in\sC^p_g$ almost surely). 
The sample paths  of the It\^o lift $\widetilde{\mbB}=(B,\widetilde{\bB})$ of $B$ (with $\widetilde{\bB}$ defined by It\^o integration  $\widetilde{\bB}_{s,t}^{ij}:=\int_s^tB_{s,u}^i\dd B_u^j$) are rough paths ($\widetilde{\mbB}(\omega)\in\sC^p$ almost surely), 
but they are not geometric. 
A key property of geometric rough paths is %
the chain rule (see Lemma \ref{lem:rough_path:ito_formula}), as a consequence of It\^o's lemma.

\begin{definition}[Controlled rough paths] 
 Let $p\in[2,3)$, $T>0$, and $\mbX=(X,\bX)\in\sC^p([0,T],\R^d)$ be a \rev{rough} path. A controlled \rev{rough} path (with respect to $X$) is a pair
$$
(Y,Y')\in\C^p([0,T],\R^n)\times\C^p([0,T],\R^{n\times d}),
$$
where $Y'$ is called the Gubinelli derivative of $Y$, such that the remainder term $R^Y:\rev{[0,T]^2}\to\R^n$ given by
\begin{equation}\label{eq:remainder}
R^Y_{s,t}:=Y_{s,t}-Y'_sX_{s,t}
\end{equation}
satisfies $\|R^Y\|_{\frac{p}{2}}<\infty$. 

The set $\sD^p_X=\sD^p_X([0,T],\R^n)$ denotes the set of controlled \rev{rough} paths with respect to $X$. 
\end{definition}
A controlled rough path $(Y,Y')\in\sD^p_X$ looks like $X$ over short intervals: $Y_t\approx Y_s+Y'_sX_{s,t}$ for small $|t-s|$. 
For example, $(f(X),\nabla f(X))$ is a controlled rough path if $f\in C^2$.  Controlled rough paths %
are sufficiently smooth (with respect to the rough path $\mbX$) to allow for a notion of \textit{rough integral} against $\mbX$, defined below.

\begin{proposition}[Rough integration {\cite[Proposition 2.6]{Friz2018}}]
\label{prop:rough_integral_welldefined:error_bound}
Let $p\in[2,3)$, $T>0$, $\mbX=(X,\bX)\in\sC^p([0,T],\R^d)$ be a rough path, and $(Y,Y')\in\sD^p_X([0,T],\R^{n\times d})$  be a controlled rough path. Then, the rough integral of $(Y,Y')$ against $\mbX$, defined as the limit over all partitions $\pi$ of $[0,T]$ with vanishing mesh size
\begin{equation}\label{eq:rough_int}
\int_0^TY_r\dd\mbX_r
:=
\lim_{|\pi|\to0}\sum_{[s,t]\in\pi}Y_sX_{s,t}+Y'_s\bX_{s,t},
\end{equation} 
exists\footnote{We use the identification $\R^{n\times d\times d}=\L(\R^{d \times d},\R^n)$ to make sense of the last term $Y'_s\bX_{s,t}$.}. Moreover, for any $0\leq s<t\leq T$, we have the estimate %
\begin{equation}\label{eq:rough_int:error_bound}
\left\|\int_s^tY_r\dd\mbX_r-Y_sX_{s,t}-Y'_s\bX_{s,t}\right\|\leq C_p(\|R^Y\|_{\frac{p}{2},[s,t]}\|X\|_{p,[s,t]}+\|Y'\|_{p,[s,t]}\|\bX\|_{\frac{p}{2},[s,t]})
\end{equation} 
 for a constant $C_p$ that only depends on $p$.
\end{proposition} 
For intuition, let $X$ be a scalar-valued path that is sufficiently smooth so that $\bX_{s,t}:=\int_s^t(X_r-X_s)\dd X_r$ is well-defined and $(Y,Y')=(f(X),\nabla f(X))$ with $f\in C^2$ is a controlled path. Then, a Taylor approximation gives
$\int_s^tf(X_r)\dd X_r\approx f(X_s)\int_s^t\dd X_r+\nabla f(X_s)\int_s^t(X_r-X_s)\dd X_r=Y_sX_{s,t}+Y'_s\bX_{s,t}$, which is like the left hand side of \eqref{eq:rough_int:error_bound}. 
The rough integral generalizes this intuition to cases where $X$ is too irregular for the integral $\int_s^tX_{s,r}\dd X_r$ to be well-defined (e.g., if $X$ is a sample of Brownian motion): 
We first define  the enhancement $\bX$ and then define the rough integral of $(Y,Y')$ against $(X,\bX)$  by \eqref{eq:rough_int}.
The term $Y'_s\bX_{s,t}$ is key to ensuring that  \eqref{eq:rough_int} is well-posed and different choices of $\bX$ (e.g., via It\^o or Stratonovich integrals) give different results.

Next, we give a few useful results about rough integration and functions of controlled rough paths. 
Let  $(X,\bX)\in\sC^p$ and $(\rev{\widetilde{X}},\rev{\widetilde{\bX}})\in\sC^p$ be two rough paths,  and  
$(Y,Y')\in\sD^p_X$ and $(\rev{\widetilde{Y}},\rev{\widetilde{Y}}')\in\sD^p_{\rev{\widetilde{X}}}$ be two controlled rough paths. 
We write $\Delta X=X-\rev{\widetilde{X}}$ and similarly for $\Delta Y,\Delta Y',\Delta R^Y$, and as in \cite{Friz2018}, we define 
\begin{align*}
M_{Y'}&:=\|Y_0'\|+\|Y'\|_p,
&& 
\hspace{-4mm}
K_Y:= \|Y_0'\|+\|Y'\|_p+\|R^Y\|_\frac{p}{2} =  M_{Y'}+\|R^Y\|_\frac{p}{2},
\\
\Delta M_{Y'} &:= \|\Delta Y_0'\|+\|\Delta Y'\|_p,
&&
\hspace{-4mm}
\Delta K_Y := \|\Delta Y_0'\|+\|\Delta Y'\|_p+\|\Delta R^Y\|_\frac{p}{2}=\Delta M_{Y'}+\|\Delta R^Y\|_\frac{p}{2}.
\end{align*}

\begin{lemma}[Stability of rough integration {\cite[Lemma 3.4]{Friz2018}}]
\label{lem:rough_int_stability}
Let $p\in[2,3)$, $T>0$,  
 $\mbX=(X,\bX)\in\sC^p$ and $\widetilde{\mbX}=(\rev{\widetilde{X}},\rev{\widetilde{\bX}})\in\sC^p$ be two \rev{rough} paths,  and  
$(Y,Y')\in\sD^p_X$ and $(\rev{\widetilde{Y}},\rev{\widetilde{Y}}')\in\sD^p_{\rev{\widetilde{X}}}$ be two controlled \rev{rough} paths. Then,
$
(Z,Z'):=\big(
\int_0^\cdot Y_s\dd\mbX_s,Y
\big) \in\sD^p_X$ and $(\tilde{Z},\tilde{Z}'):=\big(\int_0^\cdot \rev{\widetilde{Y}}_s\dd\widetilde{\mbX}_s,\rev{\widetilde{Y}}\big) \in\sD^p_{\rev{\widetilde{X}}}.
$
Moreover, %
\begin{align}
\label{lem:rough_int_stability:remainder}
\|R^{\int_0^\cdot Y_s\dd\mbX_s}-R^{\int_0^\cdot \rev{\widetilde{Y}}_s\dd\widetilde{\mbX}_s}\|_{\frac{p}{2}}
&\leq 
C_p
(1+\|\mbX\|_p+\|\widetilde{\mbX}\|_p)
\big(K_{\rev{\widetilde{Y}}}\|\Delta\mbX\|_p+\|\mbX\|_p\Delta K_Y\big).
\end{align}
\end{lemma}

\begin{lemma}[$(\sigma(\cdot,Y),\nabla \sigma(\cdot,Y)Y')$ is a controlled path]\label{lem:rough_path:sigma(.,Y):controlled}
Let $p\in[2,3)$, $T>0$, $\sigma\in C^2_b([0,T]\times\R^n,\R^{n\times d})$, 
$\mbX=(X,\bX)\in\sC^p([0,T],\R^d)$ be a \rev{rough} path, and 
$(Y,Y')\in\sD^p_X([0,T],\R^n)$ be a controlled \rev{rough} path. Then, 
$(\sigma(\cdot,Y_\cdot),\sigma(\cdot,Y_\cdot)')
:=
\big(\sigma(\cdot,Y_\cdot),\rev{\nablaof{x}\sigma}(\cdot,Y_\cdot)Y'_\cdot
\big)\in\sD^p_X\rev{([0,T]\times\R^n,\R^{n\times d})}$.  
\rev{Also}, there exists $C_p\geq 1$ such that 
\begin{align}
\label{eq:controlled_path:pvar_norm}
\|Y\|_p&\leq C_p(\|Y'\|_\infty\|X\|_p+\|R^Y\|_\frac{p}{2})
\\
\label{eq:|Y|_p_ineq}
&\leq C_p(1+\|X\|_p)K_Y,
\\
\label{eq:sigma(.,Y):pvar}
\|\sigma(\cdot,Y_\cdot)\|_p&\leq 
C_p\|\sigma\|_{C^1_b}
(
M_{Y'}\|X\|_p+\|R^Y\|_\frac{p}{2}
+
T
),
\\
\label{eq:sigma(.,Y)':pvar}
\|\sigma(\cdot,Y_\cdot)'\|_p&\leq 
C_p\|\sigma\|_{C^2_b}
K_Y(1+K_Y+T)(1+\|X\|_p),
\\
\label{eq:RY:p/2var:Y^2+RY+T}
\|R^{\sigma(\cdot,Y_\cdot)}\|_\frac{p}{2}&\leq
C_p\|\sigma\|_{C^2_b}(
\|Y\|_p^2
+
\|R^Y\|_\frac{p}{2}
+
T
)
\\
\label{eq:RY:p/2var:KYs}
&\leq
C_p\|\sigma\|_{C^2_b}
(
K_Y(1+
K_Y)(1+\|X\|_p)^2
+
T
).
\end{align}
\end{lemma}
Lemma \ref{lem:rough_path:sigma(.,Y):controlled} is a straightforward extension  of \cite[Lemmas 3.5 and 3.6]{Friz2018} to time-varying maps $\sigma(t,x)$.
\begin{remark}[Smoothness of $\sigma$]\label{remark:sigma_smoothness}
The assumption $\sigma\in C^2_b$ implies that  
{\small $\|
\sigma
\|_\infty
+
\|
\rev{\nablaof{t}\sigma}
\|_\infty
+
\|
\rev{\nablaof{x}\sigma}
\|_\infty
+
\|
\rev{\nablaof{t}^2\sigma}
\|_\infty
+
\|
\rev{\nablaof{x}^2\sigma}
\|_\infty
+
\|
\rev{\nablaof{x}\nablaof{t}\sigma}
\|_\infty
$} %
 is bounded. %
It is stronger than assumptions on $\sigma$  used for classical PMPs derived via It\^o calculus due to the use of rough path theory, see  Section \ref{sec:introduction} and \cite[Section 4.1]{Perkowski2016}. 
It can be relaxed to assuming that   $\sigma(t,\cdot)\in C_b^2$ for almost every (a.e.) $t\in[0,T]$ and that $\sigma(\cdot,x)$ is uniformly \smash{$\frac{2}{p}$}-H\"older continuous. Under this assumption and similar ones on derivatives $\rev{\nablaof{x}^k\sigma}$, most results in this work still hold by replacing terms in $T$ by terms in $T^\frac{2}{p}$. However, the assumption of $\frac{2}{p}$-H\"older continuity in $t$ is only barely weaker \rev{than} assuming Lipschitz continuity  (in particular, $\sigma(t,x)=B_t(\omega)$ with $X=B(\omega)$ a sample path of Brownian motion is not $\frac{2}{p}$-H\"older continuous), so we assume that $\sigma\in C_b^n$ to simplify notations and results. %
Assuming a controlled structure $\sigma_t(\cdot)\approx \sigma_s(\cdot)+\sigma'_s(\cdot)(X_t-X_s)$ \cite{Friz2024,Bugini2024}  could be considered in future work. See also Remark \ref{remark:sigma_smoothness:2}.
\end{remark}
The next result gives a chain rule for functions of  controlled rough paths and geometric rough paths.
\begin{lemma}[It\^o's formula for geometric rough paths {\cite[Theorem 7.7]{Friz2020}}]
\label{lem:rough_path:ito_formula}
Let $p\in[2,3)$, $T>0$, $f\in C^3$, 
$\mbX\in\sC^p_g([0,T],\R^d)$ be a geometric \rev{rough} path, and 
$(Y,Y')\in\sD^p_X([0,T],\R^n)$ be a controlled \rev{rough} path such that 
 $$
Y_t=Y_0+\int_0^tY'_s\dd\mbX_s+\Gamma_t$$ for all  
$t\in[0,T]$
for a path of finite $\frac{p}{2}$-variation $\Gamma\in\C^{\frac{p}{2}}([0,T],\R^n)$ and a controlled \rev{rough} path  $(Y',Y'')\in\sD^p_X([0,T],\R^{n\times d})$.  
Then, 
\begin{equation}\label{eq:rough_path:ito_formula}
f(Y_t)=f(Y_0)+\int_0^t\nabla f(Y_u)Y'_u\dd\mbX_u
+\int_0^t\nabla f(Y_u)\dd\Gamma_u
\end{equation}
for all $t\in[0,T]$, where the first integral is a rough integral, and the second integral is a Young integral.
\end{lemma}
\ifarxiv
Lemma \ref{lem:rough_path:ito_formula} is standard. It is a particular case of  \cite[Theorem 7.7]{Friz2020}, assuming that the rough path $\mbX$ is geometric so its bracket is zero. Although \cite[Theorem 7.7]{Friz2020} is formulated for $\frac{1}{p}$- and $\frac{2}{p}$-H\"older continuous paths, the proof follows similarly for   finite $p$- and $\frac{p}{2}$-variation paths.  
\fi

\rev{Lemma \ref{lem:pvar:RX_intervals} below is an extension of \cite[Lemma 2.3]{Allan2020} (Lemma \ref{lem:pvar:intervals}) for rough paths and remainders of controlled paths. Its proof follows the proof of \cite[Lemma 2.3]{Allan2020} with minor changes.}

\begin{lemma}[\rev{Bounds on $\bX$ and $R^Y$ when stitching a partition}]
\label{lem:pvar:RX_intervals}
\rev{Let $p\in[2,3)$, $T>0$, 
$0=t_0<t_1<\dots<t_n=T$ be a partition of $[0,T]$, 
$\mbX=(X,\bX)\in\sC^p([0,T],\R^d)$, and $(Y,Y')\in\sD^p_X([0,T],\R^{n\times d})$. Then, 
\begin{align}
\label{eq:pvar:bX_intervals}
\|\bX\|_{\frac{p}{2},[0,T]}
&\leq 
2n^2
\bigg(
\bigg(
\sum_{i=1}^n\|\bX\|^\frac{p}{2}_{\frac{p}{2},[t_{i-1},t_i]}
\bigg)^\frac{2}{p}
+
\bigg(
\sum_{i=1}^n\|X\|^p_{p,[t_{i-1},t_i]}
\bigg)^\frac{2}{p}
\bigg),
\\
\label{eq:pvar:RX_intervals}
\|R^Y\|_{\frac{p}{2},[0,T]}
&\leq 
2n^2
\bigg(
\bigg(
\sum_{i=1}^n\|R^Y\|^\frac{p}{2}_{\frac{p}{2},[t_{i-1},t_i]}
\bigg)^\frac{2}{p}
+
\bigg(
\sum_{i=1}^n\|Y'\|^p_{p,[t_{i-1},t_i]}
\bigg)^\frac{1}{p}
\bigg(
\sum_{i=1}^n\|X\|^p_{p,[t_{i-1},t_i]}
\bigg)^\frac{1}{p}
\bigg).
\end{align}
}
\end{lemma}

\subsection{The greedy partition and Gaussian rough paths}\label{sec:preliminaries:greedy_gaussian_paths}
To prove \PMP, we will use short variations  around the optimal solution %
 (see Proposition \ref{prop:linear_variation}), %
and show that the difference between the solutions to the corresponding rough differential equations is small \textit{in expectation}  over the driving signal $\mbX=\mbB(\omega)$ and initial conditions $x_0$. %
An immediate 
challenge with this approach is that classical error bounds in rough path theory depend on constants that are typically not integrable \cite{Bayer2016,Cass2017}, as they depend exponentially on the $p$-variation rough path norm  $\|\mbX\|_p$ (see for example \cite[Theorem 3.9]{Friz2018}) or on the exponential of the $\frac{1}{p}$-H\"older constant of $\mbX$ (see e.g. the proof of \cite[Theorem 7.9]{Allan2021} or \cite[Theorem 8.5]{Friz2020}). %
Thus, we cannot use classical error bounds from rough path theory %
in the proof of \pmp, as they may not be integrable.

A solution to this challenge was identified in \cite{Cass2013} and refined  in \cite{Friz2013}. %
It consists of %
using a quantity $N_\alpha(\mbX)$ that counts the number of $\alpha$-increments of the %
rough path norm of $\mbX$ over an interval (Definition \ref{def:Nalpha}). Importantly,  for particular \textit{Gaussian rough paths} $\mbB(\omega)$, the exponential of $N_\alpha(\mbB)$ is integrable (Theorem \ref{thm:gaussian_rough_paths}). 
Thus, in this work, we use ideas in \cite{Cass2013,Friz2013} (see also \cite{Bayer2016,Riedel2017,Cass2017,Cass2019}) and derive finer error bounds as a function of $N_\alpha(\mbB)$, which will give us integrable bounds that can be used to prove \pmp. %

\begin{definition}[Control \rev{and greedy partition}]%
\label{def:Nalpha}
Let $T>0$. A control is a continuous map $w:\rev{[0,T]^2}\to\R$  such that \rev{$w(s,t)\geq 0$} and $w(s,t)+w(t,u)\leq w(s,u)$ for any $\rev{s,t,u\in[0,T]}$. %

Given a control $w$, a resolution $\alpha>0$ and $[s,t]\subseteq[0,T]$, we  define the sequence
\begin{align*}
\tau_0&=s,
\qquad
\tau_{i+1}=\inf\{u : w(\tau_i,u)\geq \alpha, \tau_i<u\leq t\}\wedge t\ \, \text{for} \ \, i\in\N,
\end{align*}
with the convention $\inf\emptyset=+\infty$, and
$$
N_{\alpha,[s,t]}(w):=\sup \{n\in\N\cup \{0\}:\tau_n< t\}.
$$

The \textit{greedy partition} of the interval $[s,t]$ is defined as the partition
$$
\{\tau_i, i=0,1,\dots,N_{\alpha,[s,t]}(w)+1\}.
$$
Given $p\in[2,3)$ and $\mbX\in\sC^p$, %
the control  $w_\mbX$ %
and $N_{\alpha,[s,t]}(\mbX)$  are defined %
\rev{as}  $
w_\mbX(s,t)
:=
\|X\|_{p,[s,t]}^p+\|\bX\|_{\frac{p}{2},[s,t]}^\frac{p}{2}$ and $N_{\alpha,[s,t]}(\mbX):=N_{\alpha,[s,t]}(w_\mbX)$.
\end{definition} 
The \textit{control} $w$ in Definition \ref{def:Nalpha} should not be confused with the \textit{control input} $u$ in \ocp. The distinction should be clear from context: \rev{$w_\mbX$} bounds variations of $\mbX$, whereas $u$ steers the dynamical system in \ocp.

Our error bounds for solutions to RDEs will depend on the inhomogeneous  rough path norm $\|\mbX\|_p$ due to oscillations of the driving signal $\mbX$ (see \eqref{eq:rough_int:error_bound} and \eqref{lem:rough_int_stability:remainder} which contain terms in  $\|\mbX\|_p$), and on  terms increasing linearly over time  due to the drift $b(t,Y,u)$ and the time-varying diffusion $\sigma(t,Y)$  (see  Lemma \ref{lem:bounds_int_b_ds}, and \eqref{eq:sigma(.,Y):pvar} and \eqref{eq:RY:p/2var:Y^2+RY+T} which contain terms in \rev{$T$}). These terms relate 
to the control  $w_\mbX(\cdot)$ and to $N_\alpha(\mbX)$ via Corollary \ref{cor:Nalpha:NX_NXtilde_NT:small_intervals},  
using two  lemmas that are similar to \cite[Lemma 4.9]{Cass2013} and \cite[Lemma 5]{Bayer2016}, see also \cite[Lemma 3]{Friz2013}. 
The results below should be considered standard, and their proofs are %
\ifarxiv
in the appendix.
\fi  

\begin{lemma}\label{lem:Nalpha<=w(0,T)}
Let $T>0$, $w:\rev{[0,T]^2}\to\R$ be a control, and $\alpha>0$. Then, $\alpha N_{\alpha,[0,T]}(w)\leq w(0,T)$. 
\end{lemma}
\begin{lemma}\label{lem:sum_Nalpha}
Let $T>0$, $C\geq1$ be a constant, $w_1,\dots,w_n:\rev{[0,T]^2}\to\R$ be $n$ controls,  and  $w:\rev{[0,T]^2}\to\R$ be the control defined by $w(s,t)=C\sum_{j=1}^n w_j(s,t)$. 
Then, for any $\alpha>0$ and $[s,t]\subseteq[0,T]$,
\begin{equation}\label{eq:Nalpha(w)<=sumNalpha(wj)}
N_{\alpha,[s,t]}(w)
\leq
C\bigg(
2\sum_{j=1}^nN_{\alpha,[s,t]}(w_j)+n
\bigg).
\end{equation}
\end{lemma}

\begin{corollary}\label{cor:Nalpha:NX_NXtilde_NT:small_intervals}
Let $p\in[2,3)$, $T>0$, $\mbX,\widetilde{\mbX}\in\sC^p$, $C_p=6^p$, and define the control  $w:\rev{[0,T]^2}\to\R$  by 
$$
w(s,t)=C_p(w_\mbX(s,t)+w_{\widetilde{\mbX}}(s,t)+|t-s|).
$$
Then,  for any $\alpha>0$ and  $[s,t]\subseteq[0,T]$,
\begin{equation}\label{eq:Nalpha<=3CpNalpha_X_and_time}
N_{\alpha,[s,t]}(w)\leq 5C_p(
N_{\alpha,[s,t]}(\mbX)+N_{\alpha,[s,t]}(\widetilde{\mbX})+|t-s|/\alpha+1
).
\end{equation}
Moreover, given   $0<\alpha\leq 1$  and   any interval $[s,t]\subseteq[0,T]$ small-enough so that  $w(s,t)\leq \alpha$, we have
\begin{equation}\label{eq:|X|+|Xtilde|+|dt|<=alpha^p}
\|\mbX\|_{p,[s,t]}+\|\widetilde{\mbX}\|_{p,[s,t]}+|t-s|\leq
 \alpha^\frac{1}{p}.
\end{equation}
Finally, for $0<\alpha\leq 1$ and any interval $[s,t]\subseteq[0,T]$, with $C_{p,\alpha}=\rev{66}e\alpha^\frac{1}{p}$,
\begin{equation}\label{eq:||X||_p<=exp(N_alpha^p)}
\|\mbX\|_{p,[s,t]}+\|\widetilde{\mbX}\|_{p,[s,t]}+|t-s|
\leq 
C_{p,\alpha}\exp\left(N_{\alpha,[s,t]}(w)\right).
\end{equation}
\end{corollary}
Next, we introduce \textit{enhanced Gaussian processes}, which are stochastic processes $\mbB=(B,\bB)$ that consist 
of a Gaussian process $B$ and of its enhancement $\bB$ defined so that  the sample paths $\mbB(\omega)=(B(\omega),\bB(\omega))$ are geometric rough paths. The sample paths $\mbB(\omega)$ are  called \textit{Gaussian rough paths}. An example of enhanced Gaussian process is  the Stratonovich lift  of a Brownian motion. %
The next assumption ensures that a Gaussian process $B$ can be lifted to an enhanced Gaussian process $\mbB=(B,\bB)$, and that the exponential of $N_{\alpha}(\mbB)$ (which will appear in our analysis through applications of \eqref{eq:Nalpha<=3CpNalpha_X_and_time} and \eqref{eq:||X||_p<=exp(N_alpha^p)}) is integrable. 
This assumption is considered in \cite{Friz2016,Bayer2016}, and is verified for a large number of Gaussian processes \cite{Friz2016}. %

\begin{assumption}[Gaussian process with regular covariance {\cite[Condition 10]{Bayer2016}}]\label{assum:Gaussian_lift}
 Let $T>0$ and $B = (B_1,\dots,B_d)$ %
 be a centered, continuous, $\R^d$-valued Gaussian process with independent components. Assume that the covariance of every component has H\"older dominated finite mixed $(1, \rho)$-variation for some $\rho\in[1,2)$ on $[0,T]^2$, that is, there exists $K<\infty$ such that for $k =1,\dots,d$, 
uniformly over
$0\leq s < t\leq T$, 
$\sup_{(t_i)\in\mathcal{P}([s,t]),(t_j')\in\mathcal{P}([s,t])
}
(\sum_{t_j'}(\sum_{t_i}
\E[B_{t_i,t_{i+1}}^kB_{t_j',t_{j+1}'}^k]
)^\rho
)^\frac{1}{\rho}\leq K|t-s|^\frac{1}{\rho}.
$
\end{assumption} 
\begin{theorem}[Enhanced Gaussian process and Gaussian rough paths] %
\label{thm:gaussian_rough_paths}
Let $T>0$, $\rho\in[1,\frac{3}{2})$, $p\in(2\rho,3)$, $(\Omega,\F,\Prob)$ be a probability space, and $B$  be a centered, continuous, $\R^d$-valued Gaussian process with independent components   satisfying Assumption \ref{assum:Gaussian_lift} for $\rho$. 
Then, there exists a unique stochastic process  $\mbB=(B,\bB)$  
 that is the natural lift of $B$, and whose sample paths are geometric \rev{rough} paths, that is, $\mbB(\omega)=(B(\omega),\bB(\omega))\in\sC^p_g([0,T],\R^d)$ almost surely. 
Moreover, for any $\alpha>0$ and $D\geq 0$, 
\begin{equation}\label{eq:gaussian_rough_paths:E[exp(Nalpha)]<infty}
\E\left[\exp\left(DN_{\alpha,[0,T]}(\mbB)\right)\right]<\infty.
\end{equation}
$\mbB$ is called %
enhanced Gaussian process, and its sample paths $\mbB(\omega)$ are called Gaussian rough paths.
\end{theorem}
The first claim of Theorem \ref{thm:gaussian_rough_paths} is in \cite[Theorem
15.33]{Friz2010}  (or \cite[Theorem 10.4]{Friz2020}, see also \cite[Corollary 2.3]{Friz2016}), 
and the second follows from results in \cite{Friz2013} (see also \cite[Theorem 11]{Bayer2016}). 

In Theorem \ref{thm:gaussian_rough_paths}, the enhancement $\bB$ in the ``natural'' lift $\mbB=(B,\bB)$ can be defined in various equivalent ways, %
 see \cite[Remark 10.7]{Friz2020}. In \cite[Theorem 10.4]{Friz2020},  the diagonal elements of $\bB$ are defined as $\bB^{ii}_{s,t}=\frac{1}{2}(B_{s,t}^i)^2$, the  off-diagonal terms 
 $\bB^{ij}$ are defined as an $L^2$ limit  $\bB_{s,t}^{ij}:=\lim_{|\pi|\to 0}\sum_{[u,v]\in\pi}B^i_{s,u}B^j_{u,v}$, and the other terms $\bB^{ji}$ are defined as $\bB^{ji}_{s,t}:=-\bB^{ij}_{s,t}+B^i_{s,t}B^j_{s,t}$, so that the algebraic conditions \eqref{eq:chen's_relation} and \eqref{eq:rough_path:integration_by_parts} are satisfied. %
 See also other equivalent definitions in \cite[Exercise 10.11]{Friz2020} and in \cite[Theorem
15.33]{Friz2010}. 

\section{Rough differential equations (RDEs) and error bounds}\label{sec:rdes}
In this section, we study the RDEs%
\begin{align}
\label{eq:RDE}
Y_t&=y+\int_0^tb(s,Y_s,u_s)\dd s+\int_0^t\sigma(s,Y_s)\dd\mbX_s,
\qquad\qquad\ \ \ \ \, t\in[0,T],
\\
\label{eq:RDE:linear}
V_t &=v +  
\int_0^t
\rev{\nablaof{x}b}(s,Y_s,u_s)V_s\dd s
+  
\int_0^t
\rev{\nablaof{x}\sigma}(s,Y_s)V_s\dd\mbX_s,
\quad t\in[0,T].
\end{align}
In Section \ref{sec:rdes:calculus}, we provide additional results  for controlled rough paths and rough integrals. 
In Section \ref{sec:rdes:existence_unicity}, we show that  solutions to the RDEs \eqref{eq:RDE} and \eqref{eq:RDE:linear} exist and are unique.
In Sections \ref{sec:rdes:bounds} and \ref{sec:rdes:bounds:linear}, we derive bounds for the solutions to these RDEs  by leveraging greedy partitions. Finally, in Section \ref{sec:rdes:random_integrable}, we study  solutions to random RDEs driven by Gaussian rough paths $\mbX=\mbB(\omega)$ and derive integrable errors bounds. 
\ifarxiv
For conciseness, 
\fi 
the detailed proofs of some results in this section are in the 
\ifarxiv
appendix.
\fi
 
\subsection{Calculus with rough paths: controlled rough paths and rough integration}\label{sec:rdes:calculus}
The lemmas in this section are slight variations of results in \cite[Section 3.2]{Friz2018} to highlight constants and to support a time-varying diffusion  $\sigma$. 
Their proofs are long but straightforward\ifarxiv  \ and are in the 
appendix.
\fi

\begin{lemma}[Products of controlled paths are controlled paths]\label{lem:control_path:product}
Let $p\in[2,3)$, $T>0$, $\mbX=(X,\bX)\in\sC^p$, and $(Y,Y'),(Z,Z')\in\sD^p_X$ be two controlled rough paths. Then, $(YZ,(YZ)')\in\sD^p_X$ with Gubinelli derivative 
$(YZ)'=ZY'+YZ'$. 
Moreover,
\begin{align}
\label{lem:control_path:product:YZ_p}
\|YZ\|_p
&\leq
C_p(\|Y\|_\infty\|Z\|_p+\|Z\|_\infty\|Y\|_p),
\\
\label{lem:control_path:product:(YZ)'_p}
\|(YZ)'\|_p
&\leq
C_p(\|Z\|_\infty\|Y'\|_p+\|Y'\|_\infty\|Z\|_p+\|Y\|_\infty\|Z'\|_p+\|Z'\|_\infty\|Y\|_p),
\\
\label{lem:control_path:product:R^YZ_p/2}
\|R^{YZ}\|_\frac{p}{2}
&\leq
C_p(\|Y\|_\infty\|R^Z\|_\frac{p}{2}+\|R^Y\|_\frac{p}{2}\|Z\|_\infty+\|Y\|_p\|Z\|_p).
\end{align}
\end{lemma}

\begin{lemma}[Error bounds for controlled paths]\label{lem:rough_path:f(Y)'-f(tilde(Y))'} 
Let $p\in[2,3)$, $T>0$, $\mbX=(X,\bX)\in\sC^p$ and $\widetilde{\mbX}=(\rev{\widetilde{X}},\rev{\widetilde{\bX}})\in\sC^p$,  
 $(Y,Y')\in\sD^p_X$ and $(\rev{\widetilde{Y}},\rev{\widetilde{Y}}')\in\sD^p_{\rev{\widetilde{X}}}$,  and $\sigma\in C^3_b$. Then, 
\begin{align}\label{eq:controlled_path:Y-Ytilde_p}
\|Y-\rev{\widetilde{Y}}\|_p
&\leq 
C_p\big(\Delta M_{Y'}\|X\|_p+M_{\rev{\widetilde{Y}}'}\|\Delta X\|_p+\|\Delta R^Y\|_\frac{p}{2}\big),
\\
\label{eq:controlled_path:sigma(Y)-sigma(Ytilde)_p}
\|\sigma(\cdot,Y)'-\sigma(\cdot,\rev{\widetilde{Y}})'\|_p
&\leq
C_p\|\sigma\|_{C_b^3}
(1+
K_Y+K_{\rev{\widetilde{Y}}}
)^3
(1+\|X\|_p+\|\rev{\widetilde{X}}\|_p)^3
(
1+T
)
\rev{\,\times}
\\
\nonumber
&\qquad\qquad\qquad\qquad
\rev{\big(}
\|\Delta X\|_p
+
\|\Delta Y_0\|
+
\|\Delta R^Y\|_\frac{p}{2}
+
\|\Delta Y'_0\|
+
\|\Delta Y'\|_p
\big),
\\
\label{eq:controlled_path:R^sigma(Y)-R^sigma(Ytilde)_p/2}
\|R^{\sigma(\cdot,Y)}-R^{\sigma(\cdot,\rev{\widetilde{Y}})}\|_\frac{p}{2}
&\leq 
C_p\|\sigma\|_{C^3_b}
(1+K_Y+K_{\rev{\widetilde{Y}}})^3
(1+\|X\|_p+\|\rev{\widetilde{X}}\|_p)^2
	(
1+T
)
\rev{\,\times}
\\
\nonumber
&\qquad\qquad\qquad\qquad
\rev{\big(}
\|\Delta Y_0\|+\|\Delta R^Y\|_\frac{p}{2}+\Delta M_{Y'}\|X\|_p+\|\Delta X\|_p\big).
\end{align}
\end{lemma}
\begin{remark}[Smoothness of $\sigma$]\label{remark:sigma_smoothness:2}
As discussed in Remark \ref{remark:sigma_smoothness}, the assumption that $\sigma$ is smooth in $t$ can be relaxed to assuming that $\sigma(\cdot,x)$ is uniformly \smash{$\frac{2}{p}$}-H\"older continuous. 
However, the proof of the inequality \eqref{eq:controlled_path:R^sigma(Y)-R^sigma(Ytilde)_p/2} for $\|R^{\sigma(\cdot,Y)}-R^{\sigma(\cdot,\rev{\widetilde{Y}})}\|_\frac{p}{2}$ breaks under weaker regularity assumptions. In particular, it breaks if we only assume that $
\|
\sigma(t,x)-
\sigma(s,x)
\|
\leq C|t-s|^{\frac{1}{2}-\epsilon} 
$ (e.g., if $\sigma(t,x)=B_t(\omega)$ is a sample path of  Brownian motion).
\end{remark}

\begin{lemma}[The rough integral $\int\rev{\sigma(\cdot,Y)}\dd\mbX$ defines a controlled path]\label{lem:rough_path:intfdX:error_bounds}
Let $p\in[2,3)$, $T>0$, $\mbX=(X,\bX)\in\sC^p$, $(Y,Y')\in\sD^p_X$, and $\sigma\in C^2_b$.  Then, 
$(Z,Z'):=\left(
\int_0^\cdot \sigma(s,Y_s)\dd\mbX_s,\sigma(\cdot,Y_\cdot)
\right)\in\sD^p_X$, 
and   
\begin{align}\label{eq:rough_path:R^int_sig_dX:p/2}
\|R^{\int_0^\cdot \sigma(s,Y_s)\dd\mbX_s}\|_{\frac{p}{2}}&\leq 
C_p\|\sigma\|_{C_b^2}(1+K_Y)^2(1+\|X\|_p)^2(1+T)
\|\mbX\|_p.
\end{align}
\end{lemma}

\begin{lemma}[Error bounds for $(\int\rev{\sigma(\cdot,Y)}\dd\mbX-\int\rev{\sigma(\cdot,\rev{\widetilde{Y}})}\dd\widetilde{\mbX})$]\label{lem:rough_path:intfdX_fY:error_bounds}
Let $p\in[2,3)$, $T>0$, $\mbX=(X,\bX)\in\sC^p$ and $\widetilde{\mbX}=(\rev{\widetilde{X}},\rev{\widetilde{\bX}})\in\sC^p$,  
$(Y,Y')\in\sD^p_X$ and $(\rev{\widetilde{Y}},\rev{\widetilde{Y}}')\in\sD^p_{\rev{\widetilde{X}}}$,  and $\sigma\in C^3_b$. Then,
\begin{align}\label{eq:rough_path:intfdX_fY:delta_sigma}
\|\sigma(\cdot,Y)-\sigma(\cdot,\rev{\widetilde{Y}})\|_p
&\leq 
C_p\|\sigma\|_{C^2_b}
(1 + K_Y+K_{\rev{\widetilde{Y}}})^2
(1+\|X\|_p+\|\rev{\widetilde{X}}\|_p+T)
\rev{\,\times}
\\
\nonumber
&\hspace{1.7cm}
\rev{\big(}
\|\Delta X\|_p+\|\Delta Y_0\|+(\|\Delta Y_0'\|+\|\Delta Y'\|_p)\|X\|_p+\|\Delta R^Y\|_\frac{p}{2}\big),
\\
\label{eq:rough_path:intfdX_fY:delta_R}
\|R^{\int_0^\cdot \sigma(s,Y_s)\dd\mbX_s}-R^{\int_0^\cdot \sigma(s,\rev{\widetilde{Y}}_s)\dd\widetilde{\mbX}_s}\|_{\frac{p}{2}}
&\leq 
C_p\|\sigma\|_{C_b^3}
(1+K_Y+K_{\rev{\widetilde{Y}}}+T)^3(1+\|\mbX\|_p+\|\widetilde{\mbX}\|_p)^4
(1+T)
\rev{\,\times}
\\
\nonumber
&\hspace{-34mm}
\rev{\big(}
\|\Delta\mbX\|_p+
\|\mbX\|_p(\|\Delta Y_0\|+\|\Delta Y_0'\|+\|\Delta R^Y\|_\frac{p}{2}+\|\Delta Y'\|_p+(\|\Delta Y_0'\|+\|\Delta Y'\|_p)\|X\|_p+\|\Delta X\|_p)
\big).
\end{align}
\end{lemma}

\subsection{Rough differential equations: existence and unicity of solutions}\label{sec:rdes:existence_unicity}

In this section, we study the existence and unicity of solutions to the  RDEs \eqref{eq:RDE} and \eqref{eq:RDE:linear}. 
Since these RDEs include a drift term $\int_0^\cdot b(t,Y_t,u_t)\dd t$, we use the following assumption and lemma to bound its $\frac{p}{2}$-variation.  
\begin{assumption}[Regularity of $b$]\label{assumption:b}%
Let $T>0$ and $U\subseteq\R^m$. 
The map $b:[0,T]\times\R^n\times U\to\R^n$ satisfies:
\begin{itemize}
\item 
$
b(\cdot, x, u):[0,T]\to\R^n
$ 
is measurable for all $(x,u)\in\R^n\times U$,
\item 
$
b(t,\cdot,\cdot):\R^n\times U\to\R^n
$
is continuous for almost every $t\in[0,T]$,
\item $b$ is bounded and Lipschitz in $x$: 
There exists a constant $C_b\geq 0$ such that 
$\|b(t,x,u)\|\leq C_b$ and   
$\|b(t,x,u)-b(t,\tilde{x},u)\|\leq C_b\|x-\tilde{x}\|$ for almost every $t\in[0,T]$, all $x,\tilde{x}\in\R^n$ and all $u\in U$.
\end{itemize}
\end{assumption}

\begin{lemma}[$p$-variations of Lebesgue integrals]\label{lem:bounds_int_b_ds}
Let $p\in[2,3)$, $T>0$, $U\subseteq\R^m$, %
 $b:[0,T]\times\R^n\times U\to\R^n$ satisfy Assumption \ref{assumption:b},  
$Y,\rev{\widetilde{Y}}:[0,T]\to\R^n$ be two continuous paths, and $u,\tilde{u}\in L^\infty([0,T],U)$. Then,  
\begin{align}\label{eq:bounds_int_b_ds}
\left\|\int_0^\cdot b(s,Y_s,u_s)\dd s\right\|_\frac{p}{2} &\leq C_{p,b}T,
\  \ \text{and} \ \ \ 
\left\|\int_0^\cdot (b(s,Y_s,u_s)-b(s,\rev{\widetilde{Y}}_s,u_s))\dd s\right\|_\frac{p}{2} \leq C_{p,b}T\|Y-\rev{\widetilde{Y}}\|_\infty.
\end{align}
Moreover, if $b$ is also Lipschitz in $u$, so that  $\|b(t,x,u)-b(t,x,\tilde{u})\|\leq C_b\|u-\tilde{u}\|$ 
for almost every $t\in[0,T]$, all $x\in\R^n$, and all $u,\tilde{u}\in U$, then 
\begin{align}\label{eq:bounds_int_b_ds:dy_du}
\left\|\int_0^\cdot (b(s,Y_s,u_s)-b(s,\rev{\widetilde{Y}}_s,\tilde{u}_s))\dd s\right\|_\frac{p}{2} 
\leq 
C_{p,b}T
\big(
\|Y-\rev{\widetilde{Y}}\|_\infty
+
\|u-\tilde{u}\|_{L^\infty([0,T],U)}
\big).
\end{align}
\end{lemma}
\ifarxiv
\begin{proof}
 By our assumptions,  $\int_0^t b(s,Y_s,u_s)\dd s$ and $\int_0^t b(s,\rev{\widetilde{Y}}_s,u_s)\dd s$ are well-defined Lebesgue integrals. %
Given $s,t\in[0,T]$, %
we obtain  $\|\int_s^t b(r,Y_r,u_r)\dd r\|
\leq C_b|t-s|$  
and 
$\|\int_s^t (b(r,Y_r,u_r)-b(r,\rev{\widetilde{Y}}_r,u_r))\dd r\|
\leq C_b\|Y-\rev{\widetilde{Y}}\|_\infty|t-s|$, 
and assuming that $b$ is moreover Lipschitz in $u$, 
$\|\int_s^t (b(r,Y_r,u_r)-b(r,\rev{\widetilde{Y}}_r,\tilde{u}_r))\dd r\|
\leq C_b(\|Y-\rev{\widetilde{Y}}\|_\infty+\|u-\tilde{u}\|_\infty)|t-s|$. 
The desired inequalities then follow from \eqref{eq:path_finite_var:sum_p/2vars} in  Lemma \ref{lem:pvariation:inequalities}.
\end{proof}
\fi

\begin{remark}[On the boundedness of the drift $b$]\label{remark:b_bounded}
We assume that  $b$ is bounded, which is stronger than making a linear growth assumption (LG) in $x$ (i.e., $\|b(t,0,u)\|\leq C_b$, so that $\|b(t,x,u)\|\leq C_b(1+\|x\|)$), under which the bound in \eqref{eq:bounds_int_b_ds} becomes $\|\int_0^\cdot b(s,Y_s,u_s)\dd s\|_\frac{p}{2} \leq C_{p,b}(1+\|Y\|_\infty)T$. However, the proof of Theorem \ref{thm:rdes:existence_unicity} about existence and unicity of solutions to  nonlinear RDEs  breaks under LG, because it relies on a fixed point argument and on stitching solutions on intervals $[0,t]$ whose size is independent of the initial conditions, and the size of these intervals does depend on the initial conditions under LG. Similarly, our bounds for nonlinear RDEs (e.g., in Theorem \ref{thm:rdes:integrable})  assume that $b$ is bounded, as they rely on bounding $p$- and $\frac{p}{2}$-variations of solutions to nonlinear RDEs over short intervals, and these $p$-variations would otherwise depend on the initial conditions under LG. 
For these reasons, the case where $b$ is linear in $x$ (which does not satisfy Assumption \ref{assumption:b}) is handled by a separate analysis for linear RDEs, see  Theorems \ref{thm:rde:linear:existence_uniqueness}  and  \ref{thm:rdes:integrable}.
\end{remark}

\subsubsection{Nonlinear rough differential equations}\label{sec:rdes:existence_unicity:nonlinear}
\begin{theorem}[Nonlinear RDE \eqref{eq:RDE}: existence and unicity of solutions]
\label{thm:rdes:existence_unicity}
Let $p\in[2,3)$, $T>0$, $y\in\R^n$, $U\subseteq\R^m$, $u\in L^\infty([0,T],U)$, 
$b:[0,T]\times\R^n\times U\to\R^n$ satisfy Assumption \ref{assumption:b},  $\sigma\in C_b^3([0,T]\times\R^n,\R^{n\times d})$, and   
$\mbX=(X,\bX)\in\sC^p([0,T],\R^d)$ be a \rev{rough} path.

Then, there exists a unique %
solution $(Y,Y')\in\sD^p_X([0,T],\R^n)$ with $Y'=\sigma(\cdot,Y)$ to the nonlinear RDE $Y_t=y+\int_0^tb(s,Y_s,u_s)\dd s+\int_0^t\sigma(s,Y_s)\dd\mbX_s$ in \eqref{eq:RDE},
where $(\sigma(\cdot,Y),\sigma(\cdot,Y)')=(\sigma(\cdot,Y),\rev{\nablaof{x}\sigma}(\cdot,Y)Y')\in\sD^p_X$.  
\end{theorem}
\ifarxiv
The proof of Theorem  \ref{thm:rdes:existence_unicity} relies on a classical fixed point argument and is in the 
appendix. 
 We adapt arguments in the proof of \cite[Theorem 2.5]{Allan2020} and \cite[Theorem 3.8]{Friz2018} for our different problem setting with time-varying coefficients $(b,\sigma)$ and a control $u$ that is not smooth but  is not in the diffusion $\sigma$.
\fi

\subsubsection{Linear rough differential equations}\label{sec:rdes:existence_unicity:linear}
Existence and unicity of solutions to the linear RDE \eqref{eq:RDE:linear} does not follow from Theorem \ref{thm:rdes:existence_unicity}, since the coefficients of a linear RDE are not bounded a-priori. First, we prove the existence and unicity of solutions to generic linear RDEs with drift (Theorem \ref{thm:rde:linear:existence_uniqueness}). Then, we apply it to the linear RDE \eqref{eq:RDE:linear} (Corollary \ref{cor:rde:linear:existence_uniqueness}).  Such a result can be considered standard, %
although we could not find a result  in the literature that can handle a time-varying diffusion and an irregular drift. 

Note that Theorem \ref{thm:rde:linear:existence_uniqueness} for linear RDEs does not rely on Theorem \ref{thm:rdes:existence_unicity}, so 
Theorem \ref{thm:rde:linear:existence_uniqueness} ensures that the RDE \eqref{eq:RDE} has a unique solution if $b$ is linear but is not bounded, see Remark \ref{remark:b_bounded}.

\begin{theorem}[Linear RDEs: existence and unicity of solutions]\label{thm:rde:linear:existence_uniqueness}
Let $p\in[2,3)$,  $T>0$,  $v\in\R^n$ be an initial condition, $A\in L^\infty([0,T],\R^{n\times n})$ be an integrable map, 
$\mbX\in\rev{\sC^p}([0,T],\R^d)$ be a \rev{rough} path, 
and $(\Sigma,\Sigma')\in\sD_X^p([0,T],\R^{n\times d \times n})$   be a controlled \rev{rough} path.   
Then, there exists a unique solution $(V,V')\in\sD^p_X([0,T],\R^n)$ with $V'=\Sigma V$ to the linear RDE 
\begin{align}
\label{eq:rde_linear}
V_t &=v +  
\int_0^t
A_sV_s\dd s
+  
\int_0^t
\Sigma_sV_s\dd\mbX_s,
\qquad t\in[0,T].
\end{align}  
\end{theorem}
The proof %
of Theorem \ref{thm:rde:linear:existence_uniqueness} consists of rewriting  \eqref{eq:rde_linear} as a linear RDE with constant coefficients driven by a \rev{new rough path} defined using $(A,\Sigma,\mbX)$ 
and concluding \rev{with \cite[Theorem 2]{Lejay2009}}\ifarxiv
, see  the appendix.
\fi  

\begin{assumption}[Stronger regularity of $b$]\label{assumption:b:stronger:linear_rde}
Let $T>0$, $U\subseteq\R^m$, and  $b:[0,T]\times\R^n\times U\to\R^n$ satisfy
\begin{itemize} 
\item 
$b(\cdot, x, u):[0,T]\to\R^n$  is  measurable for all $(x,u)\in\R^n\times U$,  
\item  
$b(t,\cdot,\cdot):\R^n\times U\to\R^n$ 
is continuous for almost every $t\in[0,T]$,
\item 
$b(t,\cdot,u):\R^n\to\R^n$  is continuously differentiable  for almost every $t\in[0,T]$  and  all $u\in U$,  
\item %
$\big\|\rev{\nablaof{x}b}(t,x,u)\big\|\leq C_b$ for almost every $t\in[0,T]$ and all $(x,u)\in\R^n\times U$ for some constant $C_b\geq 0$.
\end{itemize}
\end{assumption}
\begin{corollary}[Linearized RDE \eqref{eq:RDE:linear}: existence and unicity of solutions]\label{cor:rde:linear:existence_uniqueness}
Let $p\in[2,3)$,  $T>0$,  $v\in\R^n$, $U\subseteq\R^m$, 
$u\in L^\infty([0,T],U)$, 
$b:[0,T]\times\R^n\times U\to\R^n$  satisfy  Assumption \ref{assumption:b:stronger:linear_rde}, 
$\sigma\in C_b^3([0,T]\times\R^n,\R^{n\times d})$, 
$\mbX\in\rev{\sC^p}([0,T],\R^d)$ be a \rev{rough} path, and 
$(Y,Y')\in\sD^p_X([0,T],\R^n)$ be a controlled \rev{rough} path.  
Then, there exists a unique solution $(V,V')\in\sD^p_X([0,T],\R^n)$ with $V'=\rev{\nablaof{x}\sigma}(\cdot,Y_\cdot)V_\cdot$ to the linear RDE $V_t =v +  
\int_0^t
\rev{\nablaof{x}b}(s,Y_s,u_s)V_s\dd s
+  
\int_0^t
\rev{\nablaof{x}\sigma}(s,Y_s)V_s\dd\mbX_s$ in \eqref{eq:RDE:linear}. 
\end{corollary}
\begin{proof}
By Assumption \ref{assumption:b:stronger:linear_rde}, $\rev{\nablaof{x}b}(\cdot,Y_\cdot,u_\cdot)\in L^\infty([0,T],\R^{n\times n})$. By Lemma \ref{lem:rough_path:sigma(.,Y):controlled}, $\big(\rev{\nablaof{x}\sigma}(\cdot,Y_\cdot),\rev{\nablaof{x}^2\sigma}(\cdot,Y_\cdot)Y'_\cdot\big)\\\in\sD_X^p([0,T],\R^{n\times d\times n})$. The conclusion then follows from Theorem \ref{thm:rde:linear:existence_uniqueness}.
\end{proof}

\subsection{Bounds on solutions to nonlinear RDEs}\label{sec:rdes:bounds}
In this section, we derive bounds for solutions to the nonlinear RDE \eqref{eq:RDE} by leveraging greedy partitions and the quantity $N_{\alpha,[0,T]}(\mbX)$ in Definition \ref{def:Nalpha}.  The first step consists of deriving bounds that are independent of $\|\mbX\|_p$ and $T$ over intervals short-enough. 
The size $\alpha$ of these intervals is then used to define a greedy partition and derive finer bounds over the entire interval $[0,T]$ as a function of the quantity $N_{\alpha,[0,T]}(\mbX)$. Importantly, the final bounds in Lemma \ref{lem:rdes:bound:entire_interval} and Proposition \ref{prop:rdes:error_bound:entire_interval}
are integrable for Gaussian rough paths $\mbX=\mbB(\omega)$, see Theorem \ref{thm:gaussian_rough_paths} and Theorem \ref{thm:rdes:integrable}.  

\ifarxiv
For conciseness, the proofs of multiple results in this section are in the appendix.
\else
\fi
\begin{proposition}[Bounds for solutions to RDEs on short intervals]\label{prop:bounds_pvars_solutions_RDEs}
Let $p,T,y,U,u,b,\sigma,\mbX,Y,Y'$ be as in Theorem \ref{thm:rdes:existence_unicity}, where $b$ satisfies Assumption \ref{assumption:b}, $\sigma\in C_b^3$, and $(Y,Y')\in\sD^p_X$ solves the RDE \eqref{eq:RDE}. 
Then,
\begin{align}
\label{eq:bounds_pvars_solutions_RDEs:sigma(Y)'_p}
\|\sigma(\cdot,Y)'\|_p &\leq C_{p,\sigma}(\|Y\|_p+T).
\end{align}
Moreover, there \rev{exist} two constants $C_{p,b,\sigma}\geq 1$ and $0<\alpha_{p,b,\sigma}<1$ such that 
\begin{align}
\label{eq:bounds_pvars_solutions_RDEs:small_intervals:Y}
\|Y\|_{p,I} &\leq C_{p,b,\sigma},
\\
\label{eq:bounds_pvars_solutions_RDEs:small_intervals:R}
\|R^Y\|_{\frac{p}{2},I} &\leq C_{p,b,\sigma},
\\
\label{eq:bound:KY}
K_{Y,I}= \|Y_{t_0}'\|+\|Y'\|_{p,I}+\|R^Y\|_{\frac{p}{2},I}
&\leq
C_{p,b,\sigma}.
\end{align}
for any interval $I=[t_0,t_1]\subseteq[0,T]$ small-enough so that  $
\|\mbX\|_{p,I}+|I|\leq
 \alpha_{p,b,\sigma}^\frac{1}{p}$. 
\end{proposition}
\ifarxiv
To show Proposition \ref{prop:bounds_pvars_solutions_RDEs}, we take inspiration from the proof of  \cite[Proposition 2.4]{Allan2020}. 
 The main differences are handling a time-varying diffusion $\sigma(\cdot,Y)$ and using an interval $I$ small-enough so that the quantities in   \eqref{eq:bounds_pvars_solutions_RDEs:small_intervals:Y}-\eqref{eq:bound:KY} are bounded by a constant  that only depends on $(p,b,\sigma)$ and not on $(X,Y,I)$. 
\fi

\begin{proposition}
[Error bound for solutions to RDEs on short intervals]\label{prop:rdes:error_bound}
Let $p\in[2,3)$, $T>0$, $y,\tilde{y}\in\R^n$, $U\subseteq\R^m$, 
$u,\tilde{u}\in L^\infty([0,T],U)$, 
$b:[0,T]\times\R^n\times U\to\R^n,(t,x,u)\mapsto b(t,x,u)$ satisfy Assumption \ref{assumption:b} and be Lipschitz in $u$, 
$\sigma\in C_b^3([0,T]\times\R^n,\R^{n\times d})$,  
$\mbX,\widetilde{\mbX}\in\sC^p([0,T],\R^d)$ be two \rev{rough} paths, 
 $(Y,Y')\in\sD^p_X([0,T],\R^n)$ and $(\widetilde{Y},\widetilde{Y}')\in\sD^p_{\rev{\widetilde{X}}}([0,T],\R^n)$ with $(Y',\widetilde{Y}')=(\sigma(\cdot,Y_\cdot),\sigma(\cdot,\widetilde{Y}_\cdot))$ be the solutions to the RDEs 
\begin{align*}
Y_t=y+\int_0^tb(s,Y_s,u_s)\dd s+\int_0^t\sigma(s,Y_s)\dd\mbX_s,
\quad
\widetilde{Y}_t=\tilde{y}+\int_0^tb(s,\widetilde{Y}_s,\tilde{u}_s)\dd s+\int_0^t\sigma(s,\widetilde{Y}_s)\dd\widetilde{\mbX}_s,
\quad t\in[0,T].
\end{align*}
Then, there \rev{exist} two constants $C_{p,b,\sigma}\geq 1$ and $0<\alpha_{p,b,\sigma}<1$ such that
\begin{gather}\label{eq:bounds_pvars_solutions_RDEs:small_intervals:combined}
\|Y\|_{p,I},\, 
\|\widetilde{Y}\|_{p,I},\, 
\|R^Y\|_{\frac{p}{2},I},\, 
\|R^{\widetilde{Y}}\|_{\frac{p}{2},I},\, 
K_{Y,I},\, 
K_{\widetilde{Y},I}\leq
C_{p,b,\sigma},
\quad\text{and}
\\[2mm]
\label{eq:RDE:DY'+dRY:close}
\|Y'-\widetilde{Y}'\|_{p,I}+\|R^Y-R^{\widetilde{Y}}\|_{\frac{p}{2},I}
\leq
C_{p,b,\sigma}(\|Y_{t_0}-\widetilde{Y}_{t_0}\|
+
\|\Delta\mbX\|_{p,I}
+
|I|\|u-\tilde{u}\|_{L^\infty(I,U)}),
\end{gather}
for any interval $I=[t_0,t_1]\subseteq[0,T]$ such that $\|\mbX\|_{p,I}+\|\widetilde{\mbX}\|_{p,I}+|I|\leq
 \alpha_{p,b,\sigma}^\frac{1}{p}$.
\end{proposition}
\ifarxiv
The proof of Proposition \ref{prop:rdes:error_bound} takes similar steps as in the proofs of \cite[Proposition 2.6]{Allan2020} and \cite[Theorem 3.9]{Friz2018}, and as  when proving the contractivity of $\M_t$ for Theorem \ref{thm:rdes:existence_unicity}. The main differences compared to prior work are using bounds with $\sigma(\cdot,Y)$ that 
is time-varying and working on short-enough intervals to obtain a bound with a constant $C_{p,b,\sigma}$ that only depends on $(p,b,\sigma)$, and not on $(X,\rev{\widetilde{X}},Y,\widetilde{Y},I)$. This last difference is the key to obtain integrable bounds. 
\fi

\begin{lemma}[Bounds for solutions to RDEs on long intervals]\label{lem:rdes:bound:entire_interval}
Let $p,T,y,U,u,b,\sigma,\mbX,Y,Y'$ be as in Theorem \ref{thm:rdes:existence_unicity}, where $b$ satisfies Assumption \ref{assumption:b}, $\sigma\in C_b^3$, and $(Y,Y')\in\sD^p_X$ solves the RDE \eqref{eq:RDE}.  
Define $N_{\alpha,[0,T]}(\mbX)$ as in Definition \ref{def:Nalpha}. 
Then, there exist constants $C_{p,T,b,\sigma}\geq 1$ and $0<\alpha_{p,b,\sigma}<1$ such that  
\begin{align}
\label{eq:RDE:Yp_Y'p_RYp/2:bound:full_interval}
\|Y\|_{p,[0,T]}+
\|Y'\|_{p,[0,T]}+
\|R^Y\|_{\frac{p}{2},[0,T]}
&\leq
C_{p,T,b,\sigma}\exp\left(
C_{p,T,b,\sigma}N_{\alpha_{p,b,\sigma},[0,T]}(\mbX)
\right),
\\
\label{eq:RDE:Y:bound:full_interval}
\|Y\|_{\infty,[0,T]}
&\leq
C_{p,T,b,\sigma}\exp\left(
C_{p,T,b,\sigma}N_{\alpha_{p,b,\sigma},[0,T]}(\mbX)
\right)+\|y\|.
\end{align}
\end{lemma}
\begin{proof}
Define $N_{\alpha,I}(w)$ and $w$ as in Definition \ref{def:Nalpha} and Corollary \ref{cor:Nalpha:NX_NXtilde_NT:small_intervals} \rev{(with $\widetilde{\mbX}=0$)}, respectively. 
Let $C_{p,b,\sigma}\geq 1$ and $\alpha_{p,b,\sigma}>0$ be the constants in Proposition \ref{prop:bounds_pvars_solutions_RDEs}, 
and $I=[s,t]\subseteq[0,T]$ be any interval such that $w(s,t)\leq\alpha:=\alpha_{p,b,\sigma}$. 
Then, by   Corollary \ref{cor:Nalpha:NX_NXtilde_NT:small_intervals}, $
\|\mbX\|_{p,I}+|I|
\mathop{\leq}\limits^{\eqref{eq:|X|+|Xtilde|+|dt|<=alpha^p}}
 \alpha_{p,b,\sigma}^\frac{1}{p}$, so  Proposition \ref{prop:bounds_pvars_solutions_RDEs} implies that %
$\|Y\|_{p,I} 
\mathop{\leq}\limits^{\eqref{eq:bounds_pvars_solutions_RDEs:small_intervals:Y}}
C_{p,b,\sigma}$
and $\|Y'\|_{p,I}+\|R^Y\|_{\frac{p}{2},I}
\mathop{\leq}\limits^{\eqref{eq:bound:KY}} C_{p,b,\sigma}$.
Thus, as defined in Definition \ref{def:Nalpha}, with $N=N_{\alpha,[0,T]}(w)$, the partition 
 $\{\tau_i, i=0,1,\dots,N+1\}$ 
of the interval $[0,T]$, which satisfies $w(\tau_i,\tau_{i+1})\leq \alpha$ for all $i$,  is such that 
 $$
\|Y\|_{p,[0,T]}\leq (N+1)(
\sum_{i=0}^N\|Y\|_{p,[\tau_i,\tau_{i+1}]}^p
)^\frac{1}{p}
\leq (N+1)
\sum_{i=0}^N\|Y\|_{p,[\tau_i,\tau_{i+1}]}\leq (N+1)^2C_{p,b,\sigma}
\leq
\rev{C_{p,b,\sigma}\exp(N)},
$$
where the first inequality follows from  Lemma \ref{lem:pvar:intervals}. The same inequality holds for $\|Y'\|_{p,[0,T]}$ and $\|R^Y\|_{\frac{p}{2},[0,T]}$\rev{, using also \eqref{eq:pvar:RX_intervals} in Lemma \ref{lem:pvar:RX_intervals} to bound $\|R^Y\|_{\frac{p}{2},[0,T]}$.} The desired first inequality \eqref{eq:RDE:Yp_Y'p_RYp/2:bound:full_interval} then follows from $\exp(N)\leq \exp(5C_p(
N_{\alpha,[0,T]}(\mbX)+T/\alpha+1))\leq C_{p,T,b,\sigma,\alpha}\exp(C_pN_{\alpha,[0,T]}(\mbX))$ by \eqref{eq:Nalpha<=3CpNalpha_X_and_time}. 
Finally, the desired inequality \eqref{eq:RDE:Y:bound:full_interval} follows from $\|Y\|_{\infty,[0,T]}\leq \|Y\|_{p,[0,T]}+\|Y_0\|$ by  \eqref{eq:path_finite_var:infty_ineq}, and we conclude.
\end{proof} 
Proposition \ref{prop:rdes:error_bound:entire_interval} below is the main result of this section. It is similar to \cite[Theorem 4]{Bayer2016}, but includes  
a drift term %
with a control $u$ and a time-varying diffusion $\sigma$.  
Its proof is similar to the proof of Lemma \ref{lem:rdes:bound:entire_interval} by appropriately  replacing inequalities   but is slightly longer, so we provide it in  
\ifarxiv
the appendix.
\fi
\begin{proposition}
[Error bound for solutions to RDEs on long intervals]\label{prop:rdes:error_bound:entire_interval}
Define $p,T,y,\tilde{y},U,u,\tilde{u},b,\sigma$, $\mbX,\widetilde{\mbX}$
as in Proposition \ref{prop:rdes:error_bound}, 
where $b$ satisfies Assumption \ref{assumption:b} and is Lipschitz in $u$ and $\sigma\in C_b^3$, and let $(Y,Y')\in\sD^p_X$ and $(\widetilde{Y},\widetilde{Y}')\in\sD^p_{\widetilde{X}}$  with $Y'=\sigma(\cdot,Y_\cdot)$ and $\widetilde{Y}'=\sigma(\cdot,\widetilde{Y}_\cdot)$ 
be the solutions to the RDEs
\begin{align*}
Y_t=y+\int_0^tb(s,Y_s,u_s)\dd s+\int_0^t\sigma(s,Y_s)\dd\mbX_s,
\quad
\widetilde{Y}_t=\tilde{y}+\int_0^tb(s,\widetilde{Y}_s,\tilde{u}_s)\dd s+\int_0^t\sigma(s,\widetilde{Y}_s)\dd\widetilde{\mbX}_s,
\quad t\in[0,T].
\end{align*} 
Then, there exist constants $C_{p,b,\sigma}\geq 1$ and $0<\alpha_{p,b,\sigma}<1$ such that  
for any interval $I=[t_0,t_1]\subseteq[0,T]$, 
\begin{align}\label{eq:RDE:DY'+dRY:close:full_interval}
\|Y'-\widetilde{Y}'\|_{p,I}+\|R^Y-R^{\widetilde{Y}}\|_{\frac{p}{2},I}
&\leq
C_{p,b,\sigma}\exp\left(
C_{p,b,\sigma}N_{\alpha_{p,b,\sigma},I}(w)
\right)(\|\Delta Y_{t_0}\|
+
\|\Delta \mbX\|_{p,I}
+
|I|\|\Delta u\|_{L^\infty,I}
)
\\
\label{eq:RDE:DY:close:full_interval}
\|Y-\widetilde{Y}\|_{\infty,I} 
&\leq
C_{p,b,\sigma}\exp\left(
C_{p,b,\sigma}N_{\alpha_{p,b,\sigma},I}(w)
\right)(\|\Delta Y_{t_0}\|
+
\|\Delta \mbX\|_{p,I}
+
|I|\|\Delta u\|_{L^\infty,I}
)
\end{align}
where $\|\Delta u\|_{L^\infty,I}=\|u-\tilde{u}\|_{L^\infty(I,U)}$, and $(N_{\alpha_{p,b,\sigma},I}(w),w)$ are  defined in Definition \ref{def:Nalpha} and Corollary \ref{cor:Nalpha:NX_NXtilde_NT:small_intervals}. %

Moreoever, if $u$ and $\tilde{u}$ only differ on a subinterval $J\subseteq I$, i.e., $u_t=\tilde{u}_t$ for almost all $t\in I\setminus J$, then
\begin{align}\label{eq:RDE:DY'+dRY:close:full_interval:u_subinterval}
\|Y'-\widetilde{Y}'\|_{p,I}+\|R^Y-R^{\widetilde{Y}}\|_{\frac{p}{2},I}
&\leq
C_{p,b,\sigma}\exp\left(
C_{p,b,\sigma}N_{\alpha_{p,b,\sigma},I}(w)
\right)(\|\Delta Y_{t_0}\|
+
\|\Delta \mbX\|_{p,I}
+
|J|\|\Delta u\|_{L^\infty,I}
)
\\
\label{eq:RDE:DY:close:full_interval:u_subinterval}
\|Y-\widetilde{Y}\|_{\infty,I} 
&\leq
C_{p,b,\sigma}\exp\left(
C_{p,b,\sigma}N_{\alpha_{p,b,\sigma},I}(w)
\right)(\|\Delta Y_{t_0}\|
+
\|\Delta \mbX\|_{p,I}
+
|J|\|\Delta u\|_{L^\infty,I}
)
\end{align}
for some constants $C_{p,b,\sigma}\geq 1$ and $0<\alpha_{p,b,\sigma}<1$. 
\end{proposition}

\ifarxiv
\subsection{Bounds on solutions to linear RDEs  and for the Jacobian flow}\label{sec:rdes:bounds:linear}
\fi
\ifarxiv
Next, we derive bounds on solutions to linear RDEs and on the Jacobian flow of nonlinear RDEs, stated in Lemmas \ref{lem:rde:linear:bounded_solutions} and  \ref{lem:rde:linearized:bounded_solutions}.   
For conciseness, the proofs of these results are in the 
appendix.
\fi

Lemma \ref{lem:rde:linear:bounded_solutions} below is similar to \cite[Proposition 8.13]{Friz2020}, generalizing it to RDEs with drift, control input, and time-varying coefficients, and giving bounds that are integrable   for Gaussian rough paths $\mbX=\mbB(\omega)$ (see Theorems \ref{thm:gaussian_rough_paths} and  \ref{thm:rdes:integrable}). 
The proof follows similar steps as the proof of Proposition \ref{prop:rdes:error_bound:entire_interval}
using the bounds in Proposition \ref{prop:rough_integral_welldefined:error_bound}, Corollary \ref{cor:Nalpha:NX_NXtilde_NT:small_intervals}, 
Lemma \ref{lem:control_path:product}, 
and the  Gr\"onwall Lemma for rough paths \cite[Lemma 2.12]{Deya2019}.

\begin{lemma}[Bounds for solutions to linear RDEs]\label{lem:rde:linear:bounded_solutions}
Let $p,T,y,A,\mbX,\Sigma,\Sigma'$ be as in Theorem \ref{thm:rde:linear:existence_uniqueness}, where $A\in L^\infty$, $(\Sigma,\Sigma')\in\sD_X^p$, and  $(V,V')\in\sD^p_X$ solves the linear RDE  \eqref{eq:rde_linear}.  
Define $N_{\alpha,[0,T]}(\mbX)$ as in Definition \ref{def:Nalpha}. 
Assume that  there \rev{exist} two constants $C_\Sigma\geq 1$ and $0<\alpha_\Sigma<1$ such that 
\begin{align*} 
\|\Sigma\|_{\infty,I} + \|\Sigma\|_{p,I} +\|\Sigma'\|_{p,I}+\|R^\Sigma\|_{\frac{p}{2},I}
&\leq
C_\Sigma
\end{align*}
for any interval $I\subseteq[0,T]$ such that $\|\mbX\|_{p,I}+|I|\leq
 \alpha_\Sigma^{\frac{1}{p}}$.  
Then,  
\begin{align}
\label{eq:RDE:linear:Vp_V'p_RVp/2:bound}
\|V\|_{p,[0,T]}+
\|V'\|_{p,[0,T]}+
\|R^V\|_{\frac{p}{2},[0,T]}
&\leq 
C_{p,T,A,\Sigma}\exp\left(C_{p,T,A,\Sigma}N_{\alpha_{p,A,\Sigma},[0,T]}(\mbX)\right)\|v\|,
\\
\label{eq:RDE:linear:V:bound}
\|V\|_{\infty,[0,T]}
&\leq 
C_{p,T,A,\Sigma}\exp\left(C_{p,T,A,\Sigma}N_{\alpha_{p,A,\Sigma},[0,T]}(\mbX)\right)\|v\|,
\end{align}
where the constants $C_{p,T,A,\Sigma}\geq 1$ and $0<\alpha_{p,A,\Sigma}\leq\alpha_\Sigma$    only depend on $(p,T,\|A\|_\infty,C_\Sigma,\alpha_\Sigma)$. 
\end{lemma}

Lemma \ref{lem:rde:linearized:bounded_solutions} below is similar to \cite[Theorem 6.5]{Cass2013} (see also \cite[Proposition 5]{Friz2013}), generalizing it to RDEs with  drift $b$, control input $u$, and time-varying diffusion $\sigma$. Note that the bounds \eqref{eq:rde:linearized:Vp_V'p_RVp/2:bound} and \eqref{eq:rde:linearized:bounded_solutions} do not depend on $Y$: This fact is used in the proof of Lemma \ref{lem:rde:linearized:continuous} that gives the continuity of the map $(y,\mbX,u)\mapsto V$. 
Lemma \ref{lem:rde:linearized:bounded_solutions}  follows from 
 the bounds in Lemma \ref{lem:rough_path:sigma(.,Y):controlled} and Lemma \ref{lem:rde:linear:bounded_solutions}.

\begin{lemma}[Bounds for the Jacobian flow]
\label{lem:rde:linearized:bounded_solutions}
Let $p\in[2,3)$,  $T>0$,  $y,v\in\R^n$, $U\subseteq\R^m$, 
$u\in L^\infty([0,T],U)$, 
$b:[0,T]\times\R^n\times U\to\R^n$  satisfy  Assumption \ref{assumption:b:stronger:linear_rde}, 
$\sigma\in C_b^3([0,T]\times\R^n,\R^{n\times d})$,  
$\mbX\in\rev{\sC^p}([0,T],\R^d)$, 
and $(Y,Y'),(V,V')\in\sD^p_X([0,T],\R^n)$ 
with $Y'=\sigma(\cdot,Y_\cdot)$ and $V'=\rev{\nablaof{x}\sigma}(\cdot,Y_\cdot)V_\cdot$ 
be the solutions to the  RDEs $Y_t=y+\int_0^tb(s,Y_s,u_s)\dd s+\int_0^t\sigma(s,Y_s)\dd\mbX_s$ in \eqref{eq:RDE}  and  $V_t =v +  
\int_0^t
\rev{\nablaof{x}b}(s,Y_s,u_s)V_s\dd s
+  
\int_0^t
\rev{\nablaof{x}\sigma}(s,Y_s)V_s\dd\mbX_s$ in \eqref{eq:RDE:linear}. 
Then, there \rev{exist} two constants $C_{p,T,b,\sigma}\geq 1$ and $0<\alpha_{p,b,\sigma}<1$ such that 
\begin{align}
\label{eq:rde:linearized:Vp_V'p_RVp/2:bound}
\|V\|_{p,[0,T]}+
\|V'\|_{p,[0,T]}+
\|R^V\|_{\frac{p}{2},[0,T]}
&\leq 
C_{p,T,b,\sigma}\exp\left(C_{p,T,b,\sigma}N_{\alpha_{p,b,\sigma},[0,T]}(\mbX)\right)\|v\|,
\\
\label{eq:rde:linearized:bounded_solutions}
\|V\|_{\infty,[0,T]}
&\leq 
C_{p,T,b,\sigma}\exp\left(C_{p,T,b,\sigma}N_{\alpha_{p,b,\sigma},[0,T]}(\mbX)\right)\|v\|,
\end{align} 
where 
$N_{\alpha,[0,T]}(\mbX)$ is defined in Definition \ref{def:Nalpha}.
\end{lemma}

\subsection{Integrability of solutions to nonlinear and linear RDEs}\label{sec:rdes:random_integrable}
Finally, we combine the results from the previous sections to show that pathwise solutions  to nonlinear and linear RDEs driven by Gaussian rough paths $\mbB(\omega)$ are  integrable. First, in Lemmas \ref{lem:rdes:continuity} and  \ref{lem:rde:linearized:continuous}, we show the continuity of the It\^o-Lyons map $(y,\mbX)\mapsto Y_{(y,\mbX)}$, that is, that solutions to RDEs are continuous with respect to the initial conditions and the driving rough path. With this result follows the measurability of the map $\omega\mapsto Y_{(y(\omega),\mbB(\omega))}$ that assigns the pathwise solution to a random RDE driven by a Gaussian rough path $\mbX=\mbB(\omega)$. Thanks to the favorable integrability properties of enhanced Gaussian processes $\mbB$ in Theorem \ref{thm:gaussian_rough_paths}, we conclude that such solutions are integrable in Theorem \ref{thm:rdes:integrable}. These results will be used throughout  the proof of \pmp, e.g., to make sense of the maximality condition \eqref{eq:spmp:maximizality_condition} that involves $\E\left[
H(t,x_t,v,p_t,\mathfrak{p}_0)
\right]$.

\begin{lemma}
[Continuity of solutions to RDEs]\label{lem:rdes:continuity}
Define $p,T,y,\tilde{y},U,u,\tilde{u},b,\sigma,\mbX,\widetilde{\mbX}$
as in Proposition \ref{prop:rdes:error_bound}, 
where $b$ satisfies Assumption \ref{assumption:b} and is Lipschitz in $u$ and $\sigma\in C_b^3$, 
let $(Y,Y')\in\sD^p_X$ and $(\widetilde{Y},\widetilde{Y}')\in\sD^p_{\widetilde{X}}$  with $Y'=\sigma(\cdot,Y_\cdot)$ and $\widetilde{Y}'=\sigma(\cdot,\widetilde{Y}_\cdot)$  be the solutions to the RDEs
\begin{align*}
Y_t=y+\int_0^tb(s,Y_s,u_s)\dd s+\int_0^t\sigma(s,Y_s)\dd\mbX_s,
\quad
\widetilde{Y}_t=\tilde{y}+\int_0^tb(s,\widetilde{Y}_s,\tilde{u}_s)\dd s+\int_0^t\sigma(s,\widetilde{Y}_s)\dd\widetilde{\mbX}_s,
\quad t\in[0,T],
\end{align*} and assume that there exists a constant $M\geq 0$ such that   $\|\mbX\|_p,\|\widetilde{\mbX}\|_p\leq M$. Then,%
\begin{align}
\label{eq:RDE:DY'+dRY:close:full_interval:continuity}
\|Y'-\widetilde{Y}'\|_{p,[0,T]}+
\|R^Y-R^{\widetilde{Y}}\|_{\frac{p}{2},[0,T]}
&\leq 
C_{p,T,b,\sigma,M}
\left(
\|y-\tilde{y}\| 
+
\|\Delta \mbX\|_{p,[0,T]}
 + 
 T\|u-\tilde{u}\|_{L^\infty([0,T],U)}
\right),
\\
\label{eq:RDE:DY:close:full_interval:continuity}
\|Y-\widetilde{Y}\|_{\infty,[0,T]}
&\leq 
C_{p,T,b,\sigma,M}
\left(
\|y-\tilde{y}\|
+
\|\Delta \mbX\|_{p,[0,T]}
 + 
 T\|u-\tilde{u}\|_{L^\infty([0,T],U)}
\right)
\end{align}
for a constant $C_{p,T,b,\sigma,M}\geq 0$.
\end{lemma}
\begin{proof}
First, the RDEs have unique solutions $(Y,Y')\in\sD^p_X$  and $(\widetilde{Y},\widetilde{Y}')\in\sD^p_{\widetilde{X}}$ thanks to Theorem \ref{thm:rdes:existence_unicity}. Let $\alpha=\alpha_{p,b,\sigma}$ be as in Proposition \ref{prop:rdes:error_bound:entire_interval}, and $I=[0,T]$. 
By Lemma \ref{lem:Nalpha<=w(0,T)}, $\alpha N_{\alpha,[0,T]}(\mbX)\leq w_\mbX(0,T)=\|X\|^p_p+\|\bX\|_\frac{p}{2}^\frac{p}{2}
\leq M^p+M^\frac{p}{2}$. By  \eqref{eq:Nalpha<=3CpNalpha_X_and_time} in Corollary \ref{cor:Nalpha:NX_NXtilde_NT:small_intervals}, $N_{\alpha,[0,T]}(w)\leq 5C_p(
N_{\alpha,[0,T]}(\mbX)+N_{\alpha,[0,T]}(\widetilde\mbX)+T/\alpha+1
)$. Thus, $N_{\alpha,[0,T]}(w)\leq C$ for a constant $C$ that depends only on $(p,T,b,\sigma,M)$. Finally, the inequalities \eqref{eq:RDE:DY'+dRY:close:full_interval:continuity} and \eqref{eq:RDE:DY:close:full_interval:continuity} follow from \eqref{eq:RDE:DY'+dRY:close:full_interval} and \eqref{eq:RDE:DY:close:full_interval} for a new constant $C_{p,T,b,\sigma,M}\geq 0$. 
\end{proof}

The continuity of solutions to  \rev{linear} RDEs is proved under the following stronger assumption. 
\begin{assumption}[Stronger regularity of $b$]\label{assumption:b:stronger}
Let $T>0$,  $U\subseteq\R^m$, and $b:[0,T]\times\R^n\times U\to\R^n$. The map $b$ satisfies Assumption \ref{assumption:b:stronger:linear_rde} and is such that 
$\|b(t,x,u)\|+\big\|\rev{\nablaof{x}b}(t,x,u)\big\| \leq C_b$ 
and 
$\big\|\rev{\nablaof{x}b}(t,x,u)-\rev{\nablaof{x}b}(t,\tilde{x},u)\big\|\leq C_b\|x-\tilde{x}\|$ for some constant $C_b\geq 0$, almost every $t\in[0,T]$, all $x,\tilde{x}\in\R^n$, and all $u\in U$.
\end{assumption}

\begin{lemma}[Continuity of the Jacobian flow]\label{lem:rde:linearized:continuous}
Define $p,T,y,\tilde{y},U,u,\tilde{u},b,\sigma,\mbX,\widetilde{\mbX},Y,Y',\widetilde{Y},\widetilde{Y}',M$ as in Lemma \ref{lem:rdes:continuity}, 
and assume that $b$ satisfies  Assumption \ref{assumption:b:stronger}, 
 $(t,x,u)\mapsto \rev{\nablaof{x}b}(t,x,u)$ is Lipschitz in $u$,  \rev{and} 
$\sigma\in C_b^4([0,T]\times\R^n,\R^{n\times d})$. %
Let $v\in C_b^1(\R^n,\R^n)$ and  $(V,V')\in\sD^p_X([0,T],\R^n)$ and $(\rev{\widetilde{V}},\rev{\widetilde{V}}')\in\sD^p_{\widetilde{X}}([0,T],\R^n)$ 
with $V'=\rev{\nablaof{x}\sigma}(\cdot,Y_\cdot)V_\cdot$  and $\widetilde{V}'=\rev{\nablaof{x}\sigma}(\cdot,\widetilde{Y}_\cdot)\widetilde{V}_\cdot$  
be the solutions to the linear RDEs 
{\small
\begin{align*}
V_t &= 
v(y)
+
\int_0^t
\hspace{-3pt}
\rev{\nablaof{x}b}(s,Y_s,u_s)V_s\dd s +
\int_0^t
\hspace{-3pt}
\rev{\nablaof{x}\sigma}(s,Y_s)V_s\dd\mbX_s
\  \text{and}\    
\widetilde{V}_t = 
v(\tilde{y})
+
\int_0^t
\hspace{-3pt}
\rev{\nablaof{x}b}(s,\widetilde{Y}_s,\tilde{u}_s)\widetilde{V}_s\dd s +
\int_0^t
\hspace{-3pt}
\rev{\nablaof{x}\sigma}(s,\widetilde{Y}_s)\widetilde{V}_s\dd\widetilde{\mbX}_s
\end{align*}
}%
where  $t\in[0,T]$.
Then,  for a constant $C_{p,T,b,\sigma,v,M}\geq 0$,
\begin{align*}
\|V'-\widetilde{V}'\|_{p,[0,T]}+
\|R^V-R^{\widetilde{V}}\|_{\frac{p}{2},[0,T]}
&\leq 
C_{p,T,b,\sigma,v,M}
\left(
\|y-\tilde{y}\| 
+
\|\Delta \mbX\|_{p,[0,T]}
 + 
 T\|u-\tilde{u}\|_{L^\infty([0,T],U)}
\right),
\\
\|V-\widetilde{V}\|_{\infty,[0,T]}&\leq 
C_{p,T,b,\sigma,v,M}
\left(
\|y-\tilde{y}\| 
+
\|\Delta \mbX\|_{p,[0,T]}
 + 
 T\|u-\tilde{u}\|_{L^\infty([0,T],U)}
\right).
\end{align*}
\end{lemma}
\begin{proof}
First, the linear RDEs also have unique solutions $(V,V')\in\sD^p_X$ and $(\rev{\widetilde{V}},\rev{\widetilde{V}}')\in\sD^p_{\widetilde{X}}$ thanks to Corollary  \ref{cor:rde:linear:existence_uniqueness}.  
Second, by Lemma \ref{lem:rde:linearized:bounded_solutions},  
there \rev{exist} two constants $C_{p,T,b,\sigma}\geq 1$ and $0<\alpha_{p,b,\sigma}<1$  (importantly, they do not depend on $Y,\widetilde{Y}$) such that 
$$
\|V\|_{\infty,[0,T]}+\|V'\|_{\infty,[0,T]}\leq C_{p,T,b,\sigma}\exp\left(C_{p,T,b,\sigma}N_{\alpha_{p,b,\sigma},[0,T]}(\mbX)\right)\|v(y)\|,
$$ 
and similarly for $\|\widetilde{V}\|_{\infty,[0,T]}+\|\widetilde{V}'\|_{\infty,[0,T]}$. 
By Lemma \ref{lem:Nalpha<=w(0,T)}, $\alpha_{p,b,\sigma} N_{\alpha_{p,b,\sigma},[0,T]}(\mbX)\leq w_\mbX(0,T)=\|X\|^p_p+\|\bX\|_\frac{p}{2}^\frac{p}{2}
\leq M^p+M^\frac{p}{2}$, since $\|\mbX\|_p=\|X\|_p+\|\bX\|_\frac{p}{2}\leq M$.  Thus,
$$
\|V\|_{\infty,[0,T]}+\|V'\|_{\infty,[0,T]}+\|\widetilde{V}\|_{\infty,[0,T]}+\|\widetilde{V}'\|_{\infty,[0,T]}\leq C_{p,T,b,\sigma,v,M}.
$$
Thus, we may assume that $(V,V')$ and $(\widetilde{V},\widetilde{V}')$ solve the nonlinear RDEs
\begin{align*}
V_t &=v(Y_0) +  
\int_0^t
\hat{b}(s,(Y_s,V_s),u_s)\dd s
+  
\int_0^t
\hat\sigma(s,(Y_s,V_s))\dd\mbX_s,
\quad t\in[0,T]
\end{align*}  
(and similarly for $\widetilde{V}$) for some bounded coefficients $(\hat{b},\hat\sigma)$, where $\hat{b}:[0,T]\times\R^{2n}\times U\to\R^n$ satisfies Assumption \ref{assumption:b} and is Lipschitz in $u$, and $\hat\sigma\in C_b^3([0,T]\times\R^{2n},\R^{n\times d})$. Thus, the pair $((Y,V),(Y',V'))$ solves the RDE
\begin{align*}
\begin{bmatrix}
Y_t \\ V_t
\end{bmatrix} &=
\begin{bmatrix}
y \\ v(y)
\end{bmatrix} +  
\int_0^t
\begin{bmatrix}
b(s,Y_s,u_s) \\ \hat{b}(s,(Y_s,V_s),u_s)
\end{bmatrix}
\dd s
+  
\int_0^t
\begin{bmatrix}
\sigma(s,Y_s) \\ 
\hat\sigma(s,(Y_s,V_s))
\end{bmatrix}
\dd\mbX_s,
\end{align*}  
and similarly for $(\widetilde{Y},\widetilde{V})$. Then using \eqref{eq:RDE:DY'+dRY:close:full_interval:continuity} and  \eqref{eq:RDE:DY:close:full_interval:continuity} in Lemma \ref{lem:rdes:continuity},
\begin{align*} 
\|(Y,V)-(\widetilde{Y},\widetilde{V})\|_{\infty,[0,T]} 
&\leq
C_{p,T,b,\sigma,v,M}(
\|(y,v(y))-(\tilde{y},v(\tilde{y}))\|
+
\|\Delta\mbX\|_{p,[0,T]}
+
T\|u-\tilde{u}\|_{L^\infty([0,T],U)}
),
\end{align*}
and similarly for $\|(Y',V')-(\widetilde{Y}',\widetilde{V}')\|_{p}+\|R^{(Y,V)}-R^{(\widetilde{Y},\widetilde{V})}\|_{\frac{p}{2}}$,
from which we deduce the desired result.
\end{proof}

\begin{theorem}[Integrability of pathwise solutions to random RDEs]\label{thm:rdes:integrable} 
Let $T>0$, $\rho\in[1,\frac{3}{2})$, $p\in(2\rho,3)$,  $(\Omega,\F,\Prob)$ be a probability space, $B$ be a centered, continuous, $\R^d$-valued Gaussian process with independent components satisfying Assumption \ref{assum:Gaussian_lift} for $\rho$, and $\mbB$ be the associated enhanced Gaussian process in Theorem \ref{thm:gaussian_rough_paths}. 
Let $\ell\geq1$,  $y,\tilde{y}\in L^\ell(\Omega,\R^n)$, $U\subseteq\R^m$, 
$u\in L^\infty([0,T],U)$, 
$b:[0,T]\times\R^n\times U\to\R^n$ satisfy Assumption \ref{assumption:b},  
$\sigma\in C_b^3([0,T]\times\R^n,\R^{n\times d})$,  
and for almost every $\omega\in\Omega$, define $(Y(\omega),Y'(\omega))\in\sD^p_{B(\omega)}([0,T],\R^n)$ with $Y'(\omega)=\sigma(\cdot,Y_\cdot(\omega))$ as the solution to the nonlinear RDE
\begin{align*}
Y_t(\omega)&=y(\omega)+\int_0^tb(s,Y_s(\omega),u_s)\dd s+\int_0^t\sigma(s,Y_s(\omega))\dd\mbB_s(\omega),
\quad\ \ t\in[0,T].
\end{align*} 
Then, $Y\in L^\ell(\Omega,C([0,T],\R^n))$. 

Moreover, let $v\in C_b^1(\R^n,\R^n)$,   
 assume that 
$b$  satisfies  Assumption \ref{assumption:b:stronger} and 
$\sigma\in C_b^4$, and for almost every $\omega\in\Omega$, define $(V(\omega),V'(\omega))\in\sD^p_{B(\omega)}([0,T],\R^n)$ with $V'(\omega)=\rev{\nablaof{x}\sigma}(\cdot,Y_\cdot(\omega))V_\cdot(\omega)$ as the solution to the linear RDE
\begin{align*}
V_t(\omega) &=v(\tilde{y}(\omega)) +  
\int_0^t
\rev{\nablaof{x}b}(s,Y_s(\omega),u_s)V_s(\omega)\dd s
+  
\int_0^t
\rev{\nablaof{x}\sigma}(s,Y_s(\omega))V_s(\omega)\dd\mbB_s(\omega),
\quad t\in[0,T].
\end{align*}   
Then, $V\in L^\ell(\Omega,C([0,T],\R^n))$.
\end{theorem}
\begin{proof}
We first prove that $Y\in L^\ell(\Omega,C([0,T],\R^n))$. %
Note that the assumption that $b$ is Lipschitz in the control input $u$ is not needed, as Lemma \ref{lem:rdes:continuity} and other  results used below still hold without this assumption since the control $u=\tilde{u}$ is fixed.

\textit{1) Measurability}: As is standard in rough path theory, we express the solution map $Y$ as a composition of the measurable map $(y,\mbB):\Omega\to\R^n\times\sC^p_g$ and the pathwise solution to the RDE \eqref{eq:RDE}, which is continuous with respect to the initial condition and the driving signal.  First, 
 $(\sC^p_g,d_p)$ is a (complete) metric space (with  Borel sets defined by its metric topology), where  
$d_p(\mbX,\widetilde\mbX):=\|\Delta X_0\|+\|\Delta\mbX\|_p$ denotes the inhomogeneous rough path metric. 
By Theorem \ref{thm:gaussian_rough_paths}, the map
\begin{align*}
\mbB:
\, 
(\Omega,\F,\Prob)&\to
(\sC^p_g([0,T],\R^d),d_p),
\  
\omega\mapsto \mbB(\omega)=(B(\omega),\bB(\omega))
\end{align*}
is measurable. Also, by Theorem \ref{thm:rdes:existence_unicity}, for any  initial conditions $\bar{y}\in\R^n$ and    \rev{rough} path $\mbX\in \sC_g^p$, the RDE   \eqref{eq:RDE} has a unique solution, denoted by $\big(\widehat{Y}_{(\bar{y},\mbX)},(\widehat{Y}_{(\bar{y},\mbX}))'\big)\in\sD^p_X$, so the map 
\begin{align*} 
\widehat{Y}_{(\cdot,\cdot)}:
\,
 \left(
\R^n\times \sC^p_g([0,T],\R^d),
\|\cdot\|\oplus d_p
\right)
&\to  
\left(C([0,T],\R^n),\|\cdot\|_\infty\right),
\ 
(\bar{y},\mbX)\mapsto \widehat{Y}_{(\bar{y},\mbX)}
\end{align*} 
is well-defined, where $\|\cdot\|\oplus d_p$ is the product metric.  
Moreover, if $\bar{y},\tilde{y}\in\R^n$ and $\mbX,\widetilde{\mbX}\in \sC_g^p$ satisfy  $\|\mbX\|_p,\|\widetilde{\mbX}\|_p\leq M$ for some $M\geq 0$, then  $\big\|\widehat{Y}_{(\bar{y},\mbX)}-\widehat{Y}_{(\tilde{y},\widetilde\mbX)}\big\|_\infty\leq 
C_{p,T,b,\sigma,M}
\left(
\|\bar{y}-\tilde{y}\| 
+
\|\Delta\mbX\|_p
\right)$ by \eqref{eq:RDE:DY:close:full_interval:continuity} in Lemma \ref{lem:rdes:continuity}, 
so the map $\widehat{Y}_{(\cdot,\cdot)}$ is continuous. 
Thus,  the map
\begin{align*}
Y:%
\,
(\Omega,\F,\Prob)&\to
\left(C([0,T],\R^n),\|\cdot\|_\infty\right),
\ 
\omega\mapsto \widehat{Y}_{(y(\omega),\mbB(\omega))} 
\end{align*}
is measurable, since it is the composition of the measurable map $(y(\cdot),\mbB(\cdot))$ and the  continuous map $\widehat{Y}_{(\cdot,\cdot)}$.

\textit{2) Integrability}: 
By Lemma \ref{lem:rdes:bound:entire_interval},  there exist constants $C_{p,T,b,\sigma}\geq 1$ and $0<\alpha_{p,b,\sigma}<1$ such that  
\begin{align*}
\|Y\|_{\infty,[0,T]}
&\mathop{\leq}^{\eqref{eq:RDE:Y:bound:full_interval}}
C_{p,T,b,\sigma}\exp\left(
C_{p,T,b,\sigma}N_{\alpha_{p,b,\sigma},[0,T]}(\mbB)
\right)+\|y\|
\end{align*}
almost surely, 
where $N_{\alpha,I}(\mbB)$ is defined in Definition \ref{def:Nalpha}. $\E[\|Y\|_\infty^\ell]<\infty$ then follows from Theorem \ref{thm:gaussian_rough_paths}.

Thus, $Y\in L^\ell(\Omega,C([0,T],\R^n))$. 

The proof that $V\in L^\ell(\Omega,C([0,T],\R^n))$ is identical, using Corollary \ref{cor:rde:linear:existence_uniqueness} (instead of Theorem \ref{thm:rdes:existence_unicity}), 
Lemma \ref{lem:rde:linearized:continuous} (instead of Lemma \ref{lem:rdes:continuity}) 
and Lemma \ref{lem:rde:linearized:bounded_solutions} (instead of Lemma \ref{lem:rdes:bound:entire_interval}).
\end{proof}

\begin{corollary}[Integrable error bound]\label{cor:rdes:integrable:error_bound} 
Define 
$T,p,\Omega,\F,\Prob,\mbB,\ell,y,\tilde{y},U,\sigma$
as in Theorem \ref{thm:rdes:integrable} with $\ell>1$ and $\sigma\in C_b^3$,   
let $u,\tilde{u}\in L^\infty([0,T],U)$, 
 $b:[0,T]\times\R^n\times U\to\R^n,(t,x,u)\mapsto b(t,x,u)$ satisfy Assumption \ref{assumption:b} and be Lipschitz in $u$, 
and  $Y,\widetilde{Y}\in L^\ell(\Omega,C([0,T],\R^n))$ be the pathwise solutions to the random  RDEs
\begin{align*}
Y_t&=y+\int_0^tb(s,Y_s,u_s)\dd s+\int_0^t\sigma(s,Y_s)\dd\mbB_s,
\quad
\widetilde{Y}_t=\tilde{y}+\int_0^tb(s,\widetilde{Y}_s,\tilde{u}_s)\dd s+\int_0^t\sigma(s,\widetilde{Y}_s)\dd\mbB_s,
\quad t\in[0,T].
\end{align*} 
as in Theorem \ref{thm:rdes:integrable}. 
Then, there exist constants $C:=C_{p,\ell,T,b,\sigma}\geq 1$ and  $0<\alpha:=\alpha_{p,b,\sigma}<1$ such that  
\begin{align*}
\E\big[
\|Y-\widetilde{Y}\|_{\infty,[0,T]} 
\big]
&\leq
C
\E\rev{\left[
\exp\left(
CN_{\alpha,[0,T]}(\mbB)
\right)\right]}
\big(
\E\big[\|y-\tilde{y}\|^\ell\big]^\frac{1}{\ell}
+
T\|u-\tilde{u}\|_{L^\infty([0,T],U)}
\big),
\end{align*}
where $\E[\exp(CN_{\alpha,[0,T]}(\mbB))]<\infty$. 
\end{corollary}
\begin{proof}
First, $\E[\exp(DN_{\alpha,[0,T]}(\mbB))]<\infty$ for any $D\geq 0$ and $\alpha>0$ by Theorem \ref{thm:gaussian_rough_paths}. 
Second, H\"older's inequality gives  $\E[\exp(DN_\alpha)\|\Delta y\|]\leq \E[\exp(\frac{\ell}{\ell-1} DN_\alpha)]^\frac{\ell-1}{\ell}\E[\|\Delta y\|^\ell]^\frac{1}{\ell}$ and similarly for $\E[\exp(DN_\alpha)T\|\Delta u\|_{L^\infty}]$. 
Finally, the desired inequality  %
follows from \eqref{eq:RDE:DY:close:full_interval} %
in Proposition \ref{prop:rdes:error_bound:entire_interval} and  
 $N_{\alpha,[0,T]}(w)\leq C_p(
N_{\alpha,[0,T]}(\mbB)+T/\alpha+1
)$ by \eqref{eq:Nalpha<=3CpNalpha_X_and_time} in Corollary \ref{cor:Nalpha:NX_NXtilde_NT:small_intervals}.
\end{proof}

\section{The  Pontryagin Maximum Principle (\pmp)}\label{sec:pmp}

We now prove \pmp (Theorem \ref{thm:pmp}). First, in Section \ref{sec:linearization}, we consider particular variations of solutions to RDEs, called \textit{needle-like variations} \cite{Pontryagin1986,LeeMarkus1967,Agrachev2004,Bonnard2005}, and show that these variations can be approximated well using the solution to a linearized RDE (Lemma \ref{lem:needle_like_error:etas}) along the optimal solution $(x,u)$ to \ocp. The use of needle-like variations is a standard method for deriving a PMP that can handle the control constraints $u_t\in U$. 
Finally, in Section \ref{sec:pmp_proof}, we state the main assumptions for  \pmp and prove the result.

\subsection{Needle-like variations}\label{sec:linearization}

The  needle-like variations rely on the concept of a Lebesgue point. %

\begin{definition}[Lebesgue point]\label{def:lebesgue_point}
Let %
$T>0$, \revproof{$\ell\geq 1$,} %
$U\subseteq\R^m$, 
$u\in L^\infty([0,T],U)$, 
$b:[0,T]\times\R^n\times U\to\R^n$ satisfy  Assumption \ref{assumption:b}, %
\revproof{$Y\in L^\ell(\Omega,C([0,T],\R^n))$, and $V\subset U$ be a countable dense subset}.  
We say that $t_1\in[0,T]$ is a Lebesgue point %
if 
\revproof{almost surely,}
\begin{equation}
\label{eq:lebesgue_point:lim}
\lim_{h\to0}\frac{1}{h}\int_{t_1}^{t_1+h}b(t,Y_t,u_t)\dd t= b(t_1,Y_{t_1},u_{t_1})
\ \ \revproof{\text{and}\ \  
\lim_{h\to0}\frac{1}{h}\int_{t_1}^{t_1+h}b(t,Y_t,v)\dd t= b(t_1,Y_{t_1},v)
\ \text{for all $v\in V$.}}
\end{equation}
\end{definition}

\revproof{Given a map $f:\R\to\R^n$, we write $f=o(h)$ if $\lim_{h\to 0}\tfrac1hf(h)=0$. By Definition \ref{def:lebesgue_point}, at a Lebesgue point,  a.s., $ \int_{t_1}^{t_1+h}b(t,Y_t,u_t)\dd t- h\, b(t_1,Y_{t_1},u_{t_1})=o(h)$ and $
\int_{t_1}^{t_1+h}b(t,Y_t,v)\dd t-hb(t_1,Y_{t_1},v)=o(h)$ for all $v\in V$.
}

\revproof{
Almost every $t_1\in[0,T]$ is a Lebesgue point. 
Indeed, let $V=\{v_1,v_2,\dots\}\subset U$ be a dense subset. 
Almost surely, $t\mapsto b(t,Y_t(\omega),u_t)$ and each $t\mapsto b(t,Y_t(\omega),v_n)$ is in $L^1([0,T],\R^n)$.  
By the Lebesgue differentiation theorem, there exist sets $I_u(\omega),I_1(\omega),I_2(\omega),\dots\subseteq[0,T]$ of full measure such that $\lim_{h\to0}\frac1h\int_{t_1}^{t_1+h}b(t,Y_t(\omega),u_t)\dd t= b(t_1,Y_{t_1}(\omega),u_{t_1})$  $\forall t_1\in I_u(\omega)$ and 
$\lim_{h\to0}\frac1h\int_{t_1}^{t_1+h}b(t,Y_t(\omega),v_n)\dd t= b(t_1,Y_{t_1}(\omega),v_n)$ $\forall t_1\in I_n(\omega)$. 
The set $I(\omega):=I_u(\omega)\cap I_1(\omega)\cap\dots\subseteq[0,T]$ has full measure,  
and any $t_1\in I(\omega)$ satisfies \eqref{eq:lebesgue_point:lim}. 
So, by Fubini's theorem applied to the indicator function of the set $A=\{(t_1,\omega)\in[0,T]\times\Omega:\smash{\lim\limits_{h\to 0}}\frac1h\int_{t_1}^{t_1+h}b(t,Y_t(\omega),u_t)\dd t \neq b(t_1,Y_{t_1}(\omega),u_{t_1}) 
\ \text{or} \  
\smash{\lim\limits_{h\to 0}}\frac1h\int_{t_1}^{t_1+h}b(t,Y_t(\omega),v_n)\dd t \neq b(t_1,Y_{t_1}(\omega),v_n) \ \text{for some } n=1,2,\dots\}$, almost every $t\in[0,T]$ is a Lebesgue time\footnote{Fubini's theorem gives  
$0=
\int_\Omega\int_0^T\mathbf{1}_A(t,\omega)\dd t\dd\Prob(\omega)=\int_0^T\int_\Omega\mathbf{1}_A(t,\omega)\dd\Prob(\omega)\dd t=\int_0^T\Prob(A_t)\dd t$ with $A_t(\omega)=\{\omega\in\Omega:(t,\omega)\in A\}$. Since $\Prob(A_t)\geq 0$, we get $\Prob(A_t)=0$ for almost every $t_1\in[0,T]$, so almost every $t_1\in[0,T]$ is a Lebesgue time.}. 
}

\revproof{
If $b$ is also Lipschitz in $u$, then for any Lebesgue point $t_1$ and $\bar{u}_1\in U$, almost surely,
\begin{align}\label{eq:lebesgue_point:limat_baru1}
\int_{t_1}^{t_1+h}
b(s,Y_s,\bar{u}_1)
\dd s
&=
h\, b(t_1,Y_{t_1},\bar{u}_1)+ o(h).
\end{align}
Indeed, with $a_h:=\frac1h\int_{t_1}^{t_1+h}
b(s,Y_s,\bar{u}_1)
\dd s
- b(t_1,Y_{t_1},\bar{u}_1)$ and $(v_n)\subset V$ a sequence converging to $\bar u_1$,  we have
$
\mathop{\lim}_{h\to 0}
\|
a_h
\|
\leq
\lim_{n\to\infty}
\mathop{\lim\sup}_{h\to 0}
\|
\int_{t_1}^{t_1+h}
\frac1h b(s,Y_s,v_n)
\dd s
-
b(t_1,Y_{t_1},v_n)
\|
+
2C\|\bar u_1-v_n\|
=
0.
$
}

\revproofnote{Definition \ref{def:lebesgue_point} previously defined a Lebesgue point $t_1$ to only satisfy  $
\lim\limits_{h\to0}\frac{1}{h}\int_{t_1}^{t_1+h}b(t,Y_t,u_t)\dd t= b(t_1,Y_{t_1},u_{t_1})
$. However, this property that holds only at $u$ is not sufficient for the proof of Proposition  \ref{prop:linear_variation}, since \eqref{eq:Delta(t):term1} evaluates $b$ at $\bar{u}$ inside the integral, and not at $u$. 
\\
\indent As a counter-example, consider $T=1$, $U=[0,1]$, $x_0=0$, $\sigma=0$, 
$$
b(t,x,v)=v\textbf{1}_{\{t_1\}}(t),
\qquad u\equiv 0,
$$
and the needle-like variation $\bar{u}_1=1$ over $[t_1,T]$ for some  $t_1\in[0,T]$, which would be a Lebesgue time by the old definition. 
Then, $$Y^{\pi_1}=Y=0, 
\quad\text{and}\quad 
V^{\pi_1}=V^{\pi_1}_{t_1}=b(t_1,Y_{t_1},\bar{u}_1)-b(t_1,Y_{t_1},u_{t_1})=1,
$$
but $\|Y^{\pi_1}-Y-\eta_1V^{\pi_1}\|_{\infty,[t_1+\eta_1,T]}=\eta_1$.
}
\begin{proposition}
[Needle-like variations and linearized RDEs]
\label{prop:linear_variation}  
Let $p\in[2,3)$, $T>0$, $y\in\R^n$, $U\subseteq\R^m$, 
$u\in L^\infty([0,T],U)$, 
$b:[0,T]\times\R^n\times U\to\R^n,(t,x,u)\mapsto b(t,x,u)$   satisfy Assumption \ref{assumption:b:stronger}  and be Lipschitz in $u$, $\sigma\in C_b^4([0,T]\times\R^n,\R^{n\revproof{\times d}})$, and 
$\mbX\in\rev{\sC^p}([0,T],\R^d)$ be a \rev{rough} path. 
Let $(Y,Y')\in\sD^p_X$  and $(Y^{\pi_1},(Y^{\pi_1})')\in\sD^p_X$ with $Y'=\sigma(\cdot,Y_\cdot)$ and 
 $(Y^{\pi_1})'=\sigma(\cdot,Y^{\pi_1}_\cdot)$ solve the RDEs 
\begin{align*}
Y_t&=y+\int_0^tb(s,Y_s,u_s)\dd s+\int_0^t\sigma(s,Y_s)\dd\mbX_s,
\ \, 
Y^{\pi_1}_t=y+\int_0^tb(s,Y^{\pi_1}_s,u^{\pi_1}_s)\dd s+\int_0^t\sigma(s,Y^{\pi_1}_s)\dd\mbX_s,
\ \  
t\in[0,T],
\end{align*}
\revproof{where} %
$t_1\in[0,T]$ \revproof{is a Lebesgue point}, %
$\eta_1\in[0,T-t_1]$, %
$\bar{u}_1\in U$, %
\revproof{and} the needle-like variation $\pi_1=(t_1,\eta_1,\bar{u}_1)$ of $u$ \revproof{is defined} by 
 $$
 u^{\pi_1}_t=\begin{cases}
 \bar{u}_1\quad&\text{if }t\in[t_1,t_1+\eta_1],
 \\
 u_t&\text{otherwise}.
 \end{cases}
 $$ 
 \revproof{Let }
$(V^{\pi_1},(V^{\pi_1})')\in\sD^p_X$ with $(V^{\pi_1})'=\rev{\nablaof{x}\sigma}(\cdot,Y_\cdot)V^{\pi_1}_\cdot$ solve the linear RDE
\begin{subequations}\label{eq:rde:linear_variation}
\begin{align}
V^{\pi_1}_t&=b(t_1,Y_{t_1},\bar{u}_1)-b(t_1,Y_{t_1},u_{t_1}),
&&t\in[0,t_1],
\\
V^{\pi_1}_t &= 
V^{\pi_1}_{t_1}
+
\int_{t_1}^t
\rev{\nablaof{x}b}(s,Y_s,u_s)V^{\pi_1}_s\dd s +
\int_{t_1}^t
\rev{\nablaof{x}\sigma}(s,Y_s)V^{\pi_1}_s\dd\mbX_s,
 &&t\in[t_1,T].
\end{align}
\end{subequations}
Then, there exist constants \rev{$C>0$ and $0<\alpha<1$ that depend on $(p,T,b,\sigma,u,\bar{u}_1)$} such that 
\begin{align}
\label{eq:needle_like_deltasols}
\|Y^{\pi_1}-Y\|_{\infty,[0,T]}
&\leq 
\rev{C}\exp\left(
\rev{C}N_{\rev{\alpha},[0,T]}(\mbX)
\right)
\eta_1,
\\
\label{eq:needle_like_deltasols_variation}
\|Y^{\pi_1}-Y-\eta_1V^{\pi_1}\|_{\infty,[t_1\revproof{+\eta_1},T]}&\leq \rev{C}\exp\left(
\rev{C}N_{\rev{\alpha},[0,T]}(\mbX)
\right)
\revproof{o(\eta_1)},
\end{align}
where %
$N_{\rev{\alpha},[0,T]}(\mbX)$ is defined in Definition \ref{def:Nalpha}. %
\end{proposition} 
\revproofnote{Proposition \ref{prop:linear_variation} previously stated 
$$
\|Y^{\pi_1}-Y-\eta_1V^{\pi_1}\|_{\infty,[t_1+\eta_1,T]}
\mathop{\leq}^{\text{(Incorrect)}} \rev{C}\exp\left(
\rev{C}N_{\rev{\alpha},[0,T]}(\mbX)
\right)
\eta_1^2.
$$
As a counter-example, consider $T=1$, $U=[0,1]$, $x_0=0$, $\sigma=0$, $t_1\in[0,T]$ and $$
b(t,x,v)=\sqrt{|t-t_1|}v,
\qquad u\equiv0,
\qquad \bar{u}_1=1.
$$
Then, $t_1$ is a Lebesgue point and $V^{\pi_1}=0$, but
$$
Y^{\pi_1}_{t_1+\eta_1}-Y_{t_1+\eta_1}-\eta_10=\int_{t_1}^{t_1+\eta_1}\sqrt{|t-t_1|}1\dd t=\frac23\eta^\frac32,
$$
which cannot be bounded by $C\eta_1^2$ as $\eta_1\to 0$. Yet, $\eta^\frac32=o(\eta_1)$, which is sufficient for the proof.}
\begin{proof}
There exists unique solutions $Y,Y^{\pi_1},V^{\pi_1}$ to the RDEs  
thanks to Theorem \ref{thm:rdes:existence_unicity} and  Corollary \ref{cor:rde:linear:existence_uniqueness}.

\revproof{Define $N_{\alpha,[0,T]}(w)$ and $w$ as }in Definition \ref{def:Nalpha} and Corollary \ref{cor:Nalpha:NX_NXtilde_NT:small_intervals} \rev{(with $\widetilde{\mbX}=0$), respectively.} 
The \revproof{first} inequality \eqref{eq:needle_like_deltasols} follows from \eqref{eq:RDE:DY:close:full_interval:u_subinterval} in Proposition \ref{prop:rdes:error_bound:entire_interval}  and 
 $N_{\alpha,[0,T]}(w)\leq C_p(
N_{\alpha,[0,T]}(\mbX)+T/\alpha+1
)$ by \eqref{eq:Nalpha<=3CpNalpha_X_and_time} in Corollary \ref{cor:Nalpha:NX_NXtilde_NT:small_intervals}\revproof{.}%

In the remainder of the proof, we show the inequality  \eqref{eq:needle_like_deltasols_variation}. 
For conciseness, we denote $\Delta Y=Y^{\pi_1}-Y$,  $\Delta R^Y=R^{\Delta Y}=\rev{R^{Y^{\pi_1}}-R^Y}$, $(b(x,u),\sigma(x))$ for $(b(t,x,u),\sigma(t,x))$, and similarly for derivatives, and $$
\Delta:=Y^{\pi_1}-Y
-
\eta_1V^{\pi_1}.
$$ 
Let $t\geq t_1+\eta_1$.  We have
\begin{align*}
\Delta Y_t
&=
\int_{t_1}^t
b(Y^{\pi_1}_s,u^{\pi_1}_s)
\dd s
-
\int_{t_1}^t
b(Y_s,u_s)
\dd s
+
\int_{t_1}^t
(
\sigma(Y^{\pi_1}_s)-\sigma(Y_s)
)\dd\mbX_s
\\
&=
\int_{t_1}^{t_1+\eta_1}
b(Y^{\pi_1}_s,\bar{u}_1)
\dd s
+
\int_{t_1+\eta_1}^t
b(Y^{\pi_1}_s,u_s)
\dd s
-
\int_{t_1}^t
b(Y_s,u_s)
\dd s
+
\int_{t_1}^t
(
\sigma(Y^{\pi_1}_s)-\sigma(Y_s)
)
\dd\mbX_s
\\
&=
\int_{t_1}^{t_1+\eta_1}
(
b(Y^{\pi_1}_s,\bar{u}_1)
-
b(Y^{\pi_1}_s,u_s)
)
\dd s
+
\int_{t_1}^t
(
b(Y^{\pi_1}_s,u_s)
-
b(Y_s,u_s)
)
\dd s
+
\int_{t_1}^t
(
\sigma(Y^{\pi_1}_s)-\sigma(Y_s)
)
\dd\mbX_s,
\end{align*} 
so that, 
since $V^{\pi_1}_t = V^{\pi_1}_{t_1}+ 
\int_{t_1}^t\rev{\nablaof{x}b}(Y_s,u_s)V^{\pi_1}_s\dd s+
\int_{t_1}^t\rev{\nablaof{x}\sigma}(Y_s)V^{\pi_1}_s\dd\mbX_s$,
\begin{align*}
\Delta_t
&=
\int_{t_1}^{t_1+\eta_1}
(
b(Y^{\pi_1}_s,\bar{u}_1)
-
b(Y^{\pi_1}_s,u_s)
)
\dd s
-\eta_1V^{\pi_1}_{t_1}
+ 
\int_{t_1}^t
\left(
b(Y^{\pi_1}_s,u_s)
-
b(Y_s,u_s)
-
\eta_1\rev{\nablaof{x}b}(Y_s,u_s)V^{\pi_1}_s
\right)
\dd s
\\
&\qquad+
\int_{t_1}^t
\left(
\sigma(Y^{\pi_1}_s)-\sigma(Y_s)
-
\eta_1\rev{\nablaof{x}\sigma}(Y_s)V^{\pi_1}_s
\right)\dd\mbX_s.
\end{align*}
Next, by Taylor's Theorem \cite{Folland1990},
$b(Y^{\pi_1},u)
-
b(Y,u)
=
\rev{\nablaof{x}b}(Y,u)
\Delta Y
+
\int_0^1
\big(
\rev{\nablaof{x}b}(Y+\theta\Delta Y,u)-\rev{\nablaof{x}b}(Y,u)
\big)
\dd\theta
\Delta Y$, 
and similarly for $\sigma(Y^{\pi_1})
-
\sigma(Y)$ by the mean value theorem, so 
\begin{align}\label{eq:Delta(t):full}
\Delta_t
&=
\int_{t_1}^{t_1+\eta_1}
(
b(Y^{\pi_1}_s,\bar{u}_1)
-
b(Y^{\pi_1}_s,u_s)
)
\dd s
-\eta_1V^{\pi_1}_{t_1}
+
\int_{t_1}^t
\int_0^1
\bigg(
\rev{\nablaof{x}b}(Y_s+\theta\Delta Y_s,u_s)-\rev{\nablaof{x}b}(Y_s,u_s)
\bigg)
\dd\theta
\Delta Y_s
\dd s
\\
\nonumber
&
\quad+
\int_{t_1}^t
\int_0^1\rev{\nablaof{x}^2\sigma}(Y_s+\theta\Delta Y_s)(1-\theta)\dd\theta
\Delta Y^{\otimes 2}_s
\dd\mbX_s 
+
\left(
\int_{t_1}^t
\rev{\nablaof{x}b}(Y_s,u_s)
\rev{\Delta_s}
\dd s+
\int_{t_1}^t
\rev{\nablaof{x}\sigma}(Y_s)
\Delta_s
\dd\mbX_s
\right).
\end{align}
Next, \revproof{define $w$ as in Corollary \ref{cor:Nalpha:NX_NXtilde_NT:small_intervals}. W}e show that the first three terms are $o\revproof{_w}(\eta_1)$  in the following sense:
\begin{align}\label{eq:needlelike:o(eta1)}
A=o\revproof{_w}(\eta_1)\quad\text{if}\quad 
\|A\|&\leq 
C_{p,b,\sigma\rev{,u,\bar{u}_1}}
\exp\left(
C_{p,b,\sigma}
N_{\alpha_{p,b,\sigma},[0,T]}(w)
\right)
\revproof{o(\eta_1)}.
\end{align}
First,
since $t_1$ is a Lebesgue point  and $Y^{\pi_1}_{t_1}=Y_{t_1}$,
\begin{align}
\int_{t_1}^{t_1+\eta_1}
(
b(Y^{\pi_1}_s,\bar{u}_1)
-
b(Y^{\pi_1}_s,u_s)
)
\dd s
-\eta_1V^{\pi_1}_{t_1}
&=
\eta_1(
b(Y_{t_1},\bar{u}_1)
-
b(Y_{t_1},u_{t_1})
)-\eta_1V^{\pi_1}_{t_1}
+
o\revproof{_w}(\eta_1)
=
o\revproof{_w}(\eta_1).
\label{eq:Delta(t):term1}
\end{align}

\revproof{
Indeed,
$
b(Y^{\pi_1}_s,\bar{u}_1)
-
b(Y^{\pi_1}_s,u_s)
=
b(Y^{\pi_1}_s,\bar{u}_1)
-
b(Y_s,\bar{u}_1)
+
b(Y_s,\bar{u}_1)
-
b(Y_s,u_s)
+
b(Y_s,u_s)
-
b(Y^{\pi_1}_s,u_s)$, so
\begin{align*}
\int_{t_1}^{t_1+\eta_1}
\|b(Y^{\pi_1}_s,\bar{u}_1)
-
b(Y_s,\bar{u}_1)\|
\dd s
&\leq
C
\|Y^{\pi_1}-Y\|_\infty\eta_1
\mathop{\leq}^{\eqref{eq:RDE:DY:close:full_interval:u_subinterval}} 
C\exp\left(
CN_{\rev{\alpha},[0,T]}(w)
\right)
\eta_1^2
=
o_w(\eta_1),
\\
\int_{t_1}^{t_1+\eta_1}
\big(
b(Y_s,\bar{u}_1)
-
b(Y_s,u_s)
\big)
\dd s
&\mathop{=}^{(\eqref{eq:lebesgue_point:lim},\eqref{eq:lebesgue_point:limat_baru1})}
\eta_1 b(Y_{t_1},\bar{u}_1)
-
\eta_1 b(Y_{t_1},u_{t_1})+o(\eta_1)
\\
\int_{t_1}^{t_1+\eta_1}
\|
b(Y_s,u_s)
-
b(Y^{\pi_1}_s,u_s)
\|
\dd s
&\mathop{\leq}^{\eqref{eq:RDE:DY:close:full_interval:u_subinterval}} 
o_w(\eta_1).
\end{align*}
}
\hspace{-2mm}
Second, 
$\rev{\nablaof{x}b}(x,u)$ is Lipschitz in $x$ and $\|\Delta Y\|_\infty^2=o\revproof{_w}(\eta_1)$ by \eqref{eq:RDE:DY:close:full_interval:u_subinterval}, so  
\begin{equation}
\label{eq:Delta(t):term2}
\int_{t_1}^t
\int_0^1
\bigg(
\rev{\nablaof{x}b}(Y_s+\theta\Delta Y_s,u_s)-\rev{\nablaof{x}b}(Y_s,u_s)
\bigg)
\dd\theta
\Delta Y_s
\dd s
=o\revproof{_w}(\eta_1).
\end{equation}
Third, we bound the last rough integral, noting that
\begin{equation}\label{eq:RDE:DY'+dRY:close:full_interval:needle_like}
\|\Delta Y'\|_p
+
\|\Delta R^Y\|_\frac{p}{2}
+
\|\Delta Y\|_\infty
+
\|\Delta Y\|_p
\mathop{\leq}^{\eqref{eq:RDE:DY'+dRY:close:full_interval:u_subinterval},\eqref{eq:RDE:DY:close:full_interval:u_subinterval},
\rev{\eqref{eq:controlled_path:Y-Ytilde_p}},  \eqref{eq:||X||_p<=exp(N_alpha^p)}}
C_{p,b,\sigma\rev{,u,\bar{u}_1}}
\exp\left(
C_{p,b,\sigma}N_{\alpha_{p,b,\sigma},[0,T]}(w)
\right)
\eta_1.
\end{equation} 
Next, we define  $W_t = \Delta Y_t\otimes\Delta Y_t$ and $Z_t=
\int_0^1\rev{\nablaof{x}^2\sigma}(Y_t+\theta\Delta Y_t)(1-\theta)\dd\theta$.
By Lemma \ref{lem:control_path:product}, $(W,W')\in\sD^p_X$ is a controlled path with Gubinelli derivative $ 
W'=\rev{\Delta Y\otimes \Delta Y'+\Delta Y'\otimes \Delta Y}$,
\begin{subequations}\label{eq:ineqs:W}
\begin{align}
\|W'\|_p
&\mathop{\leq}^{\eqref{lem:control_path:product:(YZ)'_p}}
C_p(\|\Delta Y\|_\infty+\|\Delta Y\|_p+\|\Delta Y'\|_p+\|\Delta Y'\|_\infty)^2
\mathop{=}^{\eqref{eq:RDE:DY'+dRY:close:full_interval:needle_like}}
o\revproof{_w}(\eta_1),
\\
\|R^{W}\|_\frac{p}{2}
&\mathop{\leq}^{\eqref{lem:control_path:product:R^YZ_p/2}}
C_p(\|\Delta Y\|_\infty\|R^{\Delta Y}\|_\frac{p}{2}+\|\Delta Y\|_p^2)
\mathop{=}^{\eqref{eq:RDE:DY'+dRY:close:full_interval:needle_like}}
o\revproof{_w}(\eta_1),
\\
\|W\|_p
&\mathop{\leq}^{\eqref{eq:controlled_path:pvar_norm}}
C_p(\|W'\|_\infty\|X\|_p+\|R^W\|_\frac{p}{2})
\mathop{=}^{\eqref{eq:RDE:DY'+dRY:close:full_interval:needle_like},\eqref{eq:||X||_p<=exp(N_alpha^p)}}o\revproof{_w}(\eta_1),
\\
\|W\|_\infty
&\mathop{=}^{\eqref{eq:RDE:DY:close:full_interval:u_subinterval}}
o\revproof{_w}(\eta_1).
\end{align}
\end{subequations}
By Lemma \ref{lem:rough_path:sigma(.,Y):controlled}, $(Z,Z')\in\sD^p_X$ is also a controlled path, since $\rev{\nablaof{x}^2\sigma}\in C^2_b$ as $\sigma\in C^4_b$. 
By Lemma \ref{lem:control_path:product}, $(ZW,(ZW)')\in\sD^p_X$ is a controlled path with Gubinelli derivative  $
(ZW)'=WZ'+ZW'$, and   
\begin{subequations}\label{eq:ineqs:ZW}
\begin{align}
\|(ZW)'\|_p
&\mathop{\leq}^{\eqref{lem:control_path:product:(YZ)'_p}}
C_p(\|W\|_\infty\|Z'\|_p+\|Z'\|_\infty\|W\|_p+\|Z\|_\infty\|W'\|_p+\|W'\|_\infty\|Z\|_p)
\mathop{=}^{\eqref{eq:ineqs:W}}
o\revproof{_w}(\eta_1),
\label{eq:ineqs:ZW:||(ZW)'||_p}
\\
\|R^{ZW}\|_\frac{p}{2}
&\mathop{\leq}^{\eqref{lem:control_path:product:R^YZ_p/2}}
C_p(\|Z\|_\infty\|R^W\|_\frac{p}{2}+\|R^Z\|_\frac{p}{2}\|W\|_\infty+\|Z\|_p\|W\|_p)
\mathop{=}^{\eqref{eq:ineqs:W}}
o\revproof{_w}(\eta_1),
\\
\|(ZW)'\|_\infty
&\mathop{\leq}^{\eqref{eq:path_finite_var:infty_ineq}}
\|(ZW)'_0\|+\|(ZW)'\|_p
\mathop{=}^{\eqref{eq:ineqs:ZW:||(ZW)'||_p}}
o\revproof{_w}(\eta_1),
\end{align}
\end{subequations}
where the quantities $\|Z\|_p,\|Z'\|_p,\|Z'\|_\infty,\|R^Z\|_{\frac{p}{2}}$ can be bounded by $C_{p,b,\sigma}\exp\left(
C_{p,b,\sigma}N_{\alpha_{p,b,\sigma},[0,T]}(w)
\right)$ using   \eqref{eq:RDE:Yp_Y'p_RYp/2:bound:full_interval} and \eqref{eq:RY:p/2var:Y^2+RY+T}.  
Thus, for any $s,t>t_1$, 
\begin{align}
\left\|\int_s^tZ_sW_s\dd\mbX_s\right\|
&\mathop{\leq}^{\eqref{eq:rough_int:error_bound}}
\|Z_sW_sX_{s,t}\|+\|(ZW)'_s\bX_{s,t}\|+ C_p(\|R^{ZW}\|_{\frac{p}{2},[s,t]}\|X\|_{p,[s,t]}+\|(ZW)'\|_{p,[s,t]}\|\bX\|_{\frac{p}{2},[s,t]})
\nonumber
\\
&\leq
C_p\|\mbX\|_p(\|Z\|_{\infty,[s,t]}\|W\|_{\infty,[s,t]}+\|(ZW)'\|_{\infty,[s,t]}+\|R^{ZW}\|_{\frac{p}{2},[s,t]}+\|(ZW)'\|_{p,[s,t]})
\nonumber
\\
&\mathop{=}^{\eqref{eq:||X||_p<=exp(N_alpha^p)},\eqref{eq:ineqs:W},\eqref{eq:ineqs:ZW}}
o\revproof{_w}(\eta_1),
\ \  \text{so that }
\nonumber
\\
\label{eq:Delta(t):term3}
&\quad
\int_{t_1}^t
\int_0^1\rev{\nablaof{x}^2\sigma}(Y_s+\theta\Delta Y_s)(1-\theta)\dd\theta
\Delta Y^{\otimes 2}_s
\dd\mbX_s
=o\revproof{_w}(\eta_1).
\end{align}
Thus, using (\eqref{eq:Delta(t):term1},\eqref{eq:Delta(t):term2},\eqref{eq:Delta(t):term3}),  \eqref{eq:Delta(t):full} can be written as
\begin{align*}
\Delta_t
&=
o\revproof{_w}(\eta_1) + 
\int_{t_1}^t
\rev{\nablaof{x}b}(s,Y_s,u_s)
\Delta_s
\dd s
+
\int_{t_1}^t
\rev{\nablaof{x}\sigma}(s,Y_s)
\Delta_s
\dd\mbX_s,
\end{align*}
\revproof{
Finally, 
\begin{align*}
\|\Delta\|_{\infty,[t_1+\eta_1,T]}
&\mathop{\leq}^{\eqref{eq:rde:linearized:bounded_solutions}}
C_{p,T,b,\sigma}
\exp(C_{p,T,b,\sigma}N_{\alpha_{p,b,\sigma},[0,T]}(\mbX))
\|o_w(\eta_1)\|
\\
&\, 
\mathop{=}\limits^{\eqref{eq:needlelike:o(eta1)}}
\,
C_{p,b,\sigma\rev{,u,\bar{u}_1}}
\exp(
C_{p,T,b,\sigma}N_{\alpha_{p,b,\sigma},[0,T]}(\mbX)
)
\exp(
C_{p,b,\sigma}N_{\alpha_{p,b,\sigma},[0,T]}(w)
)
\revproof{o(\eta_1)}
\\
&\mathop{\leq}
\limits^{\eqref{eq:Nalpha<=3CpNalpha_X_and_time}}
C_{p,b,\sigma\rev{,u,\bar{u}_1}}
\exp(
C_{p,T,b,\sigma}N_{\alpha_{p,b,\sigma},[0,T]}(\mbX)
)
\revproof{o(\eta_1)},
\end{align*} 
}
\hspace{-1mm}which concludes the proof.
\end{proof}
Proposition \ref{prop:linear_variation}  can be extended to needle-like variations with multiple spikes using a standard argument by induction,  
 \ifarxiv
see  Corollary \ref{cor:needle_like_error:etas} in the 
appendix.
\fi 
\rev{
The next result follows from Proposition \ref{prop:linear_variation},  Theorem \ref{thm:gaussian_rough_paths}, and Theorem  \ref{thm:rdes:integrable}.}

\begin{lemma}[Needle-like variation formula for random RDEs]\label{lem:needle_like_error:etas}
Let %
$\rho\in[1,\frac{3}{2})$, $p\in(2\rho,3)$, %
and define $T,U,u,b,\sigma$ as in Proposition \ref{prop:linear_variation}.
Let $\Omega,\F,\Prob,B,\mbB,\ell,y,Y,Y'$ be as in Theorem \ref{thm:rdes:integrable}, with  $\mbB=(B,\bB)$ an enhanced Gaussian process and $(Y(\omega),Y'(\omega))\in\sD^p_{B(\omega)}([0,T],\R^n)$  the pathwise solution to the RDE
\begin{align}
\tag{\ref{eq:RDE}}
Y_t(\omega)&=y(\omega)+\int_0^tb(s,Y_s(\omega),u_s)\dd s+\int_0^t\sigma(s,Y_s(\omega))\dd\mbB_s(\omega),
\quad\ \ t\in[0,T].
\end{align}  
\revproof{Given} $q\in\N$, let $0<t_1<\dots<t_q<T$ be Lebesgue points %
(Definition \ref{def:lebesgue_point}), 
$\bar{u}_1,\dots,\bar{u}_q\in U$, 
 $0\leq \eta_i< t_{i+1}-t_i$ for $i=1,\dots, q-1$ and $0\leq\eta_q< T-t_q$, and define the needle-like variation $\pi=\{t_1,\dots,t_q,\eta_1,\dots,\eta_q,\bar{u}_1,\dots,\bar{u}_q\}$ of $u$ as the control $u^\pi$ defined by%
 $$
 \revproof{
 u^\pi_t=\begin{cases}
 \bar{u}_i\quad&\text{if }t\in[t_i,t_i+\eta_i],
 \\
 u_t&\text{otherwise}.
 \end{cases}
 }
 $$
As in Theorem \ref{thm:rdes:integrable}, 
let $(Y^\pi(\omega),(Y^\pi(\omega))')\in\sD^p_{B(\omega)}([0,T],\R^n)$ be the pathwise solution to the RDE
\begin{align*} 
Y^\pi_t(\omega)=y(\omega)+\int_0^tb(s,Y^\pi_s(\omega),u^\pi_s)\dd s+\int_0^t\sigma(s,Y^\pi_s(\omega))\dd\mbB_s(\omega),
\  t\in[0,T],
\end{align*}
and for $i=1,\dots,q$, let $(V^{\pi_i}(\omega),(V^{\pi_i}(\omega))')\in\sD^p_{B(\omega)}([0,T],\R^n)$ be the pathwise solutions to the RDEs 
\begin{subequations}
\begin{align*}
V^{\pi_i}_t(\omega) &= 
V^{\pi_i}_{t_i}(\omega)
+
\int_{t_i}^t
\rev{\nablaof{x}b}(s,Y_s(\omega),u_s)V^{\pi_i}_s(\omega)\dd s +
\int_{t_i}^t
\rev{\nablaof{x}\sigma}(s,Y_s(\omega))V^{\pi_i}_s(\omega)\dd\mbB_s(\omega),
 &&t\in[t_i,T],
\\
V^{\pi_i}_t(\omega)&=b(t_i,Y_{t_i}(\omega),\bar{u}_i)-b(t_i,Y_{t_i}(\omega),u_{t_i}),
&&t\in[0,t_i],
\end{align*}
\end{subequations}
with $Y,Y^\pi,V^{\pi_1},\dots,V^{\pi_q}\in L^\ell(\Omega,C([0,T],\R^n))$.  
Let $\varphi:\R^n\to\R^{\rev{k}}$ be continuously differentiable and satisfy $\big\|\rev{\nabla\varphi}(x)\big\|
\leq C_\varphi$ 
and 
$\big\|\rev{\nabla\varphi}(x)-\rev{\nabla\varphi}(\tilde{x})\big\|\leq C_\varphi\|x-\tilde{x}\|$ for all $x,\tilde{x}\in\R^n$ for some constant $C_\varphi<\infty$. 
 Then, 
\begin{align}
\label{eq:needle_like_deltasols:B}
\E\left[
\left\|
\varphi(Y^\pi)-\varphi(Y)
\right\|_{\infty,[0,T]}
\right]
&\leq 
C\,
\E\left[\exp(CN_{\alpha,[0,T]}(\mbB))\right]
\sum_{i=1}^q\eta_i,
\\
\label{eq:needle_like_deltasols_variation:B}
\E\bigg[
\bigg\|
\varphi(Y^\pi)-\varphi(Y)-\sum_{i=1}^q\eta_i\rev{\nabla\varphi}(Y)V^{\pi_i}
\bigg\|_{\infty,[t_q\revproof{+\eta_q},T]}
\bigg]
&\leq 
C\,
\E\left[\exp(CN_{\alpha,[0,T]}(\mbB))\right]
\revproof{o\bigg(\sum_{i=1}^q\eta_i\bigg),}
\end{align}
where $0<C<\infty$ and $0<\alpha<1$ \rev{depend} on $(p,T,b,\sigma,\rev{u,\bar{u}_1,\dots,\bar{u}_q,}\varphi)$, and $\E\left[\exp(CN(\mbB))\right]<\infty$.
\end{lemma}
\begin{proof}
For conciseness, we only prove the case $q=1$ corresponding to  one variation  $\pi_1:=\{t_1,\eta_1,\bar{u}_1\}$ of the control input $u$, as the case $q\geq2$ follows by using 
\ifarxiv
Corollary \ref{cor:needle_like_error:etas} 
\fi 
instead of Proposition \ref{prop:linear_variation}. Also, note that  Theorem \ref{thm:rdes:integrable} ensures that the pathwise solutions to the RDEs are well-defined and integrable.

By Proposition \ref{prop:linear_variation}, there exist constants $C>0$ and $0<\alpha<1$ \rev{depending} on $(p,T,b,\sigma\rev{,u,\bar{u}_1})$ such that 
\begin{align*}  
\|Y^{\pi}-Y\|_{\infty,[0,T]}
&\mathop{\leq}^{\eqref{eq:needle_like_deltasols}} 
C\exp\left(
CN_{\alpha,[0,T]}(\mbB)
\right)
\eta_1,
\quad
\|Y^{\pi}-Y-\eta_1V^{\pi_1}\|_{\infty,[t_1\revproof{+\eta_1},T]}
\mathop{\leq}^{\eqref{eq:needle_like_deltasols_variation}}
C\exp\left(
CN_{\alpha,[0,T]}(\mbB)
\right)
\revproof{o(\eta_1)}
\end{align*}
almost surely, 
and where $\E\left[\exp(CN_{\alpha,[0,T]}(\mbB))\right]<\infty$ by Theorem \ref{thm:gaussian_rough_paths}. %
Also, with $\Delta Y:=Y^{\pi}-Y$, by the mean value theorem and   Taylor's Theorem \cite{Folland1990},  
\begin{align*}
\|\varphi(Y^{\pi})-\varphi(Y)\|
&=
\left\|
\int_0^1\rev{\nabla\varphi}(Y+\theta\Delta Y)\dd\theta\Delta Y
\right\|
\leq C_\varphi\|\Delta Y\|.
\\
\left\|
\varphi(Y^{\pi})-\varphi(Y)-\eta_1\rev{\nabla\varphi}(Y)V^{\pi_1}
\right\|
&=
\left\|
\rev{\nabla\varphi}(Y)
(\Delta Y-\eta_1V^{\pi_1})
+
\int_0^1
\bigg(
\rev{\nabla\varphi}(Y+\theta\Delta Y)-\rev{\nabla\varphi}(Y)
\bigg)
\dd\theta
\Delta Y
\right\|
\\
&\leq 
C_\varphi(\|\Delta Y-\eta_1V^{\pi_1}\|+\|\Delta Y\|^2).
\end{align*}
The desired inequalities \eqref{eq:needle_like_deltasols:B} and \eqref{eq:needle_like_deltasols_variation:B} follow from the previous inequalities. 
\end{proof}

\subsection{Proof of the  Pontryagin Maximum Principle}\label{sec:pmp_proof}

We prove \pmp under the assumptions below.

\begin{assumption}[Assumptions for \ocp and \pmp (Theorem \ref{thm:pmp})]
\label{assum:pmp}
Let $T>0$, $\rho\in[1,\frac{3}{2})$, $p\in(2\rho,3)$, 
$\ell\geq 2$, 
$U\subseteq\R^m$,  
and $(\Omega,\F,\Prob)$ be a probability space. 
\begin{itemize}
\item The drift   $b:[0,T]\times\R^n\times U\to\R^n$ and cost  $f:[0,T]\times\R^n\times U\to\R$  
satisfy Assumption \ref{assumption:b:stronger} and are Lipschitz in $u$ for a constant $C_{b,f}\geq 0$, that is:
\begin{itemize} 
\item 
$b(\cdot, x, u):[0,T]\to\R^n$  is  measurable for all $(x,u)\in\R^n\times U$,  
\item  
$b(t,\cdot,\cdot):\R^n\times U\to\R^n$ 
is continuous for almost every $t\in[0,T]$,
\item 
$b(t,\cdot,u):\R^n\to\R^n$  is  continuously differentiable  for almost every $t\in[0,T]$  and  all $u\in U$,  
\item 
$\|b(t,x,u)\|+\big\|\rev{\nablaof{x}b}(t,x,u)\big\| \leq C_{b,f}$ 
and 
$\big\|\rev{\nablaof{x}b}(t,x,u)-\rev{\nablaof{x}b}(t,\tilde{x},u)\big\|\leq C_{b,f}\|x-\tilde{x}\|$  
for almost every $t\in[0,T]$, all $x,\tilde{x}\in\R^n$, and all $u\in U$,
\item %
$\|b(t,x,u)-b(t,x,\tilde{u})\|\leq C_{b,f}\|u-\tilde{u}\|$ for almost every $t\in[0,T]$, all $x\in\R^n$, and all $u,\tilde{u}\in U$,
\end{itemize}
and similarly for $f$.

\item The diffusion $\sigma:[0,T]\times\R^n\to\R^{n\times d}$ satisfies $\sigma\in C^4_b([0,T]\times\R^n,\R^{n\times d})$.
\item The terminal cost $g:\R^n\to\R$ and constraints $h:\R^n\to\R^r$ are continuously differentiable and %
satisfy 
$\big\|\rev{\nabla g}(x)\big\|+\big\|\rev{\nabla h}(x)\big\|
\leq C_{g,h}$
and 
$\big\|\rev{\nabla g}(x)-\rev{\nabla g}(\tilde{x})\big\|
+
\big\|\rev{\nabla h}(x)-\rev{\nabla h}(\tilde{x})\big\|
\leq C_{g,h}\|x-\tilde{x}\|$ for all $x,\tilde{x}\in\R^n$ 
for a constant $C_{g,h}>0$.

\item The initial conditions $x_0$ satisfy $x_0\in L^\ell(\Omega,\R^n)$.
\item $B:\Omega\to C([0,T],\R^d)$  is a centered, continuous, $\R^d$-valued Gaussian process with independent components   satisfying Assumption \ref{assum:Gaussian_lift} for $\rho$, and the driving signal 
$\mbB=(B,\bB):\Omega\to\sC^p_g([0,T],\R^d)$ is the enhanced Gaussian process in Theorem \ref{thm:gaussian_rough_paths} associated  to $B$ that satisfies \eqref{eq:gaussian_rough_paths:E[exp(Nalpha)]<infty}.
\end{itemize}
\end{assumption}

These assumptions on $(b,f)$ are reasonable, see for instance \cite{Diehl2016,Allan2020,BonalliLewESAIM2022} that make similar assumptions. 
The boundedness assumption on $b$ holds for control-affine drifts $b(t,x,u)=b_0(t,x)+b_1(t,x)u$ if  $(b_0,b_1)$ satisfy Assumption \ref{assumption:b:stronger} and  $U$ is bounded, and similarly for  $f$ if it is quadratic in  $u$ and $U$ is bounded.  See also
Remark \ref{remark:b_bounded} for technical considerations.  
Assuming  boundedness is reasonable in applications, as dynamical systems are often constrained to operate in a bounded set of  conditions representing physical constraints (i.e., $x_t\in\mathcal{X}$ with $\mathcal{X}\subset\R^n$ a compact set). 
 The assumption
$\sigma\in C^4_b$ is discussed in Remarks \ref{remark:sigma_smoothness} and \ref{remark:sigma_smoothness:2}: It is stronger than assumptions used in standard PMPs using It\^o calculus due to the use of rough path theory. %
As discussed in Sections \ref{sec:introduction} and \ref{sec:preliminaries}, the class of enhanced Gaussian processes $\mbB=(B,\bB)$ from Theorem \ref{thm:gaussian_rough_paths} is large (see the many examples in \cite{Friz2016}) and covers scenarios where $B$ is not a semimartingale (such as with fractional Brownian motion) that cannot be tackled via It\^o calculus. In particular, the Stratonovich lift of Brownian motion satisfies this assumption. %

\textbf{Proof sketch}: We now  prove Theorem \ref{thm:pmp} in the following steps:
\begin{itemize}
  \setlength\itemsep{0pt}
\item \textit{Augmented state, needle-like variations, and end-point mapping}: We define an augmented state $\tilde{x}=(x,x^0)$ that contains the state $x$ and its associated cost $x^0$,  the needle-like variations $u^\pi$ around the optimal control $u$, and the  end-point mapping $F^q$ that evaluates variations of the expected  terminal cost and constraints as the needle-like variation $\pi$ changes.
\item \textit{Variational linearization and separation argument}: We evaluate the end-point mapping $F^q$ around the optimal control $u$, and argue by contradiction to deduce the existence of the non-trivial Lagrange multiplier $(\mathfrak{p}_0,\dots,\mathfrak{p}_r)$ using a hyperplane separation argument based on Brouwer's fixed point theorem.
\item \textit{Adjoint \rev{equation} \eqref{eq:spmp:pmp_equations} and transversality condition \eqref{eq:spmp:transversality_condition:pT}}: We define the adjoint vector $p$ as the solution to the random linear RDE \eqref{eq:spmp:pmp_equations}, using the backward representation of rough integrals  \cite[Proposition 5.12]{Friz2020} to ensure that the terminal value $p_T$ satisfies the   transversality condition \eqref{eq:spmp:transversality_condition:pT}.
\item \textit{Maximality condition \eqref{eq:spmp:maximizality_condition}}: We use a contradiction argument by combining the inequality from the previous hyperplane separation theorem with  It\^o's lemma to deduce the maximality condition \eqref{eq:spmp:maximizality_condition}.
\end{itemize}
The main differences with standard proofs of previous stochastic PMPs \cite{Peng1990,Yong1999,BonalliLewESAIM2022} are 
\begin{itemize}
  \setlength\itemsep{0pt}
\item using random linear RDEs to evaluate the effect of needle-like variations  (here, the integrability of our bounds using greedy partitions and the integrability properties of Gaussian rough paths is key), 
\item defining the adjoint vector using a random RDE instead of FBSDEs (switching between forward and backward integration is done via \cite[Proposition 5.12]{Friz2020}, noting that rough path theory does not rely on non-anticipativity or martingale arguments), and 
\item deducing the maximality condition using a pathwise use of It\^o's lemma  %
to   mimic the proof of the deterministic PMP (i.e., the chain rule of It\^o's lemma gives us $\tilde{p}_t^\top \tilde{v}_t=\tilde{p}_T^\top \tilde{v}_T$ for all $t\in[0,T]$ almost surely \eqref{eq:pmp:proof:pv=cst} as in the deterministic setting, e.g., see \cite[page 254]{LeeMarkus1967} or  \cite[page 161]{Bonnard2005}).  
\end{itemize}

\textbf{Preliminary step -- Augmented state, needle-like variations, and end-point mapping.}
For any control $v\in L^\infty([0,T],U)$, consider the cost-augmented state $\tilde{x}_t = (x_t,x^0_t)$ and the random RDE 
\begin{align}
\label{eq:sde:cost_augmented}
\tilde{x}_t 
&=
\tilde{x}_0 + 
\int_0^t
\tilde{b}(s,\tilde{x}_s,v_s)\dd s + 
\int_0^t
\tilde{\sigma}(s,\tilde{x}_s)\dd\mbB_s
= 
\begin{bmatrix}
x_0
\\
0
\end{bmatrix}
+
\int_0^t
\begin{bmatrix}
b(s,x_s,v_s)
\\
f(s,x_s,v_s)
\end{bmatrix}\dd s + 
\int_0^t
\begin{bmatrix}
\sigma(s,x_s) \\ 0
\end{bmatrix}\dd\mbB_s,
\ t\in[0,T],
\end{align}
which has a well-defined solution $\tilde{x}^v=(x^v,x^{0,v})\in L^\ell(\Omega,C([0,T],\R^{n+1}))$ by Theorem \ref{thm:rdes:integrable} with  almost surely 
 $(\tilde{x}^v(\omega),(\tilde{x}^v(\omega))')\in\sD^p_{B(\omega)}([0,T],\R^{n+1})$. 
The first $n$ components $x^v$ of $\tilde{x}^v$ is the state trajectory associated to the control $v$, and  the last component $x^{0,v}$ of $\tilde{x}^v$ is the accumulated cost over the trajectory $(x^v,v)$. 
Given the optimal control $u$ that solves \ocp, $\tilde{x}:=\tilde{x}^u$ denotes the optimal cost-augmented state trajectory. 

Next, we define the needle-like variations of   the optimal control $u$ as in Lemma \ref{lem:needle_like_error:etas}. Given $q\in\N$, let $0<t_1<\dots<t_q<T$ be Lebesgue points %
(Definition \ref{def:lebesgue_point} \revproof{with $\widetilde b$ instead of $b$}), 
$\bar{u}_1,\dots,\bar{u}_q\in U$, 
 $0\leq \eta_i< t_{i+1}-t_i$ for $i=1,\dots, q-1$ and $0\leq\eta_q< T-t_q$, and define the needle-like variation $\pi=\{t_1,\dots,t_q,\eta_1,\dots,\eta_q,\bar{u}_1,\dots,\bar{u}_q\}$ of $u$ as the control $u^\pi$ defined by%
 $$
 u^\pi_t=\begin{cases}
 \bar{u}_i\quad&\text{if }t\in[t_i,t_i+\eta_i],
 \\
 u_t&\text{otherwise}.
 \end{cases}
 $$
Let $\tilde{x}^\pi:=\tilde{x}^{u^\pi}$ be the solution to the random RDE \eqref{eq:sde:cost_augmented} associated to the control $u^\pi$. 
 As shown in  Lemma \ref{lem:needle_like_error:etas},  the difference $x^\pi-x$ can be approximated using the pathwise solutions $\tilde{v}^{\pi_i}\in L^\ell(\Omega,C([t_i,T],\R^{n+1}))$ to the random linear RDEs 
\begin{align}\label{eq:sde:cost_augmented:linear}
\tilde{v}^{\pi_i}_t 
&=
\tilde{b}(t_i,\tilde{x}_{t_i},\bar{u}_i)
-
\tilde{b}(t_i,\tilde{x}_{t_i},u_{t_i})
 + 
\int_{t_i}^t
\rev{\nablaof{\tilde{x}}\tilde{b}}(s,\tilde{x}_s,u_s)
\tilde{v}^{\pi_i}_s
\dd s + 
\int_{t_i}^t
\rev{\nablaof{\tilde{x}}\tilde{\sigma}}(s,\tilde{x}_s)
\tilde{v}^{\pi_i}_s
\dd\mbB_s,
\quad t\in[t_i,T],
\end{align}
where $i=1,\dots,q$ and $(\tilde{v}^{\pi_i}(\omega),(\tilde{v}^{\pi_i}(\omega))')\in\sD^p_{B(\omega)}$ almost surely.   
 Next, we define the map
$$
\widetilde\Phi:\R^{n+1}\to\R^{r+1}, \
 \tilde{x}=(x,x^0)
\mapsto
\big(
h(x), \, x^0+g(x)
\big),
$$
which evaluates the terminal constraint and total cost associated to \ocp for an augmented state $\tilde{x}$. 
Also,  with  $
\R_+^q=\left\{\eta=(\eta_1,\dots,\eta_q)\in\R^q: \eta_1\geq 0,\dots,\eta_q\geq 0\right\}$, 
$\delta=\min\{T-t_q, t_{i+1}-t_i,  i=1,\dots, \revproof{q-1}\}$ and $B^q_\delta=\{\eta\in\R^q:\|\eta\|<\delta\}$, we define the end-point mapping  $\rev{F}:B^q_\delta\cap\R^q_+\to\R^{r+1}$ by 
\begin{equation}\label{eq:spmp:needlelike_endpoint_mapping}
\rev{F}(\eta)
=
\E\left[
\widetilde\Phi\left(\tilde{x}^\pi_T\right)
-
\widetilde\Phi\left(\tilde{x}_T\right)
\right]
=
\begin{bmatrix}
\E\big[
h(x^\pi_T)-h(x_T)
\big]
\\
\E\big[x^{\pi,0}_T+g(x^\pi_T)
-
(x_T^0+g(x_T))\big]
\end{bmatrix}
\end{equation}
which satisfies $\rev{F}(0)=0$, 
and  the linear map $\dd \rev{F_0}:\R^q_+\to\R^{r+1}$ by
\begin{equation}\label{eq:end_point_mapping:differential}
\dd \rev{F_0}(\eta)=
\sum_{i=1}^q\eta_i
\E\left[
\rev{\nabla\widetilde\Phi}(\tilde{x}_T)
\tilde{v}^{\pi_i}_T
\right].
\end{equation}
The map 
$\dd \rev{F_0}$ is the Gateaux differential of $\rev{F}$ at $0$ in the direction $\eta$:%
\begin{equation}\label{eq:end_point_mapping:gateaux_differentiable}
\lim_{ \alpha>0, \, \alpha\to 0 }
\frac{\rev{F}(\alpha\eta)}{\alpha}=
\dd \rev{F}_0(\eta)
\
\text{ for any }\eta\in\R^q_+.
\end{equation}
Indeed,  by Lemma \ref{lem:needle_like_error:etas} and Jensen's inequality,  for a constant $C>0$  and $\alpha>0$ small-enough, 
$$
\big\|
\rev{F}(\alpha\eta)-\dd \rev{F_0}(\alpha\eta)
\big\|
\leq
\E\bigg[
\bigg\|
\widetilde\Phi\left(\tilde{x}^\pi_T\right)
-
\widetilde\Phi\left(\tilde{x}_T\right)
-
\sum_{i=1}^q\alpha\eta_i
\rev{\nabla\widetilde\Phi}(\tilde{x}_T)
\tilde{v}^{\pi_i}_T
\bigg\|
\bigg]
\mathop{\leq}^{\eqref{eq:needle_like_deltasols_variation:B}}
\revproof{o\bigg(\alpha\sum_{i=1}^q\eta_i\bigg)}
$$
so   \eqref{eq:end_point_mapping:gateaux_differentiable} follows after dividing by $\alpha$ and taking the limit as $\alpha\to 0$.  
 
Note also that multiple derivatives in $(\tilde{b},\tilde{\sigma},\widetilde\Phi)$ are zero as they do not depend on $x^0$, so
\begin{equation}\label{eq:tilde_b_sig_Phi}
\rev{\nablaof{\tilde{x}}\tilde{b}}(t,x,u)^\top
=
\begin{bmatrix}
\rev{\nablaof{x}b}(t,x,u) & \rev{\nablaof{x}f}(t,x,u) 
\\
0 & 0
\end{bmatrix},
\quad
\rev{\nablaof{\tilde{x}}\tilde{\sigma}}(t,x)^\top
=
\begin{bmatrix}
\rev{\nablaof{x}\sigma} & 0 
\\
0 & 0
\end{bmatrix},
\quad
\rev{\nabla\widetilde\Phi}(x)=\begin{bmatrix}
\rev{\nabla h}(x) & 0 
\\
\rev{\nabla g}(x) & 1
\end{bmatrix}.
\end{equation}

\textbf{Step 1) -- Variational linearization and separation argument.}
The first step of the proof %
proceeds by contradiction, using the end-point mapping $\rev{F}$ %
and a Brouwer fixed point argument.  We start by defining the 
closed convex cone
$$
K=
\overline{
\left\{
\sum_{i=1}^{\rev{\tilde{q}}}\alpha_i
\E\left[
\rev{\nabla\widetilde\Phi}(\tilde{x}_T)
\tilde{v}^{\pi_i}_T
\right]
\,:\,
\alpha_i\geq 0,
\
\pi_i\rev{\ =\{t_i,\eta_i,\bar{u}_i\}}\text{ is a needle-like variation of $u(\cdot)$},
\
\rev{\tilde{q}}\in\N
\right\}
},
$$
and note that 
$K\subset\R^{r+1}$. Indeed, by contradiction,  if $K=\R^{r+1}$, then \rev{there exist $q\in\N$ needle-like variations $\pi_i=\{t_i,\eta_i,\bar{u}_i\}$ such that $K=\{
\sum_{i=1}^q\alpha_i
\E[
\nabla\widetilde\Phi(\tilde{x}_T)
\tilde{v}^{\pi_i}_T
]:\alpha_i\geq 0\}$. Using these needle-like variations $\pi_i$ and $\pi=\{t_1,\dots,t_q,\eta_1,\dots,\eta_q,\bar{u}_1,\dots,\bar{u}_q\}$ of $u(\cdot)$,  we define the maps $\rev{F}$ and $\dd \rev{F_0}$, and obtain} $\dd \rev{F_0}(\R^q_+)=K=\R^{r+1}$, so $0\in\Int(\rev{F}(B^q_\delta\cap\R^q_+))$ by \cite[Lemma 12.4]{Agrachev2004} (whose proof relies on Brouwer's fixed point theorem). %
In particular, there exists another feasible trajectory $(x^\pi,u^\pi)$ with a strictly lower cost $\E[x^{\pi,0}_T+g(x^\pi_T)]<\E[x^0_T+g(x_T)]$, so $(x,u)$ is not optimal, which is a contradiction. 
Thus, $K\subset\R^{r+1}$. 

\ifarxiv
\begin{figure}[t]
\centering
\includegraphics[width=0.4\textwidth]{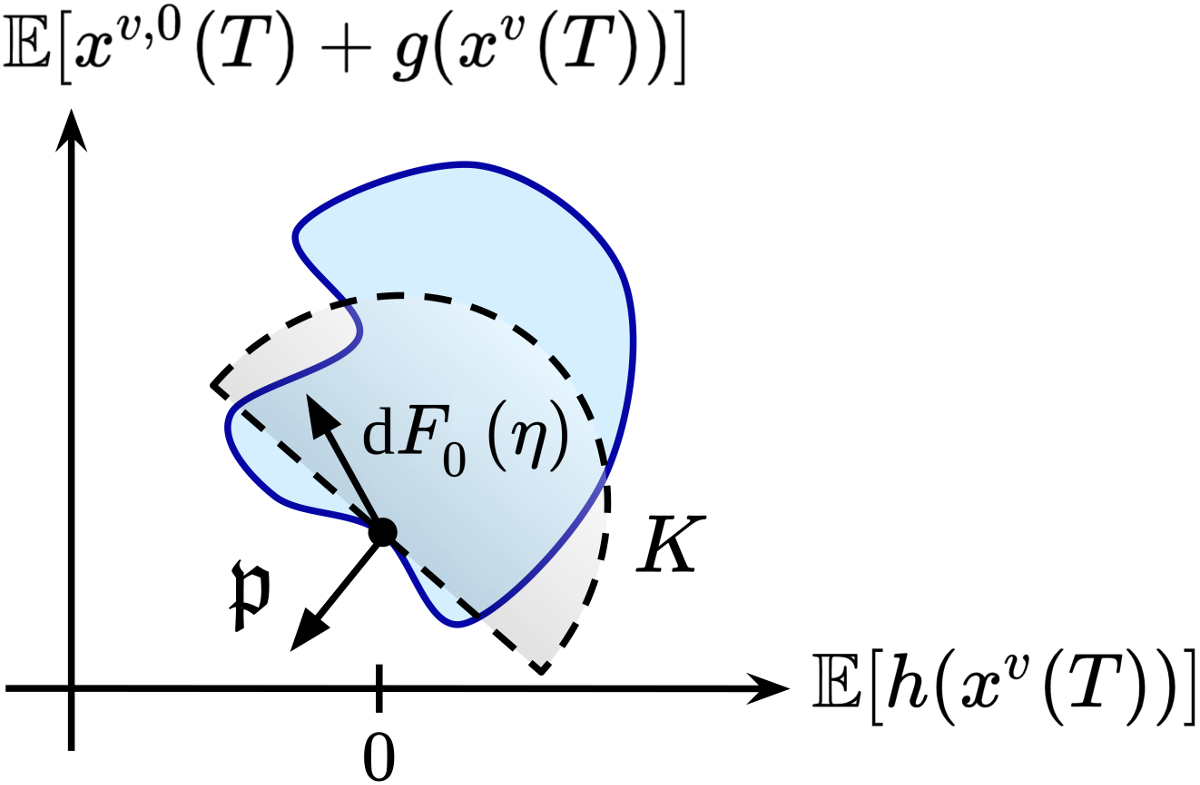}
\caption{Separation hyperplane argument at the optimal trajectory $(x,u)$. The blue region is the reachable set  $\{\E[\widetilde\Phi(\tilde{x}^v)], v\in L^\infty([0,T],U)\}$. The cone  $K$ \rev{(often called the first Pontryagin cone) satisfies $K\subset\R^{r+1}$}, otherwise we would be able to find  a feasible trajectory $(x^v,v)$ with lower cost, so $(x,u)$ would not be optimal.}
\label{fig:comparison}
\end{figure}
\fi

By the  hyperplane separation theorem, there exists a non-zero vector $\mathfrak{p}=(\mathfrak{p}_1,\dots,\mathfrak{p}_r,\mathfrak{p}_0)\in\R^{r+1}$   such that $\mathfrak{p}_0\leq 0$ and 
$\mathfrak{p}^\top z
\leq 
0$ for all $z\in K$. By renormalizing $\mathfrak{p}$, we may assume that $\mathfrak{p}_0\in\{0,-1\}$. The condition $\mathfrak{p}^\top z
\leq 
0$ for all $z\in K$ can be written as 
\begin{equation}\label{eq:spmp:proof:mu^TEgradPhiv<=0}
\mathfrak{p}^\top\,
\E\left[
\rev{\nabla\widetilde\Phi}(\tilde{x}_T)
\tilde{v}^{\pi_1}_T
\right]
\leq 
0
\ \ 
\text{for all needle-like variations $\pi_1=(t_1,\eta_1,\bar{u}_1)$ of $u$}.
\end{equation} 

\textbf{Step 2) -- Adjoint equation \eqref{eq:spmp:pmp_equations} and transversality condition \eqref{eq:spmp:transversality_condition:pT}.}
Define the reversed processes 
\rev{$\big(\overleftarrow{\ \tilde{x}_t},\overleftarrow{\ u_t^{\phantom{a}}}\big)
:=
\big(\tilde{x}_{T-t},u_{T-t}\big)$  and $\overleftarrow{\mbB}:=(\overleftarrow{B},\overleftarrow{\bB})$ with $\overleftarrow{B}_t:=B_{T-t}$ and $\overleftarrow{\bB}_{s,t}:=\bB_{T-s,T-t}$ (see \cite[Exercise 2.6]{Friz2020})}. %
Then, using Theorem \ref{thm:rdes:integrable} 
(noting with \cite[Proposition 5.12]{Friz2020} that $\overleftarrow{\tilde{x}}$ is the pathwise solution to an RDE driven by $\overleftarrow{\mbB}$, and since  including a drift term in \cite[Proposition 5.12]{Friz2020} poses no difficulty as  linear RDEs with drift can be rewritten as driftless linear RDEs along a new geometric rough path as in the proof of Theorem \ref{thm:rde:linear:existence_uniqueness}), 
we define the pathwise solution $\overleftarrow{\tilde{p}}\in L^\ell(\Omega,C([0,T],\R^{n+1}))$ to the random linear RDE
\begin{align*}
\overleftarrow{\ \tilde{p}_t}
 &=
\mathfrak{p}^\top \rev{\nabla\widetilde\Phi}(\tilde{x}_T) 
+ \int_0^t
\rev{\nablaof{\tilde{x}}\tilde{b}}
\big(
\revproof{T-}s,
\overleftarrow{\ \tilde{x}_s},
\overleftarrow{\ u_s^{\phantom{a}}}
\big)^\top
\overleftarrow{\ \tilde{p}_s}\dd s
+
\int_0^t
\rev{\nablaof{\tilde{x}}\tilde{\sigma}}\big(\revproof{T-}s,\overleftarrow{\ \tilde{x}_s}\big)^\top \overleftarrow{\ \tilde{p}_s}\dd\overleftarrow{\mbB}_s,
\quad
t\in[0,T],
\end{align*}
where  $\big(\overleftarrow{\tilde{p}}(\omega),\overleftarrow{\tilde{p}}(\omega)'\big)\in\sD_{\overleftarrow{B}(\omega)}^p$ almost surely. 
Next, we define the adjoint vector $\tilde{p}\in L^\ell(\Omega,C([0,T],\R^\revproof{n+1}))$ by $(\tilde{p}_t,\tilde{p}_t'):=\big(\overleftarrow{\tilde{p}}_{T-t},\overleftarrow{\tilde{p}}'_{T-t}\big)$ for all $t\in[0,T]$ almost surely. 
By \cite[Proposition 5.12]{Friz2020}, %
the sample paths satisfy $(\tilde{p}(\omega),\tilde{p}(\omega)')\in\sD^p_{B(\omega)}$ almost surely and solve the linear RDE 
\begin{align*}
 \tilde{p}_t(\omega)&=
\tilde{p}_0(\omega)
-\int_0^t
\rev{\nablaof{\tilde{x}}\tilde{b}}(s,\tilde{x}_s(\omega),u_s)^\top
\tilde{p}_s(\omega)\dd s
-
\int_0^t
\rev{\nablaof{\tilde{x}}\tilde{\sigma}}(s,\tilde{x}_s(\omega))^\top \tilde{p}_s(\omega)\dd\mbB_s(\omega),
\quad
t\in[0,T],
\end{align*}
with $\tilde{p}_0=\overleftarrow{\ \tilde{p}_T}\in L^\ell(\Omega,\R^{n+1})$,
 which gives  the adjoint equation \eqref{eq:spmp:pmp_equations}, since multiple derivatives in $(\tilde{b},\tilde{\sigma})$ are zero by \eqref{eq:tilde_b_sig_Phi}. Also,  $\tilde{p}_T=\overleftarrow{\ \tilde{p}_0}=\mathfrak{p}^\top
\rev{\nabla\widetilde\Phi}(\tilde{x}_T) $ gives the transversality condition \eqref{eq:spmp:transversality_condition:pT}  
\begin{align*} 
\tilde{p}_T
&= 
\mathfrak{p}^\top
\rev{\nabla\widetilde\Phi}(\tilde{x}_T) 
\mathop{=}^{\eqref{eq:tilde_b_sig_Phi}}
\left(
\sum_{i=1}^r\mathfrak{p}_i
\rev{\nabla h_i}(x_T)
+
\mathfrak{p}_0  
\rev{\nabla g}(x_T),
\ 
\mathfrak{p}_0
\right)
\ \text{ almost surely}.
\end{align*}   
Also $\dd p^0_t=0$ by \eqref{eq:tilde_b_sig_Phi}, so  $p^0_t=\mathfrak{p}_0$ for all $t\in[0,T]$.

\textbf{Step 3) -- Maximality condition \eqref{eq:spmp:maximizality_condition}.} 
We observe that $
\tilde{v}^{\pi_1}_t{}^\top
\rev{\nablaof{\tilde{x}}\tilde{b}}^\top\tilde{p}_t
=
\tilde{p}_t{}^\top
\rev{\nablaof{\tilde{x}}\tilde{b}}\tilde{v}^{\pi_1}_t$  
and  that  
$
\tilde{v}^{\pi_1}_t{}^\top
\rev{\nablaof{\tilde{x}}\tilde{\sigma}}^\top\tilde{p}_t
=
\sum_{i=1}^n
[\tilde{v}^{\pi_1}_t]^i
\sum_{k=1}^n\hspace{-2pt}
\rev{\nablaof{\tilde{x}^i}\tilde{\sigma}^{k\cdot}}[\tilde{p}_t]^k
=
\tilde{p}_t^\top
\rev{\nablaof{\tilde{x}}\tilde{\sigma}}
\tilde{v}^{\pi_1}_t
$,  
where $(\rev{\nablaof{\tilde{x}}\tilde{b}},\rev{\nablaof{\tilde{x}}\tilde{\sigma}})
=
(\rev{\nablaof{\tilde{x}}\tilde{b}}(t,\tilde{x}_t,u_t),\rev{\nablaof{\tilde{x}}\tilde{\sigma}}(t,\tilde{x}_t))$ for conciseness. 
Then, by It\^o's lemma (Lemma \ref{lem:rough_path:ito_formula}), almost surely, for any $s,t\in[\revproof{t_1},T]$, %
 \begin{align}\label{eq:pmp:proof:pv=cst}
 \hspace{-2pt}
\tilde{p}_t^\top 
\tilde{v}^{\pi_1}_t
=
\tilde{p}_s^\top
\tilde{v}^{\pi_1}_s
+
\int_s^t
\tilde{v}^{\pi_1}_r{}^\top
\bigg(-
\rev{\nablaof{\tilde{x}}\tilde{b}}^\top\tilde{p}_r\dd r - \rev{\nablaof{\tilde{x}}\tilde{\sigma}}^\top\tilde{p}_r\dd\mbB_r
\bigg)+
\int_s^t
\tilde{p}_r^\top
\bigg(
\rev{\nablaof{\tilde{x}}\tilde{b}}\tilde{v}^{\pi_1}_r\dd r + \rev{\nablaof{\tilde{x}}\tilde{\sigma}}
\tilde{v}^{\pi_1}_r\dd\mbB_r
\bigg)
=\tilde{p}_s^\top
\tilde{v}^{\pi_1}_s.
\end{align}
Thus, $\E\left[
\tilde{p}_t^\top
\tilde{v}^{\pi_1}_t
\right]
=
\E\left[
\tilde{p}_T^\top
\tilde{v}^{\pi_1}_T
\right]
=
\E\left[
\mathfrak{p}^\top
\rev{\nabla\widetilde\Phi}(\tilde{x}_T)
\tilde{v}^{\pi_1}_T
\right]$ for almost every $t\in[\revproof{t_1},T]$, with  $\E\left[
\tilde{p}_t^\top
\tilde{v}^{\pi_1}_t
\right]<\infty$ by H\"older's inequality since $\tilde{p}_t,\tilde{v}_t\in L^\ell(\Omega,\R^{n+1})$ with $\ell\geq 2$. We combine this result with \eqref{eq:spmp:proof:mu^TEgradPhiv<=0} and obtain 
\begin{equation}\label{eq:spmp:proof:pv_neq0_for_maximality}
\E\left[
\tilde{p}_t^\top
\tilde{v}^{\pi_1}_t
\right]
\leq 0
\ \ \text{for a.e. }t\in[0,T],
\   
\text{for any needle-like variations $\pi_1=(t_1,\eta_1,\bar{u}_1)$ of $u$}.
\end{equation}

Finally, by contradiction, suppose that \eqref{eq:spmp:maximizality_condition} does not hold. Then, there \rev{exist} a control $u^1(\cdot)$ and a subset of $[0,T]$ of positive measure on which
$$
\E\left[H(t,x_t,p_t,\mathfrak{p}_0,u_t)\right]
<
\E\left[H(t,x_t,p_t,\mathfrak{p}_0,u^1_t)\right].
$$
Let $t_1$ be a Lebesgue point of this subset of $[0,T]$. Then, 
$$
\E\left[
\tilde{p}_{t_1}^\top \tilde{b}(t_1,x_{t_1},u_{t_1})
\right]
<
\E\left[
\tilde{p}_{t_1}^\top \tilde{b}(t_1,x_{t_1},\bar{u}_1)
\right]
$$
for some $\bar{u}_1\in U$. Then, if we define \revproof{a} needle-like variation $\pi_1=(t_1,\revproof{\eta_1},\bar{u}_1)$ with associated variation vector $\tilde{v}^{\pi_1}$, we obtain
$$
\E\left[
\tilde{p}_{t_1}^\top 
	\left(\tilde{b}(t_1,x_{t_1},\bar{u}_1)-\tilde{b}(t_1,x_{t_1},u_{t_1})\right)
\right]
	=
\E\left[
\tilde{p}_{t_1}^\top \tilde{v}^{\pi_1}_{t_1}
\right]
>0,
$$
which contradicts \eqref{eq:spmp:proof:pv_neq0_for_maximality}. Thus, the maximality condition \eqref{eq:spmp:maximizality_condition} holds, which concludes the proof of \pmp.

\section{Numerical example}\label{sec:example}
As a brief application of \pmp, we implement the indirect shooting method presented in Section \ref{sec:introduction} %
 for a %
regulation task. %
We consider the open-loop (OL) 
optimal control problem (OCP)
\begin{equation}\label{OLOCP}%
\tag{\olocp}
\begin{cases}
\min\limits_{u\in L^\infty([0,T],\R^m)} \qquad 
	&\E\left[
\int_0^T \frac{1}{2}(x_t^\top Q x_t +u_t^\top R u_t)\dd t
\right]
\\
\ \ \ \textrm{\rev{subject to}} \ \ \  
&
x_t = x_0+\int_0^t (A(x_s)x_s+\bar{B}u_s)\dd s +\int_0^t  \sigma(x_s)\circ\dd B_s,
\quad t\in[0,T],
\end{cases}
\end{equation}
where $n=m=3$, $A(x)=-J^{-1}S(x)J$ with $ 
S(x)=
\SmallMatrix{
0&-x_3&x_2
\\
x_3&0&-x_1
\\
-x_2&x_1&0
}$ and  
$J=\text{diag}(J_1,J_2,J_3)=\SmallMatrix{
3&0&0
\\
0&2&0
\\
0&0&4
}$, 
$\bar{B}=J^{-1}$, 
$\sigma(x)=0.4\,\textrm{diag}(x)$, 
$R=\text{diag}(R_1,R_2,R_3)=3I_{3\times 3}$,  
$Q=10I_{3\times 3}$,  
$x_0=\frac{\pi}{180}(-1, -4.5, 4.5)$,   
$B$ is a  standard $n$-dimensional Brownian motion, and the SDE is a Stratonovich SDE. 
This problem may represent a stabilization task for the angular velocity of a spacecraft with nonlinear rigid body dynamics. 
By composing all functions in %
\olocp with a smooth cut-off function, %
we may assume that Assumption \ref{assum:pmp} holds, so candidate optimal solutions are described by \pmp: %
\begin{align*}
&\begin{cases}
H(x,u,p,\mathfrak{p}_0) 
=
p^\top(A(x)x+\bar{B}u)
-
\frac{1}{2}(x^\top Q x +u^\top R u),
\\
\rev{\nablaof{u}H}
=
p^\top B
-
u^\top R
\mathop{\implies}\limits^\eqref{eq:spmp:maximizality_condition} 
u_t
=
R^{-1}
\bar{B}^\top\E\left[
p_t
\right],
\\
p_T
\mathop{=}\limits^\eqref{eq:spmp:transversality_condition:pT} 
0,
\end{cases}
\end{align*}
where $\mathfrak{p}_0=-1$ since there are no final state constraints.  

\textbf{Numerical resolution:} 
We consider two algorithms that use $M\in\N$ independent samples $B^i$ of $B$. 
\begin{enumerate}[leftmargin=5mm]
\item[1)] Direct method (\texttt{Direct}): We search for the control $u\in L^\infty([0,T],\R^m)$ and the sample paths $(x^i)_{i=1}^M\in C([0,T],\R^{Mn})$ that solve the sample average approximation 
\begin{equation*}
\min\limits_u \  
\frac{1}{M}\sum_{i=1}^M
\int_0^T \frac{1}{2}(x_t(\omega^i)^\top Q x_t(\omega^i) +u_t^\top R u_t)\dd t 
\ \ \ \textrm{s.t.} \ \ \   
x_t^i = x_0+\int_0^t(A(x_s^i)x_s^i+\bar{B}u_s)\dd s + \int_0^t\sigma(x_s^i)\circ \dd B_s^i.
\end{equation*}
Numerically, we discretize the problem with $(N+1)\in\N$ nodes by optimizing over  $\hat{u}=(\hat{u}_0,\dots,\hat{u}_N)\in\R^{(N+1)m}$ and $(\hat{x}^i)_{i=1}^M=((\hat{x}^i_0,\dots,\hat{x}^i_N))_{i=1}^M\in\R^{M(N+1)n}$ and discretizing the SDE with %
a Milstein scheme. %
We solve the resulting finite-dimensional optimization problem  via sequential quadratic programming (SQP) \cite[Chapter 18]{2006}. We do not enforce trust region constraints nor use a linesearch and thus always take full steps at each SQP iteration. We return a solution once the difference between SQP iterates $\Delta:=(\Delta\hat{u},(\Delta\hat{x}^i)_{i=1}^M)$ satisfies $\|\Delta\|_\infty<\epsilon$. 
This method is  in \cite[Section 6.2]{Lew2024}.
\item[2)] Indirect shooting method (\texttt{Indirect}): %
We search for  adjoint vector initial  values $(p_0^i)_{i=1}^M\in\R^{Mn}$ 
satisfying
\begin{align}\label{eq:example:shooting_problem}
\begin{bmatrix}
p_T^1
\\
\vdots
\\
p_T^M
\end{bmatrix}
=
0,
\ 
\text{ where }
\begin{cases}
x_T^i = x_0+\int_0^T 
(A(x_t^i)x_t^i+\bar{B}u_t^M)
\dd t +\int_0^T  \sigma(x_t^i)\dd\mbB_t^i,
\\[1mm]
p_T^i
=
p_0^i
-\int_0^T\rev{\nablaof{x}H}(x_t^i,u_t^M,p_t^i,\mathfrak{p}_0)\dd t
-
\int_0^T
\rev{\nabla\sigma}(x_t^i)^\top p_t^i\dd\mbB_t^i,
\\ 
u_t^M=R^{-1}\bar{B}^\top 
\left(\frac{1}{M}\sum\limits_{i=1}^M
p_t^i\right)
\text{ for a.e. }
t\in[0,T],
\end{cases}
\end{align} 
where the $\mbB^i$'s are the Stratonovich  lifts of the $B^i$'s.  
We define the map $$F:\, \R^{Mn}\to\R^{Mn},\  
(p_0^i)_{i=1}^M
\mapsto
(p_T^i)_{i=1}^M$$ that solves the coupled RDE in \eqref{eq:example:shooting_problem} from initial  values $(p_0^i)_{i=1}^M$ and returns $(p_T^i)_{i=1}^M$.  %
Numerically, we integrate the RDE in \eqref{eq:example:shooting_problem} using the estimate \eqref{eq:rough_int:error_bound} for rough integrals, see  
\ifarxiv
the appendix
\fi 
for details. 
Then, we solve the equation $F\left( (p_0^i)_{i=1}^M \right)=0$ via Newton's method by iteratively defining 
$$
({}^{(\ell+1)}p_0^i)_{i=1}^M=
({}^{(\ell)}p_0^i)_{i=1}^M
-
\left(
\nabla F\left(({}^{(\ell)}p_0^i)_{i=1}^M\right)
\right)
^{-1}
F\left(({}^{(\ell)}p_0^i)_{i=1}^M\right)
\qquad
\ell=0,1,\dots,
$$ 
starting from an initial guess $({}^{(0)}p_0^i)_{i=1}^M$. 
We return a solution $({}^{(\ell)}p_0^i)_{i=1}^M$ once $\|F(({}^{(\ell)}p_0^i)_{i=1}^M\rev{)}\|_\infty < \epsilon$. 
\end{enumerate}
The two methods are implemented in Python using \texttt{JAX} \cite{jax2018github}. %
The quadratic programs at each SQP iteration of \texttt{Direct} are solved using \texttt{OSQP} \cite{osqp2020}. 
The tolerance threshold is set to $\epsilon=10^{-6}$. 
We use zero initial guesses for both methods. 
We checked that solutions returned by the \texttt{Direct}  and \texttt{Indirect} methods are close to each other. 
Computation times are measured on a laptop with
an 1.10GHz Intel Core i7-10710U CPU.  
Code to reproduce results is available at \blue{\url{https://github.com/ToyotaResearchInstitute/rspmp}}.

\textbf{Results}: A solution to \olocp for $(T,M,N)=(2,10,40)$ is reported in Figure \ref{fig:openloop_feedback}. The  state and control trajectories converge to zero over time (in average for the state trajectories) to minimize the cost. The adjoint vector trajectories $p^i$ start from different initial conditions $p_0^i$ and are all zero at the final time ($p_T^i=0$) to satisfy  the transversality condition  of \PMP.

We report median computation time \rev{speedups} and costs over $20$ runs of each method in Figure \ref{fig:comparison}. 
The proposed \texttt{Indirect} method is about $10\times$ faster than the \texttt{Direct} method. %
Also, the costs associated to the solutions to the sampled problems decrease %
as the sample size $M$ increases. %
Solutions are sensitive to the sample size, and computing low-cost solutions with high certainty over the sampling procedure is achievable with a reasonable sample size ($M=10$) for this problem.

 \textbf{Discussion}: %
First, results in Figure \ref{fig:comparison} (right) suggest that the proposed method may be asymptotically optimal, since the cost and its variance decrease as the sample size $M$ increases. Proving such asymptotic optimality properties of the method for certain classes of problems  is of interest for future work.  
Second, the \texttt{Indirect} method is significantly faster than the \texttt{Direct} method, thanks to leveraging the structure of the problem encoded in \pmp %
to optimize over only the $Mn$ variables $p_0^i$ for the \texttt{Indirect} method versus optimizing over the $M(N+1)n+(N+1)m$ variables $(\hat{x}^i,\hat{u}^i)$ for the \texttt{Direct} method.  
However, indirect methods  typically have higher numerical sensitivity to the choice of initial guess.  %
This tradeoff is well-known in the deterministic optimal control literature, motivating the future development of multiple shooting and homotopy methods \cite{Trelat2012,Bonalli2018} for stochastic optimal control.

\ifarxiv
\begin{figure}[t]
\centering
\includegraphics[width=0.65\textwidth,trim={4mm 5mm 15mm 5mm},clip]{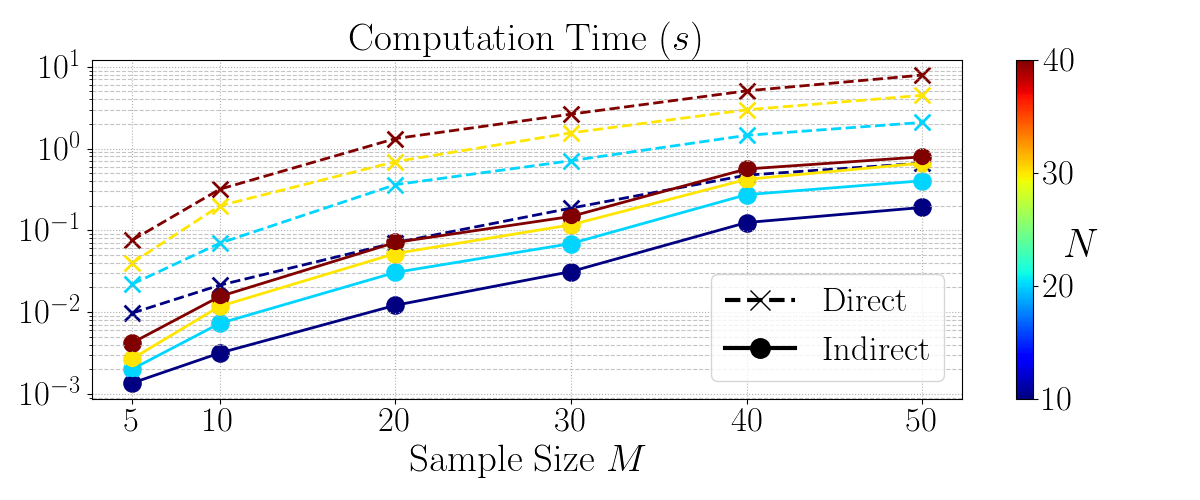}
\begin{minipage}{.49\textwidth}
\centering
{
\addtolength{\tabcolsep}{-0.3em}
\rev{
\begin{tabular}{|cc|cccccc|}
\toprule
\multicolumn{2}{|c|}{Speedup}
 & \multicolumn{6}{c|}{Sample Size $M$}
\\[1pt]
\multicolumn{2}{|c|}{\texttt{Direct} / \texttt{Indirect}}
 & $5$ & $10$ & $20$ & $30$ & $40$ & $50$
\\
\midrule
& $10 $ & $7.2$ & $6.8$ & $5.8$ & $6.0$ & $3.8$ & $3.5$
\\
Horizon & $20 $ & $11.0$ & $9.6$ & $11.8$ & $10.3$ & $5.3$ & $5.2$
\\
$N$ & $30 $ & $14.9$ & $17.0$ & $13.3$ & $13.3$ & $7.1$ & $6.7$
\\
& $40 $ & $18.0$ & $20.5$ & $18.6$ & $17.8$ & $9.0$ & $10.0$
\\
\bottomrule
\end{tabular}
}
}%
\end{minipage}%
\hspace{1mm}
\begin{minipage}{.49\textwidth}
\centering
\includegraphics[width=0.95\textwidth]{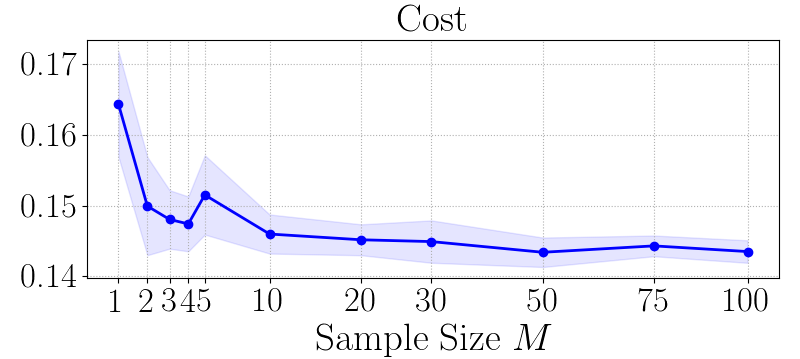}
\end{minipage}%
\caption{
\textbf{Top}: Median computation time for the \texttt{Direct} and \texttt{Indirect} methods solving \olocp for $T=3$ and varying sample sizes $M$ and horizons $N$. 
\textbf{Bottom left}: Median \rev{speedup over} the \texttt{Direct} \rev{method} \rev{(a speedup of $10$ corresponds to the \texttt{Indirect} method returning a solution $10\times$ faster than the \texttt{Direct} method)}. 
\textbf{Bottom right}: Median cost $\pm$ one median absolute deviation intervals of  solutions returned by the \texttt{Indirect} method for different samples sizes for $N=40$, evaluated with $10^4$ Monte Carlo samples.}
\label{fig:comparison}
\end{figure}

\fi
 
\textbf{Feedback optimization:} 
Next, we consider the feedback (FB) optimal control problem (OCP)
\begin{equation}
\label{FBOCP}%
\tag{\fbocp}
\begin{cases}
\min\limits_{k\in L^\infty([0,T],\R^m)} \qquad 
	&\E\left[
\int_0^T \frac{1}{2}(x_t^\top Q x_t +x_t^\top K_t^\top R K_t x_t)\dd t
\right]
\\
\ \ \ \textrm{\rev{subject to}} \ \ \  
&
x_t = x_0+\int_0^t (A(x_s)+\bar{B}K_s)x_s\dd s +\int_0^t  \sigma(x_s)\circ\dd B_s,
\quad t\in[0,T],
\end{cases}
\end{equation}
where we optimize  over the diagonal feedback gain $K_t=\text{diag}(k_t)\in\R^{m\times m}$. 
The problem \fbocp derives from \olocp by considering the feedback control $u=Kx$ and optimizing over the gains $K$. Since the gains $k$ are deterministic, the resulting feedback control law $u=Kx$ is causal and this formulation fits within our framework. 
Candidate optimal solutions are described by \pmp: 
\begin{align*} 
\fbocp:
&\begin{cases}
H(x,k,p,\mathfrak{p}_0)
=
p^\top(A(x)+\bar{B}K)x
-
\frac{1}{2}(x^\top Q x +x^\top K^\top RK x),
\\
\rev{\nablaof{k}H}
=
(1/J_j)p_jx_j
-
k_jR_jx_j^2
\mathop{\implies}\limits^\eqref{eq:spmp:maximizality_condition} 
k_{j,t}
=
(J_jR_j)^{-1}\E\left[
p_{j,t}x_{j,t}
\right]
/
\E\left[(x_{j,t})^2\right], \  j=1,2,3,
\\
p_T
\mathop{=}\limits^\eqref{eq:spmp:transversality_condition:pT} 
0.
\end{cases}
\end{align*}
Using these necessary conditions, we implement the %
\texttt{Indirect} shooting method %
presented previously  %
and solve  \fbocp for $(T,M,N)=(2,10,40)$ from a zero initial guess $({}^{(0)}p_0^i)_{i=1}^M$. 
We find that solving \fbocp is numerically sensitive to the choice of initial guess. %
Thus, using a homotopy method, we solve \fbocp for $R_j\in\{100, 99.9,\dots,3\}$, using the solution $({}^{(0)}p_0^i)_{i=1}^M$ computed for the previous value of $R$ as an initial guess for each solve via Newton's method. 

Results in Figure \ref{fig:openloop_feedback} show that  state trajectories solving \fbocp  have slightly lower variance than those solving \olocp, which is the result of optimizing over a state feedback control trajectory. %

\section{Conclusion and outlook}\label{sec:conclusion} 
The optimality conditions in \pmp provide new insights onto the structure of solutions to stochastic optimal control problems when optimizing over deterministic open-loop controls or  parameterized feedback controls. 
By using rough path theory \rev{instead of It\^o calculus}, the results and analysis rely on a pathwise \rev{approach} instead of FBSDEs %
and can handle Gaussian processes $B$ that are not semimartingales, %
at the expense of stronger regularity assumptions on the diffusion  $\sigma$. 
The main motivation for deriving \pmp is the development of new algorithms for stochastic optimal control that 
can more easily borrow ideas from the deterministic optimal control literature, %
such as indirect shooting methods.

The following directions of future research are interesting.  
First, extending our results to the case where the diffusion $\sigma$ depends on the control is non-trivial, as  it may lead to a degenerate formulation  with irregular controls as described in \cite{Diehl2016,Allan2020}, and %
rough path theory relies on coefficients that are smooth-enough in time. 
Second, extending \pmp to tackle more general settings such as risk-averse optimal control problems \cite{LewMPC2024,Bonalli2023} would  provide valuable insights and be useful in applications. 
Finally, while the proposed indirect shooting method %
works well on the example in Section \ref{sec:example}, proving theoretical properties for the method such as asymptotic optimality \cite{Phelps2016,Shapiro2021,Lew2024,Melnikov2024} or robustness to discretization \cite{Allan2023} remains an open problem. %
As %
the unknown initial value of the adjoint vector $p_0$ is a random variable (i.e., the search space is  infinite-dimensional, as opposed to the deterministic setting where $p_0\in\R^n$), future analysis of the proposed indirect method may require innovative proof techniques and inspire \rev{the development of} faster and more robust algorithms.

\subsection*{Acknowledgements}
T.L. would like to thank Riccardo Bonalli for many stimulating discussions on the deterministic and stochastic PMPs  
and \rev{for} his helpful feedback. \rev{The author  also thanks the reviewers for their insightful comments and suggestions that have helped improve the  manuscript.}

\ifarxiv

\appendix

\addtocontents{toc}{\protect\setcounter{tocdepth}{1}}

\section{Additional proofs}\label{sec:appendix}

\textbf{Contents}
\begin{enumerate}
\item[]%
\hyperref[apdx:proofs:intro]{\ref{apdx:proofs:intro}\hspace{1mm} Preliminary remarks}
\hspace*{\fill} 
\pageref{apdx:proofs:intro}
\item[]%
\hyperref[apdx:proofs:preliminary_results]{\ref{apdx:proofs:preliminary_results}\hspace{1mm} Proofs of preliminary results (Section \ref{sec:preliminaries})}
\hspace*{\fill} 
\pageref{apdx:proofs:preliminary_results}
\item[]%
\hyperref[apdx:proofs:rdes]{\ref{apdx:proofs:rdes}\hspace{1mm} Proofs of rough differential equation results (Section \ref{sec:rdes})}
\hspace*{\fill} 
\pageref{apdx:proofs:rdes}
\item[]%
\hyperref[apdx:proofs:pmp]{\ref{apdx:proofs:pmp}\hspace{1mm} Additional proofs for the PMP (Section \ref{sec:pmp})}
\hspace*{\fill} 
\pageref{apdx:proofs:pmp}
\item[]%
\hyperref[apdx:proofs:example]{\ref{apdx:proofs:example}\hspace{1mm} Additional details for the indirect shooting method (Section \ref{sec:example})}
\hspace*{\fill} 
\pageref{apdx:proofs:example}
\end{enumerate}

\subsection{Preliminary remarks}\label{apdx:proofs:intro}
This appendix contains additional proofs for  results in the main manuscript. 
Thoughout, %
we use the inequalities 
$(|a_1|+|a_2|)^\frac{1}{p}\leq |a_1|^\frac{1}{p}+|a_2|^\frac{1}{p}$ and $|\sum_{i=1}^Na_i|^p\leq N^p\sum_{i=1}^N|a_i|^p$ 
for any $p\geq 1$, $a_i\in\R$, and $N\in\N$. 
Also, any function $f\in C_b^n$ satisfies the mean value theorem $
f(x+h)=f(x)
+
\frac{1}{k!}\sum_{k=1}^{n-1}\nabla^k f(x)h^{\otimes k}
+
\frac{1}{(n-1)!}
\int_0^1
\nabla^nf(x+\theta h)
h^{\otimes n}
(1-\theta)^{n-1}\dd\theta$. 

\subsection{Proofs of preliminary results (Section \ref{sec:preliminaries})}\label{apdx:proofs:preliminary_results}

\subsubsection{Rough paths, controlled rough paths, and rough integration (Section \ref{sec:preliminaries:rough_paths})}
\begin{proof}[Proof of Lemma \ref{lem:pvariation:inequalities}]
To show \eqref{eq:path_finite_var:infty_ineq}, we write
$$
\|X\|_\infty=\sup_{t\in[0,T]}\|X_t\|
\leq\|X_0\|+\sup_{t\in[0,T]}\|X_t-X_0\|^\frac{p}{p}
\leq 
\|X_0\|+\Big(\sup_{t\in[0,T]}\|X_t-X_0\|^p\Big)^\frac{1}{p}
\leq
\|X_0\|+\|X\|_p.
$$

To show \eqref{eq:path_finite_var:sum_p/2vars}, given an arbitrary partition  $\pi$ of $[0,T]$,
 \begin{align*}
 \sum_{[s,t]\in\pi}\|X_{s,t}\|^\frac{p}{2}
 &\leq
 \sum_{[s,t]\in\pi}
 \Big(
 \sum_{i=1}^n\|Y^i_{s,t}\|\|\widetilde{Y}^i_{s,t}\|+
 \sum_{j=1}^m\|Z^j_{s,t}\|
 +
 c|t-s|
 \Big)^\frac{p}{2}
 \\
 &\leq
 (n+m+1)^\frac{p}{2}
 \Big(
 \sum_{i=1}^n
 \sum_{[s,t]\in\pi}
 \|Y^i_{s,t}\|^\frac{p}{2}\|\widetilde{Y}^i_{s,t}\|^\frac{p}{2}
 +
 \sum_{j=1}^m
 \sum_{[s,t]\in\pi}
 \|Z^j_{s,t}\|^\frac{p}{2}
 +
 c^\frac{p}{2}\sum_{[s,t]\in\pi}
 |t-s|^\frac{p}{2}
 \Big)
 \\
 &\leq
 C_{\rev{m,n,}p}
 \Big(
 \sum_{i=1}^n
 \Big(\sum_{[s,t]\in\pi}
 \|Y^i_{s,t}\|^p
 \Big)^\frac{1}{2}
  \Big(\sum_{[s,t]\in\pi}
 \|\widetilde{Y}^i_{s,t}\|^p
 \Big)^\frac{1}{2}
 +
 \sum_{j=1}^m
 \|Z^j\|_\frac{p}{2}^\frac{p}{2}
 +
 (c\,T)^\frac{p}{2}
 \Big)
 \\
 &\leq
 C_{\rev{m,n,}p}
 \Big(
 \sum_{i=1}^n
 \|Y^i\|_p^\frac{p}{2}
 \|\widetilde{Y}^i\|_p^\frac{p}{2}
 +
 \sum_{j=1}^m
 \|Z^j\|_\frac{p}{2}^\frac{p}{2}
 +
 (c\,T)^\frac{p}{2}
 \Big),
 \end{align*}
 where we used $|\sum_{i=1}^Na_i|^\frac{p}{2}\leq N^\frac{p}{2}\sum_{i=1}^N|a_i|^\frac{p}{2}$  and H\"older's inequality $\sum_{i=1}^n|a_ib_i|\leq(\sum_{i=1}^n|a_i|^p)^\frac{1}{p}(\sum_{i=1}^n|b_i|^q)^\frac{1}{q}$ for any $a,b\in\R^n$ and $\frac{1}{p}+\frac{1}{q}=1$, so that
 $\sum_{i=1}^n\left(|a_i|^\frac{p}{2}|b_i|^\frac{p}{2}\right)
 \leq
 \left(\sum_{i=1}^n|a_i|^p\right)^\frac{1}{2}\left(\sum_{i=1}^n|b_i|^p\right)^\frac{1}{2}$.
The bound above is independent of the choice of partition $\pi$. 
 Thus, \eqref{eq:path_finite_var:sum_p/2vars} follows after taking the supremum over all partitions $\pi$ of $[0,T]$ and using the inequality $(\sum_i|a_i|)^\frac{2}{p}\leq\sum_i|a_i|^\frac{2}{p}$ for $p\geq 2$.  \eqref{eq:path_finite_var:sum_pvars} is shown similarly.   

To show $\|X\|_p\leq\|X\|_\frac{p}{2}$, we  write 
$
\|X\|_p^p
=\sup_\pi\sum_{[s,t]\in\pi}\left(\|X_{s,t}\|^\frac{p}{2}\right)^2
\leq
(
\sup_\pi\sum_{[s,t]\in\pi}\|X_{s,t}\|^\frac{p}{2})^2
=
\|X\|_\frac{p}{2}^p
$ since $\sum_i a_i^2\leq\left(\sum_i a_i\right)^2$ for $a_i\geq 0$.

To show \eqref{eq:sigma(.,X):pvar}, we use the smoothness of $\sigma$: for all $s,t\in[0,T]$,
\begin{align}\label{eq:sigma(.,Y)_st_bound}
\|\sigma(\cdot,X)_{s,t}\|&\leq\|\sigma(t,X_t)-\sigma(t,X_s)\|+\|\sigma(t,X_s)-\sigma(s,X_s)\|\leq \|\sigma\|_{C_b^1}(
\|X_{s,t}\|
+
|t-s|), 
\end{align}
and we conclude with \eqref{eq:path_finite_var:sum_pvars}.  

To show \eqref{eq:Delta_sigma(.,X):pvar},  
for any $s,t\in[0,T]$, we write
\begin{align}\label{eq:deltaf_timevarying:st:0}
\|(\sigma(\cdot,X)-\sigma(\cdot,\rev{\widetilde{X}}))_{s,t}\|
&\leq
\|
(\sigma(t,X)-\sigma(t,\rev{\widetilde{X}}))_{s,t}
\| 
+
\| 
(\sigma(\cdot,X_s)-\sigma(\cdot,\rev{\widetilde{X}}_s))_{s,t}
\|.
\end{align}
We bound the first term in \eqref{eq:deltaf_timevarying:st:0} next. Denoting $\Delta X = X-\rev{\widetilde{X}}$ and for any $\theta\in\R$,
$$
\rev{\nablaof{x}\sigma}(t,\rev{\widetilde{X}}_t+\theta\Delta X_t)
=
\rev{\nablaof{x}\sigma}(t,\rev{\widetilde{X}}_s+\theta\Delta X_s)
+
\int_0^1
\rev{\nablaof{x}^2\sigma}(t,\rev{\widetilde{X}}_s+\theta\Delta X_s+\zeta(\rev{\widetilde{X}}_{s,t}+\theta\Delta X_{s,t}))\dd\zeta
((1-\theta)\rev{\widetilde{X}}_{s,t}+\theta X_{s,t})
$$
so that by denoting $Z_{\theta,\zeta}=\rev{\widetilde{X}}_s+\theta\Delta X_s+\zeta(\rev{\widetilde{X}}_{s,t}+\theta\Delta X_{s,t})$ and using the mean value theorem again,
\begin{align*}
\|\sigma(t,X)_{s,t}-\sigma(t,\rev{\widetilde{X}})_{s,t}\|
&=
\left\|
\int_0^1\nabla \sigma(t,\rev{\widetilde{X}}_t+\theta\Delta X_t)\dd\theta\Delta X_t
-
\int_0^1\nabla \sigma(t,\rev{\widetilde{X}}_s+\theta\Delta X_s)\dd\theta\Delta X_s
\right\|
\\
&\hspace{-3cm}=
\left\|
\int_0^1\nabla \sigma(t,\rev{\widetilde{X}}_s+\theta\Delta X_s)\dd\theta(\Delta X_t-\Delta X_s)
+
\int_0^1\int_0^1
\nabla^2 \sigma(t,Z_{\theta,\zeta})\dd\zeta
((1-\theta)\rev{\widetilde{X}}_{s,t}+\theta X_{s,t})
\dd\theta\Delta X_t 
\right\|
\\
&\leq
\|\sigma\|_{C_b^2}
\left(
\|\Delta X_{s,t}\|
+
(\|X_{s,t}\|+\|\rev{\widetilde{X}}_{s,t}\|)\|\Delta X\|_\infty
\right).
\end{align*}
The second term in \eqref{eq:deltaf_timevarying:st:0} is bounded similarly:
\begin{align}\label{eq:deltaf_timevarying:st:2ndterm}
\| 
(\sigma(\cdot,X_s){-}\sigma(\cdot,\rev{\widetilde{X}}_s))_{s,t}
\|
&=
\left\|
\int_0^1\hspace{-1mm}
\Big(\rev{\nablaof{x}\sigma}(t,\rev{\widetilde{X}}_s+\theta\Delta X_s){-} \rev{\nablaof{x}\sigma}(s,\rev{\widetilde{X}}_s+\theta\Delta X_s)\Big)
\dd\theta\Delta X_s
\right\|
\leq \|\sigma\|_{C_b^2}|t-s|\|\Delta X\|_\infty.
\end{align}
Thus, \eqref{eq:deltaf_timevarying:st:0} can be bounded as
\begin{align}\label{eq:deltaf_timevarying:st}
\|(\sigma(\cdot,X)-\sigma(\cdot,\rev{\widetilde{X}}))_{s,t}\|
&\leq
\|\sigma\|_{C_b^2}
\big(
\|\Delta X_{s,t}\|
+
(\|X_{s,t}\|+\|\rev{\widetilde{X}}_{s,t}\|+
|t-s|)\|\Delta X\|_\infty
\big).
\end{align}
The conclusion follows from  \eqref{eq:path_finite_var:sum_pvars} and $\|\Delta X\|_\infty\leq\|\Delta X_0\|+\|\Delta X\|_p$ in \eqref{eq:path_finite_var:infty_ineq}.
\end{proof}

\begin{proof}[Proof of Lemma \ref{lem:rough_path:sigma(.,Y):controlled}]
To show \eqref{eq:controlled_path:pvar_norm} and \eqref{eq:|Y|_p_ineq},   we write
\begin{align*}
\|Y_{s,t}\|
\mathop{\leq}^{\eqref{eq:remainder}}
\|Y'_s\|\|X_{s,t}\|+\|R^Y_{s,t}\|
\leq 
\|Y'\|_\infty\|X_{s,t}\|+\|R^Y_{s,t}\|
\mathop{\leq}^{\eqref{eq:path_finite_var:infty_ineq}}
(\|Y'_0\|+\|Y'\|_p)\|X_{s,t}\|+\|R^Y_{s,t}\|
\end{align*}
for any $s,t\in[0,T]$. %
Then,  \eqref{eq:controlled_path:pvar_norm} and \eqref{eq:|Y|_p_ineq} follow from \eqref{eq:path_finite_var:sum_pvars}  and 
$\|R^Y\|_p\leq\|R^Y\|_\frac{p}{2}$ in Lemma \ref{lem:pvariation:inequalities}.

To show \eqref{eq:sigma(.,Y):pvar}, for any $s,t\in[0,T]$, we write
\begin{align*}
\|\sigma(\cdot,Y)_{s,t}\|
\mathop{\leq}^{\eqref{eq:sigma(.,Y)_st_bound},\eqref{eq:remainder}}
 \|\sigma\|_{C_b^1}(\|Y'_s\|\|X_{s,t}\|+\|R^Y_{s,t}\|+|t-s|)
\mathop{\leq}^{\eqref{eq:path_finite_var:infty_ineq}} 
\|\sigma\|_{C_b^1}(M_{Y'}\|X_{s,t}\|+\|R^Y_{s,t}\|+|t-s|^\frac{2}{p}),
\end{align*}
where we used %
$\|Y'\|_\infty\leq\|Y'_0\|+\|Y'\|_p=M_{Y'}$ %
in the second inequality. Then, \eqref{eq:sigma(.,Y):pvar} follows from \eqref{eq:path_finite_var:sum_pvars}.

To show \eqref{eq:sigma(.,Y)':pvar}, denoting $\nablax\sigma=\rev{\nablaof{x}\sigma}$, for any $s,t\in[0,T]$, we write
\begin{align*}
&\|\sigma(\cdot,Y)'_{s,t}\|
=
\|
\nablax\sigma(t,Y_t)Y_t'
-\nablax\sigma(t,Y_t)Y_s'+\nablax\sigma(t,Y_t)Y_s'
-\nablax\sigma(s,Y_t)Y_s'+\nablax\sigma(s,Y_t)Y_s'
-\nablax\sigma(s,Y_s)Y_s'
\|
\\
&\hspace{17mm}\leq 
\|\nablax\sigma(t,Y_t)\|\|Y'_t-Y'_s\|+
\|\nablax\sigma(t,Y_t)-\nablax\sigma(s,Y_t)\|\|Y_s'\| +
\|\nablax\sigma(s,Y_t)-\nablax\sigma(s,Y_s)\|\|Y_s'\| 
\\
&\hspace{17mm}\leq 
\|\sigma\|_{C^2_b}(\|Y'_{s,t}\|
+
(|t-s|
+
\|Y_{s,t}\|)
\|Y_s'\|
)
\\
&\hspace{17mm}\leq 
\|\sigma\|_{C^2_b}(
\|Y'_{s,t}\|+
(
|t-s|
+
M_{Y'}\|X_{s,t}\|+\|R^Y_{s,t}\|)
M_{Y'}
),
\\
&\implies
\|\sigma(\cdot,Y)'\|_p
\mathop{\leq}^{\eqref{eq:path_finite_var:sum_pvars}}
C_p\|\sigma\|_{C^2_b}(
\|Y'\|_p
+
M_{Y'}(
T
+
M_{Y'}\|X\|_p+\|R^Y\|_\frac{p}{2}
)
)
\\
&\hspace{27mm}\leq
C_p\|\sigma\|_{C^2_b}M_{Y'}(
1+
T+
M_{Y'}\|X\|_p+\|R^Y\|_\frac{p}{2})
\\
&\hspace{27mm}\leq
C_p\|\sigma\|_{C^2_b}K_Y(1+K_Y+T)(1+\|X\|_p).
\hspace{3.5cm}
(K_Y =  M_{Y'}+\|R^Y\|_\frac{p}{2})
\end{align*}

To show \eqref{eq:RY:p/2var:Y^2+RY+T} and \eqref{eq:RY:p/2var:KYs}, for any $s,t\in[0,T]$,
\begin{align*}
\|R^{\sigma(\cdot,Y)}_{s,t}\|
&\mathop{=}^{\eqref{eq:remainder}}
\|
\sigma(s,Y_t)-\sigma(s,Y_s)-\nablax \sigma(s,Y_s)Y'_sX_{s,t}
+\sigma(t,Y_t)-\sigma(s,Y_t)
\|
\\
&\mathop{=}^{\eqref{eq:remainder}}
\|
\sigma(s,Y_t)-\sigma(s,Y_s)-\nablax \sigma(s,Y_s)Y_{s,t}-\nablax \sigma(s,Y_s)R^Y_{s,t}
+
\sigma(t,Y_t)-\sigma(s,Y_t)\|
\\
&\leq
\left\|
\int_0^1\nablax^2 \sigma(s,Y_s+\theta Y_{s,t})Y_{s,t}^{\otimes 2}(1-\theta)\dd\theta
\right\| 
+\|\nablax \sigma(s,Y_s)R^Y_{s,t}\|
+
\|\sigma(t,Y_t)-\sigma(s,Y_t)\|
\\
&\leq 
\|\sigma\|_{C^2_b}(\|Y_{s,t}\|^2+\|R^Y_{s,t}\|+|t-s|),
\\
\implies
\|R^{\sigma(\cdot,Y)}\|_\frac{p}{2}
&\mathop{\leq}^{\eqref{eq:path_finite_var:sum_p/2vars}}
C_p
\|\sigma\|_{C^2_b}(
\|Y\|_p^2+\|R^Y\|_\frac{p}{2}+T)
\\
&\mathop{\leq}^{\eqref{eq:|Y|_p_ineq}}
C_p\|\sigma\|_{C^2_b}(
(1+\|X\|_p)^2K_Y^2
+
\|R^Y\|_\frac{p}{2}
+
T
)
\\
&\leq
C_p\|\sigma\|_{C^2_b}(
K_Y(K_Y(1+\|X\|_p)^2+1))+T
)
\leq
\eqref{eq:RY:p/2var:KYs}.
\end{align*}

Finally, from the last three inequalities,  $(\sigma(\cdot,Y),\sigma(\cdot,Y)')\in\C^p\times\C^p$ and $\|R^{\sigma(\cdot,Y)}\|_\frac{p}{2}<\infty$, so we conclude that $(\sigma(\cdot,Y),\sigma(\cdot,Y)')\in\sD^p_X$. 
\end{proof}

\begin{proof}[Proof of Lemma \ref{lem:rough_path:ito_formula}]
This result is a particular case of  \cite[Theorem 7.7]{Friz2020}, assuming that the rough path $\mbX$ is geometric, so that its bracket is zero. The proof below is standard. We may assume that $f\in C^3_b$  since it is only evaluated on the image of the path $Y$ which is bounded.

First,   the rough integral in \eqref{eq:rough_path:ito_formula} is well-defined, since $(\nabla f(Y)Y',(\nabla f(Y)Y')')\in\sD^p_X$ 
with $(\nabla f(Y)Y')'=(\nabla f(Y)Y''+\nabla f^2(Y)(Y'\otimes Y'))$ by Lemma \ref{lem:control_path:product} (note that  Lemma \ref{lem:rough_path:ito_formula} is not used in the proof of Lemma \ref{lem:control_path:product}). The   Young integral in \eqref{eq:rough_path:ito_formula} is also well-defined  since $\nabla f(Y)\in\C^p$ and $\Gamma\in\C^{\frac{p}{2}}$  so that   $\frac{1}{p}+\frac{2}{p}>1$ (see for example  \cite[Theorem 6.8]{Friz2010} or \cite[Proposition 5.2]{Allan2021}). 

Next, we show \eqref{eq:rough_path:ito_formula}. Let $s,t\in[0,T]$. Then, %
$$
Y_{s,t}
\mathop{=}^{\eqref{eq:rough_int:error_bound}}
Y'_sX_{s,t}+Y''_s\bX_{s,t}+\Gamma_{s,t}+K_{s,t},
\ \text{with }
\|K_{s,t}\|\leq C_p(\|R^Y\|_{\frac{p}{2},[s,t]}\|X\|_{p,[s,t]}+\|Y'\|_{p,[s,t]}\|\bX\|_{\frac{p}{2},[s,t]}),
$$
and by the mean value theorem, with $Z_{s,t}:=\frac{1}{2}\int_0^1\nabla^3f(Y_s+\theta Y_{s,t})Y_{s,t}^{\otimes 3}(1-\theta)^2\dd\theta$,
\begin{align*}
f(Y)_{s,t}
&=
\nabla f(Y_s)Y_{s,t}+\frac{1}{2}\nabla^2f(Y_s)Y_{s,t}^{\otimes 2}
+
Z_{s,t}
\\
&=
\nabla f(Y_s)(Y'_sX_{s,t}+Y''_s\bX_{s,t}+\Gamma_{s,t})
+
\nabla f(Y_s)K_{s,t}
+
\frac{1}{2}\nabla^2f(Y_s)(Y'_sX_{s,t}+R^Y_{s,t})^{\otimes 2}
+
Z_{s,t}
\\
&=
\nabla f(Y_s)(Y'_sX_{s,t}+Y''_s\bX_{s,t})+\nabla f(Y_s)\Gamma_{s,t}
+
\frac{1}{2}\nabla^2f(Y_s)(Y'_s\otimes Y'_s)(X_{s,t}\otimes X_{s,t})
+
\\
&\qquad\quad
\frac{1}{2}\nabla^2f(Y_s)
(Y'_sX_{s,t}\otimes R^Y_{s,t}+R^Y_{s,t}\otimes Y'_sX_{s,t} +R^Y_{s,t}\otimes R^Y_{s,t})
+
\nabla f(Y_s)K_{s,t}
+
Z_{s,t},
\end{align*}
Next, define $\text{Sym}(\bX_{s,t})$ as  the symmetric part of $\bX_{s,t}$, with $\text{Sym}(\bX_{s,t})^{ij}=\frac{1}{2}(\bX_{s,t}^{ij}+\bX_{s,t}^{ji})$. Since $\mbX$ is geometric, $X_{s,t}^{ij}X_{s,t}^{ji}=\bX_{s,t}^{ij}+\bX_{s,t}^{ji}=2\text{Sym}(\bX_{s,t})^{ij}$, so $\frac{1}{2}(X_{s,t}\otimes X_{s,t})=\text{Sym}(\bX_{s,t})$. Also, $\nabla^2f(Y_s)(Y'_s\otimes Y'_s)$ is symmetric, so $\nabla^2f(Y_s)(Y'_s\otimes Y'_s)\bX_{s,t}=\nabla^2f(Y_s)(Y'_s\otimes Y'_s)\text{Sym}(\bX_{s,t})$, since ``the contraction of a symmetric matrix with an antisymmetric matrix is zero'' \cite[Proof of Proposition 6.9]{Allan2021}. Thus, 
\begin{align*}
f(Y)_{s,t} 
&=
\nabla f(Y_s)Y'_sX_{s,t}
+
(\nabla f(Y_s)Y''_s+\nabla f^2(Y_s)(Y'_s\otimes Y'_s))\bX_{s,t}
+
\nabla f(Y_s)\Gamma_{s,t} 
+
\\
&\qquad\quad
\frac{1}{2}\nabla^2f(Y_s)
(Y'_sX_{s,t}\otimes R^Y_{s,t}+R^Y_{s,t}\otimes Y'_sX_{s,t} +R^Y_{s,t}\otimes R^Y_{s,t})
+
\nabla f(Y_s)K_{s,t}
+
Z_{s,t}
\\
&=
(\nabla f(Y)Y')_sX_{s,t}
+
(\nabla f(Y)Y')'_s\bX_{s,t}
+
\nabla f(Y_s)\Gamma_{s,t} 
+
W_{s,t},
\end{align*}
where %
$W_{s,t}:=\frac{1}{2}\nabla^2f(Y_s)
(Y'_sX_{s,t}\otimes R^Y_{s,t}+R^Y_{s,t}\otimes Y'_sX_{s,t} +R^Y_{s,t}\otimes R^Y_{s,t})
+
\nabla f(Y_s)K_{s,t}
+
Z_{s,t}$ can be bounded as
$$
\|W_{s,t}\|\leq C_{p,f,\|Y\|_\infty,\|Y'\|_\infty}
\left(
\|X\|_{p,[s,t]}\|R^Y\|_{\frac{p}{2},[s,t]}
+
\|R^Y\|_{\frac{p}{2},[s,t]}^2
+
\|Y'\|_{p,[s,t]}\|\bX\|_{\frac{p}{2},[s,t]}
+
\|Y\|_{p,[s,t]}^3
\right),
$$
Since $X,Y,Y'$ have finite $p$-variation and $\bX,R^Y$ have finite $\frac{p}{2}$-variation, up to a time reparameterization, we may assume that $X,Y,Y'$ are $\frac{1}{p}$-H\"older continuous and $\bX,R^Y$ are $\frac{2}{p}$-H\"older continuous \cite[Proposition 5.14]{Friz2010}, so that
$$
\|W_{s,t}\|\leq C_{p,f,X,\mbX,Y,Y',R^Y}|t-s|^\frac{3}{p}.
$$
Thus, for an arbitrary partition $\pi$ of $[0,T]$, we obtain
\begin{align*}
f(Y_T)-f(Y_0)
=
\sum_{[s,t]\in\pi}
(\nabla f(Y)Y')_sX_{s,t}
+
(\nabla f(Y)Y')'_s\bX_{s,t}
+
\sum_{[s,t]\in\pi}
\nabla f(Y_s)\Gamma_{s,t}
+
\sum_{[s,t]\in\pi}
W_{s,t}.
\end{align*}
The conclusion follows after taking the limit over all partitions $\pi$ of $[0,T]$ with vanishing mesh size in the above: %
The first sum is a rough integral, 
the second is a Young integral, which is  well-defined  by \cite[Theorem 6.8]{Friz2010},  since $\nabla f(Y)\in\C^p$ and $\Gamma\in\C^{\frac{p}{2}}$  so that   $\frac{1}{p}+\frac{2}{p}>1$, and the third   is zero, since 
$$
\sum_{[s,t]\in\pi}
\|W_{s,t}\|
\leq 
C_{p,f,X,\mbX,Y,Y',R^Y} 
\sum_{[s,t]\in\pi}|t-s|^\frac{p}{p}
|\pi|^{\frac{3-p}{p}}
=
C_{p,f,X,\mbX,Y,Y',R^Y}
T
|\pi|^{\frac{3-p}{p}}
\to 0\ \text{as }|\pi|\to0,
$$
so that $\lim_{|\pi|\to0}\sum_{[s,t]\in\pi}
W_{s,t}=0$. 
This concludes the proof of Lemma \ref{lem:rough_path:ito_formula}.
\end{proof}

\begin{proof}[Proof of Lemma \ref{lem:pvar:RX_intervals}]
\rev{From Chen's relation	 
\eqref{eq:chen's_relation} and from \eqref{eq:remainder},  $\bX_{s,t}=\bX_{s,u}+\bX_{u,t}+X_{s,u}\otimes X_{u,t}$ and $R^Y_{s,t}=R^Y_{s,u}+R^Y_{u,t}+Y'_{s,u}X_{u,t}$ for any $s,u,t\in[0,T]$. Thus, for any partition $s=u_0<u_1<\dots<u_N=t$ of $[s,t]\subseteq[0,T]$,
$$
\bX_{s,t}=\sum_{i=1}^N\bX_{u_{i-1},u_i}+\sum_{i=1}^{N-1}X_{u_{i-1},u_i}\otimes X_{u_i,t},
\ \ \text{and}\ \ 
R^Y_{s,t}=\sum_{i=1}^NR^Y_{u_{i-1},u_i}+\sum_{i=1}^{N-1}Y'_{u_{i-1},u_i}X_{u_i,t}.
$$
The first equation is in \cite[Exercise 2.5]{Friz2020}, and the second  is derived similarly. 
 Next, to show \eqref{eq:pvar:RX_intervals}, we follow the proof of \cite[Lemma 2.3]{Allan2020}. Let $0=s_0<s_1<\dots<s_N=T$ be any partition of $[0,T]$. Then, the union of the partitions $\{t_i\}_{i=0}^n$ and $\{s_j\}_{j=0}^N$ of $[0,T]$ can be labeled either as
$$
s_{j-1}=t_0^j<t_1^j<\dots<t_{n_j}^j=s_j
\ \, \forall j=1,\dots,N, 
\quad \text{or as}\ \ 
t_{i-1}=s_0^i<s_1^i<\dots<s_{N_i}^i=t_i
\ \, \forall i=1,\dots,n,
$$
where $n_j\leq n$ for all $j=1,\dots,N$.} 
\rev{Next, we have
$$
\sum_{j=1}^N\sum_{i=1}^{n_j}\|R^Y_{t^j_{i-1},t^j_i}\|^\frac{p}{2}
\leq
\sum_{i=1}^n\sum_{j=1}^{N_i}\|R^Y_{s^i_{j-1},s^i_j}\|^\frac{p}{2}
\leq
\sum_{i=1}^n\|R\|^\frac{p}{2}_{\frac{p}{2},[t_{i-1},t_i]},
\ \ \text{and}\ \  
\sum_{j=1}^N\sum_{i=1}^{n_j}\|Y'_{t^j_{i-1},t^j_i}\|^p 
\leq
\sum_{i=1}^n\|Y'\|^p_{p,[t_{i-1},t_i]}.
$$
Also, using the bounds $\|X_{t^j_i,t^j_{n_j}}\|\leq\sum_{i=1}^{n_j}\|X_{t^j_{i-1},t^j_i}\|$ and $|\sum_{i=1}^na_i|^p\leq n^p\sum_{i=1}^n|a_i|^p$, we obtain
$$
\sum_{j=1}^N\sum_{i=1}^{n_j}
\|X_{t^j_i,t^j_{n_j}}\|^p 
\leq
\sum_{j=1}^N\sum_{i=1}^{n_j}\left(\sum_{i=1}^{n_j}
\|X_{t^j_{i-1},t^j_i}\| 
\right)^p
\leq
n^p
\sum_{j=1}^N\sum_{i=1}^{n_j}
\|X_{t^j_{i-1},t^j_i}\|^p
\leq
n^p
\sum_{i=1}^n 
\|X\|^p_{p,[t_{i-1},t_i]}.
$$
Finally, we combine the inequalities above and obtain
\begin{align*}
\sum_{j=1}^N\|R^Y_{s_{j-1},s_j}\|^\frac{p}{2}
&=
\sum_{j=1}^N\bigg\|
\sum_{i=1}^{n_j}R^Y_{t_{i-1}^j,t_i^j}
+
\sum_{i=1}^{n_j}Y'_{t_{i-1}^j,t_i^j}X_{t_i^j,t_{n_j}^j}
\bigg\|^\frac{p}{2}
\\
&\leq
2^\frac{p}{2}n^\frac{p}{2}
\bigg(
\sum_{j=1}^N
\sum_{i=1}^{n_j}
\left\|
R^Y_{t_{i-1}^j,t_i^j}
\right\|^\frac{p}{2}
+
\sum_{j=1}^N
\sum_{i=1}^{n_j}
\left\|
Y'_{t_{i-1}^j,t_i^j}X_{t_i^j,t_{n_j}^j}
\right\|^\frac{p}{2}
\bigg)
\\
&\leq
2^\frac{p}{2}n^\frac{p}{2}
\bigg(
\sum_{j=1}^N
\sum_{i=1}^{n_j}
\left\|
R^Y_{t_{i-1}^j,t_i^j}
\right\|^\frac{p}{2}
+
\bigg(
\sum_{j=1}^N
\sum_{i=1}^{n_j} 
\left\|
Y'_{t_{i-1}^j,t_i^j}
\right\|^p
\bigg)^\frac{1}{2}
\bigg(
\sum_{j=1}^N
\sum_{i=1}^{n_j} 
\left\|
X_{t_i^j,t_{n_j}^j}
\right\|^p
\bigg)^\frac{1}{2}
\bigg)
\\
&\leq
2^\frac{p}{2}n^\frac{p}{2}n^\frac{p}{2}
\bigg(
\sum_{i=1}^n\left\|R^Y\right\|^\frac{p}{2}_{\frac{p}{2},[t_{i-1},t_i]}
+
\bigg(
\sum_{i=1}^n 
\left\|Y'\right\|^p_{p,[t_{i-1},t_i]}
\bigg)^\frac{1}{2}
\bigg(
\sum_{i=1}^n 
\left\|X\right\|^p_{p,[t_{i-1},t_i]}
\bigg)^\frac{1}{2}
\bigg).
\end{align*}
As the choice of partition $\{s_j\}_{j=0}^N$ was arbitrary,   \eqref{eq:pvar:RX_intervals} follows after taking the supremum over all partitions of $[0,T]$. The proof of  \eqref{eq:pvar:bX_intervals} is identical.}
\end{proof}

\subsubsection{The greedy partition and Gaussian rough paths (Section \ref{sec:preliminaries:greedy_gaussian_paths})}
\begin{proof}[Proof of Lemma \ref{lem:Nalpha<=w(0,T)}] The proof follows \cite[Lemma 4.9]{Cass2013}. Define the greedy partition $\{\tau_i\}_{i=0}^{N_{\alpha,[s,t]}(w)+1}$ of the interval $[0,T]$. Since $w$ is a control, $w(\tau_i,\tau_{i+1})+w(\tau_{i+1},\tau_{i+2})\leq w(\tau_i,\tau_{i+2})$, so 
\begin{align*}
\alpha N_{\alpha,[s,t]}(w)
&=
\sum_{i=0}^{N_{\alpha,[s,t]}(w)-1}
w(\tau_i,\tau_{i+1})
\leq
w(0,\tau_{N_{\alpha,[s,t]}(w)})
\leq
w(0,T),
\end{align*}
and the conclusion follows.
\end{proof}
\begin{proof}[Proof of Lemma \ref{lem:sum_Nalpha}]
The proof is inspired from the proof of \cite[Lemma 3]{Friz2013}. 
First, we define the accumulated $\alpha$-local $w$-variations \cite[Definition 4.1]{Cass2013}
$$
w_\alpha(s,t)=
\sup_{\substack{
\pi=\{t_i\}\subset[s,t]
\\
w(t_i,t_{i+1})\leq\alpha}
}
\sum_{i=0}^{N-1}w(t_i,t_{i+1}),
\qquad
w_{j,\alpha}(s,t)=
\sup_{\substack{
\pi=\{t_i\}\subset[s,t]
\\
w_j(t_i,t_{i+1})\leq\alpha}
}
\sum_{i=0}^{N-1}w_j(t_i,t_{i+1}),
\ \, j=1,\dots,n,
$$
where the supremums are over all partitions $\pi=\{s=t_0<t_1<\dots<t_N=t\}$ of $[s,t]$ such that $w(t_i,t_{i+1})\leq\alpha$  (or $w_j(t_i,t_{i+1})\leq\alpha$, respectively) for all $i=0,\dots,N-1$. We have
\begin{align*}
\frac{w_\alpha(s,t)}{C}
&\leq
\sup_{\substack{
\pi=\{t_i\}\subset[s,t]
\\
w(t_i,t_{i+1})\leq\alpha}
}
\sum_{i=0}^{N-1}\sum_{j=1}^nw_j(t_i,t_{i+1})
\leq
\sum_{j=1}^n
\sup_{\substack{
\pi=\{t_i\}\subset[s,t]
\\
w_j(t_i,t_{i+1})\leq\alpha}
}
\sum_{i=0}^{N-1}w_j(t_i,t_{i+1})
=
\sum_{j=1}^n
w_{j,\alpha}(s,t).
\end{align*}
where the last inequality follows from the fact that $w_j(t_i,t_{i+1})\leq\alpha$ for all $j$ if $w(t_i,t_{i+1})\leq\alpha$. 
 Then,
\begin{align*}
\alpha N_{\alpha,[s,t]}(w)
&=
\sum_{i=0}^{N_{\alpha,[s,t]}(w)-1}
w(\tau_i,\tau_{i+1})
\leq
w_\alpha(s,t)
\leq
C\sum_{j=1}^n
w_{j,\alpha}(s,t).
\end{align*}
Finally, from \cite[Proposition 4.11]{Cass2013} (see also the proof of \cite[Lemma 3]{Friz2013}), $w_{j,\alpha}(s,t)\leq \alpha(2N_{\alpha,[s,t]}(w_j)+1)$ for $j=1,\dots,n$, and the conclusion follows.
\end{proof}

\begin{proof}[Proof of Corollary \ref{cor:Nalpha:NX_NXtilde_NT:small_intervals}]
First, $w$ is continuous, $w(t,t)=0$, and $w(s,t)+w(t,u)\leq w(s,u)$ for any $\rev{s,t,u\in[0,T]}$, so $w$ is a control. 
Second, consider the control $w_T$ defined by $w_T(s,t)=|t-s|$, so that $w(s,t)=C_p(w_\mbX(s,t)+w_{\widetilde{\mbX}}(s,t)+w_T(s,t))$. Any interval $[s,t]\subseteq[0,T]$ can be partitioned into intervals of size at most $\alpha$, so $N_{\alpha,[s,t]}(w_T)\leq 1+|t-s|/\alpha$. The inequality \eqref{eq:Nalpha<=3CpNalpha_X_and_time} then follows from \eqref{eq:Nalpha(w)<=sumNalpha(wj)} in Lemma \ref{lem:sum_Nalpha}.

To show \eqref{eq:|X|+|Xtilde|+|dt|<=alpha^p},  let   $0<\alpha\leq 1$  and     $[s,t]\subseteq[0,T]$ be an interval small-enough to satisfy $w(s,t)\leq \alpha$, so that
$$
\|\bX\|_{\frac{p}{2},[s,t]}^\frac{p}{2}\leq w_\mbX(s,t)\leq w(s,t)\leq\alpha\leq 1,
\qquad
|t-s|\leq w(s,t)\leq\alpha\leq 1.
$$
Then, %
$|t-s|^p\leq |t-s|$ and $\|\bX\|_{\frac{p}{2},[s,t]}^p\leq\|\bX\|_{\frac{p}{2},[s,t]}^\frac{p}{2}$, so that
{\small
\begin{align}\label{eq:mbX^p<=wX}
\|\mbX\|_{p,[s,t]}^p
&=
\big(
\|X\|_{p,[s,t]}+\|\bX\|_{\frac{p}{2},[s,t]}
\big)^p
\leq
2^p\big(
\|X\|_{p,[s,t]}^p+\|\bX\|_{\frac{p}{2},[s,t]}^p
\big)
\leq
2^p\big(
\|X\|_{p,[s,t]}^p+\|\bX\|_{\frac{p}{2},[s,t]}^\frac{p}{2}
\big)
=
2^pw_\mbX(s,t).
\end{align}
}%
Thus,
{\small
\begin{align*}
(
\|\mbX\|_{p,[s,t]}+\|\widetilde{\mbX}\|_{p,[s,t]}+|t-s|
)^p
&\leq
3^p(
\|\mbX\|_{p,[s,t]}^p+\|\widetilde{\mbX}\|_{p,[s,t]}^p+|t-s|^p
)
\leq
6^p(
w_{\mbX}(s,t)+
w_{\widetilde{\mbX}}(s,t)+
|t-s|
)
= w(s,t),
\end{align*}
}%
and \eqref{eq:|X|+|Xtilde|+|dt|<=alpha^p} follows with $w(s,t)\leq\alpha$.

To show \eqref{eq:||X||_p<=exp(N_alpha^p)}, consider the greedy partition $\{\tau_i, i=0,1,\dots,N_{\alpha,[s,t]}(w)+1\}$, which is such that \eqref{eq:|X|+|Xtilde|+|dt|<=alpha^p} holds on any subinterval $[\tau_i,\tau_{i+1}]$ of this partition since $w(\tau_i,\tau_{i+1})\leq\alpha$, so that  $\|\mbX\|_{p,[\tau_i,\tau_{i+1}]}\leq\alpha^\frac{1}{p}$, $\|\widetilde{\mbX}\|_{p,[\tau_i,\tau_{i+1}]}\leq\alpha^\frac{1}{p}$, and $|\tau_{i+1}-\tau_i|\leq\alpha^\frac{1}{p}$. Then, by Lemma \ref{lem:pvar:intervals} \rev{and Lemma \ref{lem:pvar:RX_intervals}}, with $N:=N_{\alpha,[s,t]}(w)$,
\rev{\begin{align*}
\|\mbX\|_{p,[s,t]}
&\leq
(N+1)
\bigg(
\sum_{i=0}^N
\|X\|_{p,[\tau_i,\tau_{i+1}]}^p
\bigg)^\frac{1}{p}
+
2(N+1)^2
\bigg(
\bigg(
\sum_{i=0}^N\|\bX\|^\frac{p}{2}_{\frac{p}{2},[\tau_i,\tau_{i+1}]}
\bigg)^\frac{2}{p}
+
\bigg(
\sum_{i=0}^N\|X\|^p_{p,[\tau_i,\tau_{i+1}]}
\bigg)^\frac{2}{p}
\bigg)
\\
&\leq
5(N+1)^3\alpha^\frac{1}{p},
\end{align*}
where we used $\|X\|_{p,[\tau_i,\tau_{i+1}]}\leq\alpha^\frac{1}{p}$, $\|\bX\|_{\frac{p}{2},[\tau_i,\tau_{i+1}]}\leq\alpha^\frac{1}{p}$, and $\alpha\leq 1$ so $\alpha^\frac{2}{p}\leq\alpha^\frac{1}{p}$. Then,}
\begin{align*}
\|\mbX\|_{p,[s,t]}+\|\widetilde{\mbX}\|_{p,[s,t]}+|t-s|
&\leq
\rev{5(N+1)^3\alpha^\frac{1}{p}+5(N+1)^3\alpha^\frac{1}{p}}
+
\sum_{i=0}^N|\tau_{i+1}-\tau_i|
\leq
\rev{11\alpha^\frac{1}{p}(N+1)^3.}
\end{align*}
\rev{U}sing \rev{$
11\alpha^\frac{1}{p}(N+1)^3\leq
11\alpha^\frac{1}{p}3!\exp(N+1)\leq 66e\alpha^\frac{1}{p}\exp(N)$}, 
we conclude that \eqref{eq:||X||_p<=exp(N_alpha^p)} holds  
with  $C_{p,\alpha}=\rev{66}e\alpha^\frac{1}{p}$, which concludes the proof.
\end{proof}

\begin{proof}[Proof of Theorem \ref{thm:gaussian_rough_paths}]
The first claim follows from \cite[Theorem
15.33]{Friz2010}  (or \cite[Theorem 10.4]{Friz2020}, see also \cite[Corollary 2.3]{Friz2016}). 
Then, following the proof of \cite[Theorem 11]{Bayer2016}, for $\frac{1}{q}=\frac{1}{2\rho}+\frac{1}{2}>\frac{1}{2}$, by  \cite[Lemma 5 and Corollary 2]{Friz2013}, for a constant $C>0$,  
$\Prob(N_{\alpha,[0,T]}(\mbB)\geq r)\leq \exp(
-C\alpha^\frac{2}{p}r^\frac{2}{q}
)$ 
for every $r>0$.  
Thus, for $R>0$,
\begin{align*}
\E\left[\exp\left(DN_{\alpha,[0,T]}(\mbB)\right)\right]
&=
\int_0^\infty\Prob(\exp(DN_{\alpha,[0,T]}(\mbB)) \geq s)\dd s 
\\
&\leq
R+\int_R^\infty\Prob(\exp(DN_{\alpha,[0,T]}(\mbB)) \geq s)\dd s
\\
&=
R+\int_{\frac{\log(R)}{D}}^\infty\Prob(N_{\alpha,[0,T]}(\mbB) \geq r)D\exp(Dr)\dd r
&&\text{($r=\log(s)/D$)}
\\
&\leq
R+D\int_{\frac{\log(R)}{D}}^\infty
\exp\left(
-C\alpha^\frac{2}{p}r^\frac{2}{q}
+
Dr
\right)
\dd r,
\end{align*}
which is bounded since  $\frac{1}{q}>\frac{1}{2}$.
\end{proof}

\subsection{Proofs of rough differential equation results %
(Section \ref{sec:rdes})}\label{apdx:proofs:rdes}
\subsubsection{Calculus with rough paths: controlled rough paths and rough integration (Section \ref{sec:rdes:calculus})}
\begin{proof}[Proof of Lemma \ref{lem:control_path:product}]
The first inequality \eqref{lem:control_path:product:YZ_p} follows from  
$\|(YZ)_{s,t}\|%
\leq\|Y\|_\infty\|Z_{s,t}\|+\|Z\|_\infty\|Y_{s,t}\|$ and \eqref{eq:path_finite_var:sum_pvars} in Lemma \ref{lem:pvariation:inequalities}. 
The second inequality \eqref{lem:control_path:product:(YZ)'_p} follows similary from 
\begin{align*}
\|(YZ)'_{s,t}\|
&\leq\|(ZY')_{s,t}\|+\|(YZ')_{s,t}\|
\leq
\|Z\|_\infty\|Y'_{s,t}\|+\|Y'\|_\infty\|Z_{s,t}\|
+
\|Y\|_\infty\|Z'_{s,t}\|+\|Z'\|_\infty\|Y_{s,t}\|
\end{align*}
and \eqref{eq:path_finite_var:sum_pvars} in Lemma \ref{lem:pvariation:inequalities}. 
 To show the third inequality \eqref{lem:control_path:product:R^YZ_p/2}, we write
\begin{align*}
R^{YZ}_{s,t} &= 
(YZ)_{s,t}-(YZ)'_s X_{s,t}
=
Y_tZ_t-Y_sZ_s-(Y_sZ_s'+Z_sY_s') X_{s,t}
\\
&=
Y_sZ_{s,t}+Y_{s,t}Z_s+Y_{s,t}Z_{s,t}-(Y_sZ_s'+Z_sY_s') X_{s,t}
\\
&=Y_sR^Z_{s,t}+R^Y_{s,t}Z_s+Y_{s,t}Z_{s,t}
\end{align*}
so that $\|R^{YZ}_{s,t}\| \leq
\|Y\|_\infty\|R^Z_{s,t}\|+\|R^Y_{s,t}\|\|Z\|_\infty+\|Y_{s,t}\|\|Z_{s,t}\|$, and we conclude with  \eqref{eq:path_finite_var:sum_p/2vars} in  Lemma \ref{lem:pvariation:inequalities}. 
Thus, $(YZ,(YZ)')\in\sD^p_X$ with $(YZ)'=Y'Z+YZ'$. 
\end{proof}

\begin{proof}[Proof of Lemma \ref{lem:rough_path:f(Y)'-f(tilde(Y))'}]
Throughout the proof, we write $\nabla\sigma:=\rev{\nablaof{x}\sigma}$ and $\nabla^2\sigma:=\rev{\nablaof{x}^2\sigma}$ for conciseness. 

To show \eqref{eq:controlled_path:Y-Ytilde_p}, we note that
\begin{align*}
\|\Delta Y_{s,t}\|
&=
\|Y'_sX_{s,t}+R_{s,t}^Y-\widetilde{Y}'_s\rev{\widetilde{X}}_{s,t}-R_{s,t}^{\widetilde{Y}}\| 
&&(Y_{s,t}=Y'_sX_{s,t}+R^Y_{s,t})
\\
&\leq
\|\Delta Y'\|_\infty\|X_{s,t}\|+\|\widetilde{Y}'\|_\infty\|\Delta X_{s,t}\|+\|\Delta R_{s,t}\| 
\\
&\leq
\Delta M_{Y'}\|X_{s,t}\|+M_{\widetilde{Y}'}\|\Delta X_{s,t}\|+\|\Delta R_{s,t}\|, 
&&(\|\widetilde{Y}'\|_\infty\leq\|\widetilde{Y}'_0\|+\|\widetilde{Y}'\|_p=M_{\widetilde{Y}'})
\end{align*}
so \eqref{eq:controlled_path:Y-Ytilde_p} follows from applying \eqref{eq:path_finite_var:sum_pvars} in Lemma \ref{lem:pvariation:inequalities}.

To show \eqref{eq:controlled_path:sigma(Y)-sigma(Ytilde)_p}, we note that
\begin{align*}
\|(\sigma(\cdot,Y)'-\sigma(\cdot,\widetilde{Y})')_{s,t}\|
&=
\|
(\nabla \sigma(\cdot,Y)Y'-\nabla \sigma(\cdot,\widetilde{Y})\widetilde{Y}')_{s,t}
\| 
\\
&\leq
\|
(\nabla \sigma(\cdot,Y)\Delta Y')_{s,t}
\|+
\|((\nabla \sigma(\cdot,Y)-\nabla \sigma(\cdot,\widetilde{Y}))\widetilde{Y}')_{s,t}
\| 
\\
&\hspace{-4cm}\leq
\|
\nabla\sigma(t,Y_t)\Delta Y'_{s,t}
\|
+
\|
\nabla\sigma(\cdot,Y)_{s,t}\Delta Y'_s
\|
+ 
\|(\nabla\sigma(t,Y_t)-\nabla\sigma(t,\widetilde{Y}_t))\widetilde{Y}'_{s,t}
\| 
+
\|(\nabla\sigma(\cdot,Y)-\nabla\sigma(\cdot,\widetilde{Y}))_{s,t}\widetilde{Y}'_s
\|
\\
&\hspace{-4cm}\leq
\|\sigma\|_{C_b^1}\|\Delta Y'_{s,t}\|
+
\big(
\|\sigma\|_{C_b^2}\|Y_{s,t}\|+|t-s|
\big)\|\Delta Y'\|_\infty
+
\|\sigma\|_{C_b^2}\|\Delta Y\|_\infty\|\widetilde{Y}'_{s,t}\|
+
\\
&+
\|\sigma\|_{C_b^3}
\big(
\|\Delta Y_{s,t}\|
+
(
\|Y_{s,t}\|+\|\widetilde{Y}_{s,t}\|
+
|t-s|
)
\|\Delta Y\|_\infty
\big)\|\widetilde{Y}'\|_\infty,
\end{align*}
where in the last inequality, we used  \eqref{eq:sigma(.,Y)_st_bound} to bound  the second term and \eqref{eq:deltaf_timevarying:st} to bound the fourth term. 
Thus,  by \eqref{eq:path_finite_var:sum_pvars} in Lemma \ref{lem:pvariation:inequalities},
\begin{align*}
\|\sigma(\cdot,Y)'-\sigma(\cdot,\widetilde{Y})'\|_p
&\leq
C_p
\|\sigma\|_{C_b^3}
\big(
\|\Delta Y'\|_p+\|Y\|_p\|\Delta Y'\|_\infty+
\|\Delta Y\|_\infty\|\widetilde{Y}'\|_p
\\
&\qquad
+
	(
	\|\Delta Y\|_p
	+
	(\|Y\|_p+\|\widetilde{Y}\|_p)\|\Delta Y\|_\infty
	)
	\|\widetilde{Y}'\|_\infty
+
T(
\|\Delta Y'\|_\infty+\|\Delta Y\|_\infty\|\widetilde{Y}'\|_\infty
)
\big)
\\
&\hspace{-2cm}\leq
C_p\|\sigma\|_{C_b^3}
(1+\|\widetilde{Y}'_0\|+\|Y\|_p+\|\widetilde{Y}\|_p+\|\widetilde{Y}'\|_p
)^2
(
1+T
)
\big(
\|\Delta Y\|_p
+
\|\Delta Y'\|_p
+
\|\Delta Y_0\|
+
\|\Delta Y'_0\|
\big),
\end{align*}
using $\|\Delta Y\|_\infty\leq\|\Delta Y_0\|+\|\Delta Y\|_p$  by \eqref{eq:path_finite_var:infty_ineq} in the last inequality.
Next, 
by \eqref{eq:controlled_path:Y-Ytilde_p}, 
$\|\Delta Y\|_p
\leq 
C_p(\Delta M_{Y'}\|X\|_p+M_{\widetilde{Y}'}\|\Delta X\|_p+\|\Delta R^Y\|_\frac{p}{2})\leq
C_p(1+M_{\widetilde{Y}'})(1+\|X\|_p)(\Delta M_{Y'}+\|\Delta X\|_p+\|\Delta R^Y\|_\frac{p}{2})$ with $\Delta M_{Y'}=\|\Delta Y'_0\|+\|\Delta Y'\|_p$, so 
\begin{align*}
\|\sigma(\cdot,Y)'-\sigma(\cdot,\widetilde{Y})'\|_p
&\leq
C_p\|\sigma\|_{C^3_b}
(1+\|\widetilde{Y}'_0\|+\|Y\|_p+\|\widetilde{Y}\|_p+\|\widetilde{Y}'\|_p)^2
(1+M_{\widetilde{Y}'})(1+\|X\|_p)
(1+T)
\rev{\times}
\\
&\qquad\qquad\qquad\qquad\qquad\qquad\qquad
\rev{\big(}
\|\Delta X\|_p+\|\Delta R^Y\|_\frac{p}{2}
+
\|\Delta Y_0\|
+
\|\Delta Y'_0\|
+
\|\Delta Y'\|_p
\big).
\end{align*}
By \eqref{eq:|Y|_p_ineq}, 
$\|Y\|_p\leq C_p(1+\|X\|_p)K_Y$, so 
\begin{align*}
\|\sigma(\cdot,Y)'-\sigma(\cdot,\widetilde{Y})'\|_p
&\leq
C_p\|\sigma\|_{C^3_b}
(1+K_Y+K_{\widetilde{Y}})^2
(1+M_{\widetilde{Y}'})
(1+\|X\|_p+\|\widetilde{X}\|_p)^3
(1+T)
\rev{\times}
\\
&\qquad\qquad\qquad\qquad\qquad\qquad\qquad
\rev{\big(}
\|\Delta X\|_p+\|\Delta R^Y\|_\frac{p}{2}
+
\|\Delta Y_0\|
+
\|\Delta Y'_0\|
+
\|\Delta Y'\|_p
\big),
\end{align*}
and we obtain   \eqref{eq:controlled_path:sigma(Y)-sigma(Ytilde)_p}.

To show \eqref{eq:controlled_path:R^sigma(Y)-R^sigma(Ytilde)_p/2}, we first decompose  $\|R^{\sigma(\cdot,Y)}_{s,t}-R^{\sigma(\cdot,\widetilde{Y})}_{s,t}\|$ as
\begin{align*}
\|R^{\sigma(\cdot,Y)}_{s,t}-R^{\sigma(\cdot,\widetilde{Y})}_{s,t}\|
&=
\|
\sigma(\cdot,Y)_{s,t}-\nabla \sigma(s,Y_s)Y'_sX_{s,t}-\sigma(\cdot,\widetilde{Y})_{s,t}
+
\nabla \sigma(s,\widetilde{Y}_s)\widetilde{Y}'_s\rev{\widetilde{X}}_{s,t})
\|
\\
&\hspace{-15mm}\leq
\|
\sigma(\cdot,Y)_{s,t}-\nabla \sigma(s,Y_s)Y_{s,t}
-
(\sigma(\cdot,\widetilde{Y})_{s,t}-\nabla \sigma(s,\widetilde{Y}_s)\widetilde{Y}_{s,t})\|
+
\|\nabla \sigma(s,Y_s)R^Y_{s,t}
-
\nabla \sigma(s,\widetilde{Y}_s)R^{\widetilde{Y}}_{s,t}
\|
\\
&\hspace{-15mm}\leq
\|
\sigma(s,Y)_{s,t}-\nabla \sigma(s,Y_s)Y_{s,t}
-
(\sigma(s,\widetilde{Y})_{s,t}-\nabla \sigma(s,\widetilde{Y}_s)\widetilde{Y}_{s,t})\|
+
\|\nabla \sigma(s,Y_s)R^Y_{s,t}
-
\nabla \sigma(s,\widetilde{Y}_s)R^{\widetilde{Y}}_{s,t}
\|
+
\\
&\qquad
\| 
(\sigma(\cdot,Y_t)-\sigma(\cdot,\widetilde{Y}_t))_{s,t}
\|
\\
&\hspace{-15mm}=:\|A_{s,t}\|+\|B_{s,t}\|+\|C_{s,t}\|.
\end{align*}
Next, we bound the three terms.  
First, denoting $\nabla^2\sigma(y)=\nabla^2\sigma(s,y)$ and $Y^\theta_{s,t}=\widetilde{Y}_s+\theta\widetilde{Y}_{s,t}$ for conciseness,
\begin{align*}
\|A_{s,t}\| &=
\left\|
\int_0^1 
\left(
	\nabla^2\sigma(Y_s+\theta Y_{s,t})Y_{s,t}^{\otimes 2}
	- 
	\nabla^2\sigma(\widetilde{Y}_s+\theta\widetilde{Y}_{s,t})\widetilde{Y}_{s,t}^{\otimes 2}
\right)(1-\theta)\dd\theta
\right\|
\\
&=
\left\|
\int_0^1
\left(
	\nabla^2\sigma(Y^\theta_{s,t})(Y_{s,t}^{\otimes 2}
	-
	\widetilde{Y}_{s,t}^{\otimes 2})
	+ 
	(\nabla^2\sigma(Y^\theta_{s,t})-\nabla^2\sigma(\widetilde{Y}^\theta_{s,t}))\widetilde{Y}_{s,t}^{\otimes 2}
\right)(1-\theta)\dd\theta
\right\|
\\
&=
\left\|
\int_0^1
\left(
	\nabla^2\sigma(Y^\theta_{s,t})(Y_{s,t}\otimes\Delta Y_{s,t}
	-
	\Delta Y_{s,t}\otimes\widetilde{Y}_{s,t})
	+ 
	(\nabla^2\sigma(Y^\theta_{s,t})-\nabla^2\sigma(\widetilde{Y}^\theta_{s,t}))\widetilde{Y}_{s,t}^{\otimes 2}
\right)(1-\theta)\dd\theta
\right\|
\\
&\leq
\|\sigma\|_{C^3_b}((\|Y_{s,t}\|+\|\widetilde{Y}_{s,t}\|)
\|\Delta Y_{s,t}\|
+
\|\Delta Y\|_\infty\|\widetilde{Y}_{s,t}\|^2).
\end{align*}

Second,
\begin{align*}
\|B_{s,t}\|
&=
\|(\nabla \sigma(s,Y_s)-\nabla \sigma(s,\widetilde{Y}_s))R^Y_{s,t}
+
\nabla \sigma(s,\widetilde{Y}_s)(R^Y_{s,t}-R^{\widetilde{Y}}_{s,t})
\| 
\leq 
\|\sigma\|_{C^2_b}(\|\Delta Y\|_\infty\|R_{s,t}^Y\|+\|\Delta R^Y_{s,t}\|).
\end{align*}

Third, by  \eqref{eq:deltaf_timevarying:st:2ndterm}, 
$
\| 
C_{s,t}
\|
\leq \|\sigma\|_{C_b^2}|t-s|\|\Delta Y\|_\infty.
$ 
Thus, by \eqref{eq:path_finite_var:sum_p/2vars} in Lemma \ref{lem:pvariation:inequalities},
\begin{align*}
\|R^{\sigma(\cdot,Y)}-R^{\sigma(\cdot,\widetilde{Y})}\|_\frac{p}{2}
&\leq
C_p(
\|A\|_\frac{p}{2}+\|B\|_\frac{p}{2}+\|C\|_\frac{p}{2}
)
\\
&\hspace{-3cm}\leq
C_p\|\sigma\|_{C^3_b}
\left(
(\|Y\|_p+\|\widetilde{Y}\|_p)\|\Delta Y\|_p+\|\Delta Y\|_\infty\|\widetilde{Y}\|^2_p
+
\|R^Y\|_\frac{p}{2}\|\Delta Y\|_\infty+\|\Delta R^Y\|_\frac{p}{2}
+
T\|\Delta Y\|_\infty
\right)
\\
&\hspace{-3cm}\leq
C_p\|\sigma\|_{C^3_b}
\left(
\|Y\|_p+\|\widetilde{Y}\|_p+\|\widetilde{Y}\|^2_p
+
\|R^Y\|_\frac{p}{2}
+
1
+
T
\right)
(\|\Delta Y_0\|+\|\Delta Y\|_p+\|\Delta R^Y\|_\frac{p}{2})
\\
&\hspace{-3cm}\leq
C_p\|\sigma\|_{C^3_b}
(
K_Y
+
K_{\widetilde{Y}}
+
1
)^2
(1+\|X\|_p+\|\rev{\widetilde{X}}\|_p)^2
(1+T)
(\|\Delta Y_0\|+\|\Delta Y\|_p+\|\Delta R^Y\|_\frac{p}{2}).
\end{align*}
where we used $\|Y\|_p\leq C_p(1+\|X\|_p)K_Y$ from \eqref{eq:|Y|_p_ineq} in the last inequality.  
Combining this inequality with  
$$\|\Delta Y\|_p
\mathop{\leq}^{\eqref{eq:controlled_path:Y-Ytilde_p}} 
C_p(\Delta M_{Y'}\|X\|_p+M_{\widetilde{Y}'}\|\Delta X\|_p+\|\Delta R^Y\|_\frac{p}{2})\leq 
C_p(1+K_{\widetilde{Y}})(\Delta M_{Y'}\|X\|_p+\|\Delta X\|_p+\|\Delta R^Y\|_\frac{p}{2}),$$
we obtain the desired inequality \eqref{eq:controlled_path:R^sigma(Y)-R^sigma(Ytilde)_p/2}.
\end{proof}

\begin{proof}[Proof of Lemma \ref{lem:rough_path:intfdX:error_bounds}]
First, by Lemma \ref{lem:rough_path:sigma(.,Y):controlled},  $(\sigma(\cdot,Y),\sigma(\cdot,Y)')=(\sigma(\cdot,Y),\rev{\nablaof{x}\sigma}(\cdot,Y)Y')$ is a controlled path, so by Proposition \ref{prop:rough_integral_welldefined:error_bound}, the rough integral $\int_s^t \sigma(r,Y_r)\dd\mbX_r$ is well-defined. %
Thus, for any $s,t\in[0,T]$, 
\begin{align*}
\left\|(R^{\int_0^\cdot \sigma(r,Y_r)\dd\mbX_r})_{s,t}\right\|
&=
\left\|
\int_s^t \sigma(r,Y_r)\dd\mbX_r-\sigma(s,Y_s)X_{s,t}
\right\|
\\
&\mathop{\leq}^{\eqref{eq:rough_int:error_bound}}
C_p(\|\sigma(\cdot,Y)'_s\bX_{s,t}\|+\|R^{\sigma(\cdot,Y)}\|_{\frac{p}{2},[s,t]}\|X\|_{p,[s,t]}+\|\sigma(\cdot,Y)'\|_{p,[s,t]}\|\bX\|_{\frac{p}{2},[s,t]})
\\
&\leq
C_p(\|\sigma\|_{C_b^1}\|Y'\|_\infty\|\bX\|_{\frac{p}{2},[s,t]}+\|R^{\sigma(\cdot,Y)}\|_{\frac{p}{2},[s,t]}\|X\|_{p,[0,T]}+\|\sigma(\cdot,Y)'\|_{p,[0,T]}\|\bX\|_{\frac{p}{2},[s,t]}).
\end{align*}
To continue, we need the following lemma, which is similar to \eqref{eq:path_finite_var:sum_p/2vars} in Lemma \ref{lem:pvariation:inequalities}.
\begin{lemma}\label{lem:sum_p/2vars:p/2vars_subintervals}
Let $p\geq 2$, $T>0$, $c\geq 0$, $X:[0,T]\to\R^d$, $Y^i,\widetilde{Y}^i\in\C^p$ for $i=1,\dots,n$, and $Z^j\in\C^\frac{p}{2}$ for $j=1,\dots,m$. %
Then, there exists a constant $C_p\geq 1$ such that
\begin{equation}\label{eq:path_finite_var:sum_p/2vars:p/2vars_subintervals}
\begin{split}
&\|X_{s,t}\|\leq 
\sum_{i=1}^n\|Y^i\|_{p,[s,t]}\|\widetilde{Y}^i\|_{p,[s,t]}
+
\sum_{j=1}^m\|Z^j\|_{\frac{p}{2},[s,t]}+c|t-s|
\ \forall s,t\in[0,T]
\\[-2mm]
&\hspace{1cm}\implies 
\|X\|_\frac{p}{2}\leq C_p\bigg(
\sum_{i=1}^n\|Y^i\|_p\|\widetilde{Y}^i\|_p
+
\sum_{j=1}^m\|Z^j\|_\frac{p}{2}
+
c\,T
\bigg).
\end{split}
\end{equation}
\end{lemma}
\begin{proof}[Proof of Lemma \ref{lem:sum_p/2vars:p/2vars_subintervals}]
We prove the particular case $n=0$ (without $(Y^i,\widetilde{Y}^i)$) and $c=0$. The general case follows with minor modifications, see the proof of \eqref{eq:path_finite_var:sum_p/2vars} in Lemma \ref{lem:pvariation:inequalities}. 
Given any partition  $\pi$ of $[0,T]$,
 \begin{align*}
 \sum_{[s,t]\in\pi}\|X_{s,t}\|^\frac{p}{2}
 &\leq
 \sum_{[s,t]\in\pi}
 \Big(
 \sum_{j=1}^m\|Z^j\|_{\frac{p}{2},[s,t]} 
 \Big)^\frac{p}{2}
 \leq
 m^\frac{p}{2}
 \Big( 
 \sum_{j=1}^m
 \sum_{[s,t]\in\pi}
 \|Z^j\|_{\frac{p}{2},[s,t]}^\frac{p}{2}  
 \Big)
 \leq
 C_p 
 \sum_{j=1}^m
 \|Z^j\|_{\frac{p}{2}}^\frac{p}{2}   ,
 \end{align*}
 where we used $|\sum_{i=1}^Na_i|^\frac{p}{2}\leq N^\frac{p}{2}\sum_{i=1}^N|a_i|^\frac{p}{2}$ for $p\geq 2$ in the second inequality, and $\sum_{[s,t]\in\pi}
 \|Z^j\|_{\frac{p}{2},[s,t]}^\frac{p}{2}\leq  
 \|Z^j\|_{\frac{p}{2},[0,T]}^\frac{p}{2}$ in the third inequality, since $w(s,t)=\|Z^j\|_{\frac{p}{2},[s,t]}^\frac{p}{2}$ is a control, so that $\sum_{[s,t]\in\pi}w(s,t)\leq w(0,T)$.
The bound above is independent of the choice of partition $\pi$. 
 Thus, \eqref{eq:path_finite_var:sum_p/2vars:p/2vars_subintervals} follows after taking the supremum over all partitions $\pi$ of $[0,T]$ and using the inequality $(\sum_i|a_i|)^\frac{2}{p}\leq\sum_i|a_i|^\frac{2}{p}$ for $p\geq 2$.  
\end{proof}
\textit{Proof of Lemma \ref{lem:rough_path:intfdX:error_bounds} (continued).}
Thus, by Lemma \ref{lem:sum_p/2vars:p/2vars_subintervals}, 
\begin{align*}
\|R^{\int_0^\cdot \sigma(s,Y_s)\dd\mbX_s}\|_\frac{p}{2}
&\leq
C_p(\|\sigma\|_{C_b^1}\|Y'\|_\infty\|\bX\|_\frac{p}{2}+\|R^{\sigma(\cdot,Y)}\|_\frac{p}{2}\|X\|_p+\|\sigma(\cdot,Y)'\|_p\|\bX\|_\frac{p}{2})
\\
&\hspace{-2cm}\leq
C_p\|\sigma\|_{C_b^2}
\big(
K_Y
\big(
\|\bX\|_\frac{p}{2}
+
(1+K_Y)(1+\|X\|_p)^2
(
\|X\|_p+
\|\bX\|_\frac{p}{2}
)
\big)
+
T(\|X\|_p+K_Y(1+\|X\|_p)\|\bX\|_\frac{p}{2})
\big)
\\
&\leq 
C_p\|\sigma\|_{C_b^2}(1+K_Y)^2(1+\|X\|_p)^2(1+T)(
	\|\bX\|_\frac{p}{2}+\|X\|_p),
\end{align*}
where we used $\|Y'\|_\infty\leq\|Y'_0\|+\|Y'\|_p=M_{Y'}\leq K_Y$, and the inequalities  from Lemma \ref{lem:rough_path:sigma(.,Y):controlled} 
$\|\sigma(\cdot,Y_\cdot)'\|_p
\leq 
C_p
\|\sigma\|_{C^2_b}K_Y(1+K_Y+T)
(1+\|X\|_p)$ 
and 
$\|R^{\sigma(\cdot,Y_\cdot)}\|_\frac{p}{2}
\leq
C_p
\|\sigma\|_{C^2_b}
(
K_Y(1+K_Y)(1+\|X\|_p)^2
+
T
)$. This concludes the proof of \eqref{eq:rough_path:R^int_sig_dX:p/2}. 
Note also that $\|
\int_0^\cdot \sigma(s,Y_r)\dd\mbX_s
\|_p<\infty$, since
{\small
\begin{align*}
\Big\|
\int_s^t \sigma(r,Y_r)\dd\mbX_r
\Big\|
&\leq
\Big\|
\int_s^t \sigma(r,Y_r)\dd\mbX_r-\sigma(s,Y_s)X_{s,t}
\Big\|
+
\|\sigma(s,Y_s)X_{s,t}\|
=
\big\|
R^{\int_0^\cdot \sigma(r,Y_r)\dd\mbX_r}_{s,t}
\big\|
+
\|\sigma(s,Y_s)X_{s,t}\|.
\end{align*}
}%
Also, $\|\sigma(\cdot,Y)\|_p<\infty$ by \eqref{eq:sigma(.,Y):pvar}. Together, $\|
\int_0^\cdot \sigma(s,Y_s)\dd\mbX_s
\|_p,\|\sigma(\cdot,Y)\|_p,\|R^{\int_0^\cdot \sigma(s,Y_s)\dd\mbX_s}\|_\frac{p}{2}<\infty$ imply that $(Z,Z')\in\sD^p_X$, which concludes the proof. 
\end{proof}

\begin{proof}[Proof of Lemma \ref{lem:rough_path:intfdX_fY:error_bounds}]
To show \eqref{eq:rough_path:intfdX_fY:delta_sigma},
we  combine $\|\sigma(\cdot,Y)-\sigma(\cdot,\widetilde{Y})\|_p\leq C_p\|\sigma\|_{C^2_b}(1 + \|Y\|_p+\|\widetilde{Y}\|_p+T)(\|\Delta Y_0\|+\|\Delta Y\|_p)$ in \eqref{eq:Delta_sigma(.,X):pvar},  
$\|Y\|_p\leq  C_p(1+\|X\|_p)K_Y$  in \eqref{eq:|Y|_p_ineq}, 
  $\|\Delta Y\|_p
 \leq 
C_p\big(\Delta M_{Y'}\|X\|_p+M_{\widetilde{Y}'}\|\Delta X\|_p+\|\Delta R^Y\|_\frac{p}{2}\big)$ in \eqref{eq:controlled_path:Y-Ytilde_p}, 
and  $M_{\widetilde{Y}'}\leq K_{\widetilde{Y}}$, so that
\begin{align*}
\|\sigma(\cdot,Y)-\sigma(\cdot,\widetilde{Y})\|_p
&\mathop{\leq}^{\eqref{eq:Delta_sigma(.,X):pvar}} 
C_p\|\sigma\|_{C^2_b}(1 + \|Y\|_p+\|\widetilde{Y}\|_p+T)
(\|\Delta Y_0\|+\|\Delta Y\|_p)
\\
&\hspace{-3cm}\mathop{\leq}^{\eqref{eq:|Y|_p_ineq},\eqref{eq:controlled_path:Y-Ytilde_p}} 
C_p\|\sigma\|_{C^2_b}
(1+\|X\|_p+\|\rev{\widetilde{X}}\|_p+T)
(1 + K_Y+K_{\widetilde{Y}})
(\|\Delta Y_0\|+\Delta M_{Y'}\|X\|_p+M_{\widetilde{Y}'}\|\Delta X\|_p+\|\Delta R^Y\|_\frac{p}{2})
\\
&\hspace{-30mm}\leq 
C_p\|\sigma\|_{C^2_b}(1\,{+}\,\|X\|_p\,{+}\,\|\rev{\widetilde{X}}\|_p\,{+}\,T)(1 \,{+}\, K_Y\,{+}\,K_{\widetilde{Y}})^2
(\|\Delta Y_0\|\,{+}\,
	(\|\Delta Y_0'\|\,{+}\,\|\Delta Y'\|_p)\|X\|_p\,{+}\,
	\|\Delta X\|_p\,{+}\,\|\Delta R^Y\|_\frac{p}{2})
\end{align*}
which is the desired inequality \eqref{eq:rough_path:intfdX_fY:delta_sigma}. 

To show \eqref{eq:rough_path:intfdX_fY:delta_R}, we observe from \eqref{lem:rough_int_stability:remainder} in Lemma \ref{lem:rough_int_stability} that it suffices to  bound  $K_{\sigma(\cdot,\widetilde{Y})}$ and $\Delta K_{\sigma(\cdot,Y)}$.  
First,  we bound $K_{\sigma(\cdot,\widetilde{Y})}$ as follows. By \eqref{eq:sigma(.,Y)':pvar}
and \eqref{eq:RY:p/2var:KYs} in Lemma \ref{lem:rough_path:sigma(.,Y):controlled}, $
\|\sigma(\cdot,Y)'\|_p\leq 
C_p\|\sigma\|_{C^2_b}K_Y(1+K_Y+T)(1+\|X\|_p)$ and $\|R^{\sigma(\cdot,Y)}\|_\frac{p}{2}\leq
C_p\|\sigma\|_{C^2_b}(K_Y(1+
K_Y)(1+\|X\|_p)^2+T)$. Also, $\|\sigma(\cdot,\widetilde{Y})'_0\|\leq\|\sigma\|_{C_b^2}\|\widetilde{Y}_0'\|\leq\|\sigma\|_{C_b^2}K_{\widetilde{Y}}$. By combining this inequalities, 
\begin{align}
\nonumber
K_{\sigma(\cdot,\widetilde{Y})}
&= \|\sigma(\cdot,\widetilde{Y})'_0\|+\|\sigma(\cdot,\widetilde{Y})'\|_p+\|R^{\sigma(\cdot,\widetilde{Y})}\|_\frac{p}{2}
\leq
C_p\|\sigma\|_{C^2_b}(1+
K_{\widetilde{Y}}+T)^2(1+\|\rev{\widetilde{X}}\|_p)^2
\\
\label{eq:bound:Ksigma(Y)}
&\leq
C_p\|\sigma\|_{C^2_b}
(1+K_Y+K_{\widetilde{Y}}+T)^3
(1+\|X\|_p+\|\rev{\widetilde{X}}\|_p)^3
(1+T).
\end{align}
Second,  using \eqref{eq:controlled_path:sigma(Y)-sigma(Ytilde)_p} and 
\eqref{eq:controlled_path:R^sigma(Y)-R^sigma(Ytilde)_p/2} in Lemma \ref{lem:rough_path:f(Y)'-f(tilde(Y))'}, we have
\begin{align*}
\|\sigma(\cdot,Y)'-\sigma(\cdot,\widetilde{Y})'\|_p
&\leq
C_p\|\sigma\|_{C^3_b}
(1+K_Y+K_{\widetilde{Y}'})^3(1+\|X\|_p+\|\rev{\widetilde{X}}\|_p)^3
(1+T)
\rev{\times}
\\
&\qquad\qquad\qquad
\rev{\big(}
\|\Delta X\|_p+
\|\Delta Y_0\|
+\|\Delta R^Y\|_\frac{p}{2}
+
\|\Delta Y'_0\|
+
\|\Delta Y'\|_p
\big),
\\
\|R^{\sigma(\cdot,Y)}-R^{\sigma(\cdot,\widetilde{Y})}\|_\frac{p}{2}
&\leq 
C_p\|\sigma\|_{C^3_b}(1+K_Y+K_{\widetilde{Y}})^3
	(1+\|X\|_p+\|\rev{\widetilde{X}}\|_p)^3
(1+T)
\rev{\times}
\\
&\qquad\qquad\qquad
\rev{\big(}
\|\Delta Y_0\|+\|\Delta R^Y\|_\frac{p}{2}+(\|\Delta Y_0'\|+\|\Delta Y'\|_p)\|X\|_p+\|\Delta X\|_p\big).
\end{align*}
By combining these two inequalities with $\|\Delta \sigma(\cdot,Y)'_0\|\leq\|\sigma\|_{C_b^2}(1+M_{Y'})(\|\Delta Y_0\|+\|\Delta Y_0'\|)$, we get
\begin{align}
\Delta K_{\sigma(\cdot,Y)} &=\|\Delta \sigma(\cdot,Y)'_0\|+\|\Delta \sigma(\cdot,Y)'\|_p+\|\Delta R^{\sigma(\cdot,Y)}\|_\frac{p}{2}
\nonumber
\\
&\leq
C_p\|\sigma\|_{C^3_b}(1+K_Y+K_{\widetilde{Y}})^3
	(1+\|X\|_p+\|\rev{\widetilde{X}}\|_p)^3
	(1+T)
\rev{\times}
\nonumber
\\
&\hspace{2cm}
\rev{\big(}
\|\Delta Y_0\|+\|\Delta Y_0'\|+\|\Delta R^Y\|_\frac{p}{2}+\|\Delta Y'\|_p+(\|\Delta Y_0'\|+\|\Delta Y'\|_p)\|X\|_p+\|\Delta X\|_p\big).
\label{eq:bound:DeltaKsigma(Y)}
\end{align}
Finally, combining the last inequalities \eqref{eq:bound:Ksigma(Y)} and \eqref{eq:bound:DeltaKsigma(Y)} using  \eqref{lem:rough_int_stability:remainder} in  Lemma \ref{lem:rough_int_stability}, we get
\begin{align*}
\|R^{\int_0^\cdot \sigma(s,Y_s)\dd\mbX_s}-R^{\int_0^\cdot \sigma(s,\widetilde{Y}_s)\dd\widetilde{\mbX}_s}\|_{\frac{p}{2}}
&\mathop{\leq}^{\eqref{lem:rough_int_stability:remainder} }
C_p
(1+\|\mbX\|_p+\|\widetilde{\mbX}\|_p)
\big(K_{\sigma(\cdot,\widetilde{Y})}\|\Delta\mbX\|_p+\|\mbX\|_p\Delta K_{\sigma(\cdot,Y)}\big)
\\
&\hspace{-5.1cm}\mathop{\leq}^{\eqref{eq:bound:Ksigma(Y)},\eqref{eq:bound:DeltaKsigma(Y)}}
C_p\|\sigma\|_{C_b^3}
(
1+\|\mbX\|_p+\|\widetilde{\mbX}\|_p)^4
(1+K_Y+K_{\widetilde{Y}}+T)^3
(1+T)
\rev{\times}
\\
& 
\hspace{-34mm}
\rev{\big(}
\|\Delta\mbX\|_p+
\|\mbX\|_p(\|\Delta Y_0\|+\|\Delta Y_0'\|+\|\Delta R^Y\|_\frac{p}{2}+\|\Delta Y'\|_p+(\|\Delta Y_0'\|+\|\Delta Y'\|_p)\|X\|_p+\|\Delta X\|_p)\big),
\end{align*}
which is the desired inequality \eqref{eq:rough_path:intfdX_fY:delta_R}.
\end{proof}

\subsubsection{Rough differential equations: existence and unicity  of solutions (Section \ref{sec:rdes:existence_unicity})}

\subsubsection*{Nonlinear rough differential equations (Section \ref{sec:rdes:existence_unicity:nonlinear})}

\begin{proof}[Proof of Theorem \ref{thm:rdes:existence_unicity}]
First, for any $t\in(0,T]$, we define the map
$$
\M_t:\ \sD^p_X([0,t],\R^n)\to\sD^p_X([0,t],\R^n),\ 
(Y,Y')\mapsto
\left(y+\int_0^\cdot b(s,Y_s,u_s)\dd s+\int_0^\cdot\sigma(s,Y_s)\dd\mbX_s,\sigma(\cdot,Y)\right).
$$
For $\delta\geq 1$, we define 
the ball
$$
\B_t^{(\delta)}=\left\{
(Y,Y')\in\sD^p_X([0,t],\R^n): \, Y_0=y, \, Y'_0=\sigma(0,y), 
\ 
\|Y,Y'\|_{X,p,[0,t]}^{(\delta)}\leq 1
\right\}
\subset\sD^p_X.
$$
with
\begin{equation}\label{eq:controlled_paths:norm}
\|Y,Y'\|_{X,p,[0,t]}^{(\delta)}=\|Y'\|_{p,[0,t]}+\delta\|R^Y\|_{\frac{p}{2},[0,t]}.
\end{equation}
The set $\B_t^{(\delta)}$ is a closed subset of the Banach space $\sD^p_X$ (when equiped with the metric induced by the norm $\|\cdot,\cdot\|_{X,p}^\delta$, and when we restrict  $\sD^p_X$ to controlled paths with fixed initial condition $(Y_0,Y_0')=(y,\sigma(0,y))$), so 
$\B_t^{(\delta)}$ is itself a complete metric space.  
Also, $\B_t^{(\delta)}$ is nonempty, since $s\mapsto(y+\sigma(0,y)X_{0,s},\sigma(0,y))\in\B_t^{(\delta)}$.

\textbf{Invariance.} We claim that $\M_{t_1}:\B_{t_1}^{(\delta)}\to\B_{t_1}^{(\delta)}$ for $t_1>0$ small enough, i.e., that $\B_{t_1}^{(\delta)}$ is invariant under $\M_{t_1}$. Let $(Y,Y')\in\B_{t_1}^{(\delta)}$. 
Since $\|Y,Y'\|_{X,p,[0,t]}^{(\delta)}\leq 1$,  $K_Y= \|Y_0'\|+\|Y'\|_{p,[0,t]}+\|R^Y\|_{\frac{p}{2},[0,t]}$ and $K_{\widetilde{Y}}$ satisfy $K_Y,K_{\widetilde{Y}}
\leq \|\sigma\|_{C_b^0}+1+\frac{1}{\delta}$. Thus, for $\delta \geq 1$, there exists $C_\sigma\geq 0$ such that
\begin{align}\label{eq:contraction:KY}
K_Y,K_{\widetilde{Y}}\leq \tilde{C}_\sigma(1+\delta^{-1})
\leq C_\sigma,
\   
(1+K_Y+K_{\widetilde{Y}})^2\leq C_\sigma,
\   
(1+K_Y+K_{\widetilde{Y}})^3\leq  C_\sigma.
\end{align}
Next, 
$
R^{y+\int_0^\cdot b(s,Y_s,u_s)\dd s+\int_0^\cdot\sigma(s,Y_s)\dd\mbX_s}_{s,t}
=
\int_s^t b(r,Y_r,u_r)\dd s+\int_s^t\sigma(r,Y_r)\dd\mbX_r
-
\sigma(s,Y_s)X_{s,t}
=
\int_s^t b(r,Y_r,u_r)\dd r
+
R^{\int_0^\cdot\sigma(s,Y_s)\dd\mbX_s}_{s,t}$, so
\begin{align*}
\|\M_{t_1}(Y,Y')\|_{X,p}^{(\delta)}
&=
\|\sigma(\cdot,Y)\|_p+\delta\|R^{y+\int_0^\cdot b(s,Y_s,u_s)\dd s+\int_0^\cdot\sigma(s,Y_s)\dd\mbX_s}\|_{\frac{p}{2}}
\\
&
\hspace{-26mm}\leq
\|\sigma(\cdot,Y)\|_p+C_p\delta\bigg(
\left\|\int_0^\cdot b(s,Y_s,u_s)\dd s\right\|_\frac{p}{2} +
\left\|R^{\int_0^\cdot\sigma(s,Y_s)\dd\mbX_s}\right\|_\frac{p}{2}
\bigg)
&&
\eqref{eq:path_finite_var:sum_p/2vars}
\\
&\hspace{-26mm}\leq\ 
C_p\big(
\|\sigma\|_{C^1_b}
(\|Y\|_p+t_1)
+ 
\delta
(
 C_{p,b}t_1
+
C_p\|\sigma\|_{C_b^2}(1+K_Y)^2(1+\|X\|_p)^2
(1+t_1)\|\mbX\|_p
)
\big)
&&
\eqref{eq:sigma(.,X):pvar},
\eqref{eq:bounds_int_b_ds},
\eqref{eq:rough_path:R^int_sig_dX:p/2}
\\
&
\hspace{-26mm}\leq
C_p
\big(\|\sigma\|_{C^1_b}(\|R^Y\|_\frac{p}{2}+K_Y\|X\|_p+t_1)
+
\delta C_{p,b,\sigma,\|X\|_p,T}
(1+K_Y)^2 
\left(
t_1
+
\|\mbX\|_p
\right)
\big)
&&
\eqref{eq:controlled_path:pvar_norm},t_1\leq T
\\
&\hspace{-26mm}\leq 
C_1
\Big(
\frac{1}{\delta}
+
\|X\|_p
 +t_1+\delta
  \left(
t_1
+
\|\mbX\|_p
\right)
\Big),
&&\hspace{-8mm}\eqref{eq:contraction:KY},\|Y,Y'\|^{(\delta)}_{X,p}\leq 1
\end{align*} 
for a constant $C_1:=C_{p,b,\sigma,\|\mbX\|_p,T}>\frac{1}{2}$. 
Let $\delta=\delta_1:=2C_1\geq 1$. Then, 
$$
\|\M_{t_1}(Y,Y')\|_{X,p}^{(\delta_1)}
\leq 
\frac{1}{2}
+
C_1\left(
\|X\|_{p,[0,t_1]}+t_1+2C_1(t_1+\|\mbX\|_{p,[0,t_1]})
\right).
$$
Then, by taking $t_1$ small-enough,  $\|\M_{t_1}(Y,Y')\|_{X,p}^{(\delta_1)}\leq 1$, so $\M_{t_1}(Y,Y')\in\B_{t_1}^{(\delta_1)}$. Invariance is proved. 

\textbf{Contraction.}  Let $(Y,Y'),(\widetilde{Y},\widetilde{Y}')\in\B_t^{(\delta_1)}$ for some $t\in(0,t_1]$. For any $\delta\geq 1$, by \eqref{eq:path_finite_var:sum_p/2vars},
\begin{align*}
\|\M_t(Y,Y')-\M_t(\widetilde{Y},\widetilde{Y}')\|_{X,p}^{(\delta)}
&\leq 
\\
&\hspace{-37mm}
\|\sigma(\cdot,Y)-\sigma(\cdot,\widetilde{Y})\|_p+
C_p\delta\bigg(
	\left\|\int_0^\cdot (b(s,Y_s,u_s) - b(s,\widetilde{Y}_s,u_s))\dd s\right\|_\frac{p}{2}
	+
\left\|R^{\int_0^\cdot\sigma(s,Y_s)\dd\mbX_s}-R^{\int_0^\cdot\sigma(s,\widetilde{Y}_s)\dd\mbX_s}\right\|_\frac{p}{2}
	\bigg).
\end{align*}
Then, %
since $(\Delta X,\Delta Y_0,\Delta Y'_0)=(0,0,0)$, 
\begin{align*}
\|\sigma(\cdot,Y)-\sigma(\cdot,\widetilde{Y})\|_p
&\mathop{\leq}^{\eqref{eq:rough_path:intfdX_fY:delta_sigma}}
C_p\|\sigma\|_{C^2_b}(1 + K_Y+K_{\widetilde{Y}})^2(1+\|X\|_p+t)(
\|\Delta Y'\|_p\|X\|_p+\|\Delta R^Y\|_\frac{p}{2})
\\
&\mathop{\leq}^{\eqref{eq:contraction:KY}}
C_{p,\sigma,\|X\|_p,T} 
(\|\Delta Y'\|_p\|X\|_p+\|\Delta R^Y\|_\frac{p}{2}),
\\
\|R^{\int_0^\cdot \sigma(s,Y_s)\dd\mbX_s}-R^{\int_0^\cdot \sigma(s,\widetilde{Y}_s)\dd\mbX_s}\|_{\frac{p}{2}}
&\mathop{\leq}^{\eqref{eq:rough_path:intfdX_fY:delta_R}} 
C_p\|\sigma\|_{C_b^3}
(1+K_Y+K_{\widetilde{Y}}+t)^3
\rev{\times}
\\
&\hspace{3cm}
\rev{\big(}
(1+\|\mbX\|_p)^5(1+t)
\|\mbX\|_p(\|\Delta Y'\|_p
+
\|\Delta R^Y\|_\frac{p}{2})
\big)
\\
&\mathop{\leq}^{\eqref{eq:contraction:KY}}
C_{p,\sigma,\|\mbX\|_p,T}
\|\mbX\|_p
(\|\Delta Y'\|_p+\|\Delta R^Y\|_\frac{p}{2}),
\\
\left\|\int_0^\cdot (b(s,Y_s,u_s)-b(s,\widetilde{Y}_s,u_s))\dd s\right\|_\frac{p}{2} 
&\mathop{\leq}^{\eqref{eq:bounds_int_b_ds}} 
C_{p,b}t\|\Delta Y\|_\infty
\mathop{\leq}^{\eqref{eq:path_finite_var:infty_ineq},\eqref{eq:controlled_path:Y-Ytilde_p}}
C_{p,b,\|X\|_p}t(\|\Delta Y'\|_p+\|\Delta R^Y\|_\frac{p}{2}).
\end{align*}
Thus, for $\delta\geq1$ and a constant $C_2:=C_{p,b,\sigma,\|\mbX\|_p,T}>\frac{1}{2}$, we obtain 
\begin{align*}
\|\M_t(Y,Y')-\M_t(\widetilde{Y},\widetilde{Y}')\|_{X,p}^{(\delta)}
&\leq 
C_2
\left(
\|\Delta Y'\|_p\|X\|_p+\|\Delta R^Y\|_\frac{p}{2}
+
\delta(\|\mbX\|_p+t)
\big(
\|\Delta Y'\|_p+\|\Delta R^Y\|_\frac{p}{2}
\big)
\right).
\end{align*}
Next, choose $\delta=\delta_2:=2C_2>1$, so that $C_2=\frac{\delta}{2}$, and choose $t=t_2\leq t_1$ small-enough so that $C_2\delta(\|\mbX\|_p+t)\leq \frac{1}{2}$, so that in particular $C_2\|X\|_p\leq\frac{1}{2\delta}$. Then, 
\begin{align*}
\|\M_{t_2}(Y,Y')-\M_{t_2}(\widetilde{Y},\widetilde{Y}')\|_{X,p}^{(\delta_2)}
&\leq 
\frac{1}{2\delta_2}
\|\Delta Y'\|_p+\frac{\delta_2}{2}\|\Delta R^Y\|_\frac{p}{2}
+
\frac{1}{2}\left(
\|\Delta Y'\|_p+\|\Delta R^Y\|_\frac{p}{2}
\right)
\\
&=
\frac{1+\delta_2}{2\delta_2}
\left(
\|\Delta Y'\|_p+\delta_2\|\Delta R^Y\|_\frac{p}{2}
\right).
\end{align*}
Since $\delta_2>1$, we obtain  
$\|\M_{t_2}(Y,Y')-\M_{t_2}(\widetilde{Y},\widetilde{Y}')\|_{X,p}^{(\delta_2)}<\|\Delta Y'\|_p+\delta_2\|\Delta R^Y\|_\frac{p}{2}=\|Y-\widetilde{Y},Y'-\widetilde{Y}'\|_{X,p}^{(\delta_2)}$, so $\M_{t_2}$ is a contraction.

 To conclude, $\M_{t_2}:\B_{t_2}^{(\delta_2)}\to\B_{t_2}^{(\delta_2)}$ is invariant and a contraction. Thus, there exists a unique fixed point $(Y,Y')\in\sD^p_X([0,t_2],\R^n)$ of the map $\M_{t_2}$, which is the solution to the RDE \eqref{eq:RDE} satisfying $Y'=\sigma(\cdot,Y)$ over the time interval $[0,t_2]$.  Since $t_2$ was chosen independently of the initial conditions, we can stitch together solutions over time intervals $[kt_2,(k+1)t_2]_{k\in\N}$, and deduce that the RDE \eqref{eq:RDE} admits a unique solution over the entire interval $[0,T]$.
\end{proof}

\subsubsection*{Linear rough differential equations (Section \ref{sec:rdes:existence_unicity:linear})}
In this section, we prove existence and unicity of solutions to linear RDEs  with drift \eqref{eq:rde_linear} in Theorem \ref{thm:rde:linear:existence_uniqueness}. 
The proof consists of rewriting the linear RDE with drift \eqref{eq:rde_linear} as a driftless linear RDE with constant coefficients driven by a \rev{new rough path} by also interpreting the Lebesgue integral $\int A_sV_s\dd s$  as a rough integral, and concluding with \rev{\cite[Theorem 2]{Lejay2009}}.  These results %
may be considered standard, although we could not find them     in the literature.

\subsubsection*{Preliminaries for the proof of Theorem \ref{thm:rde:linear:existence_uniqueness}}
\begin{lemma}[Connections between the Lebesgue 
and rough integrals]\label{lem:rough_integral_lebesgue}
Let $p\in[2,3)$, $T>0$, 
$b\in L^\infty([0,T],\R)$, 
$\mbX\in\sC^p([0,T],\R^d)$, %
and $(Y,Y')\in\sD^p_X([0,T],\R^{1\times d})$. %
\begin{itemize} 
\item[1)]  
Define $\cT:[0,T]\to\R$ and $\bT:\rev{[0,T]^2}\to\R$   by  the Lebesgue integrals
$$
\cT_t:=
\int_0^t
b_r\dd r,
\quad
\bT_{s,t}:=
\int_s^t\cT_{s,r}b_r\dd r=
\int_s^t\int_s^rb_v\dd v\,b_r\dd r.
$$
Then, $\mbT=(\cT,\bT)\in\sC^1_g([0,T],\R)$, i.e., $\mbT$ is a geometric $1$-rough path. 

Moreover, for any Gubinelli derivative $\hat{Y}'\in \C^p([0,T],\R^{n\times d\times d})$ 
 (in particular, for $\hat{Y}'=0$), the rough integral of $(Y,\hat{Y}')$ %
against $\mbT$, defined as the limit in \eqref{eq:rough_int}, 
is well-defined, and is equal to the Lebesgue integral $\int Yb\dd t$: for any $t\in[0,T]$,
\begin{equation}\label{eq:rough_integral_lebesgue:dT_bdt}
\int_0^tY_r\dd\mbT_r
= 
\int_0^tY_rb_r\dd r.
\end{equation}
\item[2)] Define $Z:[0,T]\to\R$ and $\bZ:\rev{[0,T]^2}\to\R$  by  the rough integrals
$$
Z_t:=
\int_0^t
Y_r\dd\mbX_r,
\quad
\bZ_{s,t}
:=
\int_s^tZ_{s,r}
Y_r
\dd\mbX_r.
$$
Then, $\mbZ=(Z,\bZ)\in\sC^p([0,T],\R)$, i.e., $\mbZ$ is a \rev{rough} path.

Moreover, $\mbZ$ is geometric if  $\mbX$ is geometric, i.e., $\mbX\in\sC^p_g\implies \mbZ\in\sC^p_g$. 

\end{itemize}
\end{lemma}
The pair $(Z,Z'):=(\int_0^\cdot
Y\dd\mbX, Y)$ is a controlled path by Lemma \ref{lem:rough_int_stability}. The rough path $\mbZ$ as defined in Lemma \ref{lem:rough_integral_lebesgue} is called  the canonical rough path lift of $(Z,Z')$, see \cite[Sections 6.1-6.2]{Allan2021}.  
Also, the equality $
\int Yb\dd t
=
\int Y\dd\mbT
$ in \eqref{eq:rough_integral_lebesgue:dT_bdt} (where, on the right hand side,  the Gubinelli derivative $\hat{Y}$  of $Y$ can be defined arbitrarily since $\mbT$ is smooth-enough for $\int Y\dd\mbT$ to be equivalent to a Young integral) will allow us to rewrite the linear RDE with drift \eqref{eq:rde_linear}  as a driftless linear RDE driven by a new rough path. 
\begin{proof}[Proof of Lemma \ref{lem:rough_integral_lebesgue}] 
1) First, 
$
\left|\cT_{s,t}\right|
=
\left|\int_s^tb_r\dd r\right|
\leq 
\left|\int_s^t\|b\|_\infty\dd r\right|
=\|b\|_\infty|t-s|,
$
so $\|\cT\|_1\leq\|b\|_\infty T<\infty$, and 
{\small
\begin{equation}\label{eq:bT_st<=1/2|Z|_infty(t-s)^2}
|\bT_{s,t}|=\left|\int_s^t\int_s^rb_v\dd vb_r\dd r\right|
\leq
\left|\int_s^t\left\|\int_s^\cdot b_v\dd v\right\|_{\infty,[s,r]} b_r\dd r\right|
\leq
\|b\|_\infty
\left|\int_s^t(r-s) b_r\dd r\right|
\leq
\frac{1}{2}\|b\|_\infty^2(t-s)^2,
\end{equation}
}%
so $\|\bT\|_\frac{1}{2}\leq\|b\|_\infty^2T^2/2<\infty$, and $\|\mbT\|_1<\infty$. Moreover, 
 for any $\rev{s,r,t\in[0,T]}$, 
 {\small
\begin{align*}
\bT_{s,r}
&=
\int_s^r\int_s^ub_v\dd vb_u\dd u
=
\int_s^t\int_s^ub_v\dd vb_u\dd u-
\int_r^t\int_s^ub_v\dd vb_u\dd u
=
\bT_{s,t}-
\left(
\int_r^t
\left(
\int_s^rb_v\dd vb_u\dd u
+
\int_r^ub_v\dd vb_u\dd u
\right)
\right)
\\
&=
\bT_{s,t}-
\left(
\int_s^rb_v\dd v
\int_r^tb_u\dd u
+
\int_r^t
\int_r^ub_v\dd vb_u\dd u
\right)
=
\bT_{s,t}-\bT_{r,t}-\cT_{s,r}\cT_{r,t},
\end{align*}
}%
so $\mbT$ satisfies Chen's relation \eqref{eq:chen's_relation}. The  condition  in  \eqref{eq:rough_path:integration_by_parts} follows the integration by parts formula
\begin{align*}
\bT_{s,t}
&=
\int_s^t\int_s^rb_v\dd vb_r\dd r
=
\left[
\int_s^\cdot b_v\dd v
\int_s^\cdot b_r\dd r
\right]_s^t
-
\int_s^tb_v\int_s^vb_r\dd r\dd v
=
\cT_{s,t}^2-\bT_{s,t}
\implies\bT_{s,t}=
\frac{1}{2}\cT_{s,t}^2,
\end{align*} 
so we  conclude that $\mbT$ is a geometric $1$-rough path. 

Moreover,  
given any $\hat{Y}'\in \C^p([0,T],\R^{n\times d\times d})$   and any partition $\pi$ of $[0,T]$,
{\small
\begin{equation}\label{eq:smooth_rough_intgral_no_effect}
\Big\|
\sum_{[s,t]\in\pi}\hat{Y}'_s\bT_{s,t}
\Big\|
\mathop{\leq}^{\eqref{eq:bT_st<=1/2|Z|_infty(t-s)^2}}
\frac{1}{2}\|b\|_\infty^2\|\hat{Y}'\|_\infty
\sum_{[s,t]\in\pi}(t-s)^2
\leq
\frac{1}{2}\|b\|_\infty^2\|\hat{Y}'\|_\infty\sup_{[s,t]\in\pi}|t-s|
\sum_{[s,t]\in\pi}|t-s|
\leq
\frac{1}{2}\|b\|_\infty^2\|\hat{Y}'\|_\infty T|\pi|\to 0
\end{equation}
}%
 as $|\pi|\to0$.  
Thus, 
$$
\int_0^TY_t\dd\mbT_t
\mathop{=}^{\eqref{eq:rough_int}}
\lim_{|\pi|\to0}\sum_{[s,t]\in\pi}Y_s\cT_{s,t}+
\hat{Y}'_s\bT_{s,t}
=
\lim_{|\pi|\to0}\sum_{[s,t]\in\pi}Y_s\cT_{s,t}
=
\int_0^TY_t\dd\cT_t,
$$
where the last integral is a well-defined Young integral since $Y\in\C^p$ and $\cT\in\C^1$ with $\frac{1}{p}+1>1$ \cite[Theorem 6.8]{Friz2010}. 
Thus, the rough integral of $(Y,\hat{Y}')$ against $\mbT$ is well-defined.   
The last integral is  also the Lebesgue integral of $Y$ with respect to  $\cT$. Finally, since
$\cT_t=\cT_s+\int_s^tb_r\dd r$,   $\cT$ is absolutely continuous with respect to the Lebesgue measure, so
$\int_0^TY_t\dd\cT_t
=
\int_0^TY_tb_t\dd t$, 
and we obtain the desired result \eqref{eq:rough_integral_lebesgue:dT_bdt}.

2) Second, $(Z,Z'):=(\int_0^\cdot
Y_s\dd\mbX_s, Y)$ is a controlled path by Lemma \ref{lem:rough_int_stability}. 
$\bZ$ is also well-defined, since $ZY$ is a  controlled path  by Lemma \ref{lem:control_path:product}, and $\|\bZ\|_{\frac{p}{2}}<\infty$ thanks to \eqref{eq:rough_int:error_bound}. 
The pair $\mbZ=(Z,\bZ)$ satisfies Chen's relation \eqref{eq:chen's_relation}, which can be shown via identical computations as for the proof that $\mbT$ satisfies Chen's relation. Thus, $\mbZ$ is a \rev{rough} path ($\mbZ\in\sC^p$). 

Moreover, the statement that  $\mbZ$ is geometric if and only if $\mbX$ is geometric 
is a consequence of $[\mbZ]=\int_0^\cdot(Y_r\otimes Y_r)\dd[\mbX]_r$ from \cite[Lemma 6.8]{Allan2021}, where $[\mbZ]_t:=Z_{0,t}^2-2\bZ_{0,t}$  and  $[\mbX]^{ij}_t:=X^i_{0,t} X^j_{0,t}-(\bX^{ij}_{0,t}+\bX^{ji}_{0,t})$ denote the brackets of $\mbZ$ and of $\mbX$ and satisfy $[\mbZ]_{s,t}=Z_{s,t}^2-2\bZ_{s,t}$ and $[\mbX]^{ij}_{s,t}=X^i_{s,t} X^j_{s,t}-(\bX^{ij}_{s,t}+\bX^{ji}_{s,t})$ \cite[Lemma 6.5]{Allan2021},
so that $\{[\mbZ]=0\iff\mbZ\in\sC^p_g\}$ if  $\{[\mbX]=0\iff \mbX\in\sC^p_g\}$. 
\end{proof}

The next result allows us to combine multiple rough paths as a joint rough path. %
\begin{lemma}[Joint geometric rough path]\label{lem:rough_path:joint_geometric}
Let $p\in[2,3)$, $T>0$, $\mbX\in\sC^p_g([0,T],\R^d)$ be a \rev{rough} path, $b^i\in L^\infty([0,T],\R)$ for $i=1,\dots,m$, and $(Y^j,(Y^j)')\in\sD^p_X([0,T],\R^d)$ %
for $j=1,\dots,n$.  
For $i=1,\dots,m$, define $\mbT^i=(\cT^i,\bT^i)\in\sC^1_g([0,T],\R)$  with $\cT_t^i=
\int_0^t
b_s^i\dd s$ and $\bT_{s,t}^i=
\int_s^t\cT_{s,r}^ib_r^i\dd r$, and %
for $j=1,\dots,n$, define $\mbZ^j=(Z^j,\bZ^j)\in\rev{\sC^p}([0,T],\R)$ with $Z_t^j=
\int_0^t
Y_s^j\dd\mbX_s$ and $\bZ_{s,t}^j=
\int_s^tZ_{s,r}^j Y_r^j\dd\mbX_r$,  as in Lemma \ref{lem:rough_integral_lebesgue}. 

Then, 
the pair $\mbJ := (J,\bJ)$ defined with 
$J_t:=(\cT_t^1,\dots,\cT_t^m,Z_t^1,\dots,Z_t^n)$ 
and 
\begin{equation}\label{eq:rough_path:joint_geometric}
\bJ_{s,t}:=
\begin{bmatrix}
\bT^1_{s,t}
& \dots &
\int_s^t \cT^1_{s,r}\dd \mbT^m_r &
\cT^1_{s,t}Z^1_{s,t}-\int_s^t Z^1_{s,r}\dd \mbT^1_r & \dots &
\cT^1_{s,t}Z^n_{s,t}-\int_s^t Z^n_{s,r}\dd \mbT^1_r
\\
\vdots & \ddots & \vdots & 
\vdots & \ddots & \vdots
\\
\int_s^t \cT^m_{s,r}\dd \mbT^1_r & \dots & \bT^m_{s,t}
&
\cT^m_{s,t}Z^1_{s,t}-\int_s^t Z^1_{s,r}\dd \mbT^m_r & \dots & \cT^m_{s,t}Z^n_{s,t}-\int_s^t Z^n_{s,r}\dd \mbT^m_r
\\
\int_s^t Z^1_{s,r}\dd \mbT^1_r & \dots & \int_s^t Z^1_{s,r}\dd \mbT^m_r &  \bZ^1_{s,t} & \dots &  \rev{\int_s^t Z^1_{s,r}Y^n_r\dd\mbX_r}
\\
\vdots & \ddots & \vdots & \vdots & \ddots & \vdots & 
\\
\int_s^t Z^n_{s,r}\dd \mbT^1_r & \dots & \int_s^t Z^n_{s,r}\dd \mbT^m_r &  
\int_s^t Z^n_{s,r}Y^1_r\dd\mbX_r
& \dots & \bZ^n_{s,t}
\end{bmatrix}
\end{equation}
is a (joint) \rev{rough} path, i.e., $\mbJ=(J,\bJ)\in\rev{\sC^p}([0,T],\R^{m+n})$.

\rev{Moreover, if $\mbX$ is geometric ($\mbX\in\sC^p_g$), then $\mbJ$ is also geometric ($\mbJ\in\sC^p_g$).}
\end{lemma}
\begin{proof}
Assume that $m=n=2$ without loss of generality, and denote $J=(\cT^1,\cT^2,Z^1,Z^2)=(J^1,J^2,J^3,J^4)$, and similarly for $\bJ$. 
First, $\|\mbJ\|_p<\infty$ follows after bounding  each term in $J_{s,t}$ and $\bJ_{s,t}$ and using Lemma \ref{lem:pvariation:inequalities}. 
Next, %
we show Chen's relation \eqref{eq:chen's_relation}:
{\small
\begin{align*}
\bJ_{s,r}^{12}
&=
\int_s^r \cT^1_{s,u}\dd \mbT^2_u
\mathop{=}^{\eqref{eq:rough_integral_lebesgue:dT_bdt}}
\int_s^r\int_s^ub_v^1\dd vb_u^2\dd u
=
\int_s^t\int_s^ub_v^1\dd vb_u^2\dd u-
\int_r^t\int_s^ub_v^1\dd vb_u^2\dd u
\\
&=
\bJ_{s,t}^{12}
-
\left(
\int_r^t
\left(
\int_s^rb_v^1\dd vb_u^2\dd u
+
\int_r^ub_v^1\dd vb_u^2\dd u
\right)
\right)
=
\bJ_{s,t}^{12}
-
\left(
\int_s^rb_v^1\dd v
\int_r^tb_u^2\dd u
+
\int_r^t
\int_r^ub_v^1\dd vb_u^2\dd u
\right)
\\
&=
\bJ_{s,t}^{12}-\bJ^{12}_{r,t}-J_{s,r}^1J_{r,t}^2,
\\
\bJ_{s,r}^{31}
&=
\int_s^r Z^1_{s,u}\dd \mbT^1_u
\mathop{=}^{\eqref{eq:rough_integral_lebesgue:dT_bdt}}
\int_s^r\int_s^uY_v^1\dd\mbX_vb_u^1\dd u
=
\int_s^t\int_s^uY_v^1\dd\mbX_vb_u^1\dd u-
\int_r^t\int_s^uY_v^1\dd\mbX_vb_u^1\dd u
\\
&=
\bJ_{s,t}^{31}
-
\left(
\int_r^t
\left(
\int_s^rY_v^1\dd\mbX_vb_u^1\dd u
+
\int_r^uY_v^1\dd\mbX_vb_u^1\dd u
\right)
\right)
=
\bJ_{s,t}^{31}
-
\left(
\int_s^rY_v^1\dd\mbX_v 
\int_r^tb_u^1\dd u
+
\int_r^t
\int_r^uY_v^1\dd\mbX_vb_u^1\dd u
\right)
\\
&=
\bJ_{s,t}^{31}-\bJ^{31}_{r,t}-Z_{s,r}^1\cT_{r,t}^1
\\
&=
\bJ_{s,t}^{31}-\bJ^{31}_{r,t}-J_{s,r}^3J_{r,t}^1,
\\
\bJ_{s,r}^{13}
&=
\cT^1_{s,r}Z^1_{s,r}
-
\int_s^r Z^1_{s,u}\dd \mbT^1_u
=
\cT^1_{s,r}Z^1_{s,r}
-
\bJ_{s,r}^{31}
=
\cT^1_{s,r}Z^1_{s,r}
-
(\bJ_{s,t}^{31}-\bJ_{r,t}^{31}-\cT_{r,t}^1Z_{s,r}^1)
\\
&=
(\cT^1_{s,t}-\cT^1_{r,t})(Z^1_{s,t}-Z^1_{r,t})
-
\bJ_{s,t}^{31}+\bJ_{r,t}^{31}
+
\cT_{r,t}^1Z_{s,r}^1
\\
&=
(\cT^1_{s,t}Z^1_{s,t}
-
\bJ_{s,t}^{31})
+
\cT^1_{r,t}Z^1_{r,t}
+
\bJ_{r,t}^{31}
-\cT^1_{s,t}Z^1_{r,t}
-\cT^1_{r,t}Z^1_{s,t}
+
\cT_{r,t}^1Z_{s,r}^1
\\
&=
(\cT^1_{s,t}Z^1_{s,t}
-
\bJ_{s,t}^{31})
+
\cT^1_{r,t}Z^1_{r,t}
+
\bJ_{r,t}^{31}
-\cT^1_{s,t}Z^1_{r,t}
-\cT^1_{r,t}(Z^1_{s,t}
-
Z_{s,r}^1)
\\
&=
(\cT^1_{s,t}Z^1_{s,t}
-
\bJ_{s,t}^{31})
+
\cT^1_{r,t}Z^1_{r,t}
+
\bJ_{r,t}^{31}
-\cT^1_{s,t}Z^1_{r,t}
-\cT^1_{r,t}Z^1_{r,t}
\\
&=
(\cT^1_{s,t}Z^1_{s,t}
-
\bJ_{s,t}^{31})
+
\bJ_{r,t}^{31}
-(\cT^1_{s,t}-\cT^1_{r,t})Z^1_{r,t}
-\cT^1_{r,t}Z^1_{r,t}
\\
&=
(\cT^1_{s,t}Z^1_{s,t}
-
\bJ_{s,t}^{31})
+
\bJ_{r,t}^{31}
-
\cT_{s,r}^1Z_{r,t}^1
-\cT^1_{r,t}Z^1_{r,t}
\\
&=
(\cT^1_{s,t}Z^1_{s,t}
-
\bJ_{s,t}^{31})
-
(\cT^1_{r,t}Z^1_{r,t}
-
\bJ_{r,t}^{31})
-
\cT_{s,r}^1Z_{r,t}^1
\\
&=
\bJ_{s,t}^{13}-\bJ^{13}_{r,t}-\cT_{s,r}^1Z_{r,t}^1
\\
&=
\bJ_{s,t}^{13}-\bJ^{13}_{r,t}-J_{s,r}^1J_{r,t}^3.
\end{align*}
Also,
\begin{align*}
\rev{
\bJ_{s,t}^{34}
}
&=
\rev{
\int_s^t Z^1_{s,u}Y^2_u\dd \mbX_u
=
\int_s^r Z^1_{s,u}Y^2_u\dd \mbX_u
+
\int_r^t Z^1_{s,u}Y^2_u\dd \mbX_u
=
\bJ_{s,r}^{34}
+
\int_r^t Z^1_{r,u}Y^2_u\dd \mbX_u
+
\int_r^t Z^1_{s,r}Y^2_u\dd \mbX_u,
}
\\
&=
\rev{
\bJ_{s,r}^{34}
+
\bJ_{r,t}^{34}
+
Z^1_{s,r}Z^2_{r,t}
=
\bJ_{s,r}^{34}
+
\bJ_{r,t}^{34}
+
J^3_{s,r}J^4_{r,t}
}
\end{align*}
}%
and similar derivations show that Chen's relation \eqref{eq:chen's_relation} hold for other pairs of indices $(i,j)$. Thus, $\mbJ\in\sC^p$, i.e., $\mbJ$ is a \rev{rough} path.

Finally, \rev{if $\mbX$ is geometric}, to show that $\mbJ$ is geometric, we need to show that $\bJ^{ij}_{s,t}+\bJ^{ji}_{s,t}=J_{s,t}^iJ_{s,t}^j$ in \eqref{eq:rough_path:integration_by_parts} holds for any $i,j\in\{1,2,3,4\}$. For $i=j$, \eqref{eq:rough_path:integration_by_parts} clearly holds since $\mbT$ and $\mbZ$ are geometric. For $1\leq i\neq j\leq 2$, \eqref{eq:rough_path:integration_by_parts} follows from the integration by parts formula that holds for the Lebesgue integral. 
\rev{
For $(i,j)=(3,4)$, %
{\small
\begin{align*}
\bJ_{s,t}^{34}+
\bJ_{s,t}^{43} 
\mathop{\approx}^{\eqref{eq:rough_path:joint_geometric:bJ34}}
(Y_s^1\otimes Y_s^2)\bX_{s,t}+(Y_s^2\otimes Y_s^1)\bX_{s,t}
=
(Y_s^1\otimes Y_s^2)(\bX_{s,t}+\bX_{s,t}^\top)
\mathop{=}^{\eqref{eq:rough_path:integration_by_parts}}
(Y_s^1\otimes Y_s^2)(X_{s,t}\otimes X_{s,t})
\mathop{\approx}^{\eqref{eq:rough_path:joint_geometric:Z1Z2}} J_{s,t}^{3}J_{s,t}^{4},
\end{align*}
}%
since $\mbX$ is geometric,  
where $a_{s,t}\approx b_{s,t}$ means $a_{s,t}=b_{s,t}+o(|t-s|^{\frac{3}{p}})$, and %
since
{\small
\begin{align}
\nonumber
\bJ_{s,t}^{34}
&=
\int_s^t Z^1_{s,u}Y^2_u\dd \mbX_u
=\int_s^t Z^1_uY^2_u\dd \mbX_u- Z^1_s\int_s^tY^2_u\dd \mbX_u
\\
\label{eq:rough_path:joint_geometric:bJ34}
&\approx
Z^1_sY^2_sX_{s,t}+((Z^1)'_sY^2_s+Z^1_s(Y^2)'_s)\bX_{s,t}-Z^1_s(Y^2_sX_{s,t}+(Y^2)'_s\bX_{s,t})
=
(Y_s^1\otimes Y_s^2)\bX_{s,t},
\\
\nonumber
J_{s,t}^{3}J_{s,t}^{4}
&=
\int_s^tY_u^1\dd\mbX_u\int_s^tY_u^2\dd\mbX_u
\approx
(Y^1_sX_{s,t}+(Y^1)'_s\bX_{s,t})
(Y^2_sX_{s,t}+(Y^2)'_s\bX_{s,t})
\\
\label{eq:rough_path:joint_geometric:Z1Z2}
&\approx
(Y_s^1\otimes Y_s^2)(X_{s,t}\otimes X_{s,t}),
\end{align}
}%
where the estimates use \eqref{eq:rough_int:error_bound}  in Proposition \ref{prop:rough_integral_welldefined:error_bound}. 
Taking the sum over any partition $\pi\in\mathcal{P}([s,t])$ with vanishing meshsize, we obtain that \eqref{eq:rough_path:integration_by_parts} holds for $(i,j)=(3,4)$.} 
For other pairs of indices $(i,j)$, \eqref{eq:rough_path:integration_by_parts} holds by definition, e.g., $\bJ^{13}_{s,t}+\bJ^{31}_{s,t}=\cT^1_{s,t}Z_{s,t}^1=J_{s,t}^1J_{s,t}^3$. To conclude, $\mbJ\in\sC^p_g$. %
\end{proof}

 \begin{proof}[Proof of Theorem \ref{thm:rde:linear:existence_uniqueness}]
We rewrite the linear RDE \eqref{eq:rde_linear} as a linear RDE with constant coefficients driven by a \rev{new rough path} and conclude \rev{with \cite[Theorem 2]{Lejay2009}}. 
First, for $1\leq i,j\leq n$, let
\begin{align*}
&\mbT^{ij}=(\cT^{ij},\bT^{ij}),
\ \text{ where}
&&\hspace{-15mm}
\cT^{ij}_t:=
\int_0^t
A_s^{ij}\dd s,
&&&&\hspace{-35mm}
\bT^{ij}_{s,t}:=
\int_s^t
\cT^{ij}_{s,r}
A_r^{ij}
\dd r,
\\
&\mbZ^{ij}=(Z^{ij},\bZ^{ij}), \ 
\ \text{ where}
&&\hspace{-15mm}
Z^{ij}_t:=
	\int_0^t
\Sigma_s^{i\cdot j}
\dd\mbX_s,
&&&&\hspace{-35mm}
\bZ^{ij}_{s,t}
:=
\int_s^tZ_{s,r}^{ij} 
\Sigma_r^{i\cdot j}
\dd\mbX_r.
\end{align*} 
By Lemma \ref{lem:rough_integral_lebesgue}, each $\mbT^{ij}$  and $\mbZ^{ij}$ is a \rev{rough} path ($\mbT^{ij},\mbZ^{ij}\in\rev{\sC^p}$). 
As in Lemma \ref{lem:rough_path:joint_geometric}, define the joint \rev{rough} path $\mbJ=(J,\bJ)\in\rev{\sC^p}([0,T],\R^{2n^2})$  with %
\begin{align*}
J
&=
\big(
J^1,
\dots,
J^{n^2},
J^{n^2+1},\dots,
J^{2n^2}
\big)
=
\left(
\cT^{11},
\dots,
\cT^{1n},
\dots,
\cT^{n1},
\dots,
\cT^{nn},
Z^{11},
\dots,
Z^{1n},
\dots,
Z^{n1},
\dots,
Z^{nn}
\right),
\end{align*}
so that  $J^{(i-1)n+j}=\cT^{ij}$ and $J^{n^2+(i-1)n+j}=Z^{ij}$ for any $i,j=1,\dots,n$, and  
define $\bJ:\rev{[0,T]^2}\to\R^{2n^2\times 2n^2}$ as in \eqref{eq:rough_path:joint_geometric} in Lemma \ref{lem:rough_path:joint_geometric}, with in particular $\bJ^{ii}=\bT^{ii}$ for $i\leq n^2$ and $\bJ^{ii}=\bZ^{ii}$ for $i> n^2$. 
 
Then, we define the tensor $F\in\R^{n\times 2n^2\times n}$ with
{\small
$$
F^{i\cdot\cdot}
=
\begin{bmatrix}
0_{(i-1)n\times n} \\
I_{n\times n} \\
0_{(n-1)n\times n} \\
I_{n\times n} \\
0_{(n-i)n \times n}
\end{bmatrix}\hspace{-2pt},
\text{ so that } 
FV_t = \underbrace{
\left[
\begin{array}{cccccccccccccc}
V^1_t&\dots&V^n_t & \dots & 0 &\dots & 0  &
V^1_t&\dots&V^n_t & \dots & 0 &\dots & 0
\\
\vdots &       & \vdots  &\ddots &  
\vdots &       & \vdots   &
\vdots &       & \vdots  &\ddots &  
\vdots &       & \vdots   
\\
0 &\dots & 0  &  \dots & V^1_t&\dots&V^n_t &
0 &\dots & 0  &  \dots & V^1_t&\dots&V^n_t
\end{array}
\right]
}_{
n\times 2n^2
}
$$
}%
for any $V_t\in\R^n$, and the linear RDE
$V_t =
v
+
\int_0^t
F V_r\dd\mbJ_r$, which is a linear RDE with constant-in-time linear vector fields $F^\ell(V)=F^\ell V$ driven by the \rev{rough} path $\mbJ$. 
Thus, \rev{by \cite[Theorem 2]{Lejay2009}}, there exists a unique solution $(V,FV)\in\sD_J([0,T],\R^n)$ to the linear RDE $V_t =
v
+
\int_0^t
F V_r\dd\mbJ_r$.

This concludes the proof, because $(V,\Sigma V)$ is also a solution to the original linear RDE \eqref{eq:rde_linear}, since
\begin{align}\label{eq:rde_linear:dJ_dX}
\int_s^t
F V_r\dd\mbJ_r
=
\int_s^t
A_rV_r\dd r
+  
\int_s^t
\Sigma_rV_r\dd\mbX_r
\end{align} 
for any $s,t\in[0,T]$. The last result \eqref{eq:rde_linear:dJ_dX} follows from long but straightforward computations with  integrals, noting that 
each component on the left hand side of \eqref{eq:rde_linear:dJ_dX} satisfies
{\small
\begin{align}\label{eq:rde_linear:intFVdJ}
 \bigg(
\int_s^tFV_r\dd\mbJ_r\bigg)^i
&\approx
(FV_s)^iJ_{s,t}+(FFV_s)^i\bJ_{s,t},
\end{align}
}%
for $i=1,\dots,n$, 
where  $a_{s,t}\approx b_{s,t}$ means $a_{s,t}=b_{s,t}+o(|t-s|^{\frac{3}{p}})$ and the estimate comes from \eqref{eq:rough_int:error_bound} in Proposition \ref{prop:rough_integral_welldefined:error_bound} (up to a time reparameterization \cite[Proposition 5.14]{Friz2010}, we may assume that $J,FV,FFV$ are $\frac{1}{p}$-H\"older continuous and that $\bJ,R^{FV}$ are $\frac{2}{p}$-H\"older continuous, and similarly for $X,\bX$). 

Also, for $i=1,\dots,n$,

{\small
\begin{align} 
(FV_s)^iJ_{s,t}
&=
\sum_{j=1}^j
V^j_s
\left(
J^{(i-1)n+j}_{s,t}
+
J^{n^2+(i-1)n+j}_{s,t}
\right)
=
\Big[
\underbrace{
\ 0 \ \dots \  0_{\phantom{s}} 
}_{(i-1)n} \ 
\underbrace{V^1_{\rev{s}} \ \dots \ V^n_{\rev{s}}}_{n}  \ 
\underbrace{
\ 0 \ \dots \  0_{\phantom{s}} 
}_{(n-1)n} \ 
\underbrace{V^1_{\rev{s}} \ \dots \ V^n_{\rev{s}}}_{n}  \ 
\underbrace{
\ 0 \ \dots \  0_{\phantom{s}} 
}_{(n-i)n}
\Big]
J_{s,t},
\nonumber
\\
(FFV_s)^i\bJ_{s,t}
&=
\sum_{k=1}^n
\sum_{j=1}^n
V^j_s
\Big(
\bJ^{    (i-1)n+k,     (k-1)n+j}_{s,t}+
\bJ^{    (i-1)n+k, n^2+(k-1)n+j}_{s,t}+
\bJ^{n^2+(i-1)n+k,     (k-1)n+j}_{s,t}+
\nonumber
\\[-2mm]
&\hspace{22mm}
+
\bJ^{n^2+(i-1)n+k, n^2+(k-1)n+j}_{s,t}
\Big)
\nonumber
\\
&\approx
\sum_{k=1}^n
\sum_{j=1}^n
V^j_s
\bJ^{n^2+(i-1)n+k, n^2+(k-1)n+j}_{s,t}
\label{eq:rde_linear:FVJ}
\end{align} 
}%
and, for $i,j,k=1,\dots,n$,
\begin{align}\label{eq:rde_linear:J}
J^{n^2+(i-1)n+j}_{s,t}
&=
\int_s^t
\Sigma^{i\cdot j}_r
\dd\mbX_r
\approx
\Sigma^{i\cdot j}_sX_{s,t}
+
(\Sigma^{i\cdot j})'_s\bX_{s,t},
\quad
\bJ^{n^2+(i-1)n+k,n^2+(k-1)n+j}_{s,t}
\rev{\mathop{\approx}^{\eqref{eq:rough_path:joint_geometric:bJ34}}}
\Sigma^{i\cdot k}_s\Sigma^{k\cdot j}_s\bX_{s,t}.
\end{align}
Additional details on the computation of \rev{$(FFV_s)^i\bJ_{s,t}$} %
are provided at the end of this section.
Then, writing $V'=\Sigma V$, the right  hand side of \eqref{eq:rde_linear:dJ_dX} satisfies
{\small
\begin{align}
\int_s^t
\left(A_rV_r\right)^i\dd r
&=
\sum_{j=1}^n
\int_s^t 
A_r^{ij}V_r^j\dd r
\mathop{=}^{\eqref{eq:rough_integral_lebesgue:dT_bdt}}
\sum_{j=1}^n
\int_s^t
V_t^j\dd\boldsymbol{\cT}_t^{ij}
\approx
\sum_{j=1}^n 
\big(
V^j_sJ^{(i-1)n+j}_{s,t}+(V^j)'_s\bJ^{(i-1)n+j,(i-1)n+j}_{s,t}
\big)
\nonumber
\\
\label{eq:rde_linear:int_AVdt}
&\approx
\sum_{j=1}^n 
V^j_sJ^{(i-1)n+j}_{s,t},
\\
\int_s^t
\left(\Sigma_rV_r\right)^i\dd\mbX_r
&\approx
\sum_{j=1}^n
\left(
V^j_s\Sigma^{i\cdot j}_sX_{s,t}+
(\Sigma^{i\cdot j}_s(V^j)'_s+V^j_s(\Sigma^{i\cdot j})'_s)\bX_{s,t}
\right)
\nonumber
\\
&=
\sum_{j=1}^n
\bigg(
V^j_s
\left(
\Sigma^{i\cdot j}_sX_{s,t}+(\Sigma^{i\cdot j})'_s\bX_{s,t}
\right)+
\bigg(
\sum_{k=1}^n
\Sigma^{i\cdot j}_s
\Sigma^{j\cdot k}_sV^k_s\bigg)\bX_{s,t}
\bigg)
\nonumber
\\
&\approx
\sum_{j=1}^n
\bigg(
V^j_s
Z^{ij}_{s,t}+
V^j_s
\bigg(
\sum_{k=1}^n
\Sigma^{i\cdot k}_s\Sigma^{k\cdot j}_s
\bigg)\bX_{s,t}
\bigg)
\nonumber
\\
&\mathop{\rev{\approx}}^{\eqref{eq:rde_linear:J}}
\sum_{j=1}^n
\bigg(
V^j_s
J^{n^2+(i-1)n+j}_{s,t}
+
V^j_s
\sum_{k=1}^n
\bJ^{n^2+(i-1)n+k,n^2+(k-1)n+j}_{s,t}
\bigg).
\label{eq:rde_linear:int_SigmaVdX}
\end{align}
}%
Thus, by comparing \eqref{eq:rde_linear:int_AVdt}+\eqref{eq:rde_linear:int_SigmaVdX}  with \eqref{eq:rde_linear:intFVdJ} using \eqref{eq:rde_linear:FVJ}, and after taking the sum over any partition $\pi\in\mathcal{P}([s,t])$ with vanishing meshsize, we conclude that  \eqref{eq:rde_linear:dJ_dX} holds. 
We conclude using \eqref{eq:rde_linear:dJ_dX} that the RDE  \eqref{eq:rde_linear} can be written as
$
V_t =v +  
\int_0^t
A_tV_t\dd t
+  
\int_0^t
\Sigma_tV_t\dd\mbX_t
=
v
+
\int_0^t
F V_t\dd\mbJ_t$  and the RDE on the right hand side  has a unique solution \rev{by \cite[Theorem 2]{Lejay2009}}, so the linear RDE \eqref{eq:rde_linear} has a unique solution.
\end{proof}

\subsubsection*{Additional details on computing \rev{$(FFV_s)^i  \bJ_{s,t}$}}%
For any $V_s\in\R^n$, %
$F\in\R^{n\times 2n^2\times n}$, $FV_s\in\R^{n\times 2n^2}$, and $FFV_s\in\R^{n\times 2n^2\times 2n^2}$. Then,  denoting by $\bJ^{a:b,c:d}_{s,t}\in\R^{(b-a+1)\times (d-c+1)}$ the block of the matrix  $\bJ_{s,t}\in\R^{2n^2\times 2n^2}$ containing the rows $a$ to $b$ and columns $c$ to $d$, and 
using the identification $(FFV_s)^i\in\R^{2n^2\times 2n^2}\cong\L(\R^{2n^2\times 2n^2},\R)$ with $(FFV_s)^i\bJ=\sum_{a=1}^{2n^2}\sum_{b=1}^{2n^2}((FFV_s)^i)^{a,b}\bJ^{a,b}_{s,t}$, we have 
{\small
\allowdisplaybreaks
\begin{align*}
(FFV_s)^i  \bJ_{s,t}
&= 
F^i FV_s
\bJ_{s,t}
\\
&\hspace{-15mm}=
\begin{bmatrix}
0_{(i-1)n\times n} \\
I_{n\times n} \\
0_{(n-1)n\times n} \\
I_{n\times n} \\
0_{(n-i)n \times n}
\end{bmatrix}
\underbrace{
\left[
\begin{array}{cccccccccccccc}
V^1_s&\dots&V^n_s & \dots & 0 &\dots & 0  &
V^1_s&\dots&V^n_s & \dots & 0 &\dots & 0
\\
\vdots &       & \vdots  &\ddots &  
\vdots &       & \vdots   &
\vdots &       & \vdots  &\ddots &  
\vdots &       & \vdots   
\\
0 &\dots & 0  &  \dots & V^1_s&\dots&V^n_s &
0 &\dots & 0  &  \dots & V^1_s&\dots&V^n_s
\end{array}
\right]
}_{
=FV_s \ (n\times 2n^2 \text{ matrix})
}
\bJ_{s,t}
\\
&\hspace{-15mm}=
\underbrace{
\begin{bmatrix}
0_{(i-1)n\times n^2}
\\
FV_s
\\
0_{(n-1)n\times 2n^2}
\\
	FV_s
\\
0_{(n-i)n \times 2n^2} 
\end{bmatrix}
}_{=F^iFV_s}
\underbrace{
\begin{bmatrix}
\star & \star
\\
\bJ^{(i-1)n+1:in, 1:n^2}_{s,t} 
& 
\bJ^{(i-1)n+1:in, n^2+1:2n^2}_{s,t}
\\
\star & \star
\\
\bJ^{n^2+(i-1)n+1:n^2+in, 1:n^2}_{s,t} 
& 
\bJ^{n^2+(i-1)n+1:n^2+in, n^2+1:2n^2}_{s,t}
\\
\star & \star
\end{bmatrix}
}_{=\bJ_{s,t}}
\\
&\hspace{-15mm}=
	\left[
	\begin{array}{cccccccc}
	V^1_s&\dots&V^n_s & \dots & 0 &\dots & 0
	\\
	\vdots &       & \vdots  &\ddots &  
	\vdots &       & \vdots   
	\\
	0 &\dots & 0  &  \dots & V^1_s&\dots&V^n_s
	\end{array}
	\right]
	\Bigg(
	\bJ^{(i-1)n+1:in, 1:n^2}_{s,t} 
+ 
\bJ^{(i-1)n+1:in, n^2+1:2n^2}_{s,t}
+
\bJ^{n^2+(i-1)n+1:n^2+in, 1:n^2}_{s,t} 
+
	\\
	&\hspace{-5mm}
	\begin{bmatrix}	\bJ^{n^2+(i-1)n+1,n^2+1}_{s,t}&\dots&\bJ^{n^2+(i-1)n+1,n^2+n}_{s,t} & \dots & \dots & \dots & \dots
	\\
	\vdots&\vdots&\vdots& \ddots & \vdots & \vdots & \vdots
	\\
	\dots & \dots & \dots & \dots 
	&
	\bJ^{n^2+(i-1)n+n,n^2+(n-1)n+1}_{s,t}&\dots&\bJ^{n^2+(i-1)n+n,2n^2}_{s,t} 
	\end{bmatrix}
	\Bigg)
\\
&\hspace{-15mm}=
\sum_{k=1}^n
\sum_{j=1}^n
V^j_s
\left(
\bJ^{    (i-1)n+k,     (k-1)n+j}_{s,t}+
\bJ^{    (i-1)n+k, n^2+(k-1)n+j}_{s,t}+
\bJ^{n^2+(i-1)n+k,     (k-1)n+j}_{s,t}+
\bJ^{n^2+(i-1)n+k, n^2+(k-1)n+j}_{s,t}
\right),
\end{align*}
}%
where  the first three terms involve the Lebesgue integral  (thus, they are $o(|t-s|^{\frac{3}{p}})$) and do not play a role after summing over all partitions $\pi\in\mathcal{P}[s,t]$ and taking the limit as $|\pi|\to 0$ (see for example the computation in  \eqref{eq:smooth_rough_intgral_no_effect}), so we obtain \eqref{eq:rde_linear:FVJ}.

\subsubsection{Bounds on solutions to nonlinear RDEs (Section \ref{sec:rdes:bounds})}
\begin{proof}[Proof of Proposition \ref{prop:bounds_pvars_solutions_RDEs}] 
To show \eqref{eq:bounds_pvars_solutions_RDEs:sigma(Y)'_p}, note that $\sigma(\cdot,Y)'=\rev{\nablaof{x}\sigma}(\cdot,Y)Y'
=\rev{\nablaof{x}\sigma}(\cdot,Y)\sigma(\cdot,Y)$, so we obtain
$\|\sigma(\cdot,Y)'\|_p\leq C_p\|\rev{\nablaof{x}\sigma}\sigma\|_{C_b^1}(\|Y\|_p+T)$ by \eqref{eq:sigma(.,X):pvar}. 
 Next, for $s,t\in I$ with $I\subseteq[0,T]$ an arbitrary interval, 
 {\small
\begin{align*}
\|R_{s,t}^Y\|
&=
\|Y_{s,t}-Y_s'X_{s,t}\| 
\leq
\left\|
\int_s^t\sigma(r,Y_r)\dd\mbX_r-\sigma(s,Y_s)X_{s,t}-\sigma(\cdot,Y)'_s\bX_{s,t}
\right\|
+
\left\|
\int_s^tb(r,Y_r,u_r)\dd r
\right\|
+
\left\|\sigma(\cdot,Y)'_s\bX_{s,t}\right\|
\\
&\mathop{\leq}^{\eqref{eq:rough_int:error_bound},\eqref{eq:bounds_int_b_ds}}
C_p\big(
\|R^{\sigma(\cdot,Y)}\|_{\frac{p}{2},[s,t]}\|X\|_{p,I}+\|\sigma(\cdot,Y)'\|_{p,I}\|\bX\|_{\frac{p}{2},[s,t]}
\big)
+
C_{p,b}|t-s|
+
C_\sigma\|\bX\|_{\frac{p}{2},[s,t]}.
\end{align*}
}%
where we used $\|\sigma(\cdot,Y)'\|_\infty
=\|\rev{\nablaof{x}\sigma}(\cdot,Y)\sigma(\cdot,Y)\|_\infty\leq C_\sigma$ in the last inequality. 
Thus, by \eqref{eq:path_finite_var:sum_p/2vars:p/2vars_subintervals} in Lemma \ref{lem:sum_p/2vars:p/2vars_subintervals} (which is similar to \eqref{eq:path_finite_var:sum_p/2vars} in Lemma \ref{lem:pvariation:inequalities}),
\begin{align*}
\|R^Y\|_{\frac{p}{2},I}
&\leq
C_{p,b,\sigma}\left(
\|R^{\sigma(\cdot,Y)}\|_{\frac{p}{2},I}
\|X\|_{p,I}+
\|\sigma(\cdot,Y)'\|_{p,I}\|\bX\|_{\frac{p}{2},I}
+
|I|
+
\|\bX\|_{\frac{p}{2},I}
\right)
\\
&\hspace{-12mm}\mathop{\leq}^{\eqref{eq:RY:p/2var:Y^2+RY+T},\eqref{eq:bounds_pvars_solutions_RDEs:sigma(Y)'_p}
}
C_{p,b,\sigma}
\left(
(\|Y\|_{p,I}^2+\|R^Y\|_{\frac{p}{2},I}
+|I|
)
\|X\|_{p,I}+
(\|Y\|_{p,I}+|I|)\|\bX\|_{\frac{p}{2},I}
+
|I|
+
\|\bX\|_{\frac{p}{2},I}
\right)
\\
&\hspace{-12mm}\leq
C_{p,b,\sigma}
\left(
(\|Y\|_{p,I}^2+\|R^Y\|_{\frac{p}{2},I})
\|X\|_{p,I}+
(1+\|Y\|_{p,I})\|\bX\|_{\frac{p}{2},I}
+
|I|
+
|I|(\|X\|_{p,I}+\|\bX\|_{\frac{p}{2},I})
\right)
\\
&\hspace{-12mm}\leq
C_{p,b,\sigma}
\left(
(\|Y\|_{p,I}^2+\|R^Y\|_{\frac{p}{2},I})
\|X\|_{p,I}+
(1+\|Y\|_{p,I}^2)\|\bX\|_{\frac{p}{2},I}
+
|I|
+
|I|(\|X\|_{p,I}+\|\bX\|_{\frac{p}{2},I})
\right)
\quad\text{($|x|\leq 1+x^2$)}
\\
&\hspace{-12mm}\leq
C_{p,b,\sigma}
\left(
\|R^Y\|_{\frac{p}{2},I}
\|X\|_{p,I}+
\|\bX\|_{\frac{p}{2},I}
+
|I|
+
|I|(\|X\|_{p,I}+\|\bX\|_{\frac{p}{2},I})
+
\|Y\|_{p,I}^2(\|X\|_{p,I}+\|\bX\|_{\frac{p}{2},I})
\right),
\end{align*}
where $C_1:=C_{p,b,\sigma}\geq 1$. 
Let $\alpha_1:=\left(\frac{1}{2C_1}\right)^p$. Then, for $I$ small-enough so that $\|\mbX\|_{p,I}+|I|\leq\alpha_1^\frac{1}{p}$, we obtain 
$$
\|\mbX\|_{p,I}+|I|\leq\alpha_1^\frac{1}{p}=\frac{1}{2C_{p,b,\sigma}}< 1,
$$ 
so that $C_{p,b,\sigma}\|X\|_{p,I}\leq \frac{1}{2}$, 
$\|X\|_{p,I}+\|\bX\|_{\frac{p}{2},I}< 1$, and $|I|+|I|(\|X\|_{p,I}+\|\bX\|_{\frac{p}{2},I})
<2|I|$. Thus, 
\begin{align}
\label{eq:bounds_pvars_solutions_RDEs:R^Y_in_proof}
\|R^Y\|_{\frac{p}{2},I}
&\leq
4C_{p,b,\sigma}
\left(
\|\bX\|_{\frac{p}{2},I}
+
|I|
+
\|Y\|_{p,I}^2
\right).
\end{align}
Thus, since $\|Y\|_p\leq C_{p,\sigma}(\|X\|_p+\|R^Y\|_\frac{p}{2})$ by \eqref{eq:controlled_path:pvar_norm}, we obtain 
\begin{align*}
\|Y\|_{p,I}
&\leq
C_2
\left(
\|X\|_{p,I}
+
\|\bX\|_{\frac{p}{2},I}
+
|I|
+
\|Y\|_{p,I}^2
\right)
\end{align*} 
for a new constant $C_2:=C_{p,b,\sigma}>4C_1>1$. 
Note that if $C_2\|Y\|_{p,I}\leq\frac{1}{2}$, then
\begin{align}\label{eq:||Y||p<=2C||X||}
\|Y\|_{p,I}
&\leq
2C_2
\left(
\|X\|_{p,I}
+
\|\bX\|_{\frac{p}{2},I}
+
|I|
\right)
\leq
 C_{p,b,\sigma},
\end{align}
and \eqref{eq:bounds_pvars_solutions_RDEs:small_intervals:Y}
 is proven. If not and $C_2\|Y\|_{p,I}>\frac{1}{2}$, let $|I|$ be smaller so that  $\|Y\|_{p,I}=\frac{1}{2C_2}$.  
Then,
$$
\frac{1}{2C_2}=\|Y\|_{p,I}
\mathop{\leq}^\eqref{eq:||Y||p<=2C||X||}
2C_2
\left(
\|X\|_{p,I}
+
\|\bX\|_{\frac{p}{2},I}
+
|I|
\right)
\leq
2C_2
\alpha_1^\frac{1}{p}.
$$
Then, we obtain $\alpha_1\geq \left(\frac{1}{4C_2^2}\right)^p=:\alpha_2$. Thus,  it suffices to make $I$ smaller so that 
$
\|\mbX\|_{p,I}+|I|\leq
 \alpha_2^\frac{1}{p},
$ 
and  $\|Y\|_{p,I}=\frac{1}{2C_2}=:C_{p,b,\sigma}$, which concludes the proof of  \eqref{eq:bounds_pvars_solutions_RDEs:small_intervals:Y}.
The inequality 
\eqref{eq:bounds_pvars_solutions_RDEs:small_intervals:R}
follows from  \eqref{eq:bounds_pvars_solutions_RDEs:R^Y_in_proof} and \eqref{eq:bounds_pvars_solutions_RDEs:small_intervals:Y}.
The inequality \eqref{eq:bound:KY}  follows from the previous inequalities. 
\end{proof}

\begin{proof}[Proof of Proposition \ref{prop:rdes:error_bound}]  
First, Theorem \ref{thm:rdes:existence_unicity} ensures that there exists two unique solutions $(Y,Y')\in\sD^p_X$ and $(\widetilde{Y},\widetilde{Y}')\in\sD^p_{\rev{\widetilde{X}}}$  to the RDEs. Second, %
let $I=[t_0,t_1]\subseteq[0,T]$ be an interval, and $C_{p,b,\sigma}\geq1$ and $\alpha_1:=\alpha_{p,b,\sigma}$ be two constants from Proposition \ref{prop:bounds_pvars_solutions_RDEs} 
such that the inequalities in  \eqref{eq:bounds_pvars_solutions_RDEs:small_intervals:combined}
hold and
\begin{equation}\label{eq:rdes:error_bound:small_interval_alpha:proof}
\|\mbX\|_{p,I}+|I|\leq
 \alpha_1^\frac{1}{p},
 \quad
 \|\widetilde{\mbX}\|_{p,I}+|I|\leq
 \alpha_1^\frac{1}{p}.
\end{equation}
By choosing $I$ small-enough so that $\|\mbX\|_{p,I}+\|\widetilde{\mbX}\|_{p,I}+|I|\leq
 \alpha_1^\frac{1}{p}$ %
holds, we have shown \eqref{eq:bounds_pvars_solutions_RDEs:small_intervals:combined}. 

Next, we show \eqref{eq:RDE:DY'+dRY:close} (for an interval $I$ that is perhaps shorter).  
First, 
{\small
\begin{align*}
\|\sigma(\cdot,Y)-\sigma(\cdot,\widetilde{Y})\|_{p,I}
&\mathop{\leq}^{\eqref{eq:rough_path:intfdX_fY:delta_sigma}} 
C_p\|\sigma\|_{C^2_b}(1 + K_{Y,I}+K_{\widetilde{Y},I})^2
(1+\|X\|_{p,I}+\|\rev{\widetilde{X}}\|_{p,I}+|I|)
\rev{\times}
\\
&\hspace{1cm}
\rev{\big(}
\|\Delta X\|_{p,I}+\|\Delta Y_{t_0}\|+(\|\Delta Y_{t_0}'\|+\|\Delta Y'\|_{p,I})\|X\|_{p,I}+\|\Delta R^Y\|_{\frac{p}{2},I}\big)
\\
&\leq
C_{p,b,\sigma}\big(
\|\Delta\mbX\|_{p,I}+\|\Delta Y_{t_0}\|+\|\Delta Y'\|_{p,I}\|X\|_{p,I}+\|\Delta R^Y\|_{\frac{p}{2},I}\big),
\\
\|R^{\int_0^\cdot \sigma(s,Y_s)\dd\mbX_s}-R^{\int_0^\cdot \sigma(s,\widetilde{Y}_s)\dd\widetilde{\mbX}_s}\|_{\frac{p}{2},I}
&\mathop{\leq}^{\eqref{eq:rough_path:intfdX_fY:delta_R}} 
C_p\|\sigma\|_{C_b^3}
(1+K_{Y,I}+K_{\widetilde{Y},I}+|I|)^3
(1+\|\mbX\|_{p,I}+\|\widetilde{\mbX}\|_{p,I})^4
(1+|I|)
\rev{\times}
\\
&\hspace{-4.2cm}
\rev{\Big(}
\|\Delta\mbX\|_{p,I}+
\|\mbX\|_{p,I}\big(\|\Delta Y_{t_0}\|+
\|\Delta Y_{t_0}'\|
+
\|\Delta R^Y\|_{\frac{p}{2},I}+\|\Delta Y'\|_{p,I}+
(\|\Delta Y_{t_0}'\|+\|\Delta Y'\|_{p,I})\|X\|_{p,I}+\|\Delta X\|_{p,I}\big)\Big)
\\
&\leq
C_{p,b,\sigma}
\Big(
\|\Delta\mbX\|_{p,I}+\|\Delta Y_{t_0}\|
+
\|\mbX\|_{p,I}\big(\|\Delta Y'\|_{p,I}+\|\Delta R^Y\|_{\frac{p}{2},I}\big)\Big),
\\
\left\|\int_0^\cdot (b(s,Y_s,u_s)-b(s,\widetilde{Y}_s,u_s))\dd s\right\|_{\frac{p}{2},I} 
&\mathop{\leq}^{\eqref{eq:bounds_int_b_ds:dy_du}} 
C_{p,b}
|I|(\|\Delta Y\|_{\infty,I}+\|\Delta u\|_{L^\infty,I}
)
\\[-4mm]
&\mathop{\leq}^{\eqref{eq:path_finite_var:infty_ineq}}
C_{p,b}
|I|(\|\Delta Y_{t_0}\|+\|\Delta Y\|_{p,I}
+
\|\Delta u\|_{L^\infty,I}
)
\\[-1mm]
&\hspace{-5mm}\mathop{\leq}^{\eqref{eq:controlled_path:Y-Ytilde_p}}
C_{p,b}
|I|\big(
\|\Delta Y_{t_0}\|
+
\Delta M_{Y'}\|X\|_{p,I}+M_{\widetilde{Y}'}\|\Delta X\|_{p,I}+\|\Delta R^Y\|_{\frac{p}{2},I}
+
\|\Delta u\|_{L^\infty,I}
\big)
\\
&\hspace{-5mm}\leq
C_{p,b,\sigma}
|I|\big(
\|\Delta Y_{t_0}\|
+
\|\Delta Y_{t_0}'\|+\|\Delta Y'\|_{p,I}+
\|\Delta X\|_{p,I}+\|\Delta R^Y\|_{\frac{p}{2},I}
+
\|\Delta u\|_{L^\infty,I}
\big)
\\
&\hspace{-5mm}\leq
C_{p,b,\sigma}
\big(
\|\Delta\mbX\|_{p,I}
+
\|\Delta Y_{t_0}\|
+
|I|\|\Delta u\|_{L^\infty,I}
+
|I|
(
\|\Delta Y'\|_{p,I}+\|\Delta R^Y\|_{\frac{p}{2},I}
)
\big),
\end{align*}
}%
where we also used $\|\Delta Y'_{t_0}\|\leq\|\sigma\|_{C_b^1}\|\Delta Y_{t_0}\|$ and $\|\Delta Y'\|_{p,I}=\|\sigma(\cdot,Y)-\sigma(\cdot,\widetilde{Y})\|_{p,I}$. %

Then, for a constant $C=C_{p,b,\sigma}>1$ and any $\delta>1$, 
\begin{align*}
&\|\Delta Y'\|_{p,I}+\delta\|\Delta R^Y\|_{\frac{p}{2},I} 
\\
&\quad\leq 
\|\sigma(\cdot,Y)-\sigma(\cdot,\widetilde{Y})\|_{p,I}+
\delta
\left\|\int_0^\cdot (b(s,Y_s,u_s) - b(s,\widetilde{Y}_s,\tilde{u}_s))\dd s\right\|_{\frac{p}{2},I}
+
\delta
\left\|R^{\int_0^\cdot\sigma(s,Y_s)\dd\mbX_s}-R^{\int_0^\cdot\sigma(s,\widetilde{Y}_s)\dd\widetilde{\mbX}_s}\right\|_{\frac{p}{2},I}
\\
&\quad\leq 
C
\Big[
\|\Delta\mbX\|_{p,I}+\|\Delta Y_{t_0}\|+
\|\Delta Y'\|_{p,I}\|X\|_{p,I}+\|\Delta R^Y\|_{\frac{p}{2},I}
\,
+
\\
&\quad\qquad\qquad
\delta\left(
|I|\|\Delta u\|_{L^\infty,I}+
\|\Delta\mbX\|_{p,I}+\|\Delta Y_{t_0}\|+
(\|\mbX\|_{p,I}+|I|)\big(
\|\Delta Y'\|_{p,I}+\|\Delta R^Y\|_{\frac{p}{2},I}
\big)
\right)
\Big]
\\
&\quad\leq 
C
\Big[
(1+\delta)(\|\Delta\mbX\|_{p,I}+\|\Delta Y_{t_0}\|+|I|\|\Delta u\|_{L^\infty,I})
+
\\
&\quad\qquad\qquad
\big(\|X\|_{p,I}+\delta(\|\mbX\|_{p,I}+|I|)\big)
\|\Delta Y'\|_{p,I}
+
\big(
1+\delta(\|\mbX\|_{p,I}+|I|)
\big)
\|\Delta R^Y\|_{\frac{p}{2},I}
\Big],
\end{align*}
Next, we  choose $\delta:=2C>2$, so that $C=\frac{\delta}{2}$, and $I$ small-enough so that 
\begin{equation}\label{eq:rdes:error_bound:small_interval_alpha:proof:2}
\|\mbX\|_{p,I}+|I|\leq \frac{1}{4C\delta},
\end{equation}
so that
$C\delta(\|\mbX\|_{p,I}+|I|)\leq\frac{1}{4}$  
and $C\|X\|_p\leq\frac{1}{4\delta}\leq\frac{1}{4}$. Then, after rearranging, the previous inequality becomes
\begin{align*}
\left(1-\frac{1}{2}\right)\|\Delta Y'\|_{p,I}+\left(\delta-\frac{\delta}{2}-\frac{1}{4}\right)\|\Delta R^Y\|_{\frac{p}{2},I}
&\leq 
C
(1+\delta)(\|\Delta\mbX\|_{p,I}+\|\Delta Y_{t_0}\|+|I|\|\Delta u\|_{L^\infty,I}),
\end{align*}
Thus, since $\delta>2$ and $(1+\delta)\leq 2\delta=4C$, we obtain
\begin{align*}
\|\Delta Y'\|_{p,I}+\|\Delta R^Y\|_{\frac{p}{2},I}
&\leq 
\|\Delta Y'\|_{p,I}+\left(\delta-\frac{1}{2}\right)\|\Delta R^Y\|_{\frac{p}{2},I}
\leq 
2C(1+\delta)
(\|\Delta\mbX\|_{p,I}+\|\Delta Y_{t_0}\|+
|I|\|\Delta u\|_{L^\infty,I}
).
\\
&\leq 
8 C^2
(\|\Delta\mbX\|_{p,I}+\|\Delta Y_{t_0}\|+
|I|\|\Delta u\|_{L^\infty,I}
).
\end{align*}
Finally, by choosing $\alpha_{p,b,\sigma}=\min(\alpha_1,1/(4C\delta)^p)$ and $I$ small-enough so that $\|\mbX\|_{p,I}+\|\widetilde{\mbX}\|_{p,I}+|I|\leq
 \alpha_{p,b,\sigma}^\frac{1}{p}$ %
holds (and in particular, \eqref{eq:rdes:error_bound:small_interval_alpha:proof}
and  \eqref{eq:rdes:error_bound:small_interval_alpha:proof:2} hold), we obtain \eqref{eq:RDE:DY'+dRY:close} and conclude the proof. 
\end{proof}

\begin{proof}[Proof of Proposition \ref{prop:rdes:error_bound:entire_interval}] 
By \eqref{eq:|X|+|Xtilde|+|dt|<=alpha^p} in  Corollary \ref{cor:Nalpha:NX_NXtilde_NT:small_intervals} and by Proposition \ref{prop:rdes:error_bound}, there exists constants $C_{p,b,\sigma}\geq 1$ and $\alpha=\alpha_{p,b,\sigma}>0$ such that   \eqref{eq:bounds_pvars_solutions_RDEs:small_intervals:combined} holds and
\begin{align*}
\|Y'-\widetilde{Y}'\|_{p,[s,t]}+\|R^Y-R^{\widetilde{Y}}\|_{\frac{p}{2},[s,t]}
\mathop{\leq}^{\eqref{eq:RDE:DY'+dRY:close}}
C_{p,b,\sigma}(\|\Delta Y_s\|
+
\|\Delta\mbX\|_{p,[s,t]}
+
|t-s|\|\Delta u\|_{L^\infty,[s,t]}
)
\end{align*}
for any $[s,t]\subseteq[0,T]$ such that $w(s,t)\leq
 \alpha$. Thus, as defined in Definition \ref{def:Nalpha}, the greedy partition 
 $\{\tau_i, i=0,1,\dots,N_{\alpha,I}(w)+1\}$ 
of the interval $I$, which satisfies $w(\tau_i,\tau_{i+1})\leq \alpha$ for all $i$, is such that
\begin{align}
\nonumber
\|\Delta Y'\|_{p,[\tau_i,\tau_{i+1}]}+\|\Delta R^Y\|_{\frac{p}{2},[\tau_i,\tau_{i+1}]}
&\leq 
C_{p,b,\sigma}
(
\|\Delta Y_{\tau_i}\|
+
\|\Delta\mbX\|_{p,[\tau_i,\tau_{i+1}]}
+
|\tau_{i+1}-\tau_i|\|\Delta u\|_{L^\infty,[\tau_i,\tau_{i+1}]}
)
\\
\label{eq:DY'_p_ti_ti+1}
&\leq
C_{p,b,\sigma}
(
\|\Delta Y_{\tau_i}\|
+
\|\Delta\mbX\|_{p,I}
+
|I|\|\Delta u\|_{L^\infty,I}
)
\end{align}
for all $i$, and $\|X\|_{p,[\tau_i,\tau_{i+1}]}\mathop{\leq}\limits^{\eqref{eq:|X|+|Xtilde|+|dt|<=alpha^p}} C_{p,b,\sigma}$, and $M_{Y',[\tau_i,\tau_{i+1}]}=\|Y_{\tau_i}'\|+\|Y'\|_{p,[\tau_i,\tau_{i+1}]}\mathop{\leq}\limits^{\eqref{eq:bounds_pvars_solutions_RDEs:small_intervals:combined}} K_{Y,[\tau_i,\tau_{i+1}]}\leq C_{p,b,\sigma}$.

Next, with $\Delta M_{Y',[\tau_i,\tau_{i+1}]}=\|\Delta Y_{\tau_i}'\|+\|\Delta Y'\|_{p,[\tau_i,\tau_{i+1}]}$, %
and since $\|\Delta Y'_{\tau_i}\|\leq\|\sigma\|_{C_b^1}\|\Delta Y_{\tau_i}\|$, 
\begin{align}
\nonumber
\|\Delta Y\|_{p,[\tau_i,\tau_{i+1}]}
&\mathop{\leq}^\eqref{eq:controlled_path:Y-Ytilde_p} 
C_p\big(\Delta M_{Y',[\tau_i,\tau_{i+1}]}\|X\|_{p,[\tau_i,\tau_{i+1}]}+M_{\widetilde{Y}',[\tau_i,\tau_{i+1}]}\|\Delta X\|_{p,[\tau_i,\tau_{i+1}]}+\|\Delta R^Y\|_{\frac{p}{2},[\tau_i,\tau_{i+1}]}\big)
\\
\nonumber
&\mathop{\leq}^{\eqref{eq:bounds_pvars_solutions_RDEs:small_intervals:combined}}C_{p,b,\sigma}\big(
\|\Delta Y_{\tau_i}'\|+\|\Delta Y'\|_{p,[\tau_i,\tau_{i+1}]}
+ \|\Delta X\|_{p,[\tau_i,\tau_{i+1}]}+\|\Delta R^Y\|_{\frac{p}{2},[\tau_i,\tau_{i+1}]}\big)
\\
\label{eq:DY_p_ti_ti+1}
&\leq
C_{p,b,\sigma}
\big(
\|\Delta Y_{\tau_i}\|
+ \|\Delta X\|_{p,I}+\|\Delta Y'\|_{p,[\tau_i,\tau_{i+1}]}+\|\Delta R^Y\|_{\frac{p}{2},[\tau_i,\tau_{i+1}]}\big),
\end{align} 
so that, since $\|\Delta Y_{\tau_{i+1}}\|
\leq
\|\Delta Y_{\tau_i}\|+\|\Delta Y\|_{p,[\tau_i,\tau_{i+1}]}$ by \eqref{eq:path_finite_var:infty_ineq},
\begin{align}\label{eq:DeltaYti+1}
\|\Delta Y_{\tau_{i+1}}\|
&\mathop{\leq}^{
\eqref{eq:path_finite_var:infty_ineq},
\eqref{eq:DY_p_ti_ti+1}
}
C_{p,b,\sigma}
\big(
\|\Delta Y_{\tau_i}\|
+ \|\Delta X\|_{p,I}+\|\Delta Y'\|_{p,[\tau_i,\tau_{i+1}]}+\|\Delta R^Y\|_{\frac{p}{2},[\tau_i,\tau_{i+1}]}\big),
\\
\nonumber
\|\Delta Y\|_{\infty,[\tau_i,\tau_{i+1}]}
&\mathop{\leq}^{
\eqref{eq:path_finite_var:infty_ineq},
\eqref{eq:DY_p_ti_ti+1}
}
C_{p,b,\sigma}
\big(
\|\Delta Y_{\tau_i}\|
+ \|\Delta X\|_{p,I}+\|\Delta Y'\|_{p,[\tau_i,\tau_{i+1}]}+\|\Delta R^Y\|_{\frac{p}{2},[\tau_i,\tau_{i+1}]}\big)
\\
\label{eq:DY_inf_ti_ti+1}
&\mathop{\leq}^{\eqref{eq:DY'_p_ti_ti+1}}
C_{p,b,\sigma}
\big(
\|\Delta Y_{\tau_i}\|
+
\|\Delta\mbX\|_{p,I}
+
|I|\|\Delta u\|_{L^\infty,I}
\big).
\end{align} 
By sequentially combining \eqref{eq:DY'_p_ti_ti+1},
\eqref{eq:DeltaYti+1} and \eqref{eq:DY_inf_ti_ti+1} over the intervals $[\tau_i,\tau_{i+1}]$, we obtain 
for all $i$,
\begin{align*}
\|\Delta Y'\|_{p,[\tau_i,\tau_{i+1}]}+\|\Delta R^Y\|_{\frac{p}{2},[\tau_i,\tau_{i+1}]}
&\mathop{\leq}^{
\eqref{eq:DY'_p_ti_ti+1},
\eqref{eq:DeltaYti+1}
} 
(C_{p,b,\sigma})^{i+1}(i+1)
\left(
\|\Delta Y_{t_0}\|
+
\|\Delta\mbX\|_{p,I}
+
|I|\|\Delta u\|_{L^\infty,I}
\right),
\\
\|\Delta Y\|_{\infty,[\tau_i,\tau_{i+1}]}
&\mathop{\leq}^{
\ \ \ \, 
\eqref{eq:DY_inf_ti_ti+1}
\ \ \ \, } 
(C_{p,b,\sigma})^{i+1}(i+1)
\left(
\|\Delta Y_{t_0}\|
+
\|\Delta\mbX\|_{p,I}
+
|I|\|\Delta u\|_{L^\infty,I}
\right).
\end{align*}\rev{Next, using $\Delta R^Y_{s,t}=\Delta R^Y_{s,u}+\Delta R^Y_{u,t}+\Delta Y'_{s,u}X_{u,t}+\widetilde{Y}'_{s,u}\Delta X_{u,t}$, the inequality \eqref{eq:pvar:RX_intervals} in Lemma \ref{lem:pvar:RX_intervals} can be slightly extended to obtain
{\small
$$\|\Delta R^Y\|_{\frac{p}{2},[0,T]}\leq 
3n^2
\bigg(
\hspace{-1pt}
\sum_{i=1}^n\|\Delta R^Y\|_{\frac{p}{2},[\Delta t_i]}
+
\bigg(
\hspace{-1pt}
\sum_{i=1}^n\|\Delta Y'\|_{p,[\Delta t_i]}
\hspace{-1pt}
\bigg)
\bigg(
\hspace{-1pt}
\sum_{i=1}^n\|X\|_{p,[\Delta t_i]}
\hspace{-1pt}
\bigg)
+
\bigg(
\hspace{-1pt}
\sum_{i=1}^n\|\widetilde{Y}'\|_{p,[\Delta t_i]}
\hspace{-1pt}
\bigg)
\bigg(
\hspace{-1pt}
\sum_{i=1}^n\|\Delta X\|_{p,[\Delta t_i]}
\hspace{-1pt}
\bigg)
\hspace{-1pt}
\bigg)
\hspace{-5pt}
$$
}%
for any partition $0=t_0<t_1<\dots<t_n=T$ of $[0,T]$ and $[\Delta t_i]:=[t_{i-1},t_i]$.} 
Thus, for $N:=N_{\alpha,I}(w)$, using Lemma \ref{lem:pvar:intervals}\rev{, the inequality above, $\|X\|_{p,[\tau_j,\tau_{j+1}]}\mathop{\leq}\limits^{\eqref{eq:|X|+|Xtilde|+|dt|<=alpha^p}} C_{p,b,\sigma}$ and $\|\widetilde{Y}'\|_{p,[\tau_j,\tau_{j+1}]}\mathop{\leq}\limits^{\eqref{eq:bounds_pvars_solutions_RDEs:small_intervals:combined}} C_{p,b,\sigma}$,}
\rev{
{\small
\begin{align*}
\|\Delta Y'\|_{p,I} 
&\leq 
(N+1)
\bigg(
\sum_{j=0}^N\|\Delta Y'\|_{p,[\tau_j,\tau_{j+1}]}^p
\bigg)^\frac{1}{p}
\leq
(N+1)
\bigg(
\sum_{j=0}^N\|\Delta Y'\|_{p,[\tau_j,\tau_{j+1}]}
+
\|\Delta R^Y\|_{\frac{p}{2},[\tau_j,\tau_{j+1}]}
\bigg)
\\
&\leq
(C_{p,b,\sigma})^{N+1}(N+1)^3
\left(
\|\Delta Y_{t_0}\|
+
\|\Delta\mbX\|_{p,I}
+
|I|\|\Delta u\|_{L^\infty,I}
\right),
\\
\|\Delta R^Y\|_{\frac{p}{2},I}
&\leq 
3(N+1)^2
\bigg(
\sum_{j=0}^N\|\Delta R^Y\|_{\frac{p}{2},[\tau_j,\tau_{j+1}]}
+
\sum_{j=0}^N\|\Delta Y'\|_{p,[\tau_j,\tau_{j+1}]}
\sum_{j=0}^N\|X\|_{p,[\tau_j,\tau_{j+1}]}
\\[-3mm]
&\qquad\qquad\qquad\qquad\qquad\qquad\qquad\qquad\qquad+
\sum_{j=0}^N\|\widetilde{Y}'\|_{p,[\tau_j,\tau_{j+1}]}
\sum_{j=0}^N\|\Delta X\|_{p,[\tau_j,\tau_{j+1}]}
\bigg)
\\[-2mm]
&\leq 
3(N+1)^2
\bigg(
\sum_{j=0}^N\|\Delta R^Y\|_{\frac{p}{2},[\tau_j,\tau_{j+1}]}
+
(N+1)C_{p,b,\sigma}
\bigg(
\sum_{j=0}^N\|\Delta Y'\|_{p,[\tau_j,\tau_{j+1}]}
+
\sum_{j=0}^N\|\Delta X\|_{p,[\tau_j,\tau_{j+1}]}
\bigg)
\bigg)
\\[-2mm]
&\leq
3C_{p,b,\sigma}
(N+1)^3
\sum_{j=0}^N
\bigg(
\|\Delta R^Y\|_{\frac{p}{2},[\tau_j,\tau_{j+1}]}
+
\sum_{j=0}^N\|\Delta Y'\|_{p,[\tau_j,\tau_{j+1}]}
+
\sum_{j=0}^N\|\Delta X\|_{p,[\tau_j,\tau_{j+1}]}
\bigg)
\\
&\leq
3(C_{p,b,\sigma})^{N+2}(N+1)^5
\left(
\|\Delta Y_{t_0}\|
+
\|\Delta\mbX\|_{p,I}
+
|I|\|\Delta u\|_{L^\infty,I}
\right),
\\
\|\Delta Y\|_{\infty,I}
&\leq
(C_{p,b,\sigma})^{N+1}(N+1)^3
\left(
\|\Delta Y_{t_0}\|
+
\|\Delta\mbX\|_{p,I}
+
|I|\|\Delta u\|_{L^\infty,I}
\right),
\end{align*}
}%
} 
where we used $\|\Delta Y\|_{\infty,I}\leq \max_i\|\Delta Y\|_{\infty,[\tau_i,\tau_{i+1}]}\leq (N+1)\max_i\|\Delta Y\|_{\infty,[\tau_i,\tau_{i+1}]}$ in the last inequality. 
The desired inequalities  \eqref{eq:RDE:DY'+dRY:close:full_interval} and \eqref{eq:RDE:DY:close:full_interval}  follow using $(C_{p,b,\sigma})^{N+\rev{2}}(N+\rev{2})^{\rev{5}}\leq \exp((N+\rev{2})\log(C_{p,b,\sigma}))\rev{5}!\exp(N+\rev{2})=\rev{5!}e^{\rev{4}}C_{p,b,\sigma}\exp((\log(C_{p,b,\sigma})+1)N)$.

Finally, if $u$ and $\tilde{u}$ only differ on an interval $J=[s_0,s_1]\subseteq I$, then    \eqref{eq:RDE:DY'+dRY:close:full_interval:u_subinterval} and \eqref{eq:RDE:DY:close:full_interval:u_subinterval} follow  from sequentially applying \eqref{eq:RDE:DY'+dRY:close:full_interval}
 and 
\eqref{eq:RDE:DY:close:full_interval}
on  
$[t_0,s_0]
\cup
J
\cup
[s_1,t_1]=I$, noting that $\Delta u_t=0$ for almost every 
$t\in 
[t_0,s_0]
\cup
[s_1,t_1]$.
\end{proof}

\subsubsection{Bounds on  solutions to linear RDEs and on the Jacobian flow (Section \ref{sec:rdes:bounds:linear})}

The proof of Lemma \ref{lem:rde:linear:bounded_solutions} follows similar steps as  the proof of Proposition \ref{prop:rdes:error_bound:entire_interval} and relies on the following   Gr\"onwall Lemma for rough paths.
\begin{lemma}[Rough Gr\"onwall Lemma {\cite[Lemma 2.12]{Deya2019}}]
\label{lem:rough_gronwall}
Let $p\geq 1$, $T>0$, $C_1,\alpha>0$, $Y\in C([0,T],\R^n)$, and $w_1,w_2$  be two controls on $[0,T]$ (see Definition \ref{def:Nalpha}) such that
$$
\|Y_{s,t}\|\leq C_1\|Y\|_{\infty,[0,t]}w_1(s,t)^\frac{1}{p}+w_2(s,t)
\ \ \text{for any } [s,t]\subseteq[0,T] \text{ such that }w_1(s,t)\leq \alpha,
$$
and define $C_2=\min(1, 1/(\alpha(2C_1\exp(2))^p))$. 
Then,
$$
\|Y\|_{\infty,[0,T]}
\leq
2\exp\left(\frac{w_1(0,T)}{C_2\alpha}\right)
\left(
\|Y_0\|
+
\left\|w_2(0,\cdot)\exp\left(-w_1(0,\cdot)/(C_2 \alpha)\right)
\right\|_{\infty,[0,T]}
\right).
$$
\end{lemma}

\begin{proof}[Proof of Lemma \ref{lem:rde:linear:bounded_solutions}]
 By Theorem \ref{thm:rde:linear:existence_uniqueness}, there exists a unique solution $(V,V')\in\sD^p_X$ with $V'=\Sigma V$ to the linear RDE. Let $I=[s,t]\subseteq[0,T]$ be an interval.  By Lemma \ref{lem:control_path:product}, with $(\Sigma V)'=\Sigma V'+V\Sigma'$, %
\begin{align}
\label{eq:V'=SigmaV_p<=sum|Sigma|p(|V|inf+|V|p}
\|V'\|_{p,I}%
&\mathop{\leq}^{\eqref{lem:control_path:product:YZ_p}}
C_p(\|\Sigma\|_{\infty,I}\|V\|_{p,I}+\|V\|_{\infty,I}\|\Sigma \|_{p,I})
\leq
C_p(\|\Sigma \|_{\infty,I}+\|\Sigma \|_{p,I})(\|V\|_{\infty,I}+\|V\|_{p,I}),
\\
\nonumber
\|(\Sigma V)'\|_{p,I}
&\mathop{\leq}^{\eqref{lem:control_path:product:(YZ)'_p}}
C_p(\|V\|_{\infty,I}\|\Sigma '\|_{p,I}+\|\Sigma '\|_{\infty,I}\|V\|_{p,I}+\|\Sigma \|_{\infty,I}\|V'\|_{p,I}+\|V'\|_{\infty,I}\|\Sigma \|_{p,I})
\\
\nonumber
&\mathop{\leq}^{\eqref{eq:V'=SigmaV_p<=sum|Sigma|p(|V|inf+|V|p}}
C_p(\|\Sigma \|_{\infty,I}+\|\Sigma \|_{p,I}+\|\Sigma '\|_{\infty,I}+\|\Sigma '\|_{p,I})^2(\|V\|_{\infty,I}+\|V\|_{p,I}),
\\ 
\nonumber
\|R^{\Sigma V}\|_{\frac{p}{2},I}
&\mathop{\leq}^{\eqref{lem:control_path:product:R^YZ_p/2}}
C_p(\|\Sigma \|_{\infty,I}\|R^V\|_{\frac{p}{2},I}+\|R^\Sigma \|_{\frac{p}{2},I}\|V\|_{\infty,I}+\|\Sigma \|_{p,I}\|V\|_{p,I})
\\
\nonumber
&\leq
C_p(\|\Sigma \|_{\infty,I}+\|\Sigma \|_{p,I}+\|R^\Sigma \|_{\frac{p}{2},I})
(\|V\|_{\infty,I}+\|V\|_{p,I}+\|R^V\|_{\frac{p}{2},I}).
\end{align}
Let $I$ be small-enough to satisfy $\|\mbX\|_{p,I}+|I|\leq
 \alpha_\Sigma ^\frac{1}{p}$. Then, by Proposition \ref{prop:rough_integral_welldefined:error_bound},
\begin{align*}
\|R_{s,t}^V\| 
&=
\|V_{s,t}-\Sigma _sV_sX_{s,t}\|
\leq
\left\|\int_s^tA _rV_r\dd r \right\|
+
\left\|\int_s^t\Sigma _rV_r\dd\mbX_r-\Sigma _sV_sX_{s,t}-(\Sigma V)'_s\bX_{s,t}\right\|+\|(\Sigma V)'_s\bX_{s,t}\|
\\
&\mathop{\leq}^{\eqref{eq:rough_int:error_bound}}
\|A\|_{\infty,I}\|V\|_{\infty,I}|I|
+
C_p(\|R^{\Sigma V}\|_{\frac{p}{2},I}\|X\|_{p,I}+\|(\Sigma V)'\|_{p,I}\|\bX\|_{\frac{p}{2},I})+\|(\Sigma V)'\|_{\infty,I}\|\bX\|_{\frac{p}{2},I}.
\end{align*}
Thus, using  \eqref{eq:path_finite_var:sum_p/2vars:p/2vars_subintervals} in Lemma \ref{lem:sum_p/2vars:p/2vars_subintervals} (which is similar to \eqref{eq:path_finite_var:sum_p/2vars} in Lemma \ref{lem:pvariation:inequalities}) and the previous inequalities,
\begin{align*}
\|R^V\|_{\frac{p}{2},I}  
&\leq 
C_p
(
\|A\|_{\infty,I}+
\|\Sigma \|_{\infty,I}+\|\Sigma \|_{p,I}+\|\Sigma '\|_{\infty,I}+\|\Sigma '\|_{p,I}+\|R^\Sigma \|_{\frac{p}{2},I}
)^2
\rev{\times}
\\
&\hspace{2cm}
\rev{\big(}
\|V\|_{\infty,I}|I|
+
(\|V\|_{\infty,I}+\|V\|_{p,I}+\|R^V\|_{\frac{p}{2},I})
\|X\|_{p,I}
+
(\|V\|_{\infty,I}+\|V\|_{p,I})\|\bX\|_{\frac{p}{2},I}
\big)
\\
&\leq
C_{p,A,\Sigma}\big(
\|V\|_{\infty,I}|I|
+
(\|V\|_{\infty,I}+\|V\|_{p,I}+\|R^V\|_{\frac{p}{2},I})
\|X\|_{p,I}+
(\|V\|_{\infty,I}+\|V\|_{p,I})\|\bX\|_{\frac{p}{2},I}
\big).
\end{align*}
Next, we choose $I$ small-enough to remove $\|R^V\|_{\frac{p}{2},I}$ in the inequality above.  
Let $\alpha_1:=\left(1/(2C_{p,A,\Sigma })\right)^p \leq 1$ and choose  $I$ small-enough so that $\|X\|_{p,I}^p\leq\alpha_1$. Then,  $C_{p,A,\Sigma }\|X\|_{p,I}\leq\frac{1}{2}$, 
and we obtain
\begin{align}
\|R^V\|_{\frac{p}{2},I}  
&\leq
C_{p,A,\Sigma } 
(\|V\|_{p,I}+\|V\|_{\infty,I})
(|I|+\|\mbX\|_{p,I})
\label{eq:R^Vp/2_I<=C_Vp_Vinf_I+Xp}
\end{align}
for a new constant $C_{p,A,\Sigma }$. 
Using the equation above, $V_{s,t}=\Sigma _sV_sX_{s,t}+R^V_{s,t}$, and $\|\Sigma \|_{\infty,I}\leq C_\Sigma $, we obtain
$$
\|V\|_{p,I}
\mathop{\leq}^{\eqref{eq:path_finite_var:sum_pvars}}
C_p( 
\|\Sigma\|_{\infty,I}\|V\|_{\infty,I}\|X\|_{p,I}+\|R^V\|_{\frac{p}{2},I}
)
\mathop{\leq}^{\eqref{eq:R^Vp/2_I<=C_Vp_Vinf_I+Xp}} 
C_{p,A,\Sigma } 
(
\|V\|_{p,I}
+
\|V\|_{\infty,I}
)
(
|I|+\|\mbX\|_{p,I}
).
$$
Define the control $w$ by $w(s,t)=C_p(|t-s|+w_{\mbX}(s,t))$ with $w_{\mbX}(s,t)=\|X\|_{p,[s,t]}^p+\|\bX\|_{\frac{p}{2},[s,t]}^{\frac{p}{2}}$ and $C_p=6^p$ as in Corollary \ref{cor:Nalpha:NX_NXtilde_NT:small_intervals} \rev{(with $\widetilde{\mbX}=0$)}, and let $\alpha:=\left(1/(2C_{p,A,\Sigma})\right)^p\leq \alpha_1\leq 1$. Then, by choosing $I$ small-enough so that $w(s,t)\leq\alpha$, we have $|I|+\|\mbX\|_{p,I}\mathop{\leq}\limits^{\eqref{eq:|X|+|Xtilde|+|dt|<=alpha^p}}\alpha^\frac{1}{p}$, so that    
$C_{p,A,\Sigma }(|I|+\|\mbX\|_{p,I})\leq\frac{1}{2}$. Then,  we obtain
\begin{align}
\nonumber
\|V\|_{p,I}
&\leq  
C_{p,A,\Sigma } 
\|V\|_{\infty,I}
(|I|+\|\mbX\|_{p,I})
\mathop{\leq}^{\eqref{eq:mbX^p<=wX}} 
C_{p,A,\Sigma } 
\|V\|_{\infty,I}
(|I|+w_{\mbX}(s,t)^\frac{1}{p})
\\
\label{eq:Vp_I<=C_Vinf_w(s,t)^1/p}
&\leq
C_{p,A,\Sigma } 
\|V\|_{\infty,I}
(|I|^\frac{1}{p}+w_{\mbX}(s,t)^\frac{1}{p})
\leq
C_{p,A,\Sigma } 
\|V\|_{\infty,I}
w(s,t)^\frac{1}{p}
\end{align}
for a new constant $C_{p,A,\Sigma }\geq 0$, 
where we used $|I|\leq|I|^\frac{1}{p}$ since $|I|\leq 1$ and $p\geq 1$, and $(|a|+|b|)^p\leq 2^p(|a|^p+|b|^p)$ in the last line. 
Then, by applying the rough Gr\"onwall Lemma (Lemma \ref{lem:rough_gronwall}), we obtain 
\begin{equation}\label{eq:Vinf<=CVs:gronwall}
\|V\|_{\infty,I}
\leq
2\exp\left(w(s,t)/(C\alpha)\right)
\|V_s\|
\leq
2\exp\left(1/C\right)
\|V_s\|,
\end{equation}
for any interval $I\subseteq[0,T]$ such that $w(s,t)\leq\alpha$, and where $C>0$ and $\alpha>0$ only depend on $(p,A,\Sigma)$. 
By sequentially using this inequality on the greedy partition $\{\tau_i\}_{i=1}^{N_{\alpha,[0,T]}(w)+1}$ of $[0,T]$ (which satisfies $w(\tau_i,\tau_{i+1})\leq\alpha$ for all $i$, see Definition \ref{def:Nalpha}) as in the end of the proof of  Proposition \ref{prop:rdes:error_bound:entire_interval}, we obtain
$$
\|V\|_{\infty,[0,T]}
\leq 
\big(2\exp(1/C)\big)^{N_{\alpha,[0,T]}(w)+1}
\|v\|
\leq 
\tilde{C}\exp\big(\tilde{C}N_{\alpha,[0,T]}(w)\big)\|v\|
$$
for a new constant $\tilde{C}>0$. 
Finally,  $N_{\alpha,[s,t]}(w)\leq C_p(
N_{\alpha,[s,t]}(\mbX)+T/\alpha+1
)$ by \eqref{eq:Nalpha<=3CpNalpha_X_and_time} in Corollary \ref{cor:Nalpha:NX_NXtilde_NT:small_intervals}, and the desired inequality \eqref{eq:RDE:linear:V:bound} follows.

To show the inequality \eqref{eq:RDE:linear:Vp_V'p_RVp/2:bound}, from the previous inequalities, %
we have
\begin{align*}
\|V\|_{p,I}
+
\|V'\|_{p,I}
+
\|R^V\|_{\frac{p}{2},I}
&\mathop{\leq}^{
\eqref{eq:Vp_I<=C_Vinf_w(s,t)^1/p},\eqref{eq:V'=SigmaV_p<=sum|Sigma|p(|V|inf+|V|p},\eqref{eq:R^Vp/2_I<=C_Vp_Vinf_I+Xp}}
C\|V\|_{\infty,I}
\big(
w(s,t)^\frac{1}{p}
+
(1+w(s,t)^\frac{1}{p})
+
(w(s,t)^\frac{1}{p}+1)
\big)
\\
&\mathop{\leq}^{\eqref{eq:Vinf<=CVs:gronwall}}
C\|V_s\|
\end{align*}
for any interval $I\subseteq[0,T]$ with $w(s,t)\leq\alpha$. Again, by sequentially using this inequality on the greedy partition $\{\tau_i\}_{i=1}^{N_{\alpha,[0,T]}(w)+1}$  as in the end of the proof of  Proposition \ref{prop:rdes:error_bound:entire_interval},  the  inequality \eqref{eq:RDE:linear:Vp_V'p_RVp/2:bound}  follows.
\end{proof}

\begin{proof}[Proof of Lemma \ref{lem:rde:linearized:bounded_solutions}]
By Theorem \ref{thm:rdes:existence_unicity} and  Corollary  \ref{cor:rde:linear:existence_uniqueness}, there exists unique solutions to the two RDEs.  
First, by Proposition \ref{prop:bounds_pvars_solutions_RDEs}, there exists two constants $C_{p,b,\sigma}\geq 1$ and $0<\alpha_{p,b,\sigma}<1$ such that 
\begin{align*} 
\|Y\|_{p,I}+\|R^Y\|_{\frac{p}{2},I} + 
\|Y'\|_{\infty,I}+\|Y'\|_{p,I}+\|R^Y\|_{\frac{p}{2},I}
&\leq
C_{p,b,\sigma}.
\end{align*}
for any interval $I=[t_0,t_1]\subseteq[0,T]$ such that $\|\mbX\|_{p,I}+|I|\leq\alpha_{p,b,\sigma}^\frac{1}{p}$, noting that $Y'=\sigma(\cdot,Y_\cdot)$ with $\sigma\in C^3_b$. 

Second, let $A=\rev{\nablaof{x}b}(\cdot,Y_\cdot,u_\cdot)$, $\Sigma=\rev{\nablaof{x}\sigma}(\cdot,Y_\cdot)$ and $\Sigma'=\rev{\nablaof{x}^2\sigma}(\cdot,Y_\cdot)Y'_\cdot$. By Lemma \ref{lem:rough_path:sigma(.,Y):controlled}, since $\sigma_x:=\rev{\nablaof{x}\sigma}\in C^2_b$,  
$(\Sigma,\Sigma')\in\sD^p_X$,  $K_{Y,I}=\|Y'_{t_0}\|+\|Y'\|_{p,I}+\|R^Y\|_{\frac{p}{2},I}\leq C_{p,b,\sigma}$,  $M_{Y',I}\leq K_{Y,I}\leq C_{p,b,\sigma}$, and $\|X\|_{p,I}+|I|\leq\alpha_{p,b,\sigma}^\frac{1}{p}$,
\begin{align*}   
\|\Sigma\|_{p,I}
&\mathop{\leq}^{\eqref{eq:sigma(.,Y):pvar}} 
C_p\|\sigma_x\|_{C^1_b}
(
M_{Y',I}\|X\|_{p,I}+\|R^Y\|_{\frac{p}{2},I}
+
|I|
)
\leq C_{p,b,\sigma},
\\ 
\|\Sigma'\|_{p,I}
&\mathop{\leq}^{\eqref{eq:sigma(.,Y)':pvar}} 
C_p\|\sigma_x\|_{C^2_b}
K_{Y,I}(1+K_{Y,I}+|I|)(1+\|X\|_{p,I})
\leq C_{p,b,\sigma},
\\
\|R^{\Sigma}\|_{\frac{p}{2},I}
&\mathop{\leq}^{\eqref{eq:RY:p/2var:Y^2+RY+T}} 
C_p\|\sigma_x\|_{C^2_b}(
\|Y\|_{p,I}^2
+
\|R^Y\|_{\frac{p}{2},I}
+
|I|
)
\leq C_{p,b,\sigma}.
\end{align*}
Also,  $\|A\|_\infty\leq C_b$ and $\|\Sigma\|_{\infty,I}\leq C_\sigma$. 
Thus, for a constant $C_{p,b,\sigma}\geq 1$,  
\begin{align*} 
\|\Sigma\|_{\infty,I} + \|\Sigma\|_{p,I} +\|\Sigma'\|_{p,I}+\|R^{\Sigma}\|_{\frac{p}{2},I}
&\leq
C_{p,b,\sigma}
\end{align*}
for any interval $I\subseteq[0,T]$ such that $\|\mbX\|_{p,I}+|I|\leq
 \alpha_{p,b,\sigma}^\frac{1}{p}$.   
The conclusion follows from Lemma \ref{lem:rde:linear:bounded_solutions}.
\end{proof}

\subsection{Additional proofs for the  PMP  (Section \ref{sec:pmp})}\label{apdx:proofs:pmp}
\subsubsection{Needle-like variations (Section \ref{sec:linearization})}

Corollary \ref{cor:needle_like_error:etas} below extends Proposition \ref{prop:linear_variation} to needle-like variations with multiple spikes. 
As in the deterministic case, the proof proceeds by  induction, using Proposition \ref{prop:linear_variation} for the case with one variation $\pi_1=\{t_1,\eta_1,\bar{u}_1\}$.
\begin{corollary}[Needle-like variations]\label{cor:needle_like_error:etas}
Define $p,T,y,U,u,b,\sigma,\mbX,N_{\alpha,[0,T]}(\mbX),Y,Y'$  as in Proposition \ref{prop:linear_variation}. Given $q\in\N$, let $0<t_1<\dots<t_q<T$ be Lebesgue points 
(Definition \ref{def:lebesgue_point}), 
$\bar{u}_1,\dots,\bar{u}_q\in U$, 
 $0\leq \eta_i< t_{i+1}-t_i$ for $i=1,\dots, q-1$ and $0\leq\eta_q< T-t_q$, and define the needle-like variation $\pi=\{t_1,\dots,t_q,\eta_1,\dots,\eta_q,\bar{u}_1,\dots,\bar{u}_q\}$ of $u$ as the control $u^\pi$ defined by%
 $$
 u^\pi_t=\begin{cases}
 \bar{u}_i\quad&\text{if }t\in[t_i,t_i+\eta_i],
 \\
 u_t&\text{otherwise}.
 \end{cases}
 $$ 
Let $(Y^\pi,(Y^\pi)')\in\sD^p_X$ and  $(V^{\pi_i},(V^{\pi_i})')\in\sD^p_X$ for $i=1,\dots,q$ be the unique solutions to the RDEs
\begin{align*}
Y^\pi_t&=y+\int_0^tb(s,Y^\pi_s,u^\pi_s)\dd s+\int_0^t\sigma(s,Y^\pi_s)\dd\mbX_s,
 &&\hspace{-10mm}
 t\in[0,T], 
\\
V^{\pi_i}_t &= 
V^{\pi_i}_{t_i}
+
\int_{t_i}^t
\rev{\nablaof{x}b}(s,Y_s,u_s)V^{\pi_i}_s\dd s +
\int_{t_i}^t
\rev{\nablaof{x}\sigma}(s,Y_s)V^{\pi_i}_s\dd\mbX_s,
 &&\hspace{-10mm}
 t\in[t_i,T],
\\
V^{\pi_i}_t&=b(t_i,Y_{t_i},\bar{u}_i)-b(t_i,Y_{t_i},u_{t_i}),
 &&\hspace{-10mm}
 t\in[0,t_i].
\end{align*}
with  
 $(Y^\revproof{\pi})'=\sigma(\cdot,Y^\revproof{\pi}_\cdot)$  and   $(V^{\pi_\revproof{i}})'=\rev{\nablaof{x}\sigma}(\cdot,Y_\cdot)V^{\pi_\revproof{i}}_\cdot$. 
Then, there exists constants $C_{p,T,b,\sigma}>0$ and $0<\alpha_{p,b,\sigma}<1$ such that %
\begin{align}
\label{eq:needle_like_deltasols:multi}
\|Y^\pi-Y\|_{\infty,[0,T]}
&\leq 
C_{p,T,b,\sigma}\exp\left(
C_{p,T,b,\sigma}N_{\alpha_{p,b,\sigma},[0,T]}(\mbX)
\right)
\sum_{i=1}^q\eta_i,
\\
\label{eq:needle_like_deltasols_variation:multi}
\Big\|Y^\pi-Y-\sum_{i=1}^q\eta_iV^{\pi_i}\Big\|_{\infty,[t_\revproof{q}\revproof{+\eta_q},T]}&\leq 
C_{p,T,b,\sigma}\exp\left(
C_{p,T,b,\sigma}N_{\alpha_{p,b,\sigma},[0,T]}(\mbX)
\right)
\revproof{o\bigg(\sum_{i=1}^q\eta_i\bigg)}.
\end{align}
\end{corollary}
\begin{proof}[Proof of Corollary \ref{cor:needle_like_error:etas}] 
As in Proposition \ref{prop:linear_variation}, the first inequality  \eqref{eq:needle_like_deltasols:multi}  follows 
from sequentially applying 
\eqref{eq:RDE:DY'+dRY:close:full_interval:u_subinterval}
and 
\eqref{eq:RDE:DY:close:full_interval:u_subinterval} 
on  
$[0,t_1]
\cup
[t_1,t_1+\eta_1]
\cup
[t_1+\eta_1,t_2]
\cup
\dots
\cup
[t_q,t_q+\eta_q]
\cup
[t_q+\eta_q,T]=[0,T]$, and concluding with 
\eqref{eq:Nalpha<=3CpNalpha_X_and_time} in Corollary \ref{cor:Nalpha:NX_NXtilde_NT:small_intervals}.

 \noindent\begin{minipage}[c]{0.55\textwidth}
Next, we show the second inequality \eqref{eq:needle_like_deltasols_variation:multi}. 
The case $q=1$ is in Proposition \ref{prop:linear_variation}, so we prove \eqref{eq:needle_like_deltasols_variation:multi}  for $q\geq 2$ by induction,   assuming that it holds for $q$ and proving it for $q+1$ for a needle-like variation $\pi=\pi_{q+1}=\{t_1,\dots,t_{q+1},\eta_1,\dots,\eta_{q+1},\bar{u}_1,\dots,\bar{u}_{q+1}\}$ with associated  control $u^\pi$. Define the needle-like variation 
$$
\pi_q=\{t_1,\dots,t_q,\eta_1,\dots,\eta_q,\bar{u}_1,\dots,\bar{u}_q\}\subset\pi_{q+1}=\pi,
$$  of $u$ as the controls $u^{\pi_q}$  defined by $u^{\pi_q}_t=\bar{u}_i$ if $t\in [t_i,t_i+\eta_i]$ for $i=1,\dots,q$, and $u^{\pi_q}_t=u_t$ otherwise,   and 
let $(Y^{\pi_q},(Y^{\pi_q})')\in\sD^p_X$ be the unique solution to the RDE
{\small
\begin{align*} 
Y^{\pi_q}_t=y+\int_0^tb(s,Y^{\pi_q}_s,u^{\pi_q}_s)\dd s+\int_0^t\sigma(s,Y^{\pi_q}_s)\dd\mbX_s,
\  t\in[0,T].
\end{align*}
}%
Let $(\widetilde{V},(\widetilde{V})')\in\sD^p_X$ be the unique solution to the RDE
\end{minipage}
\hspace{0.007\textwidth}
\begin{minipage}[c]{0.43\textwidth}
\includegraphics[width=1\textwidth]{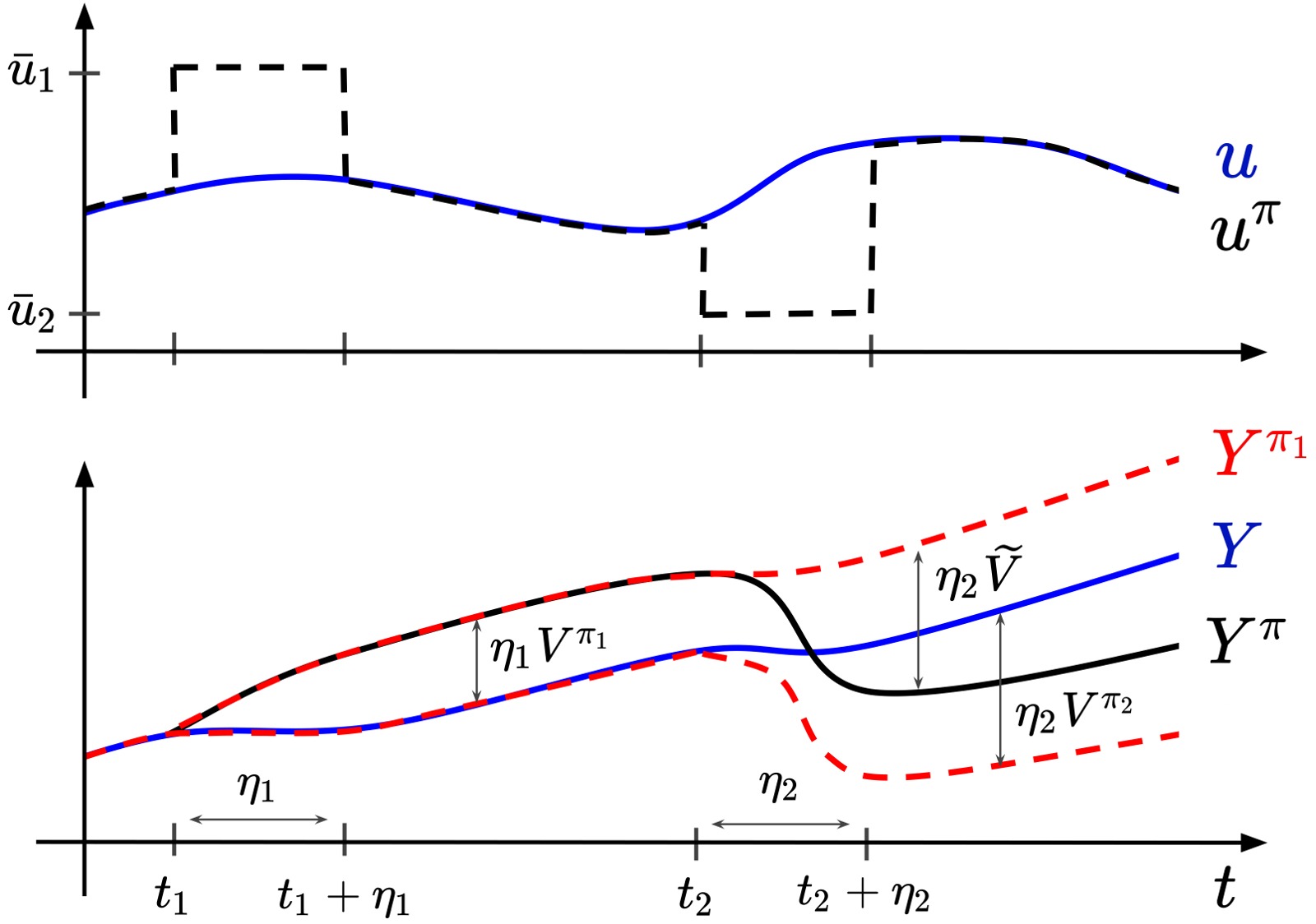}
\captionof{figure}{Needle-like variations with two spikes.}
\label{fig:needlelike}
\end{minipage}
\begin{align*}
\widetilde{V}_t &= 
\widetilde{V}_{t_{q+1}}
+
\int_{t_{q+1}}^t
\rev{\nablaof{x}b}(s,Y^{\pi_q}_s,u_s)\widetilde{V}_s\dd s +
\int_{t_{q+1}}^t
\rev{\nablaof{x}\sigma}(s,Y^{\pi_q}_s)\widetilde{V}_s\dd\mbX_s,
 &&t\in[t_{q+1},T],
\\
\widetilde{V}_t
&=
b(t_{q+1},Y^{\pi_q}_{t_{q+1}},\bar{u}_{q+1})
-
b(t_{q+1},Y^{\pi_q}_{t_{q+1}},u_{t_{q+1}}),
&&t\in[0,t_{q+1}],
\end{align*}
see Figure \ref{fig:needlelike} (corresponding to $q+1=2$).
Then,
the error $Y^\pi-Y-\sum_{i=1}^{q+1}\eta_iV^{\pi_i}$ can   be decomposed %
\revproof{into three terms:}
 \begin{align}\label{eq:needle_like_deltasols_variation:multi:threeterms}
 Y^{\pi}-Y-\sum_{i=1}^{q+1}\eta_iV^{\pi_i}
 &=
 (Y^{\pi}-Y^{\pi_q}-\eta_{q+1}\widetilde{V})
 + 
 \eta_{q+1}(\widetilde{V} -V^{\pi_{q+1}})
 +
 \Big(
 Y^{\pi_q}-Y-\sum_{i=1}^q\eta_iV^{\pi_i}
 \Big).
\end{align} 
\revproof{Let }$C=C_{p,T,b,\sigma}$ and $\alpha=\alpha_{p,b,\sigma}$ \revproof{denote appropriate constants.}

\revproofnote{Here, it is not enough to use the bound
\begin{align*} 
\Big\|Y^\pi-Y^{\pi_q}-\eta_{q+1}\widetilde{V}\Big\|_{\infty,[t_{q+1}\revproof{+\eta_{q+1}},T]}
&\mathop{\leq}^{\eqref{eq:needle_like_deltasols_variation}}
C\exp\left(
CN_{\alpha,[0,T]}(\mbX)
\right)
\revproof{o(\eta_{q+1})}.
\quad &&(\text{base case }q=1)
\end{align*}
Since this bound is evaluated around $Y^{\pi_q}$, its constant $C$ depends on  $\eta_1,\dots,\eta_q$.
}

\revproof{The first term of \eqref{eq:needle_like_deltasols_variation:multi:threeterms} can be bounded as
\begin{align} 
\label{eq:eq:needle_like_deltasols_variation:multi:dYpi:bound}
\Big\|Y^\pi-Y^{\pi_q}-\eta_{q+1}\widetilde{V}\Big\|_{\infty,[t_{q+1}\revproof{+\eta_{q+1}},T]}
&\leq
C\exp(
CN_{\alpha,[0,T]}(\mathbf X)
)
\bigg(
o(\eta_{q+1})
+
\eta_{q+1}\sum_{i=1}^q\eta_i
+
\eta_{q+1}^2
\bigg).
\end{align} 
Indeed,
\begin{align}
\label{eq:eq:needle_like_deltasols_variation:multi:Ypi}
\|Y^{\pi_q}\|_p+
\|(Y^{\pi_q})'\|_p+
\|R^{Y^{\pi_q}}\|_{\frac{p}{2}}
&\mathop{\leq}^{\eqref{eq:RDE:Yp_Y'p_RYp/2:bound:full_interval}
} 
C\exp(
CN_{\alpha,[0,T]}(\mathbf X)
),
\\
\label{eq:eq:needle_like_deltasols_variation:multi:dYpi}
\|Y^\pi{-}Y^{\pi_q}\|_\infty+
\|Y^\pi{-}Y^{\pi_q}\|_p+
\|(Y^\pi)'{-}(Y^{\pi_q})'\|_p+
\|R^{Y^\pi}{-}R^{Y^{\pi_q}}\|_{\frac{p}{2}}
&\mathop{\leq}^{\eqref{eq:RDE:DY'+dRY:close:full_interval:needle_like}
} 
C\exp(
CN_{\alpha,[0,T]}(\mathbf X)
)
\eta_{q+1},
\end{align}
independently of $\eta_1,\dots,\eta_{q+1}$. Then, we repeat the computation in the proof of the single needle-like variation case in Proposition \ref{prop:linear_variation} with 
$$
\Delta Y:=Y^\pi-Y^{\pi_q},
\qquad
\Delta := Y^\pi-Y^{\pi_q}-\eta_{q+1}\widetilde V,
$$
so
\begin{small}
\begin{align*}
\Delta_t
\mathop{=}^{\eqref{eq:Delta(t):full}}&
\int_{t_{q+1}}^{t_{q+1}+\eta_{q+1}}
\left(
b(s,Y^\pi_s,\bar u_{q+1})
-b(s,Y^\pi_s,u_s)
\right)
\dd s
-\eta_{q+1}\widetilde{V}_{t_{q+1}}
\\
&\hspace{-1cm}+
\int_{t_{q+1}}^t\int_0^1
\left(
\nabla_xb
(s,Y^{\pi_q}_s+\theta\Delta Y_s,u_s)
-\nabla_xb(s,Y^{\pi_q}_s,u_s)
\right)
\dd\theta\Delta Y_s\dd s
+
\int_{t_{q+1}}^t\int_0^1
\nabla_x^2\sigma
(s,Y^{\pi_q}_s+\theta\Delta Y_s)
(1-\theta)\dd\theta\Delta Y_s^{\otimes2}\dd\mbX_s
\\
&\hspace{-1cm}+
\int_{t_{q+1}}^t
\nabla_xb(s,Y^{\pi_q}_s,u_s)\Delta_s\dd s
+
\int_{t_{q+1}}^t
\nabla_x\sigma(s,Y^{\pi_q}_s)\Delta_s\dd\mbX_s .
\end{align*}
\end{small}
\hspace{-1mm}The first term above can be bounded as
\begin{small}
\begin{align*}
&\bigg\|
\int_{t_{q+1}}^{t_{q+1}+\eta_{q+1}}
\left(
b(s,Y^\pi_s,\bar u_{q+1})
-b(s,Y^\pi_s,u_s)
\right)
\dd s
-\eta_{q+1}
(
b(t_{q+1},Y^{\pi_q}_{t_{q+1}},\bar u_{q+1})
-b(t_{q+1},Y^{\pi_q}_{t_{q+1}},u_{t_{q+1}})
)
\bigg\|
\\
&\leq
\bigg\|
\int_{t_{q+1}}^{t_{q+1}+\eta_{q+1}}
\left(
b(s,Y_s,\bar u_{q+1})-b(s,Y_s,u_s)
\right)
\dd s 
-
\eta_{q+1}(
b(t_{q+1},Y_{t_{q+1}},\bar u_{q+1})
-b(t_{q+1},Y_{t_{q+1}},u_{t_{q+1}})
)
\bigg\|
\\
&
\qquad+
C\eta_{q+1}\|Y^\pi-Y\|_\infty
+
C\eta_{q+1}\|Y^{\pi_q}-Y\|_\infty
\\
&\mathop{\leq}^{\eqref{eq:Delta(t):term1},\eqref{eq:needle_like_deltasols:multi},\eqref{eq:eq:needle_like_deltasols_variation:multi:dYpi}}
C\exp(CN_{\alpha,[0,T]}(\mathbf X))
\left(
o(\eta_{q+1})
+
\eta_{q+1}\sum_{i=1}^q\eta_i
+
\eta_{q+1}^2
\right),
\end{align*}
\end{small}
\hspace{-3mm}
and the second and third can be bounded by $o(\eta_{q+1})$ using \eqref{eq:eq:needle_like_deltasols_variation:multi:Ypi} and \eqref{eq:eq:needle_like_deltasols_variation:multi:dYpi}. 
Thus, we obtained the desired bound \eqref{eq:eq:needle_like_deltasols_variation:multi:dYpi:bound}, 
with $o(\eta_{q+1})$ is independent of \(\eta_1,\ldots,\eta_q\).  
}

\revproof{The second term of \eqref{eq:needle_like_deltasols_variation:multi:threeterms} can be bounded as} $\|\widetilde{V} -V^{\pi_{q+1}}\|_{\infty,[t_{q+1},T]}\leq C\exp\left(
CN_{\alpha,[0,T]}(\mbX)
\right)
\sum_{i=1}^q\eta_i
$.   
Indeed,  for any $t\geq t_{q+1}$, the error can be decomposed as
{\small
\begin{align*}
\widetilde{V}_t-V^{\pi_{q+1}}_t
&=
\left(
\widetilde{V}_{t_{q+1}}-V^{\pi_{q+1}}_{t_{q+1}}
\right)
+
\int_{t_{q+1}}^t
\left(
\rev{\nablaof{x}b}(s,Y^{\pi_q}_s,u_s)\widetilde{V}_s
-
\rev{\nablaof{x}b}(s,Y_s,u_s)V^{\pi_{q+1}}_s
\right)
\dd s 
\\
&\quad+
\int_{t_{q+1}}^t
\left(
\rev{\nablaof{x}\sigma}(s,Y^{\pi_q}_s)\widetilde{V}_s
-
\rev{\nablaof{x}\sigma}(s,Y_s)V^{\pi_{q+1}}_s
\right)
\dd\mbX_s
\\
&=
\left(\widetilde{V}_{t_{q+1}}-V^{\pi_{q+1}}_{t_{q+1}}\right)
+
\int_{t_{q+1}}^t
\left(
\rev{\nablaof{x}b}(s,Y^{\pi_q}_s,u_s)(
\widetilde{V}_s-V^{\pi_{q+1}}_s
)
+
\left(\rev{\nablaof{x}b}(s,Y^{\pi_q}_s,u_s)
-
\rev{\nablaof{x}b}(s,Y_s,u_s)
\right)
V^{\pi_{q+1}}_s
\right)
\dd s 
\\
&\quad+
\int_{t_{q+1}}^t
\left(
\rev{\nablaof{x}\sigma}(s,Y^{\pi_q}_s)(\widetilde{V}_s-V^{\pi_{q+1}}_s)
+
\left(
\rev{\nablaof{x}\sigma}(s,Y^{\pi_q}_s)
-
\rev{\nablaof{x}\sigma}(s,Y_s)\right)
V^{\pi_{q+1}}_s
\right)
\dd\mbX_s
\\
&=
\left(\widetilde{V}_{t_{q+1}}-V^{\pi_{q+1}}_{t_{q+1}}\right)
+
\int_{t_{q+1}}^t
\left(\rev{\nablaof{x}b}(s,Y^{\pi_q}_s,u_s)
-
\rev{\nablaof{x}b}(s,Y_s,u_s)
\right)
V^{\pi_{q+1}}_s
\dd s 
\\
&\quad+ 
\int_{t_{q+1}}^t 
\left(
\rev{\nablaof{x}\sigma}(s,Y^{\pi_q}_s)
-
\rev{\nablaof{x}\sigma}(s,Y_s)\right)
V^{\pi_{q+1}}_s 
\dd\mbX_s
\\
&\quad
+
\int_{t_{q+1}}^t
\rev{\nablaof{x}b}(s,Y^{\pi_q}_s,u_s)(
\widetilde{V}_s-V^{\pi_{q+1}}_s
)
+
\int_{t_{q+1}}^t 
\rev{\nablaof{x}\sigma}(s,Y^{\pi_q}_s)(\widetilde{V}_s-V^{\pi_{q+1}}_s)
\dd\mbX_s.
\end{align*}
}%
By following similar arguments as  in  Proposition \ref{prop:linear_variation}, %
the first three terms in the equation above 
can be bounded by $C\exp\left(
CN_{\alpha,[0,T]}(\mbX)
\right)
\sum_{i=1}^q\eta_i
$, for example,
\begin{align*}
\|\widetilde{V}_{t_{q+1}}-V^{\pi_{q+1}}_{t_{q+1}}\|
&=
\|
b(t_{q+1},Y^{\pi_q}_{t_{q+1}},\bar{u}_{q+1})
-
b(t_{q+1},Y^{\pi_q}_{t_{q+1}},u_{t_{q+1}})
-
(b(t_{q+1},Y_{t_{q+1}},\bar{u}_{q+1})
-
b(t_{q+1},Y_{t_{q+1}},u_{t_{q+1}}))
\|
\\
&\leq
\|b(t_{q+1},Y^{\pi_q}_{t_{q+1}},\bar{u}_{q+1})
-
b(t_{q+1},Y_{t_{q+1}},\bar{u}_{q+1})\|
+
\|
b(t_{q+1},Y^{\pi_q}_{t_{q+1}},u_{t_{q+1}})
-
b(t_{q+1},Y_{t_{q+1}},u_{t_{q+1}})\|
\\
&\leq C_b\|Y^{\pi_q}_{t_{q+1}}-Y_{t_{q+1}}\|
\mathop{\leq}^{\eqref{eq:needle_like_deltasols:multi}}
C\exp\left(
CN_{\alpha,[0,T]}(\mbX)
\right)
\sum_{i=1}^q\eta_i,
\end{align*} 
and the other two integrals can be bounded as in the proof of Proposition \ref{prop:linear_variation} using the bounds \eqref{eq:rde:linearized:Vp_V'p_RVp/2:bound}
and \eqref{eq:rde:linearized:bounded_solutions}   for  $\|V\|_\infty,\|V\|_p,\|V'\|_p,\|R^V\|_\frac{p}{2}$ in Lemma \ref{lem:rde:linearized:bounded_solutions}. 
Finally, by defining $\Delta:=\widetilde{V}-V^{\pi_{q+1}}$ and 
looking at the linear RDE $\Delta_t=\Delta_{t_{q+1}}+
\int_{t_{q+1}}^t
\rev{\nablaof{x}b}(s,Y^{\pi_q}_s,u_s)\Delta_s
+
\int_{t_{q+1}}^t 
\rev{\nablaof{x}\sigma}(s,Y^{\pi_q}_s)\Delta_s
\dd\mbX_s$ whose initial value satisfies $\|\Delta_{t_{q+1}}\|\leq C\exp\left(
CN_{\alpha,[0,T]}(\mbX)
\right)
\sum_{i=1}^q\eta_i$, we conclude  that  
$\|\Delta\|_{\infty,[t_{q+1},T]}=\|\widetilde{V} -V^{\pi_{q+1}}\|_{\infty,[t_{q+1},T]}\leq C\exp\left(
CN_{\alpha,[0,T]}(\mbX)
\right)
\sum_{i=1}^q\eta_i$ using \eqref{eq:rde:linearized:bounded_solutions} in Lemma  \ref{lem:rde:linearized:bounded_solutions}. %

\revproof{The third term of \eqref{eq:needle_like_deltasols_variation:multi:threeterms} can be bounded as}
\begin{align*} 
\Big\|Y^{\pi_q}-Y-\sum_{i=1}^q\eta_iV^{\pi_i}\Big\|_{\infty,[t_{q+1}\revproof{+\eta_{q+1}},T]}
&\mathop{\leq}^{\eqref{eq:needle_like_deltasols_variation:multi}}
C\exp\left(
CN_{\alpha,[0,T]}(\mbX)
\right)
\revproof{o\bigg(\sum_{i=1}^q\eta_i\bigg)}
&&(\text{induction step for }q), 
\end{align*} 

\revproof{Thus,  by combining the last three terms in \eqref{eq:needle_like_deltasols_variation:multi:threeterms}, we obtain 
\begin{align*}
\bigg\|
Y^{\pi}-Y-\sum_{i=1}^{q+1}\eta_iV^{\pi_i}
\bigg\|_{\infty,[t_{q+1}\revproof{+\eta_{q+1}},T]}
 &\leq
 C\exp\left(
CN_{\alpha,[0,T]}(\mbX)
\right)\left(
o(\eta_{q+1})
+
\eta_{q+1}\sum_{i=1}^q\eta_i
+
\eta_{q+1}^2
+
o\left(\sum_{i=1}^q\eta_i\right)
\right).
\end{align*}
Dividing by $\sum_{i=1}^{q+1}\eta_i$ and taking the limit as $\sum_{i=1}^{q+1}\eta_i\to 0$  concludes the proof.} 
\end{proof}

\subsection{Additional details for the indirect shooting method (Section \ref{sec:example})}\label{apdx:proofs:example}
We provide additional details for the integration schemes used for the Stratonovich SDE and the coupled RDE used in the \texttt{Direct} method and the \texttt{Indirect} shooting method in Section \ref{sec:example}. We only describe the case for the open-loop problem $\olocp$, as the feedback problem $\fbocp$ only has a different drift term.

1) \texttt{Direct} method: We discretize the Stratonovich SDE in the \texttt{Direct} problem using a Milstein scheme. 
Since $\sigma$ is diagonal  and each $\frac{\partial\sigma^{jj}}{\partial x^\ell}=0$ for $\ell\neq j$,   each $j$-th component $[\hat{x}^i]^j$ of $\hat{x}^i$ is approximated as
\begin{align}\label{eq:milstein}
\hspace{-2mm}
[\hat{x}_{k+1}^i]^j
= 
[
\hat{x}_k+ 
(A(\hat{x}_k^i)\hat{x}_k^i+\bar{B}\hat{u}_k)\Delta t + \sigma(\hat{x}_k^i) \Delta B_k^i
]^j + 
\frac{1}{2}
\frac{\partial\sigma^{jj}}{\partial x^j}(\hat{x}_k^i)\sigma^{jj}(\hat{x}_k^i)
([\Delta B_k^i]^j)^2,
\    
j=1,\dots,n,
\end{align}
where  $k=0,\dots,N-1$, 
$\Delta t=\frac{T}{N}$, and 
$\Delta B_k^i=B_{(k+1)\Delta t}^i-B_{k\Delta t}^i$, see \cite[equation (3.12), Chapter 10.3]{Kloeden1992}.

2) \texttt{Indirect} method: 
To implement the map $F:\R^{Mn}\to\R^{Mn},\  
(p_0^i)_{i=1}^M
\mapsto
(p_T^i)_{i=1}^M$, we numerically integrate  the RDE in \eqref{eq:example:shooting_problem} using the estimate \eqref{eq:rough_int:error_bound} for rough integrals as%
\begin{align}\label{eq:example:rde_discretization}
\bigg[\begin{matrix}
\hat{x}_{k+1}^i
\\
\hat{p}_{k+1}^i
\end{matrix}
\bigg]
&=
\begin{bmatrix}
\hat{x}_k^i
\\
\hat{p}_k^i
\end{bmatrix}
+ 
\bar{b}\bigg(\hspace{-1mm}\begin{bmatrix}
\hat{x}_k^i
\\
\hat{p}_k^i
\end{bmatrix}\hspace{-1mm},u_k^M
\hspace{-1mm}\bigg)
\Delta t 
+ 
\bar\sigma\bigg(\hspace{-1mm}\begin{bmatrix}
\hat{x}_k^i
\\
\hat{p}_k^i
\end{bmatrix}\hspace{-1mm}\bigg)
B_{k\Delta t,(k+1)\Delta t}^i
+
\nabla\bar\sigma\bigg(\hspace{-1mm}\begin{bmatrix}
\hat{x}_k^i
\\
\hat{p}_k^i
\end{bmatrix}\hspace{-1mm}\bigg) \bar\sigma\bigg(\hspace{-1mm}\begin{bmatrix}
\hat{x}_k^i
\\
\hat{p}_k^i
\end{bmatrix}\hspace{-1mm}\bigg)
\bB^i_{k\Delta t,(k+1)\Delta t},
\end{align}
\revproofnote{Note the sign correction for $\bar\sigma$:}
where %
$\bar{b}((x,\rev{p}),u)=(b(x,u),-\rev{\nablaof{x}H}(x,u,p,\mathfrak{p}_0))$  and  %
$\bar\sigma(x,p)=(\sigma(x),\revproof{-}\rev{\nablaof{x}\sigma}(x)^\top p)$ denote the augmented drift and diffusion, $k=0,\dots,N-1$, and $\Delta t=\frac{T}{N}$. 
Because $\sigma(x)\propto\textrm{diag}(x)$ is diagonal, only the diagonal elements $[\bB^i_{k\Delta t,(k+1)\Delta t}]^{jj}=\frac{1}{2}([B_{k\Delta t,(k+1)\Delta t}^i]^j)^2$ are required to evaluate  \eqref{eq:example:rde_discretization}, and one observes that the integration rule \eqref{eq:example:rde_discretization} for $\hat{x}^i$ coincides with   
\eqref{eq:milstein}, as we show next. 
The tensor  $\nabla\bar\sigma(x,p) \bar\sigma(x,p) \in\R^{2n\times n\times n}$ is given by
\begin{align*}
\nabla\bar\sigma(x,p) \bar\sigma(x,p) 
=
\begin{bmatrix}
\rev{\nablaof{x}\sigma}(x)\sigma(x)
\revproof{-}
\rev{\nablaof{p}\sigma}(x)\rev{\nablaof{x}\sigma}(x)^\top p
\\
-\rev{\nablaof{x}}\left(
\rev{\nablaof{x}\sigma}(x)^\top p
\right)
\sigma(x)
\revproof{+}\rev{\nablaof{p}}\left(
\rev{\nablaof{x}\sigma}(x)^\top p
\right)
\rev{\nablaof{x}\sigma}(x)^\top p
\end{bmatrix}.
\end{align*}
Since $\sigma(x)\propto\textrm{diag}(x)$, 
{\small\allowdisplaybreaks
\begin{align*}
\left[
\rev{\nablaof{x}\sigma}(x)\sigma(x)
\right]^{ijk}
&=
\sum_{\ell=1}^n
\frac{\partial\sigma^{ij}}{\partial x^\ell}(x)\sigma^{\ell k}(x)
=
\begin{cases}
\frac{\partial\sigma^{ii}}{\partial x^i}(x)\sigma^{ii}(x) & \text{if }i=j=k,
\\
0&\text{otherwise},
\end{cases}
\\
\rev{\nablaof{p}\sigma}(x)\rev{\nablaof{x}\sigma}(x)^\top p
&=0.
\\
\left[
\rev{\nablaof{x}}\left(
\rev{\nablaof{x}\sigma}(x)^\top p
\right)
\sigma(x)
\right]^{ijk}
&=
\sum_{\ell=1}^n
\frac{\partial}{\partial x^\ell}
\left(
\sum_{q=1}^n
\frac{\partial\sigma^{qj}}{\partial x^i}(x) 
p^q
\right)
\sigma^{\ell k}(x)
= 
\begin{cases}
\frac{\partial^2\sigma^{ii}}{\partial^2 x^i}(x)
p^i
\sigma^{ii}(x) & \text{if }i=j=k,
\\
0&\text{otherwise},
\end{cases},
\\
\left[
\rev{\nablaof{p}}\left(
\rev{\nablaof{x}\sigma}(x)^\top p
\right)
\rev{\nablaof{x}\sigma}(x)^\top p
\right]^{ijk}
&=
\sum_{\ell=1}^n
\frac{\partial}{\partial p^\ell}
\left[
\left(
\rev{\nablaof{x}\sigma}(x)^\top p
\right)
\right]^{ij}
\left[
\rev{\nablaof{x}\sigma}(x)^\top p
\right]^{\ell k}
\\
&=
\sum_{\ell=1}^n
\frac{\partial}{\partial p^\ell}
\left(
\sum_{q=1}^n
\frac{\partial\sigma^{qj}}{\partial x^i}(x)  
p^q
\right)
\left(
\sum_{q=1}^n
\frac{\partial\sigma^{qk}}{\partial x^\ell}(x)
p^q
\right) 
\\
&=
\sum_{\ell=1}^n
\frac{\partial}{\partial p^\ell}
\left( 
\frac{\partial\sigma^{jj}}{\partial x^i}(x)
p^j
\right)
\left( 
\frac{\partial\sigma^{kk}}{\partial x^\ell}(x)
p^k
\right) 
= 
\frac{\partial\sigma^{jj}}{\partial x^i}(x)
\frac{\partial\sigma^{kk}}{\partial x^j}(x)
p^k
\\
&=
\begin{cases}
\left(
\frac{\partial\sigma^{ii}}{\partial x^i}(x)
\right)^2
p^i 	& \text{if }i=j=k,
\\
0&\text{otherwise}.
\end{cases}
\end{align*}
}%
Thus, the tensor  $\nabla\bar\sigma(x,p) \bar\sigma(x,p)$ is such that only the entries $\left[
\nabla\bar\sigma(x,p) \bar\sigma(x,p) \right]^{j,j,j}$ and $\left[
\nabla\bar\sigma(x,p) \bar\sigma(x,p) \right]^{n+j,j,j}$ are non-zero, and each $j$th index of the right hand side of \eqref{eq:example:rde_discretization} is given by
{\small
\begin{align}
\label{eq:example:rough_integral_same_as_milstein}
\left[\nabla\bar\sigma(x,p)\bar\sigma(x,p)\bB^i_{s,t}\right]^j
&=
\left[
\nabla\bar\sigma(x,p) \bar\sigma(x,p)\right]^{j,j,j}
[\bB^i_{s,t}]^{j,j}
=
\frac{1}{2}
\frac{\partial\sigma^{jj}}{\partial x^j}(x) \sigma^{jj}(x)
([B_{s,t}^i]^j)^2,
\\
\left[\nabla\bar\sigma(x,p)\bar\sigma(x,p)\bB^i_{s,t}\right]^{n+j}
&=
\left[
\nabla\bar\sigma(x,p) \bar\sigma(x,p)\right]^{n+j,j,j}
[\bB^i_{s,t}]^{j,j}
=
-\frac{1}{2}
\left(
\frac{\partial^2\sigma^{jj}}{\partial^2 x^j}(x)
\sigma^{jj}(x)
\revproof{-}
\left(
\frac{\partial\sigma^{jj}}{\partial x^j}(x)\right)^2
\right)
p^j\
([B_{s,t}^i]^j)^2,
\nonumber
\end{align}
}%
where we used $[\bB^i_{s,t}]^{jj}=(\int_s^tB_{s,u}^j\circ\dd B_u^j)^i=\frac{1}{2}([B_{s,t}^i]^j)^2$, where $\int_s^tB_{s,u}^j\circ\dd B_u^j$ is the Stratonovich integral of the $B^j_{s,\cdot}$ against $B^j$.

Thus, only the squared increments $[B_{s,t}^i]^j$ are required to evaluate  \eqref{eq:example:rde_discretization}. Using 
\eqref{eq:example:rough_integral_same_as_milstein}, we also observe that the integration rule \eqref{eq:example:rde_discretization}  for $x$ coincides with the Milstein integration scheme in \eqref{eq:milstein}.

\renewcommand{\baselinestretch}{0.95}
\bibliographystyle{IEEEtran}
\bibliography{main}

\fi
\end{document}